\documentclass[11pt,letterpaper]{amsart}

\usepackage{combelow} 

\usepackage{xspace,amssymb,amsfonts,euscript}
\usepackage{amsthm,amsmath,amscd,stmaryrd,latexsym,youngtab,verbatim,mathrsfs}
\usepackage{palatino, euscript}

\newcommand\blfootnote[1]{
\begingroup
\renewcommand\thefootnote{}\footnote{#1}%
\addtocounter{footnote}{-1}%
\endgroup
}

\usepackage{graphicx}
\usepackage[all]{xy}
\SelectTips{cm}{}

\usepackage{tikz, tikz-cd}
\usetikzlibrary{matrix, arrows, calc}
\tikzset{frontline/.style={preaction={draw=white,-,line width=6pt}},}  

\usepackage[dvips]{epsfig}
\usepackage{pinlabel}
\usepackage{psfrag}

\newcommand{\ig}[2]{\vcenter{\xy (0,0)*{\includegraphics[scale=#1]{fig/#2}} \endxy}}
\newcommand{\eqig}[2]{\vcenter{\xy (0,0)*{\includegraphics[scale=#1]{fig/#2}} \endxy}}

\IfFileExists{comments.sty}{\usepackage{comments}}

\usepackage{enumitem}

\RequirePackage{color}
\definecolor{myred}{rgb}{0.75,0,0}
\definecolor{mygreen}{rgb}{0,0.5,0}
\definecolor{myblue}{rgb}{0,0.25,0.65}
\definecolor{references}{rgb}{0,0,1}

\RequirePackage[pdftex,
colorlinks = true,
urlcolor = references, 
citecolor = references, 
linkcolor = references, 
]
 {hyperref}

\newtheorem{thm}{Theorem}[section]
\newtheorem{lemma}[thm]{Lemma}
\newtheorem{theorem}[thm]{Theorem}

\newtheorem{prop}[thm]{Proposition}
\newtheorem{proposition}[thm]{Proposition}
\newtheorem{cor}[thm]{Corollary}
\newtheorem{corollary}[thm]{Corollary}

\newtheorem{conj}[thm]{Conjecture}
\newtheorem{conjecture}[thm]{Conjecture}

\newtheorem{hypothesis}[thm]{Hypothesis}

\newtheorem*{prop*}{Proposition}
\newtheorem*{lemma*}{Lemma}

\theoremstyle{definition}
\newtheorem{defn}[thm]{Definition}
\newtheorem{definition}[thm]{Definition}

\newtheorem{notation}[thm]{Notation}
\newtheorem{ex}[thm]{Example}
\newtheorem{example}[thm]{Example}

\theoremstyle{remark}
\newtheorem{remark}[thm]{Remark}

\numberwithin{equation}{section}

    \def\DM{{\mathbb{D}}}

  \def\hg{{\mathfrak h}}

    \def\KM{{\mathbb{K}}}
    
  \def\mg{{\mathfrak m}}  
    \def\NM{{\mathbb{N}}}

    \def\RM{{\mathbb{R}}}
    \def\SM{{\mathbb{S}}}

    \def\ZM{{\mathbb{Z}}}

  \def\ab{{\mathbf a}}  \def\AC{{\mathcal{A}}}
    \def\BC{{\mathcal{B}}}
\def\CB{{\mathbf C}}  

    \def\DC{{\mathcal{D}}}
\def\EB{{\mathbf E}}  \def\eb{{\mathbf e}}  
    \def\FC{{\mathcal{F}}}
\def\GB{{\mathbf G}}    
\def\HB{{\mathbf H}}    
    \def\IC{{\mathcal{I}}}
    \def\JC{{\mathcal{J}}}
\def\KB{{\mathbf K}}    \def\KC{{\mathcal{K}}} 
    
  \def\mb{{\mathbf m}}  \def\MC{{\mathcal{M}}}
  \def\nb{{\mathbf n}}  
    \def\OC{{\mathcal{O}}}
\def\PB{{\mathbf P}}  \def\pb{{\mathbf p}}  \def\PC{{\mathcal{P}}}
    \def\QC{{\mathcal{Q}}}
  \def\rb{{\mathbf r}}  

  \def\vb{{\mathbf v}}  \def\VC{{\mathcal{V}}}
    
\def\XB{{\mathbf X}}    \def\XC{{\mathcal{X}}}
    \def\YC{{\mathcal{Y}}}

\def\a{\alpha}
\def\b{\beta}
\def\g{\gamma}

\def\d{\delta}
\def\D{\Delta}
\def\e{\varepsilon}
\def\k{\kappa}
\def\l{\lambda}

\def\s{\sigma}

\let\phi=\varphi
\let\tilde=\widetilde

\usepackage{bbm}
\def\C{{\mathbbm C}}
\def\N{{\mathbbm N}}
\def\R{{\mathbbm R}}
\def\Z{{\mathbbm Z}}
\def\Q{{\mathbbm Q}}
\def\1{\mathbbm{1}}
\newcommand{\one}{\1}

\newcommand{\bbeta}{\boldsymbol{\b}}
\newcommand{\shape}{\operatorname{sh}}
\newcommand{\sh}{\operatorname{sh}}
\newcommand{\kbbm}{\mathbbm{k}}

\newcommand{\tbarb}[1]{f_{#1}}
\newcommand{\cbb}{\mathbf{c}}
\newcommand{\abb}{\mathbf{a}}
\newcommand{\rbb}{\mathbf{r}}
\newcommand{\xbb}{\mathbf{x}}
\newcommand{\ind}{\operatorname{ind}}
\newcommand{\shift}[2]{\Sigma_{#1}(#2)}

\newcommand{\colsum}{\mathbf{c}}
\newcommand{\rowsum}{\mathbf{r}}

\newcommand{\ess}{\operatorname{ess}}

\newcommand{\FHilb}{\operatorname{FHilb}}
\newcommand{\ft}{\operatorname{ft}}

\newcommand{\supp}{\operatorname{supp}}

\newcommand{\HHH}{\operatorname{HHH}}

\newcommand{\Hecke}{\HB}

\renewcommand{\setminus}{\smallsetminus}

\newcommand{\Homg}{\operatorname{Hom^{\Z}}}
\newcommand{\Homgg}{\operatorname{Hom^{\Z\times \Z}}}
\newcommand{\Endgg}{\operatorname{End^{\Z\times \Z}}}
\newcommand{\Homc}{\underline{\operatorname{Hom}}}

\newcommand{\Endg}{\operatorname{End}^{\Z}}
\newcommand{\smMatrix}[1]{\left[\begin{smallmatrix}#1\end{smallmatrix}\right]}

\newcommand{\sqmatrix}[1]{\left[\begin{matrix} #1\end{matrix}\right]}

\newcommand{\un}{\underline}

\newcommand{\ot}{\otimes}
\newcommand{\pa}{\partial}
\newcommand{\co}{\colon}
\newcommand{\ip}[1]{\langle #1\rangle}

\renewcommand{\to}{\rightarrow}

\newcommand{\babysumset}{\stackrel{\scriptscriptstyle \sum}{\scriptstyle{\subset}}}

\renewcommand{\sl}{\mathfrak{sl}}

\newcommand{\refequal}[1]{\xy {\ar@{=}^{#1}
(-1,0)*{};(1,0)*{}};
\endxy}

\newcommand{\Hom}{\operatorname{Hom}}

\newcommand{\End}{\operatorname{End}}

\newcommand{\Ext}{\operatorname{Ext}}

\newcommand{\Ind}{\operatorname{Ind}}
\newcommand{\op}{\operatorname{op}}
\newcommand{\id}{\operatorname{id}}

\newcommand{\Id}{\operatorname{Id}}
\newcommand{\Tot}{\operatorname{Tot}}

\newcommand{\inv}{^{-1}}
\newcommand{\SYT}{\operatorname{SYT}}

\newcommand{\Ch}{\operatorname{Ch}}
\newcommand{\Cone}{\operatorname{Cone}}

\newcommand{\Bim}{{\rm Bim }}

\newcommand{\SBim}{\SM\Bim}

\newcommand{\BS}{BS}

\newcommand{\FT}{\operatorname{FT}}
\newcommand{\HT}{\operatorname{HT}}
\newcommand{\fT}{\operatorname{ft}} 
\newcommand{\hT}{\operatorname{ht}}
\newcommand{\Br}{\operatorname{Br}}

\renewcommand{\min}{\operatorname{min}}

\DeclareMathOperator{\col}{col}
\DeclareMathOperator{\row}{row}
\DeclareMathOperator{\Schu}{Schu}

\newcommand{\yoo}{\:  \Yboxdim{5pt}\Yvcentermath1 \yng(1,1)}
\newcommand{\yt}{\:  \Yboxdim{5pt}\Yvcentermath1 \yng(2)}
\newcommand{\yooo}{\:  \Yboxdim{5pt}\Yvcentermath1 \yng(1,1,1)}
\newcommand{\yto}{\:  \Yboxdim{5pt} \Yvcentermath1 \yng(2,1)}
\newcommand{\yh}{\:  \Yboxdim{5pt} \Yvcentermath1 \yng(3)}
\newcommand{\yho}{\:  \Yboxdim{5pt} \Yvcentermath1 \yng(3,1)}
\newcommand{\ytoo}{\:  \Yboxdim{5pt} \Yvcentermath1 \yng(2,1,1)}
\newcommand{\sytoth}{\: \Yvcentermath1 \young(123)}
\newcommand{\sytotch}{\: \Yvcentermath1 \young(12,3)}
\newcommand{\sytohct}{\: \Yvcentermath1 \young(13,2)}
\newcommand{\sytoctch}{\: \Yvcentermath1 \young(1,2,3)}
\newcommand{\sytot}{\: \Yvcentermath1 \young(12)}
\newcommand{\sytoct}{\: \Yvcentermath1 \young(1,2)}

\colorlet{green}{black!30!green}

\definecolor{CQG}{RGB}{0,153,76}
\definecolor{FS}{RGB}{0,76,153} 

\newcommand{\revise}[1]{{#1}}
\newcommand{\revadd}[1]{{#1}}
\newcommand{\revcomment}[1]{}

\newcommand{\tw}{\operatorname{tw}}

\begin{document}

\begin{abstract} We conjecture that the complex of Soergel bimodules associated with the full twist braid is categorically diagonalizable, for any finite Coxeter group. This utilizes the theory of categorical diagonalization introduced in an earlier paper \cite{ElHog17a}. We prove our conjecture in type $A$, and as a result we obtain a categorification of the Young idempotents.
\end{abstract}

\title{Categorical diagonalization of full twists}

\author{Ben Elias} \address{University of Oregon, Eugene}

\author{Matthew Hogancamp} \address{Northeastern University, Boston}

\maketitle

\setcounter{tocdepth}{1}
\tableofcontents

\section{Introduction}
\label{sec:intro}
{}\blfootnote{2000 Mathematics Subject Classification. Primary 20F55, 18E30, 17B10.}

In this paper we categorify some of the most important objects in the representation theory of symmetric groups: the Young idempotents. The classical Young idempotents are a family of
idempotent elements $p_T\in \Q[S_n]$, indexed by standard Young tableaux $T$ with $n$ boxes. These idempotents feature prominently in the representation theory of $S_n$, as they
decompose the regular representation into irreducible representations. In addition, these idempotents have deformations (also denoted $p_T$) in the Hecke algebra which are essential in
the construction of the colored HOMFLY-PT polynomial for links.

We categorify $p_T$ using categorical linear algebra. Namely, we apply the theory of categorical diagonalization, developed in \cite{ElHog17a}, to the full twist Rouquier complexes in
the homotopy category of Soergel bimodules.

Let us now motivate and explain these results, and their conjectural generalizations to other Coxeter groups.

\subsection{Representation theory of symmetric groups}

Fix a natural number $n$ \revise{and let $\PC(n)$ denote the set of partitions of $n$. We write $\l\vdash n$ if $\l\in \PC(n)$.  For $\l\vdash n$,} let $\SYT(\l)$ denote the set of standard Young tableaux of shape $\l$. The irreducible representations $\SM^\l$ of the symmetric
group $S_n$ over $\Q$ are indexed by $\PC(n)$, and the dimension of $\SM^\l$ is equal to the size of $\SYT(\l)$. These irreducibles $\SM^\l$ are often called \emph{Specht modules}, after
Specht \cite{Specht35} (see also Peel \cite{Peel75}) who constructed integral forms of these representations (i.e. modules over $\Z[S_n]$ which become irreducible after base change to
$\Q$).

Young constructed combinatorially a family of operators $k_T \in \Z[S_n]$ indexed by standard Young tableaux $T$ with $n$ boxes. These operators are quasi-idempotent, in that $k_T^2 =
\g_T k_T$ for some scalar $\g_T \in \Z$. Dividing by $\g_T$, one obtains idempotents $p_T \in \Q[S_n]$ which are called \emph{Young idempotents} or \emph{Young symmetrizers}, and satisfy
$\Q[S_n] p_T \cong \SM^\l$ when $T \in \SYT(\l)$. These form a complete family of primitive orthogonal idempotents in $\Q[S_n]$. See \cite{FultonTab} for more details.

The presence of Young idempotents indicates some special feature of the symmetric groups. Using the theory of semisimple algebras, one has canonical central idempotents $p_\l \in \Q[S_n]$ which project to the isotypic component of $\SM^\l$ in $\Q[S_n]$.  That is, the image of $p_\l$ is $\#\SYT(\l)$ copies of $\SM^\l$; said another way,
\begin{equation}
p_\l = \sum_{T \in \SYT(\l)} p_T.
\end{equation}
Projection to a single irreducible component within this isotypic component is not canonical, and in the context of semi-simple algebras it is atypical to have a preferred choice of such projections (e.g. $p_T$).  However, symmetric groups are peculiar in that they fit in a tower
\begin{equation} S_1 \subset S_2 \subset \cdots \subset S_n \end{equation}
of groups, from which one obtains inclusions of group algebras.  A tableau $T \in \SYT(\l)$ corresponds to a tower $\l^1 \subset \cdots \subset \l^n = \l$ of partitions, with $\l^k$ a partition of $k$ for each $k < n$. Then
\begin{equation} p_T = \prod_{k=1}^n p_{\l^k}. \end{equation}

Although the ideas are older (see e.g. \cite{Jucys, Murphy,DipperJamesIdemp}), an influential modern perspective on these projections can be found in the work of Okounkov-Vershik \cite{OV96}.
The centralizer of $\Q[S_k]$ inside $\Q[S_{k+1}]$ is generated by a single element, the \emph{Young-Jucys-Murphy (YJM) element}\footnote{\revise{Explicitly, $j_{k+1}$ is a sum of transpositions}
\begin{equation} j_{k+1} = \sum_{i=1}^k (i \; (k+1)).
\end{equation}} $j_{k+1}$.
It follows that, when inducing an irreducible
representation of $S_k$ to $S_{k+1}$, projection to its isotypic components agrees with projection to the eigenspaces of $j_{k+1}$ (and moreover the isotypic components are irreducible).
By simultaneously diagonalizing the commuting family $\{j_1, j_2, \ldots, j_n\}$, one can construct the projectors $p_T$. The Okounkov-Vershik paper gives a particularly nice proof that the simultaneous eigenvalues of the YJM operators are in bijection with standard tableaux. The eigenvalue of $j_k$ on $p_T$ is the \emph{content number} (i.e. column
number minus row number) of the $k$-th box of $T$.

Let us quickly restate the construction of $p_T$ using the previous paragraph. Suppose one has constructed the idempotent $p_U$ for a standard tableau $U$ with $n-1$ boxes. Then letting $p_U \sqcup 1 \in \Q[S_n]$ denote its image under the inclusion $\Q[S_{n-1}] \to \Q[S_n]$, we have
\begin{equation} p_U \sqcup 1 = \sum_{U \subset T} p_T \end{equation}
where the sum is over tableaux obtained from $U$ by adding a box. This identity is a concrete realization of the ``branching rule," the decomposition of $\Q[S_n](p_U \sqcup 1)\cong \Ind_{S_{n-1}}^{S_n}(\SM^{\mu})$ into isotypic components, where $\mu$ is the shape of $U$. The individual idempotents $p_T$ are projections to the eigenspaces of $j_{n}$ on the image of $p_U \sqcup 1$.


These ideas reduce the representation theory of symmetric groups to \revcomment{removed (relatively fancy)} linear algebra; the projectors $p_\l$ and $p_T$ can be constructed by diagonalizing
certain operators, whose set of joint eigenvalues are well-understood.

\revcomment{next paragraph altered significantly} 

In this approach one constructs the central idempotent $p_\l$ by first constructing non-central idempotents $p_T$. One might hope instead to
construct $p_\l$ directly by diagonalizing an element of the center of $\Q[S_n]$. The most accessible element of the center (the easiest to categorify, in some sense) is the sum
$e_1 = \sum_{k=1}^n j_k$. It acts on $\SM^\l$ with eigenvalue given by the \emph{(total) content number} of $\l$, the sum of the contents of each box, which we denote in this paper
by $\xbb(\l)$. Unfortunately, there are distinct \revise{partitions} $\l \ne \mu$ with $\xbb(\l) = \xbb(\mu)$, meaning that $e_1$ can not tell these irreducible representations
apart.

\subsection{Lifting to the Hecke algebra}

The Hecke algebra $\Hecke_n$ (in type $A$) is a $\Z[v,v\inv]$-deformation of $\Z[S_n]$, for a formal parameter $v$. By setting $v=1$ there is an isomorphism $\Hecke_n / (v-1) \cong
\Z[S_n]$. There are quasi-idempotents $k_T \in \Hecke_n$ and idempotents $p_T \in \Hecke_n \ot_{Z[v,v\inv]} \Q(v)$ analogous to those constructed by Young, which we continue (abusively) to denote $k_T$ and $p_T$. See e.g. \cite[\S 5]{DipperJamesIdemp}.

The YJM elements $\{j_k\}$ of the symmetric group deform in a naive way to the \emph{additive Jucys-Murphy elements} of $\Hecke_n$.  Dipper and James used analogous techniques to construct Young idempotents in $\Hecke_n$ by projecting to joint eigenspaces of additive Jucys-Murphy elements, see \cite[p75]{DipperJamesIdemp}. \revise{For other constructions of idempotents and their relationship to additive Jucys-Murphy elements, see \cite{AistonMorton98}.}

The Hecke algebra also contains more subtle analogs of the $\{j_k\}$ called \emph{multiplicative Jucys-Murphy (JM) elements}, denoted $\{y_k\}$. \revise{These are the images in the Hecke algebra of elements in the braid group called \emph{Jucys-Murphy braids}, also denoted $\{y_k\}$. The braid $y_k$ wraps the $k$-th strand around the first $k-1$ strands, see Notation \ref{not:JMbraid} for details. We emphasize that these braids $\{y_k\}$ commute with each other.} We learned of \revise{the multiplicative JM elements} from work of Isaev-Molev-Oskin \cite{IMO}, which was inspired by ideas of Cherednik, see the references therein. \revise{They appear earlier in work of Morton \cite{Morton02}.} If one sets $v=1$ then $y_k$ is sent to $1 \in \Z[S_n]$. Instead, the relationship \revise{between $y_k$ and $j_k$} is better encoded as a derivative at $v=1$: more precisely (see \cite[p2]{IMO})
\begin{equation} j_k = \frac{y_k - 1}{v - v^{-1}}|_{v = 1}.\end{equation}
One can also imitate the proofs of Okounkov-Vershik \cite{OV96} for these multiplicative JM elements\footnote{This is a straightforward and very worthwhile exercise. We would love to know of good references for this material.}. The eigenvalue of $y_k$ on $p_T$ is $v^{2x}$, where $x$ is the content number of the $k$-th box in $T$.

There are many reasons to prefer multiplicative JM elements over additive ones. For one, the elements $y_k$ are invertible, being in the image of the map from the braid group to $\Hecke_n$. Others will become evident soon. We will only be considering multiplicative JM elements in this paper.

The operator $\fT_n = \prod_{k=1}^n y_k \in \Hecke_n$, called the \emph{full twist}, is the multiplicative analog of the element $e_1 = \sum_{k=1}^n j_k \in \Z[S_n]$, in the sense that its ``derivative at $v=1$'' is $e_1$. It is central in $\Hecke_n$ (being the image of a central braid), and acts diagonalizably with eigenvalue on $\SM^\l$ given by $v^{2\xbb(\l)}$. \revadd{Like $e_1$, it can not distinguish between all partitions.}

So, one can construct projections $p_T$ and $p_\l$ by simultaneously diagonalizing the commuting family $\{y_1, \ldots, y_n\}$. Alternatively, one can simultaneously diagonalize the
commuting family of \emph{(partial) full twists} $\{\fT_1, \ldots, \fT_n\}$, where $\fT_i = \prod_{k=1}^i y_k$.

\subsection{The categorification}

The category $\SBim_n$ of Soergel bimodules (for $S_n$) is a full subcategory of the category of graded bimodules for the ring $R = \RM[x_1, \ldots, x_n]$. It is closed under direct sums
and direct summands, grading shifts, and tensor products, so its split Grothendieck group $[\SBim_n]$ is naturally a free $\Z[v,v\inv]$-algebra with basis $\{[B_w]\}$ indexed by the
symbols of the indecomposable Soergel bimodules (up to isomorphism and grading shift). Soergel proved that the indecomposables are indexed by $w \in S_n$, and that there is a
natural isomorphism $\Hecke_n \to [\SBim_n]$. See \S\ref{subsec:SoergelCat} for additional background on this topic.

To lift elements of the braid group (or their images in $\Hecke_n$) to this categorical world, one should work with complexes of Soergel bimodules. The homotopy category (i.e.
complexes modulo nulhomotopic chain maps) of Soergel bimodules is often called the \emph{Hecke category}, and is a triangulated categorification of $\Hecke_n$. Rouquier constructed
a bounded complex of Soergel bimodules for each braid. Thus, the bounded homotopy category $\KC^b(\SBim_n)$ contains \revise{complexes} $Y_k$ associated to JM operators, and
$\FT_k$ associated to the partial full twists.

\begin{remark} It is not known if there is a reasonable categorification of the additive Jucys-Murphy elements. \end{remark}

\subsection{Categorical linear algebra}

In previous work \cite{ElHog17a}, the authors introduced a notion of categorical diagonalization, the diagonalization of functors. Let us give a brief synopsis.

%

\revcomment{this section is shortened and significantly revised, including the statement of the theorem.}

Let $\AC$ be an additive monoidal category, and fix a complex $F\in \KC^b(\AC)$.  A  \emph{scalar object} $\Sigma \in \KC^b(\AC)$ is any finite direct sum of shifts of the monoidal identity $\one$. Here, shifts can be homological, or may also involve some internal grading with which $\AC$ is equipped.  See \S
\ref{subsec:celltriprelude} and \cite{ElHog17a} for more details. Let $\VC$ be an additive category on which $\AC$ acts.

A \emph{weak eigenobject} is nonzero object $M \in \KC^b(\VC)$ such that $F\otimes M \simeq \Sigma \ot M$. One should think that $F$ categorifies a linear operator $f$, that $\Sigma$ categorifies a scalar $s$, and that a weak eigenobject categorifies an $s$-eigenvector, as $fm = s m$. However, weak eigenobjects do not form a sufficiently nice subcategory, suitable to categorify an eigenspace. Instead, a better notion of an eigenobject
would give more control over the isomorphism $F\otimes M \simeq \Sigma \ot M$.

The key notion in \cite{ElHog17a} is that of an \emph{eigenmap} of $F$. Given a scalar object $\Sigma$ and a chain map $\a\colon \Sigma\rightarrow F$, we define the
\emph{$\a$-eigencategory} to be the the annihilator of $\Cone(\a)$ in $\KC^b(\VC)$. If this annihilator is nonzero, then we say that $\a$ is an \emph{eigenmap} for the action of
$F$ on $\KC^b(\VC)$. A nonzero object of the $\a$-eigencategory is an \emph{$\a$-eigenobject}. One should think that $\Cone(\a)$ categorifies the linear operator $f - s$, which
kills any $s$-eigenvector. Eigenmaps are a new categorical structure which have no analogue in linear algebra: they represent the ``relationship'' between the operator $f$ and its
eigenvalues.

Let $f\in A$ be an element of an algebra, and suppose there are scalars $\sigma_i$ such that $\prod_i (f-\sigma_i)=0$ and $\sigma_i-\sigma_j$ is invertible for $i\neq j$.  Then one can construct a complete collection of mutually orthogonal idempotents $p_i\in A$ which diagonalize the action of $f$ on any $A$-module, meaning that $p_i$ projects to the $\sigma_i$-eigenspace. The following \emph{Diagonalization Theorem} is a categorical analogue of this classic result. Stated loosely, it asserts that if $F$ is invertible and has ``enough'' eigenmaps, then it is categorically diagonalizable.

\begin{theorem} \label{thm:DiagThmIntro} (See \cite[Theorem 1.22]{ElHog17a} for a precise version) Suppose that $F\in\KC^b(\AC)$ is invertible.  Let $I$ be a finite set and suppose there is an $I$-indexed family of invertible scalar objects $\Sigma_\l$ with \textbf{distinct homological shifts}. Suppose there are chain maps $\a_\l\colon \Sigma_\l\rightarrow F$ such that
\begin{equation}
\label{eq:introBigOT}
\bigotimes_{\l\in I} \Cone\left(\Sigma_\l \buildrel \a_\l\over \longrightarrow F \right) \simeq 0.
\end{equation}
Then under additional technical hypotheses\footnote{\revadd{These hypotheses imply that the cones $\Cone(\a_\l)$ tensor commute with each other.}}, for each $\l\in I$ there is an explicit construction of an idempotent complex $\PB_\l$ inside the
bounded-above homotopy category $\KC^-(\AC)$, which satisfies $\Cone(\a_\l) \ot \PB_\l \simeq 0$. Hence the image of $\PB_\l$ is contained in the $\a_\l$-eigencategory. Further, the complexes $\PB_\l$ form an idempotent decomposition of identity in the following sense:
\begin{enumerate}
\item $\PB_\l\otimes \PB_\mu\simeq 0$ for $\l\neq \mu$.
\item The monoidal identity $\one\in \KC^b(\AC)$ is homotopy equivalent to a filtered complex whose subquotients are the $\PB_\l$.
\end{enumerate}
\end{theorem}

If $\VC$ is an additive category on which $\AC$ acts, then the expression for $\one\in \KC^b(\AC)$ in terms of the idempotents $\PB_\l$ induces a filtration\footnote{precisely, a semi-orthogonal decomposition} of $\KC^b(\VC)$ into its eigencategories for the action of $F$.

It is the existence of the projectors $\PB_\l$ which most
interests us, so when we can construct them we say that we have \emph{diagonalized} the functor $F$. One should think of the passage from $\KC^b(\AC)$ to $\KC^-(\AC)$ as analogous base change from $\Z[v,v\inv]$ \revise{to Laurent power series $\Z(( v^{-1} ))$, as infinite complexes allow for infinite sums.} \revadd{In this context, the requirement that $\Sigma_{\l}$ have distinct homological shifts is analogous to the classical requirement that $\sigma_i - \sigma_j$ is invertible for $i \ne j$.}


\begin{remark}
We are careful above to avoid saying explicitly that Theorem \ref{thm:DiagThmIntro} \emph{categorifies} the well known statement in linear algebra. This is because the Grothendieck group of $\KC^-(\AC)$ is zero by the Eilenberg swindle (see \cite[Remark 4.9]{ElHog17a} for further discussion). \revadd{To make an interesting statement about the Grothendieck group}, one must show that all the complexes involved
live in a suitable full subcategory of $\KC^-(\AC)$ whose Grothendieck group is a nonzero extension of $[\AC]$, as was done e.g.~in \cite{CK12a,AchStr}. Choosing this subcategory is context-dependent.  In \S \ref{subsec:groth} we define a category $\KC^{\angle}(\SBim_n)$ which contains all the complexes considered here and has a well-behaved Grothendieck group.  This justifies statements like ``$\PB_\l$ categorifies $p_\l$'' below.
\end{remark}

%

\begin{remark} In fact, our Diagonalization Theorem also proves the existence of \emph{bounded} quasi-idempotent complexes $\KB_\l$, such that $\PB_\l$ is built from of infinitely many copies of $\KB_\l$. The $\KB_\l$ enjoy a number of desirable properties, and there is an easily expressed relationship between $\KB_\l$ and $\PB_\l$ which mirrors the relationship (Koszul duality) between exterior and polynomial algebras, see \cite[\S 7.3]{ElHog17a}. \end{remark}

\subsection{The main theorem}

Now let us apply the Diagonalization Theorem \ref{thm:DiagThmIntro} to our current situation. Let $\AC = \SBim_n$, and $F$ be the full twist Rouquier complex $\FT_n\in \KC^b(\SBim_n)$. Our first major result is this. 

\begin{theorem} \label{thm:lambdaMapsIntro} (See Theorem \ref{thm:lambdaMaps}) The complex $\FT_n$ \revise{satisfies a relation of the form \eqref{eq:introBigOT},} with one eigenmap $\a_\l \co \Sigma_\l \to \FT_n$ for each partition of $n$.  The scalar object $\Sigma_\l$ is
$\one(2\xbb(\l))[-2\cbb(\l)]$, where $\xbb(\l)$ is the total content number as above, and $\cbb(\l)$ is the total column number, obtained by adding together the column number of each box
in $\l$, see Definition \ref{defn:xcr}. \end{theorem}

Note that the grading shift $2\xbb(\l)$ corresponded to the eigenvalue $v^{2\xbb(\l)}$ of $\fT_n$ on $\SM^\l$. Meanwhile, the homological shift $-2\cbb(\l)$ is invisible in the
Grothendieck group, but plays a very significant role in this paper.

\begin{remark} \revcomment{Revised this remark.} As noted above, historical approaches to the representation theory of the Hecke algebra have diagonalized the (multiplicative) JM elements $y_k$. However, the corresponding Rouquier complexes $Y_k$ are \textbf{not} categorically diagonalizable. They do not have ``enough'' eigenmaps, i.e. they do not satisfy a relation of the form \eqref{eq:introBigOT}. Henceforth, we abandon JM elements and focus on full twists. \end{remark}

Unfortunately, Theorem \ref{thm:lambdaMapsIntro} does not allow one to apply the Diagonalization Theorem. This is because the scalar objects $\Sigma_\l$ do not have distinct homological shifts.
For any $n \ge 6$ there are distinct partitions $\l \ne \mu$ with $\cbb(\l) = \cbb(\mu)$, and even $\xbb(\l) = \xbb(\mu)$. Thus we again meet the problem that $\fT_n$ can not
distinguish between distinct irreducibles, and again we solve this problem by studying the entire family of full twists. In both cases, the saving grace is the observation that, once one
fixes a partition $\l$ with $n-1$ boxes, and restricts one's attention to partitions of $n$ which contain $\l$, then these partitions have distinct values of $\xbb$ and $\cbb$. Ultimately,
this will allow us to use our Relative Diagonalization Theorem \cite[Theorem 8.2]{ElHog17a} to diagonalize $\FT_n$ if we have already diagonalized $\FT_{n-1}$.

Thus our main theorem, see Theorems \ref{thm:typeAdiag} and \ref{thm:PTprops}, states that one can simultaneously diagonalize the full twists $\{\FT_1, \ldots,
\FT_n\}$, producing complexes $\PB_\l$ for each partition of $n$, and $\PB_T$ for each standard tableaux with $n$ boxes, which categorify $p_\l$ and $p_T$ respectively.

\begin{theorem}\label{thm:introFTdiag}
Let $\FT_n\in \KC^b(\SBim_n)$ be the full twist Rouquier complex.  There exists a diagonalization $\{(\PB_\l, \a_\l)\}$ of $\FT_n$, indexed by partitions $\l$ of $n$, in the sense of \cite[Definition 6.16]{ElHog17a}.  The idempotents $\PB_\l$ are central in the homotopy category $\KC^-(\SBim_n)$, and they categorify the central idempotents in $\HB_n$.  If $\l^{(k)}$ is a partition of $k$ for $1 \le k \le n$, then
\[
\PB_{(\l^{(1)},\ldots,\l^{(n)})}:=\bigotimes_{k=1}^n \PB_{\l^{(k)}}
\]
is contractible unless $T=(\l^{(1)},\ldots,\l^{(n)})$ is a standard Young tableaux (that is, $\l^{(k)}\subset \l^{(k+1)}$ for all $1\leq k\leq n-1$).  The idempotents $\PB_T$ describe a simultaneous diagonalization of the complexes $\FT_k$, $1\leq k\leq n$; they categorify the primitive Young idempotents $p_T$ in $\Hecke_n$.
\end{theorem}

Examples of these projectors $\PB_T$ and related finite complexes $\KB_T$ can be found in \S \ref{sec:examples}. \revadd{The finite quasi-idempotents $\KB_T$ categorify Young's quasi-idempotents $k_T$.}

\begin{remark} \revcomment{new remark on centrality} In the theorem above, an object $F \in \KC(\SBim_n)$ is called \emph{central} if there exists a functorial isomorphism $X \ot F \to F \ot X$ for all objects $X$ in $\SBim_n$. As a consequence, $F$ will categorify a central element of the Hecke algebra. One could desire the stronger statement that $X \ot F \to F \ot X$ for all objects $X \in \KC^b(\SBim_n)$, and this is not even proven for $F = \FT_n$ in the literature (though it is folklore). Yet more is desired: that $\FT_n$ and $\PB_{\l}$ can be equipped with the structure of objects in the Drinfeld center of $\KC^b(\SBim_n)$, and that the eigenmaps $\a_\l$ are morphisms in the Drinfeld center. We hope to address these desires in future work. \end{remark}

\subsection{Representations and cells}

\revcomment{This section rewritten. Note: replaced the word plethysm, this could be too confusing to a portion of the audience.}

In order to explain how our theorem is proven, we must approach the representation theory of the Hecke algebra from an entirely different perspective, that of Kazhdan-Lusztig theory and cell theory. In the next section, we discuss the connection between cell theory and diagonalization.

Soergel \cite{Soer90} proved that under the isomorphism $\Hecke_n \to [\SBim_n]$, the indecomposable Soergel bimodules $\{B_w\}$ are sent to the so-called \emph{Kazhdan-Lusztig (KL)
basis} $\{b_w\}$, defined in \cite{KazLus79} by Kazhdan and Lusztig. This is a difficult theorem, and relies in an essential way on the fact that the base ring $R = \RM[x_1,
\ldots, x_n]$ is over a field of characteristic zero. In finite characteristic, the indecomposable Soergel bimodules are sent to a different basis of $\Hecke_n$, and most of the
arguments in this paper will fail.

A consequence is that the decomposition of tensor products into indecomposable objects in $\SBim_n$ is determined by multiplication of KL basis
elements. That is, if \begin{equation} \label{eq:structurecoeffsintro} b_w b_x = \sum c^y_{w,x} b_y \end{equation} for coefficients $c^y_{w,x} \in \Z[v,v\inv]$, then \begin{equation} \label{eq:structurecoeffscategorified} B_w
\ot B_x \cong \bigoplus B_y^{\oplus c^y_{w,x}},\end{equation} and in particular the coefficients $c^y_{w,x}$ are non-negative. Here, $B^{\oplus v^k}$ is shorthand for the grading shift
$B(k)$.

\begin{remark} In contrast, tensor product decomposition in $\KC^b(\SBim_n)$ is not determined in the Hecke algebra! The class of a complex in the Grothendieck group only remembers its Euler characteristic (and certainly forgets the all-important differential).  Much of this paper consists of lifting formulas in the Hecke algebra to $\KC^b(\SBim_n)$ in a particularly nice way. \end{remark}

We now describe cell theory in a fashion which is tailored to the symmetric group. While cell theory generalizes to all Coxeter groups, type $A$ has many special features which may mislead the novice reader. We refer the reader to \S\ref{sec:cells} for an exposition of cell theory in general, and to \S\ref{sec:typeAcells} for more details on type $A$. Let $\Hecke = \Hecke_n$.

To a partition $\l\vdash n$ and two tableaux $S, T \in \SYT(\l)$, the Robinson-Schensted correspondence associates a permutation $w = w(S,T) \in S_n$. We say
that $w$ is in (two-sided) cell $\l$. A crucial feature observed already in \cite{KazLus79} is that the Robinson-Schensted correspondence is compatible with the structure of ideals
in the Hecke algebra. Let $\Hecke_{\le \l}$ denote the $\Z[v,v\inv]$-span of $\{b_w\}$ for all $w$ in cells $\mu$ with $\mu \le \l$ in the dominance order on partitions. Then
$\Hecke_{\le \l}$ is a two-sided ideal, as is $\Hecke_{< \l}$ which is defined similarly. The subquotient $\Hecke_{\le \l}/\Hecke_{< \l}$, viewed as a left module for $\Hecke$, splits further
into \emph{cell modules} indexed by $T \in \SYT(\l)$. The cell module for $T$ is spanned by $\{b_{w(S,T)}\}$ as $S$ varies, and its isomorphism type is independent of $T$; indeed,
each cell module for $T \in \SYT(\l)$ is isomorphic to the Specht module $\SM^\l$. This gives an alternate construction of the irreducible representations of $\Hecke$.

Now set $\SBim := \SBim_n$. Let $w$ be in cell $\l$, and $M \in \SBim$ be arbitrary. By \eqref{eq:structurecoeffscategorified} and the fact that $\Hecke_{\le \l}$ is an ideal,
every direct summand of $M \ot B_w$ or $B_w \ot M$ has the form $B_x(k)$ where $x$ is in cells $\le \l$. Define the full \revise{additive} subcategory $\SBim_{\le \l}$ with \revise{indecomposable}
objects $\{B_w(k)\}$ for $k \in \Z$ and $w$ in cells $\le \l$. Then $\SBim_{\le \l}$ is closed under multiplication by arbitrary objects of $\SBim$, making it a \emph{thick
(two-sided) tensor ideal}. The subquotients $\SBim_{\le \l}/\SBim_{< \l}$ split into left module categories of $\SBim$, one for each $T \in \SYT(\l)$. In this way, one categorifies
the Specht modules.

\subsection{Cell filtrations and diagonalization}

Let us examine how the cell theory construction of $\SM^\l$ interfaces with the eigentheory of the full twist. There is an element $\hT_n$ in the Hecke algebra, known as the \emph{half twist} and also the image of a braid, such that $\hT_n^2 = \fT_n$. By a result of Graham (see Mathas \cite{Mathas96}), for $w \in S_n$ in cell $\l$ one has
\begin{equation} \label{eq:htactionintro} \hT_n b_w \equiv (-1)^{\cbb(\l)} v^{\xbb(\l)} b_{\Schu_L(w)} + \HB_{<\l}. \end{equation} 
Here $\Schu_L$, the left Sch\"utzenberger dual, is an involution on $S_n$ that preserves each left cell. Consequently
\begin{equation} \label{eq:ftactionintro} \fT_n b_w \equiv (-1)^{2 \cbb(\l)} v^{2 \xbb(\l)} b_{w} + \HB_{<\l}. \end{equation}
	
\revcomment{Substantial revisions from here.}

This illustrates a crucial point. It is difficult to explicitly find any eigenvectors for $\fT_n$. However, it is easy to find eigenvectors modulo lower terms in the cell
filtration. The same is true in the categorification: an eigenobject for $\FT_n$ is typically an interesting complex, while any indecomposable object $B_w$ is an ``eigenobject modulo lower
cells'' (as we prove). This motivates us to shift our goalposts: instead of looking for eigenmaps, we look for ``eigenmaps modulo lower cells,'' a concept new to this paper. Let us state the definition as it applies to \revadd{the full twist} in $\SBim$.  For greater generality, see Definition \ref{def:lequiv}.
 
\revcomment{The below definition was clumsy and clashed with Definition \ref{def:lequiv}.  We have decided to limit its scope to just $\FT$, referring to Definition \ref{def:lequiv} for the general definition.}
\begin{defn} \label{defn:introlequiv}
\revadd{Let $\Sigma_\l=\one(2\xbb(\l))[-2\cbb(\l)]$ be the scalar objects from Theorem \ref{thm:lambdaMapsIntro}.  We call a map $\a_\l \co \Sigma_\l \to \FT_n$ a \emph{$\l$-equivalence} if it the induced map $\Sigma_\l B_w\to \FT_n\ot B_w$ is a homotopy equivalence in the cellular subquotient $\KC^b(\SBim_{\le \l}/\SBim_{< \l})$.}
%
\end{defn}

\revadd{Implicit in the existence of a $\l$-equivalence is the statement that $\FT \ot B_w$ is homotopy equivalent to a} complex which is built from $\Sigma_\l B_w$ and various indecomposables in cells $< \l$. In this paper we will construct $\l$-equivalences $\Sigma_\l \to \FT_n$
with $\Sigma_\l = \one(2 \xbb(\l))[-2\cbb(\l)]$ as in Theorem \ref{thm:lambdaMapsIntro}, thus categorifying \eqref{eq:ftactionintro}.

In linear algebra, we know that an upper triangular matrix with distinct values on the diagonal is diagonalizable. A slightly stronger statement which allows for some repeated
eigenvalues is as follows. Suppose $f$ is an operator acting on a vector space $V$. Suppose $V$ is filtered by a poset $\PC$, and $f$ acts by a scalar on $V_{\le \l}/V_{< \l}$ for all $\l \in \PC$. Suppose further that the eigenvalue of $f$ on $V_{\le \l}/V_{< \l}$ and the eigenvalue on $V_{\le \mu}/V_{< \mu}$ are distinct whenever $\l$ and $\mu$ are comparable. Then $f$ is diagonalizable. This statement suffices to prove that $\fT_n$ is diagonalizable, thanks to \eqref{eq:ftactionintro} and the fact that $\xbb(\l) < \xbb(\mu)$ whenever $\l < \mu$.

One might hope for a categorical version of this statement similar to Theorem \ref{thm:DiagThmIntro}: if $F$ is equipped with a $\l$-equivalence for each $\l$, and the scalars $\Sigma_\l$ have distinct homological degrees for comparable $\l$, then $F$ is diagonalizable. We are not sure if such a statement is true or not. However, in Proposition \ref{prop:lequivmeansdiag} we do prove such a theorem in the presence of one major additional ingredient: (homological) cell triangularity.

\subsection{Homological triangularity}
\label{subsec:boundingintro}

\revcomment{Substantial revisions continue throughout this section.}

In this paper we make a new observation, which we believe is of independent interest. For $w$ in cell $\l$, not only is it true that $\FT_n \ot B_w$ is built from one copy of
$\Sigma_\l B_w$ and various indecomposables in cells $< \l$, but those indecomposables in cells $< \l$ appear in strictly lower homological degree! We say that $\FT_n$ is \emph{cell
triangular}, see Definition \ref{def:twistedunitri}. Here is a precise statement.

\begin{prop} \label{prop:sharpIntro} Let $\FT = \FT_n$, and let $\HT = \HT_n$ be the Rouquier complex for the half twist $\hT_n$. For $w$ in cell $\l$,
\begin{enumerate}
\item $\HT\otimes B_w$ is homotopy equivalent to a complex $X^0\rightarrow \cdots \rightarrow X^{\cbb(\l)}$ where $X^{\cbb(\l)}=B_{\Schu_L(w)}(\xbb(\l))$ and $X^k$ is in cells strictly less than $\l$ for $k<\cbb(\l)$.
\item $\FT\otimes B_w$ is homotopy equivalent to a complex $Y^0\rightarrow \cdots \rightarrow Y^{2\cbb(\l)}$ where $Y^{2\cbb(\l)}=B_{w}(2\xbb(\l))$ and $Y^k$ is in cells strictly less than $\l$ for $k<2\cbb(\l)$.
\end{enumerate} \end{prop}

The reader is welcome to peruse \S\ref{subsec:HT3} now, for examples of Proposition \ref{prop:sharpIntro} when $n=3$.

Cell triangularity is invisible in the Grothendieck group; we think of Proposition \ref{prop:sharpIntro} as saying that \eqref{eq:htactionintro} and \eqref{eq:ftactionintro} are
categorified in the best possible way. Proposition \ref{prop:sharpIntro} does not make it any easier to construct a chain map $\a_\l \co \Sigma_\l \to \FT_n$, but it does make it
easier to see how $\a_\l \ot B_w$ should behave for a $\l$-equivalence: this must be the inclusion of $B_w(2\xbb(\l))$ as the last term in the complex of Proposition \ref{prop:sharpIntro} (2).


In this paper, cell triangularity will be paired with a requirement of homological distinctness: that $\cbb(\l) < \cbb(\mu)$ whenever $\l < \mu$. When this happens, we say that $\FT_n$
is \emph{sharp}, see Definition \ref{def:lsharp}. In Proposition \ref{prop:lequivmeansdiag} we prove that if a complex $F$ is sharp and cell triangular and admits a $\l$-equivalence
for each $\l$, then \eqref{eq:introBigOT} holds. As a consequence, the $\l$-equivalences involved will actually be eigenmaps. We omit some technical hypotheses in this discussion. 

Our proof that $\FT_n$ is diagonalizable can now be divided into three major steps. The first is proving Proposition \ref{prop:sharpIntro}, from which one deduces that $\FT_n$ is sharp
and cell triangular. The second is constructing a family of $\l$-equivalences. The third is congealing these results into a form suitable for the Relative Diagonalization Theorem of
\cite{ElHog17a}, thereby simultaneously diagonalizing the family of full twists.

\subsection{Bounding the action of the full twist}
\label{subsec:boundingintro2}

\revcomment{We return to how thinsg used to be}

In \S\ref{subsec:twists} and \S\ref{subsec:sharpconj}, Proposition \ref{prop:sharpIntro} is reduced to a simpler statement, which is proven in Theorem \ref{thm:HTactionA}. The proof of
Theorem \ref{thm:HTactionA}, and the exposition to make it possible, comprises a significant portion of this paper. It requires a fairly thorough knowledge of Lusztig's asymptotic Hecke
algebra.

When one speaks of asymptotic data in the Hecke algebra, one refers to the lowest degree power of $v$ with a nonzero coefficient in some formula; categorically,
this is the minimal degree of a morphism between Soergel bimodules. For example, letting $c_{x,y}^z$ be as in \eqref{eq:structurecoeffsintro}, one can define\footnote{This is denoted
$\abb(z)$ by Lusztig.} $\rbb(z)$ for $z \in S_n$ to be the integer such that $v^{-\rbb(z)}$ is the minimal power of $v$ appearing in $c_{x,y}^z$ for all $x$ and $y$. It turns out that
$\rbb(z)$ only depends on the two-sided cell $\l$ of $z$, and is equal to the row number $\rbb(\l)$ of the corresponding partition.

\begin{remark} We have now defined three numerical statistics associated to a partition: the row number $\rbb(\l)$, the column number $\cbb(\l)$, and the content number $\xbb(\l)$. One
has $\xbb(\l) = \cbb(\l) - \rbb(\l)$, and $\cbb(\l) = \rbb(\l^t)$, where $\l^t$ is the transpose partition, also obtained from the two-sided cell $\l$ by multiplication by the longest
element $w_0 \in S_n$. \end{remark}

Lusztig \cite[Chapter 14]{LuszUnequal14} has a list of famous properties (P1-P15) of this statistic $\rbb(z)$ and its interaction with cells, which he proves for the symmetric group
using the fact that all Kazhdan-Lusztig polynomials are positive. We use most of these properties in the setup of the proof of Theorem \ref{thm:HTactionA}. Said another way, our results
in this paper put a homological spin on Lusztig's asymptotic Hecke algebra (c.f. Remark \ref{rmk:syzygy}).

We use these Grothendieck group considerations to reduce the proof of Proposition \ref{prop:sharpIntro} to the case when $w$ is the longest element of a parabolic subgroup $S_{k_1}
\times \ldots \times S_{k_r} \subset S_n$. We are able to prove this directly in \S\ref{subsec:bounding} and \S\ref{subsec:boundinglemma}, by a \revise{computation involving Gaussian
elimination of complexes. In the process, we compute idempotent decompositions for certain Soergel bimodules, which should be of independent interest. Our computation holds over the integers, making it useful for diagonalization in finite characteristic as well.}

\subsection{Constructing $\l$-equivalences}
\label{subsec:constructionIntro}

Two ingredients go into our construction of $\l$-equivalences. The first is a criterion for the existence of a $\l$-equivalence.

Recall that $R = \RM[x_1, \ldots, x_n]$ is the base ring of the construction of Soergel bimodules, and also is the endomorphism ring $\End(\one)$ of the monoidal identity. Using
asymptotic data once more, one can construct homogeneous polynomials $\tbarb{d} \in R$ (well-defined up to a scalar in $\RM_{>0}$) associated to involutions $d \in S_n$, by composing the
morphisms of minimal degree \[\one \to B_d \to \one.\] The span of these polynomials, as $d$ ranges over all involutions in a given cell $\l$, is proven to be an $S_n$-representation
inside $R$. In fact, we prove it is isomorphic to the Specht module $\SM^\l$! So long as there exists some involution $d$ in cell $\l$ and some map $\a$ which is not annihilated by
precomposition with $\tbarb{d}$, we are able to bootstrap the existence of some $\l$-equivalence (not necessarily $\a$ itself). This proof, accomplished in
\S\ref{subsec:lambdaEquivsufficient}, uses the irreducibility of the Specht module.

The second key ingredient in our construction of $\l$-equivalences comes from recent progress in computing the triply graded homology of torus links.

Given complexes $C,D\in \KC(\SBim_n)$, let $\Hom_{\KC(\SBim_n)}^{\Z\times \Z}(C,D)$ denote the bigraded space of homogeneous chain maps modulo homotopy:
\[
\Hom_{\KC(\SBim_n)}^{\Z\times \Z}(C,D) = \bigoplus_{i,j\in \Z}\Hom_{\KC(\SBim_n)}(C,D(i)[j]).
\]
\revcomment{starting here, several edits for clarity and length.}
In \cite{ElHog16a} the authors introduced a technique for the computation of $\Hom^{\Z\times \Z}(\one,C)$ when $C$ is a Rouquier complex of a particular class of braids. As an application, we compute $\Hom^{\Z\times \Z}(\one,\FT_n)$ as a bigraded vector space. More details are given in \S\ref{subsec:FTHHH}. In \S\ref{subsec:eigenmaptheorem} we justify that this
description of $\Hom^{\Z\times \Z}(\one,\FT_n)$ produces certain specific maps $\a_\l \co \one(2\xbb(\l))[-2\cbb(\l)] \to \FT_n$, which are candidates for $\l$-equivalences. 

In fact, $\Hom^{\Z\times \Z}(\one,\FT_n)$ is an $R$-bimodule on which the right and left actions agree. Our computation from \cite{ElHog16a} does not make the $R$-module structure clear, only
the action of $R$ in the associated graded for some filtration. This is because we describe $\Hom^{\Z\times \Z}(\one,\FT_n)$ using a (degenerate) spectral sequence. Despite this, we are able to examine the  action of $R$ on $\a_\l$ in the associated graded, and deduce that a particular polynomial $\tbarb{d} \in R$ will not kill $\a_\l$. This is sufficient to deduce the existence of some $\l$-equivalence, by the criterion mentioned above.

\begin{remark} In fact, using recent work of the second author and Gorsky \cite{GorHog17}, we can prove that $\a_\l$ itself is a $\l$-equivalence. We include this as an addendum, see \S\ref{subsec:symgrp}. \end{remark}

\subsection{Ensuring relative diagonalization}

Finally, we need to prove that $\FT_n$ can be diagonalized, relative to an inductively defined diagonalization of $\FT_{n-1}$. Fix a partition $\mu$ of size $n-1$, and let $\PB_\mu$ be the corresponding central idempotent complex in $\KC^-(\SBim_{n-1})$, which we can also view as a non-central idempotent complex in $\KC^-(\SBim_n)$. Let $\{\l^1, \ldots, \l^r\}$ denote the partitions of $n$ which are obtained from $\mu$ by adding a single box. In order to use our Relative Diagonalization Theorem \cite[Theorem 8.2]{ElHog17a}, we need to prove that
\begin{equation}\label{eq:relBigOTintro} \PB_\mu \ot \bigotimes \Cone(\a_{\l^i}) \simeq 0. \end{equation}

The corresponding statement in the Hecke algebra says that the eigenvalues of $\fT_n$ on the induction of $\SM^\mu$ to $\Hecke_n$ are those associated to $\l^i$. This is a key computation in Okounkov-Vershik \cite{OV96}. A categorical proof is slightly more difficult, as we must prove this result using cell theory instead in order to use the techniques (e.g.~$\l$-equivalences) that we have developed. To this end, we use several results of Meinolf Geck \cite{GeckInduction,GeckRelative} on induction of Kazhdan-Lusztig cells and on relative Kazhdan-Lusztig cells, which we recall in \S\ref{subsec:relativecells}.

Beyond this, the major technical difficulties involved are mostly in the Relative Diagonalization Theorem itself, and this justifies a great deal of the work done in \cite{ElHog17a}.  The resulting idempotents $\PB_T$ satisfy an inductive relationship which categorifies \cite[Equation (11)]{IMO}.

\subsection{Tightness} \label{subsec:tightness}

We have noted that there are distinct partitions $\l \ne \mu$ with the same eigenvalue of $\fT_n$, so that $\fT_n$ can not tell their irreducibles apart. Yet $\FT_n$ has distinct (and in fact, linearly independent) eigenmaps $\a_\l$ and $\a_\mu$, so it seems that $\FT_n$ should be able to tell their categorical cell modules apart, which would be an upgrade. However, this is a very subtle point which we do not prove here.

In our definition of diagonalization \cite[Definition 6.16]{ElHog17a}, the essential image of the projection $\PB_\l\otimes (-)$ must lie inside the eigencategory of $\a_\l$, but it
need not be the entire eigencategory! When the essential image of $\PB_\l\otimes (-)$ and the eigencategory of $\a_\l$ agree for all $\l$, the diagonalization is called \emph{tight}.
\revcomment{Edits in this paragraph for length.} When distinct eigenmaps $\a_\l \ne \a_\mu$ have the same scalar object $\Sigma_\l \cong \Sigma_\mu$, it is possible for there to be a nonzero object
$M$ which is an eigenobject for both $\a_\l$ and $\a_\mu$. This possibility is precluded when the diagonalization were tight. After all, if $M$ is preserved by both $\PB_\l$ and $\PB_\mu$, then $M \cong \PB_\l \PB_\mu M
\cong 0$ by orthogonality of projectors.


To illustrate the subtlety, consider the incomparable partitions $\l=(3,1,1,1)$ and $\mu=(2,2,2)$, which have the same values of $\rbb=6$, $\cbb=3$, and $\xbb=-3$. Both $\a_\l$ and
$\a_\mu$ are elements of $V = \Hom(\one(-6)[6], \FT_6)$, and $\{\a_\l,\a_\mu\}$ is a basis for this two-dimensional space. However, as observed in \cite[\S 2]{ElHog17a} or in the proof
of Lemma \ref{lem:criterion}, being a $\l$-equivalence is an open condition in this vector space $V$, the complement of a union $V_\l$ of hyperplanes. Hence almost every linear
combination of $\a_\l$ and $\a_\mu$ is both a $\l$-equivalence and a $\mu$-equivalence. If $\a_\mu$ happens to be a $\l$-equivalence (the generic situation) then there is a joint eigenobject for both $\a_\l$ and $\a_\mu$.  The question of distinguishing the $\a_\mu$ and $\a_\l$ eigencategories essentially boils down to whether it is possible to construct a $\l$-equivalence $\a_\l$ which is not a $\mu$-equivalence, and vice versa. It certainly would be quite special for $V_\mu$ to be contained in $V_\l$ for some $\l\neq \mu$, so generically one would expect the following to be true.

\begin{conjecture}\label{conj:tightness0}
The $\l$-equivalences $\a_\l$ can be chosen so that $\a_\l$ is not a $\mu$-equivalence for any $\mu\neq \l$.
\end{conjecture}

This is certainly a prerequisite for tightness.

\begin{conjecture}\label{conj:tightness2}
There exist $\l$-equivalences $\a_\l$ for which the diagonalization of $\FT_n$ from Theorem \ref{thm:typeAdiag} is tight.
\end{conjecture}

It is beyond the scope of this paper to tackle either of these conjectures. We expect Conjecture \ref{conj:tightness2} to be a straightforward consequence of Conjecture \ref{conj:tightness0}, but a proof of Conjecture \ref{conj:tightness0} will require a deeper understanding of the $R$-module structure on the bigraded hom space $\Homgg(\one,\FT_n)$. Our main tool for studying $\Hom^{\Z \times \Z}(\one,\FT_n)$ is the spectral sequence from \S\ref{subsec:FTHHH}, which only encodes the $R$-action on the associated graded, while our main tool for proving that a map is a $\l$-equivalence uses the true $R$-action (see \S\ref{subsec:constructionIntro}). One can use this to prove maps are $\l$-equivalences, because being not being annihilated by an element of $R$ in the associated graded implies not being annihilated by the same polynomial in the whole space. The converse of this is false, hence we do not show that our $\l$-equivalence is not also a $\mu$-equivalence for $\mu\neq \l$.

Finally we mention a closely related conjecture regarding the categorical analogue of the minimal polynomial for $\FT_n$.

\begin{conjecture}\label{conj:tightness1}
The $\l$-equivalences $\a_\l$ can be chosen such that
\begin{equation} \bigotimes_{\l\in \PC(n)} \Cone(\a_\l) \simeq 0, \end{equation}
with none of the factors being redundant.
\end{conjecture}
If the $\mu$-equivalence $\a_\mu$ were also a $\l$-equivalence for some $\l\neq \mu$, then the factor $\Cone(\a_\l)$ would be redundant.  Thus, Conjecture \ref{conj:tightness1} implies Conjecture \ref{conj:tightness0}.

\subsection{Other Coxeter groups}

The symmetric group is one family (type $A$) inside a broad class of groups generated by reflections, known as \emph{Coxeter groups}. Given a Coxeter group $W$ with a set of simple
reflections $S$, one has an associated Hecke algebra $\HB(W)$, a category of Soergel bimodules $\SBim(W)$, and Rouquier complexes associated with elements of the braid group of $W$. When
$W$ is finite, it has a longest element $w_0$, which lifts to an element $\hT$ inside the braid group called the \emph{half twist}. This in turn squares to a central element $\fT$ which
we call the \emph{full twist}. Mathas \cite{Mathas96} has proven an analog of \eqref{eq:htactionintro}. All the pieces are in place for a potential diagonalization of $\FT$, the Rouquier
complex of $\fT$.

Below, let $\rbb(\l)$ be defined as in \S\ref{subsec:boundingintro}; this number is commonly called Lusztig's $\abb$-function. Let $\cbb(\l) := \rbb(\l^t)$, where $\l^t$ is the two-sided cell $w_0 \l$, set $\xbb(\l) := \cbb(\l) - \rbb(\l)$.

\begin{conj} \label{conj:coxIntro} The Rouquier complex $\FT$ associated with the full twist $\fT$ is categorically diagonalizable, for any finite Coxeter group $W$. It has one eigenmap
$\a_\l \co \Sigma_\l \to \FT$ for each two-sided cell $\l$. The scalar object $\Sigma_\l$ is given by $\one(2\xbb(\l))[-2\cbb(\l)]$.\end{conj}

In a sequel to this paper, we will prove this conjecture for dihedral groups (which have only three two-sided cells) by direct computation.  

To facilitate this general conjecture, we give our exposition in terms of a general Coxeter group, until \S\ref{sec:typeAcells} when we focus on the case of the symmetric group. It should not make the exposition significantly harder to read, as various constructions (Lusztig's asymptotic Hecke algebra, Sch\"utzenberger duality) are tricky anyway, and not significantly less tricky in type $A$. The main simplification in type $A$ is that every involution is distinguished. The reader unfamiliar with Coxeter groups beyond the symmetric group can just assume $W = S_n$ throughout.

Let us discuss which arguments in this paper go through for general Coxeter groups, and which (enormous) gaps remain in the proof.

Our first step was proving Proposition \ref{prop:sharpIntro}, which said that \eqref{eq:htactionintro} was categorified in the best possible way. In \S\ref{subsec:twists} we develop some general results, which reduce this question to proving that $\HT \ot B_w$ is concentrated in homological degrees $\le \cbb(\l)$ for a single $w$ in each cell $\l$. This is the content of Conjecture \ref{conj:HTaction}. \revise{Our proof in type $A$ is the content of \S\ref{subsec:bounding}, and should not generalize in any naive way beyond type $A$.}

The second step was to construct a family of $\l$-equivalences; that they exist is Conjecture \ref{conj:eigenmap}. This in turn had two ingredients. The first was a sufficient condition for a $\l$-equivalence to exist using non-annihilation by polynomials in the Specht module. This analog of the Specht module does exist for arbitrary Coxeter groups, and it seems not to have been studied previously! We do not know what its properties might be, but it might be irreducible in general, which would enable us to reproduce the same sufficient condition.

The second ingredient was an explicit construction of particular chain maps via our computation of triply graded link homology, and in particular, our computation of $\Hom^{\Z\times \Z}(\one,\FT_n)$.  Any computation of $\Hom^{\Z\times \Z}(\one,\FT)$ for other finite Coxeter groups will likely require new ideas (or inhuman tenacity), though we will compute it for the dihedral group in the sequel.

The final step was to apply the Relative Diagonalization Theorem. This step was required because $\FT_n$ might have distinct cells with the same scalar object, and smaller full twists were required to tell these cells apart. Essentially, this step (as well as our computation of $\Hom^{\Z\times \Z}(\one,\FT_n)$ ) relies on the existence of a well-understood tower of Coxeter groups $S_1 \subset S_2 \subset \ldots \subset S_n$, with well-understood relative cell theory. We believe this is a topic which, already in the Grothendieck group, deserves far more study than it has received.

\subsection{Relation to other work}
\label{subsec:otherwork}

The idempotents $\PB_T$ and $\PB_\l$ constructed here recover several other idempotents constructed already in other contexts.  We keep the details to a minimum, but hope this will be a useful reference for experts. The categorified Jones-Wenzl idempotents \cite{CK12a,Roz10a,Rose14,Hog15} are obtained from $\PB_T$ when $T$ is the one-row tableau.  The idempotent $\PB_n$ from \cite{Hog15} is in fact isomorphic to $\PB_T$, while the others are obtained from $\PB_T$ by applying a monoidal functor from $\SBim_n$ to an appropriate additive monoidal category (such as $\sl_N$ matrix factorizations or foams).  A category $\OC$ version of the categorified Jones-Wenzl idempotent \cite{FSS12} is conjecturally related to these by a version of Koszul duality \cite{SS12}. 

In \cite{CH12} the second author and Ben Cooper categorified a complete collection of idempotents in the Temperley-Lieb algebra.  The central idempotents $P_{n,k}$ can be obtained from $\PB_{\l}$ for $\l$ a 2-row partition, while the primitive idempotents $P_\e$ can be obtained from $\PB_T$ with $T$ a 2-row standard tableau.  As further special cases, we obtain Rozansky's ``bottom projectors'' $P_{2n,0}$ \cite{Roz10b}.  Applying functors to categories of $\sl_N$ matrix factorizations, our work here generates $\sl_N$ versions of these, indexed by $N$-row partitions / tableau.  These should be closely related to Cautis' idempotents \cite{CautisClasp}.  This connection deserves further exploration and will certainly involve an adaptation of the techniques here to the setting of singular Soergel bimodules.

If $T$ is the one-column partition, then $\PB_T$ is isomorphic to $P_{1^n}$ constructed by the second author and Michael Abel \cite{AbHog17}.

Finally, let us point out the beautiful work \cite{GNR16}.  After sharing early details of this work with Gorsky, Negu\cb{t}, and Rasmussen, it was realized that categorical diagonalization provides a natural home for an emerging connection between knot homology and Hilbert schemes.  In particular, the endomorphism rings $\Endgg(\PB_T)$ are conjectured to be isomorphic to some explicit quotients of polynomial rings, which are shown to be isomorphic to the ring of functions on an open chart of the flag Hilbert scheme $\FHilb^n(\C^2)$ in \cite{GNR16}.  The most important ingredient in our story, the eigenmaps, do indeed yield an action of a polynomial ring on $\PB_T$ in which the generators have the predicted degrees.  We refer to \cite{GNR16} for details, and also to \cite{Hog15,AbHog17} for proofs of this conjecture for the one-row and one-column tableaux.  The polynomial action induced by eigenmaps follows from general principles (see for instance the proof of Lemma \ref{lem:grothPvK}).

\subsection{Structure of the paper}

\revise{The reader may have noticed that we prove our theorem by throwing the whole book at it. We use specialized results in Kazhdan-Lusztig theory, representation theory, and
homological algebra. Our positive spin is that we're not really throwing the book, we're instead reinterpreting the book homologically, taking properties of the Hecke algebra and
upgrading them to properties of the full twist complex.}

In \S\ref{sec:heckeandsbim} we give background information on Coxeter groups, the Hecke algebra, the Kazhdan-Lusztig basis, Soergel bimodules, and Rouquier complexes. We also describe
minimal complexes and Gaussian elimination, two key concepts in the homological algebra of additive categories. In addition to the basic definitions and familiar properties, we emphasize
several points which may not be familiar to the expert: \begin{enumerate} \item Rouquier's theorem about the canonicity of Rouquier complexes, and its implications for conjugation by
braids. \S\ref{subsec:RouqCanon}. \item How braids act by conjugation on polynomials. \S\ref{subsec:conjugate}. \item Minimal complexes of Rouquier complexes, and in particular, the half twist. \S\ref{subsec:minimalRouq}.
\end{enumerate}

In \S\ref{sec:cells} we give an introduction to the abstract theory of cells in monoidal categories. Then, starting in \S\ref{subsec:twists}, we study the theory of twists, which are invertible complexes sharing some of the key properties of $\HT$ and $\FT$. Sections \S\ref{subsec:twists} and \S\ref{subsec:celltri} are entirely new material.

In \S\ref{sec:cellsplus} we introduce the advanced and asymptotic cell theory of the Hecke algebra. We recall the conjectures (P1-P15) of Lusztig, and we also recall the aforementioned
result of Mathas on the action of the half twist and Sch\"utzenberger duality. Then we discuss some categorical lifts of these ideas. In \S\ref{subsec:sharpconj} we state our main
conjectures for arbitrary finite Coxeter groups: Conjecture \ref{conj:HTaction} stating that the half twist is sharp, implying that \eqref{eq:htactionintro} is categorified in the best
possible way, and Conjecture \ref{conj:eigenmap} stating that there exists a family of $\l$-equivalences. We prove some consequences of these conjectures, in particular that $\FT$
\revise{ satisfies \eqref{eq:introBigOT}.} In \S\ref{subsec:dots} we discuss the existence of certain dot maps in $\SBim$, which are maps of minimal degree from $\one$ to $B_d$ for a
distinguished involution $d$. \revadd{We suspect our conjectures imply that $B_d$ is a Frobenius algebra object.} In \S\ref{subsec:dots}, we prove the unit axiom.

In \S\ref{sec:typeAcells} we restrict our attention to type $A$, and discuss once again the various facts proven about cells and asymptotics, as they apply to the symmetric group. In \S\ref{sec:typeAtwists} we prove Conjecture \ref{conj:HTaction} in type $A$, by an explicit computation.

In \S\ref{sec:constructing} we prove Conjecture \ref{conj:eigenmap} in type $A$. After some reminders and elementary results in \S\ref{subsec:lambdaEquiv}, we use the dot maps to give an
alternate construction of the Specht module in \S\ref{subsec:specht} in type $A$, and of a ``new" representation of $W$ in other types. We use the irreducibility of the Specht module to
prove a criterion for the existence of a $\l$-equivalence in \S\ref{subsec:lambdaEquivsufficient}. In \S\ref{subsec:FTHHH} we recall our earlier results from \cite{ElHog16a} on the triply
graded homology of the full twist, and in \S\ref{subsec:eigenmaptheorem} we combine these results to prove the existence of $\l$-equivalences.

In \S\ref{sec:diag} we finally prove Theorem \ref{thm:introFTdiag}. After discussing some general consequences of diagonalization, we focus in \S\ref{subsec:subtableaux} and \S\ref{subsec:relativecells} on the interactions between $\HB(S_k)$ and $\HB(S_n)$ for $k < n$, recalling results of Geck and how they apply in type $A$. In \S\ref{subsec:implications} we use these results to prove \eqref{eq:relBigOTintro}, from which the proof of our main theorem follows relatively easily in \S\ref{subsec:proof}. 

\subsection{Acknowledgments}

We are indebted to Meinolf Geck for answering many questions by email, and for having answered so many more years ago in his papers! The first author would also like to thank Victor
Ostrik, Sasha Kleshchev, and Geordie Williamson, for a variety of useful conversations and pointers to the literature. The first author would also like to thank Benjamin Young, whose
computations ruled out some overly optimistic ideas, see Remark \ref{rmk:notametric}.

Both authors thank Eugene Gorsky, Andrei Negu\cb{t}, and Jacob Rasmussen for many enlightening conversations, and we thank Mikhail Khovanov for getting our collaboration started.  We would also like to thank Weiqiang Wang for pointing out the very nice proof of Lemma \ref{lem:ftIsCentral_decat}.

We would truly like to thank the anonymous referees for many helpful comments! We also thank Cailan Li for noticing an error in an earlier version of \eqref{thosearesomedifferentials}.

The first author was supported by NSF CAREER grant DMS-1553032, and by the Sloan Foundation. During revision, the first author was supported by DMS-2201387. The second author was supported by NSF grants DMS-1702274, DMS-1255334, and DMS-1664240.

\section{The Hecke algebra and Soergel bimodules}
\label{sec:heckeandsbim}





In this chapter, we provide the requisite background information on the Hecke algebra and Hecke category attached to a Coxeter system $(W,S)$. In the next chapter, we provide background information on cell theory. This chapter is largely review, though we call the expert reader's attention to \S\ref{subsec:conjugate} and to the examples in \S\ref{subsec:minimalRouq}.

\subsection{Notation and the basics}
\label{subsec:notation}

Throughout this chapter, fix a Coxeter group $W$ and a set of simple reflections $S\subset W$. In this paper we often restrict to \emph{type $A$} or \emph{type $A_{n-1}$}: this means that $W$ is the symmetric group $S_n$, and $S = \{s_i\}_{i=1}^{n-1}$ are the adjacent transpositions $s_i = (i\ i+1)$.

An \emph{expression} for an element $w \in W$ is a sequence $\un{w} = (s_{i_1}, s_{i_2}, \ldots, s_{i_d})$ with $s_{i_k} \in S$, such that $s_{i_1} s_{i_2} \cdots s_{i_d} = w$. When there is no ambiguity we assume that an underlined letter $\un{x}$ represents an expression for the element $x \in W$. A \emph{reduced expression} or \emph{rex} is an expression for $w$ of minimal length $d$, and this number $d$ is the \emph{length} of $w$, denoted $\ell(w)$.

If $W$ is finite, then there is a unique element of longest length, denoted $w_0$. In type $A_{n-1}$, the longest element is the permutation sending $1\leftrightarrow n$,
$2\leftrightarrow n-1$, etcetera. Its length is $\binom{n}{2}$, having a reduced expression \[w_0=s_1(s_2s_1)\cdots (s_{n-1}\cdots s_2s_1).\]

An \emph{involution} in $W$ is an element such that $w^2=1$. Longest elements are always involutions.

The set $W$ has a \emph{Bruhat order}, a partial order respecting the length function, where $x \le w$ if one can obtain a rex for $x$ as a subsequence of a rex for $w$. The identity is the unique minimal element in the Bruhat order, and $w_0$ is the unique maximal element.

A \emph{parabolic subgroup} is a subgroup $W_I$ generated by a subset $I \subset S$. Then $W_I$ is a Coxeter group with simple reflections $I$. Any rex in $S$ for an element $w \in W_I$
will only use simple reflections in $I$, and thus the length function and Bruhat order on $W_I$ are the restrictions of the corresponding notions on $W$. When $W_I$ is finite, its
longest element will be denoted $w_I$. In type $A_{n-1}$, all parabolic subgroups have the form $S_{k_1} \times S_{k_2} \times \cdots \times S_{k_r}$ with $\sum k_i = n$.

Typically the letter $s$ denotes an arbitrary simple reflection. When discussing $S_n$ with $n \le 4$ we let $s = s_1$, $t = s_2$, and $u = s_3$, by convention.

\subsection{The Hecke algebra}
\label{subsec:Hecke}

The \emph{Hecke algebra} $\HB = \HB(W)$ is a deformation of the group algebra of $W$ over the base ring $\Z[v,v\inv]$, where $v$ is a formal indeterminate. As a $\Z[v,v\inv]$-module, $\HB$ is free, with a \emph{standard basis} $\{H_w\}_{w \in W}$.  The multiplication is determined by
\begin{itemize}
\item $H_w H_v = H_{wv}$ if $\ell(wv)=\ell(w)+\ell(v)$,
\item $(H_s+v)(H_s-v\inv)=0$ for each simple reflection $s\in S$.
\end{itemize}
Here $H_1 = 1$ is the unit.

A map $f$ of $\Z[v,v\inv]$-modules is \emph{antilinear} if $f(vm) = v^{-1} f(m)$. The Hecke algebra is equipped with an antilinear automorphism called the \emph{bar involution}, uniquely
specified by $\overline{H_s}=H_s\inv$. It is an algebra homomorphism, so that $\overline{ab}=\overline{a}\overline{b}$. In terms of the standard basis, we have $\overline{H_w} =
H_{w\inv}\inv$. The following theorem is due to Kazhdan-Lusztig \cite{KazLus79}.

\begin{theorem}\label{thm:KLbasis}
There are unique elements $b_w\in \HB$ for each $w \in W$, such that $\overline{b_w}=b_w$ and
\[
b_w = H_w + \sum_{y< w} h_{y,w}(v) H_y
\]
where $h_{y,w}(v)\in v\Z[v]$. The elements $\{b_w\}_{w \in W}$ form a basis of $\HB$, which we call the \emph{Kazhdan-Lusztig basis} or \emph{KL basis}. The polynomials $h_{y,w}(v)$ are called \emph{Kazhdan-Lusztig polynomials} or \emph{KL polynomials}.
\end{theorem}

\begin{remark}
We are mostly following the notation of Soergel \cite{Soer97}, and recommend his exposition in \cite[\S 2]{Soer97} for the basics of the Hecke algebra. However, we write $b_w$ where Soergel writes $\un{H}_w$. To match the conventions of Kazhdan-Lusztig \cite{KazLus79}, set $H_w = v^{\ell(w)} T_w$, $q=v^{-2}$, $b_w = C_w'$, and $h_{y,w} = v^{\ell(w) - \ell(y)} P_{y,w}(v^{-2})$.
\end{remark}

Some additional properties of KL polynomials are recorded in this proposition.

\begin{prop} \label{prop:KLpolyprops} Fix $y < w$. The coefficient of $v^{\ell(w) - \ell(y)}$ in $h_{y,w}$ is $1$. Moreover, the coefficient of $v^k$ is zero unless $k \le \ell(w) -
\ell(y)$ and unless $k$ as the same parity as $\ell(w) - \ell(y)$. \end{prop}

By convention, one sets $h_{w,w}(v) = 1$ and $h_{y,w}(v) = 0$ if $y \nleq w$, so that $b_w = \sum_y h_{y,w}(v)H_y$. Thus, the KL polynomials give the change of basis matrix between the standard basis and the KL basis. Note also that $b_1 = H_1 = 1$.
	
\begin{ex} Suppose that $W$ is finite, with longest element $w_0$. Then \[b_{w_0} = \sum_{y \in W} v^{\ell(w_0) - \ell(y)} H_y.\]  \end{ex}
	
When $h_{y,w} = v^{\ell(w)-\ell(y)}$ this KL polynomial is called \emph{trivial}, because it has the simplest possible behavior. We say that $w$ is \emph{\revadd{(rationally)} smooth} if all Kazhdan-Lusztig
polynomials $h_{y,w}$ are trivial. For instance, the longest word of any parabolic subgroup is smooth. Any element of a dihedral group (e.g. $S_3$) is smooth.

Starting by giving examples of smooth elements is misleading, as smoothness is the exception in larger Coxeter groups. Polo \cite{Polo} proved that KL polynomials become arbitrarily
complicated, in that any polynomial whatsoever satisfying the properties of Proposition \ref{prop:KLpolyprops} is a KL polynomial for some pair of elements in $S_n$, for some $n$. There
are no known closed formulas for KL polynomials in general.

\begin{ex} \label{ex:notsmooth} For symmetric groups, the first non-smooth elements appear in $S_4$: the permutations $tsut$ and $sutsu$. Here are two nontrivial KL polynomials:
$h_{1,tsut} = {\color{red} v^2} + v^4$ and $h_{1,sutsu} = {\color{red} v^3} + v^5$ (the red terms make the KL polynomial nontrivial). \end{ex}

To conclude this section we state the Kazhdan-Lusztig inversion formula, which expresses the standard basis in terms of the KL basis.

\begin{prop}\label{prop:KLinversion} (\cite[Theorem 3.1]{KazLus79}, see also \cite[Remark 3.10]{Soer97})
If $W$ is a finite Coxeter group, then
\begin{equation} \label{eq:KLinversion} H_y = \sum_{x\leq y} (-1)^{\ell(y)-\ell(x)}h_{w_0y,w_0x}(v) b_x, \end{equation}
where $w_0\in W$ denotes the longest element.
\end{prop}

\subsection{The half twist}
\label{subsec:Hw0}

If $W$ is a finite Coxeter group and $w_0\in W$ is the longest element, then we refer to the element $H_{w_0}$ as the \emph{half twist} and $H_{w_0}^2$ as the \emph{full twist}.  These operators will be the major players in this paper. \revise{Here is a consequence of Proposition \ref{prop:KLinversion}.}

\begin{cor} \label{cor:KLhalftwist}
We have 
\begin{equation} \label{eq:KLhalftwist} H_{w_0} = \sum_{x \in W} (-1)^{\ell(w_0) - \ell(x)} h_{1, w_0 x}(v) b_x.\end{equation}
\end{cor}

Due to the importance of the half and full twists, the KL polynomials $h_{1,w}$, which are as mysterious as any, will appear often throughout this paper.

Now we discuss the operation of conjugation by $H_{w_0}$.

\begin{definition}
If $W$ is a finite Coxeter group, let $\tau:W\rightarrow W$ denote the group automorphism $\tau(x)=w_0xw_0$.
\end{definition}

Clearly $\tau^2 = \Id$. One can show that $\tau$ preserves the set of simple reflections, so it corresponds to an automorphism of the Coxeter-Dynkin diagram.  Consequently, if there are no nontrivial diagram automorphisms then $\tau=\Id_W$, and $w_0$ is central in $W$. 


\begin{lemma}\label{lem:ftIsCentral_decat}
Let $W$ be a finite Coxeter group. Then 
\begin{equation} H_{w_0} H_x H_{w_0}\inv = H_{\tau(x)}.\end{equation}
Thus the full twist $H_{w_0}^2$ is central in $\HB$, and if $W$ has no diagam automorphisms then $H_{w_0}$ is also central.
\end{lemma}

\begin{proof}
It is enough to show the result when $x = s$ is a simple reflection. We have
\[
\revadd{H_{w_0s} H_s = H_{w_0}}  = H_{\tau(s)}H_{\tau(s)w_0}  = H_{\tau(s)} H_{w_0 s}
\]
since $\tau(s)w_0 = w_0s$.  It follows that  $H_{w_0} H_s = H_{\tau(s)} H_{w_0}$.
\end{proof}

Lusztig proved \cite{Lusz81} that, after extending scalars from $\Z[v,v\inv]$ to $\Q(v)$, the Hecke algebra $\HB(W)$ of a finite Coxeter group is isomorphic to the group algebra $\Q(v)[W]$, and hence is semisimple. A central element in a semisimple algebra acts diagonalizably on any representation, so it follows that left multiplication by the full twist is a diagonalizable operator on $\HB$ (after extending scalars).

\subsection{Soergel's categorification of the Hecke algebra}
\label{subsec:SoergelCat}

If $\AC$ is an additive category, then the (split) Grothendieck group $[\AC]$ is defined to be the abelian group generated by symbols $[A]$ (called the \emph{class} of $A$) for objects $A$ of $\AC$, modulo $[X]=[A]+[B]$ whenever $X$ is isomorphic to $A \oplus B$. In particular $A\cong A'$ implies $[A]=[A']$.  Let $\AC$ be an idempotent complete additive category such that each object is a direct sum of finitely many indecomposable objects, and these indecomposable summands are unique up to isomorphism and reordering.  For instance, Krull-Schmidt categories have this property.  In this case, $[\AC]$ is a free abelian group with basis given by the isomorphism
classes of indecomposable objects. An object of $\AC$ is thus determined (up to isomorphism) by its class $[A]$ in $[\AC]$, because this symbol determines the indecomposable summands of $A$ with their multiplicities. We say that $\AC$ categorifies $[\AC]$, and that an object $A$ categorifies its class $[A]\in [\AC]$, and so on.

If $\AC$ is graded with grading shift $(1)$, then $[\AC]$ is a $\Z[v,v\inv]$-module via $v[M] = [M(1)]$. If $\AC$ is monoidal with tensor product $\ot$, then $[\AC]$ is a ring via
$[M]\cdot [N] = [M \ot N]$.

\revcomment{Removed a rambly remark.}



Now we discuss Soergel's categorification of the Hecke algebra. Let $\hg$ denote the reflection representation\footnote{More generally, $\hg$ is allowed to be any \emph{realization} of $W$ (see \cite{EWsoergelCalc}).} over $\R$. Let $R=\R[\hg^\ast]$ denote its polynomial ring, viewed as a $\Z$-graded ring where the linear functionals $\hg^*$ live in degree $2$. Then $R$ also has an action of $W$. In type $A_{n-1}$, the reader is welcome to assume that $R = \RM[x_1, \ldots, x_n]$, thought of as a graded
ring with $\deg(x_i)=2$, equipped with its natural action of $S_n$.

Let $R\mathrm{-gbimod}$ denote the category of finitely generated $\Z$-graded $(R,R)$-bimodules with degree zero $R$-bilinear maps.

The category $R\mathrm{-gbimod}$ is monoidal
 with tensor product $\ot=\otimes_R$ and monoidal identity $\one=R$, the trivial bimodule. We often abuse notation and write $MN$ instead of $M \ot N$ for the tensor product of
two $R$-bimodules.

For each subset $I \subset S$, let $R^I\subset R$ denote the $W_I$-invariant subalgebra, consisting of polynomials $f\in R$ such that $s(f)=f$ for all $s\in I$. If $I=\{s\}$ is a
singleton, then $R^{\{s\}}$ will also be denoted $R^s$. For each simple reflection $s$, let $B_s$ denote the graded $(R,R)$ bimodule defined by
\[
B_s = R\otimes_{R^s} R(1).
\]
In our convention for grading shifts, the element $1\otimes 1 \in B_s$ lies in degree $-1$.
In type $A$, we may use the slightly less cumbersome notation $B_i$ for $B_{s_i}$.

For an expression $\un{w}$, i.e. a sequence $(s_{i_1},s_{i_2},\ldots,s_{i_d})$ of simple reflections, let $\BS(\un{w})$ denote the tensor product $B_{s_{i_1}} B_{s_{i_2}} \cdots B_{s_{i_d}}$. A bimodule obtained in this way is called a \emph{Bott-Samelson bimodule}.

\begin{definition}\label{def:SBim_n}
Let $\SBim\subset R\mathrm{-gbimod}$, the category of \emph{Soergel bimodules}, denote the smallest full subcategory containing $R$ and $B_s$ ($s\in S$), and closed under direct sums, direct summands, tensor products over $R$, and grading shifts.  \revadd{The grading shift by $k\in \Z$ of $\SBim$ will be denoted by $B(-k)$.}
\end{definition}


\revise{The category $\SBim$ has the graded Krull-Schmidt property, see \cite[\S 11]{EMTW}.}


\revadd{Note that our convention is that morphisms in $\SBim$ are degree zero maps of bimodules.  Of course, we may also consider the graded hom space between Soergel bimodules, denoted $\Hom^\Z(M,N)$ and identified with the graded $\R$-vector space $\bigoplus_{n \in \ZM} \Hom_{\SBim}(M,N(n))$. We discuss this further in \S \ref{subsec:celltriprelude}.}  Since $R$ is commutative, the graded hom space $\Homg(M,N)$ is a graded $R$-bimodule in the obvious way. Note that $\Homg(R,R) \cong R$. For $f \in R$ and $g \in \Homg(M,N)$ we may write $f \ot g$ for the left action of $f$ on $g$, and $g \ot f$ for its right action.

\revadd{Given an $R$-bimodule $M$, one can define $\DM(M)$ as the graded vector space of left $R$-module morphisms $M \to R$ of all degrees. Then $\DM(M)$ inherits the structure of an $R$-bimodule, and $\DM$ is a contravariant functor. Note that $\DM(M(1)) \cong \DM(M)(-1)$. One can prove that $\SBim$ is preserved by $\DM$, and that $\DM(M \ot N) \cong \DM(M) \ot \DM(N)$. Since $\DM(B_s) \cong B_s$, Bott-Samelson bimodules are self-dual. The functor $\DM$ will categorify the bar involution on $\HB$.}

The following is an amalgamation of some major results of Soergel.

\begin{thm}[Soergel Categorification Theorem] \label{thm:SCT} (\cite[Lemma 6.13 and Satz 6.14]{Soer07}) There is an indecomposable object $B_w$ in $\SBim$ for each $w \in W$, which
appears as a summand inside $\BS(\un{w})$ for any reduced expression $\un{w}$ of $w$, and does not appear as a summand in $\BS(\un{x})$ for any shorter expression $\un{x}$. The
isomorphism classes of indecomposable objects in $\SBim$ are parametrized, up to grading shift, by $\{B_w\}_{w \in W}$. \revadd{One has $\DM(B_w) \cong B_w$.}

There is a $\ZM[v,v\inv]$-algebra isomorphism from the split Grothendieck group $[\SBim]$ to the Hecke algebra $\HB$, sending \[ [B_s] \mapsto b_s. \] Thus, for any Soergel bimodule
$B$, we let $[B]$ represent the corresponding element in $\HB$.

For any two Soergel bimodules $B$ and $B'$, the space of $R$-bimodule morphisms $\Homg(B,B')$ is a free left (resp. right) $R$-module with graded rank equal to $([B], [B'])$. Here
$(-, -) \co \HB \times \HB \to \ZM[v,v\inv]$ is the standard pairing in the Hecke algebra (see \cite[\S 2.4]{EWsoergelCalc}). This fact is known as the \emph{Soergel Hom Formula}. \end{thm}

\begin{ex} The bimodule $B_1$ is just the monoidal identity $\one$. The bimodule $B_{w_0}$ can be independently described as $R \ot_{R^{W}} R (\ell(w_0))$. \end{ex}

The Soergel categorification theorem can be proven with comparatively elementary techniques, but the following theorem, often known as Soergel's conjecture, is highly nontrivial.  It was proven for Weyl groups and dihedral groups by Soergel \cite{Soer90, Soer07} and for for general Coxeter groups by the first author and Geordie Williamson \cite{EWHodge}.  Furthermore, it relies on the fact that the base field for $\SBim$ is $\R$, and may fail for realizations defined over a field of finite characteristic.

\begin{thm} \label{thm:SC} The isomorphism $[\SBim] \to \HB$ of the Soergel Categorification Theorem sends \[ [B_w] \mapsto b_w. \] \end{thm}

Just as the KL basis $\{b_w\}$ is mysterious, so too are the bimodules $B_w$, and there are no known direct descriptions of $B_w$ in general. However, the bimodules $B_s$ are
straightforward, which is why one often restricts one's attention to Bott-Samelson bimodules.

To give an example of how the Soergel Hom Formula can be used, it implies that the KL polynomial $h_{1w}$ encodes the graded rank of $\Hom(R,B_w)$ as a left $R$-module.

A major implication of Theorem \ref{thm:SC} and the Soergel Hom Formula is that $\SBim$ is \emph{mixed} (see \cite{WebsCan15}), which is summarized by the following proposition. Note
that the indecomposable Soergel bimodules have the form $B_x(k)$ for various $x \in S_n$ and $k \in \Z$, while the particular bimodules $B_x$ where $k=0$ are special \revadd{as they are
self-dual.}

\begin{prop} \label{prop:soergelmixed} If $w \ne x$ then $\Hom(B_w, B_x(k))$ is zero unless $k \ge 1$. Also, $\Hom(B_w,B_w(k))$ is zero unless $k \ge 0$. Moreover, when $k=0$,
$\Hom(B_w,B_w)$ is one-dimensional, spanned by the identity map. \end{prop}

In other words, if one restricts to degree zero maps between self-dual indecomposable Soergel bimodules then the category looks like it is semisimple: there are no nonzero morphisms in
$\Hom(B_w, B_x)$ unless $w=x$, and $\End(B_w)$ is one-dimensional, spanned by the identity map.

In particular, the graded Jacobson radical of $\Endg(\bigoplus_{w \in W} B_w)$ is precisely the space spanned by positive degree morphisms. We will be using this property to understand
the minimal forms of complexes in $\SBim$.

\begin{remark} \label{remark:adjunction} \revadd{It is also true that for each Soergel bimodule $B$ there is another bimodule $B^\vee$, called its \emph{adjoint}, for which the functors $B \ot (-)$ and $B^\vee \ot (-)$ are biadjoint. Furthermore, there is a (contravariant, monoidally-contravariant) functor sending each object to its adjoint. The adjoint of $B_w$ is $B_{w^{-1}}$. One should not confuse adjunction with duality $\DM$, which is contravariant but monoidally covariant.} \end{remark}

\subsection{Notation for complexes} 

\revcomment{Moved this chapter here.}

\revadd{As we are about to introduce complexes of Soergel bimodules, we pause to state our general conventions for homological algebra.}

Let $\AC$ be an additive category.  We let $\KC(\AC)$ denote the homotopy category of complexes over $\AC$.  We prefer the cohomological convention for differentials, hence complexes will be denoted
\[
\cdots \buildrel d\over \rightarrow C^k \buildrel d\over \rightarrow C^{k+1}\buildrel d\over \rightarrow \cdots.
\]
We call $C^k$ the $k$-th chain object of $C$. We write $C \simeq D$ when $C$ and $D$ are homotopy equivalent complexes. Let $\KC^b(\AC),\KC^-(\AC),\KC^+(\AC)\subset \KC(\AC)$ denote the full subcategories of complexes which are bounded, respectively bounded from the right, respectively bounded from the left.

For $i \in \Z$ we let $[i]:\KC(\AC)\rightarrow \KC(\AC)$ denote the functor which shifts complexes $i$ units ``to the left.''  More precisely, $C[i]$ denotes the complex with $C[i]^k=C^{k+i}$, $d_{C[i]}=(-1)^i d_{C}$.  By a morphism $C \to D$ of degree $i$, we mean a chain map $f:C[-i]\to D$, or the equivalent data of a chain map $C \to D[i]$.

If $C,D\in \KC(\AC)$ are complexes, then we let $\Homc(C,D)$ denote the hom complex.  The $k$-th chain group of $\Homc(C,D)$ is the set of linear maps (not necessarily chain maps) from $C$ to $D$ of degree $k$, and the differential is given by the super-commutator
\[
f\mapsto d_D\circ f  - (-1)^k f\circ d_C
\]
when $f$ has degree $k$. The cohomology of $\Homc(C,D)$ yields the morphisms in the homotopy category: we have isomorphisms of graded abelian groups
\[
H^k(\Homc(C,D))\cong \Hom_{\KC(\AC)}(C[-k],D).
\]

\revadd{When $\AC$ is a graded category (like $\SBim$) and $C, D \in \KC(\AC)$, then one can construct the bigraded Hom space
\[
\Homgg(C,D) := \bigoplus_{n,m \in \Z} \Hom(C,D(n)[m]).
\]
See \S \ref{subsec:celltri} for more details.}

\subsection{Braids and Rouquier complexes}
\label{subsec:Rouq}

Theorem \ref{thm:SC} says that the indecomposable bimodules $B_w$ categorify the KL basis. It is natural to wonder how to categorify the standard basis using Soergel's theory. Given that expressing the standard basis in terms of the KL basis requires signs (see \eqref{eq:KLinversion}), one should expect that the standard basis can be expressed using complexes of Soergel bimodules. The standard basis is the image of the braid group under its map to the Hecke algebra, so we will instead associate complexes of Soergel bimodules to braids.

Associated to an arbitrary Coxeter group $W$ we have the braid group $\Br(W)$, which is generated by invertible elements $\s_s$ for $s \in S$ called \emph{crossings}.  The Hecke algebra is a quotient of the group algebra $\ZM[v,v\inv] [\Br]$, and the quotient map sends $\s_s$ to $H_{s}$.

The definition below is due to Rouquier \cite{Rou04}.

\begin{definition}\label{def:Rouquier}
For each simple reflection $s\in W$, with corresponding braid generator $\s = \s_s$, define the following two complexes of Soergel bimodules, living in the bounded homotopy category $\KC^b(\SBim)$.
\[ F(\s) = \Big(0 \to \un{B_{s}}(0) \longrightarrow R(1) \to 0\Big),\]
\[ F(\s\inv) = \Big(0 \to R(-1) \longrightarrow \un{B_{s}}(0) \to 0 \Big). \]
We have underlined the terms in homological degree zero. If $\bbeta$ is a braid word (a word in the generators $\sigma_s^\pm$), then we define $F(\bbeta)$ to be the corresponding tensor product of complexes $F(\sigma_s^\pm)$. Such a complex $F(\bbeta)$ is called a \emph{Rouquier complex}.
\end{definition}

When $\bbeta$ and $\bbeta'$ are two braid words for the same element $\b$ of the braid group, Rouquier proved that $F(\bbeta)$ and $F(\bbeta')$ are homotopy equivalent, and in fact this equivalence is canonical.

\begin{remark} \label{rmk:dotsdefn} The differential in $F(\s)$ is the multiplication map $f \ot g \mapsto fg$, which is a bimodule map $B_s \to R(1)$. The differential in
$F(\s\inv)$ is the coproduct\footnote{\revise{It is the coproduct for the graded Frobenius extension $R^s \subset R$.}} map $1 \mapsto \frac{1}{2}(\a_s \ot 1 + 1 \ot \a_s)$, which is a bimodule map $R \to B_s(1)$.  We note for future use that these maps can be composed to give a bimodule map $R \to R(2)$, equal to multiplication by the simple root $\a_s$. In type $A$, one has $\a_{s_i} = x_i - x_{i+1}$.
\end{remark}

\revcomment{Removed statements about dots above, added remark below.}

\begin{remark} \revise{For the reader familiar with the diagrammatic Hecke category from \cite{EKho} or \cite{EWsoergelCalc}, the differentials in $F(\s)$ and $F(\s\inv)$ are depicted as \emph{dots},
and we often call them by this name. Multiplication by $\a_s$ is depicted as a \emph{barbell}, and we use this name when we wish to emphasize that it appears as the composition of two dots.}
\end{remark}

If $\AC$ is an idempotent-complete additive category and $\KC^b(\AC)$ is its bounded homotopy category, then the \revise{inclusion of $\AC$ into $\KC^b(\AC)$ (as complexes supported in degree zero) induces an isomorphism between the split Grothendieck group of $\AC$ and the triangulated Grothendieck group of $\KC^b(\AC)$. The inverse isomorphism} sends a complex $C^\bullet$ to its Euler characteristic $\sum_{i}(-1)^i[C^i]$. Under the
isomorphism $[\KC^b(\SBim)] \cong [\SBim] \cong \HB$, it is clear that $[F(\s)] = b_s - v = H_s$ and $[F(\s\inv)] = b_s - v^{-1} = H_s\inv$. From this it follows that the class of any
Rouquier complex $F(\bbeta)$ in the Grothendieck group agrees with the image of its braid $\b$ under the quotient map to the Hecke algebra.

We conclude this section with the following fact (see \cite[Lemma 6.5]{EWHodge}).

\begin{lemma}\label{lem:BabsorbsF}
Let $s$ be a simple reflection and $x \in W$ an element such that $sx<x$. In the Hecke algebra, we have $H_s b_x = v^{-1} b_x$. In the categorification, we have
\begin{equation} \label{eq:sdown}
F(\s) B_x \simeq B_x(-1).
\end{equation}
\end{lemma}

\subsection{Canonicity of Rouquier complexes}
\label{subsec:RouqCanon}

Let us rephrase Rouquier's result from \cite{Rou04} on the canonical homotopy equivalence between Rouquier complexes for the same braid.

\begin{prop}(Rouquier canonicity) \label{prop:rouquiercanonicity}  For each pair $\bbeta$, $\bbeta'$ of braid words representing a braid $\b$, there exists a homotopy equivalence $\phi_{\bbeta,\bbeta'} \co F(\bbeta) \to F(\bbeta')$.  These homotopy equivalences are transitive: $\phi_{\bbeta',\bbeta''} \circ \phi_{\bbeta,\bbeta'} \simeq \phi_{\bbeta,\bbeta''}$.  Moreover, if $\bbeta'$ and $\bbeta'$ represent the same braid and $\s$ is a crossing, then $\phi_{\bbeta \s,\bbeta' \s} \simeq \phi_{\bbeta,\bbeta'} \ot \Id_{F(\s)}$. \end{prop}


\begin{remark} An explicit description of these homotopy equivalences in type $A$ can be found in \cite{EKra}. There, a homotopy equivalence is defined for each
braid relation (and for the relation $\s \s^{-1} = 1$), and the movie moves (certain coherence relations) are checked between these homotopy equivalences. This gives another proof of
Rouquier's result in type $A$. \end{remark}

\begin{remark} Because $F = F(\bbeta)$ is invertible, \revise{tensoring with $F^{-1}$ induces an isomorphism} $\Homgg(F,F) \cong \Homgg(\one,\one) \cong R$. Hence, the space of homotopy equivalences from $F(\bbeta)$ to $F(\bbeta')$ is one dimensional, and proving Rouquier canonicity amounts to keeping track of certain scalars.

It is a simple consequence of cell theory (see \S\ref{sec:cells}) that no shift of the identity bimodule occurs as a direct summand of $B_x\otimes B_y$ for nontrivial $x,y$, hence each Rouquier complex has a unique summand which is a shift of $R$. A homotopy equivalence $F(\bbeta) \to F(\bbeta')$ is homotopic to $\phi_{\bbeta,\bbeta'}$ if and only if it induces
the identity map on the \revadd{unique} $R$ summands. \end{remark}

Henceforth, whenever we work in the homotopy category $\KC^b(\SBim_n)$, we write $F(\b)$ to represent $F(\bbeta)$ for any choice of braid word $\bbeta$ representing $\b$, the choice
being irrelevant up to canonical isomorphism.

We will use Rouquier canonicity to prove certain statements about tensor-commutativity. Suppose that the braids $\b$ and $\g$ commute. Then 
\[ F(\b) F(\g) F(\b\inv) \simeq F(\g) \] because
they represent the same braid. When identifying these complexes, we always use the canonical homotopy equivalence. If both $\b_1$ and $\b_2$ commute with $\g$, then the composition
\begin{equation} \label{eq:conjblah1} F(\g) \to F(\b_2) F(\g) F(\b_2\inv) \to F(\b_1) F(\b_2) F(\g) F(\b_2\inv) F(\b_1\inv) \to F(\b_1 \b_2) F(\g) F((\b_1 \b_2)\inv)\end{equation} agrees
with the map \begin{equation} \label{eq:conjblah2} F(\g) \to F(\b_1 \b_2) F(\g) F((\b_1\b_2)\inv) \end{equation} up to homotopy. This is because every arrow is the canonical homotopy
equivalence. 

For future reference, we record a consequence for the center of the braid group. \revcomment{Made the following into a lemma, as suggested.}

\begin{lemma} \label{lem:centerofbraidgroup} If $\g$ is in the center of the braid group, and $\beta$ is any braid, then the canonical homotopy equivalence induces an isomorphism
\[ F(\g) F(\b) \to F(\b) F(\g) \]
in $\KC^b(\SBim_n)$. These isomorphisms are compatible with composition of braids, in that the two compositions
\[ F(\g) F(\b_1) F(\b_2) \to F(\b_1) F(\g) F(\b_2) \to F(\b_1) F(\b_2) F(\g) \]
and
\[ F(\g) F(\b_1) F(\b_2) \to F(\g) F(\b_1 \b_2) \to F(\b_1 \b_2) F(\g) \to F(\b_1) F(\b_2) F(\g) \]
agree up to homotopy. \end{lemma}

\begin{proof} Every arrow is a canonical homotopy equivalence as in Proposition \ref{prop:rouquiercanonicity}. \end{proof}
	
\subsection{Trivialities about the polynomial action}
\label{subsec:trivial}

Let $C \in \KC^b(\SBim_n)$ be an arbitrary complex. Then $\Homgg(\one,C)$ is naturally a bigraded $R$-bimodule, but it is one on which the right and left actions of $R$ agree. More precisely, observe that $R=\Endgg(\one)$.  For each $f\in R$, the following morphisms are all equal for any $\a \in \Homgg(\one,C)$. 
\begin{equation} \label{eq:silly} f \cdot \a = \a \cdot f = \a \circ f = \a \ot f = f \ot \a. \end{equation}
The first two terms are the left and right action. The next term is composition. The last two terms implicitly use the isomorphisms \[\Homgg(\one, C) \cong \Homgg(\one \ot \one, \one \ot C) \cong \Homgg(\one \ot \one, C \ot \one).\]

\subsection{Braid conjugation acting on morphisms}
\label{subsec:conjugate}

Given a braid $\b$, we let 
\[ \Psi_\b \co \KC^b(\SBim) \to \KC^b(\SBim)\] denote the functor which, on objects, sends $M \mapsto F(\b) \ot M \ot F(\b)^{-1}$, and on morphisms, sends $f
\mapsto 1 \ot f \ot 1$. We refer to this functor as \emph{conjugation} by the Rouquier complex of a braid. It is an invertible functor, with inverse $\Psi_{\b^{-1}}$. \revise{All Hom spaces in this section are taken within the category $\KC^b(\SBim)$.}

Let $\g$ be a braid that commutes with $\b$. Then, as noted in \S\ref{subsec:RouqCanon}, there is a canonical homotopy equivalence $F(\g) \to \Psi_{\b}(F(\g))$, which we temporarily denote $\rho_{\g,\b}$. If $\g$ and $\g'$ are both braids which commute with $\b$, and $f \in \Homg(F(\g),F(\g'))$, we define $\psi_\b(f) \in \Homg(F(\g),F(\g'))$ as the composition
\begin{equation} F(\g) \buildrel {\rho_{\g,\b}} \over \longrightarrow F(\b) \ot F(\g) \ot F(\b\inv) \buildrel {1 \ot f \ot 1} \over \longrightarrow F(\b) \ot F(\g') \ot F(\b\inv) \buildrel {\rho_{\g',\b}\inv} \over \longrightarrow F(\g'). \end{equation}
We refer to $\psi_\b$ acting on $\Homg(F(\g),F(\g'))$ as the \emph{conjugation action of braids on morphisms}.

It is clear that \begin{equation} \label{eq:conjfunctorial} \psi_\b(f \circ g) = \psi_\b(f) \circ \psi_\b(g)\end{equation} for morphisms $f$ and $g$ where this makes sense. It is also straightforward that, for two braids $\b_1$ and $\b_2$ which both commute with $\g$ and $\g'$, one has
one has \begin{equation} \label{eq:conjisaction} \psi_{\b_1 \b_2}(f) \simeq \psi_{\b_1}(\psi_{\b_2}(f)). \end{equation}
This follows from the equality of \eqref{eq:conjblah1} and \eqref{eq:conjblah2}.

The following lemma is not immediately obvious, and describes how multiplication by polynomials interacts with the conjugation action.

\begin{lemma} \label{lem:howtoconj} Let $f \in R$ be a homogeneous polynomial, let $\b$ and $\g$ braids with $\g$ in the center of the braid group, and let $\a \in \Homg(\one,F(\g))$. Then the following two maps are homotopic:
\begin{equation}\label{eq:conjugation1}
\psi_\b(f \cdot \a) \simeq w(f) \cdot \psi_\b(\a),
\end{equation}
where $w \in W$ is the image of $\b$ under the standard map from the braid group to $W$. \end{lemma}

\begin{proof} Recall from \S\ref{subsec:trivial} that $f \cdot \a$ agrees with $\a \circ f$, viewing $f$ as a chain map in $\Homg(\one,\one)$. Then, by functoriality
\eqref{eq:conjfunctorial}, it is enough to prove the result for $\g = 1$, and $\a$ the identity map. By \eqref{eq:conjisaction}, it is enough to prove this result when $\b$ is either
$\s$ or $\s\inv$ for a simple reflection $s$. Thus it is enough to show that $s(f) \cdot \rho_{1,\s}$ is homotopic to $(1 \ot f \ot 1) \circ \rho_{1,\s}$, and a similar statement
for $\s\inv$.

\revcomment{Reorganization, not an exercise any more, tried to eliminate diagrammatic language.} The homotopy between these two maps is the map $R \to B_s$ obtained by multiplying by $\pa_s(f)$, and then applying the coproduct map. Here, $\pa_s(f) := \frac{f - s(f)}{\a_s} \in R$ is a Demazure operator applied to $f$. The homotopy equivalence $\rho_\s$ can be found in diagrammatic language on
\cite[page 18, move 1a]{EKra}. That this homotopy suffices follows quickly from \cite[Equation 5.2]{EWsoergelCalc}. \end{proof}

We use this lemma in \S\ref{subsec:lambdaEquivsufficient}.

\subsection{Minimal complexes}
\label{subsec:minimal}

This section is an aside on homological algebra. We return to the applications to Soergel bimodules in the next section.  Many statements regarding chain complexes are greatly simplified by the existence of minimal complexes, which are nice representatives of complexes up to homotopy equivalence.  In turn, minimal complexes mostly make sense in the setting of Krull-Schmidt categories.

Recall that an idempotent-complete additive category is called \emph{Krull-Schmidt} if each object is a direct sum of finitely many indecomposables, and the endomorphism ring of every indecomposable object is a local  ring. \revcomment{removed distracting sentence on a non-krull-shmidt category} \revise{When $\AC$ is Krull-Schmidt, let $\JC(\AC)$ denote the ideal generated by all non-isomorphisms between indecomposable objects. This agrees with the Jacobson radical of the category.} \revcomment{several rewordings below}

\begin{definition}\label{def:minimalCx} Let $\AC$ be a additive \revise{Krull-Schmidt} category. A complex $D\in \KC(\AC)$ is called \emph{minimal} if the differential
$d:D^i\rightarrow D^{i+1}$ lies in $\JC(\AC)$ for all $i$. \end{definition}

When two indecomposable objects $X$ appear as summands of chain objects in adjacent homological degrees, and the part of the differential between them is an isomorphism, this pair can be
cancelled. There is a process known as Gaussian elimination of complexes \cite{DBNfast} which produces a new complex, whose chain objects agree with the original complex except with the
two copies of $X$ removed (by taking a complementary direct summand), and which is homotopy equivalent to the original complex. We think of Gaussian elimination as a ``deformation
retract'' of complexes. Repeating this process one obtains a complex where no summand of any differential is an isomorphism.

\begin{proposition}\label{prop:minimalCx1}
If $\AC$ is an \revise{Krull-Schmidt} additive category, then every $D\in \KC(\AC)$ deformation retracts onto a minimal complex $D_{\min}$. \qed
\end{proposition}


Since the differential \revise{in a minimal complex} is in $\JC$, and $\JC$ is an ideal, then any nulhomotopic chain map between minimal complexes is also in $\JC$.

\begin{lemma}\label{lem:minimalCxIso}
If $\AC$ is Krull-Schmidt and $D_1,D_2\in \KC(\AC)$ are minimal, then any homotopy equivalence $\phi:D_1\rightarrow D_2$ is an isomorphism of complexes.
\end{lemma}

\begin{corollary}\label{cor:minimalCxUnique}
If $\AC$ is Krull-Schmidt, then the minimal complex $D_{\min}$ from Proposition \ref{prop:minimalCx1} is unique up to isomorphism (not merely homotopy equivalence). \qed
\end{corollary}

This discussion applies, mutatis mutandis, to the graded context. \revise{The category $\SBim$ is graded Krull-Schmidt.} \revise{For information on} the graded Jacobson radical $\JC(\AC)$, a graded Krull-Schmidt category, etcetera \revise{ see \cite[Chapter 11 appendices]{EMTW}}.

\subsection{Minimal complexes, perversity, the half twist}
\label{subsec:minimalRouq}

\begin{definition}\label{def:rouquierMinCx}
For each $w\in W$, let $F_w$ denote the minimal complex of the Rouquier complex of a positive braid lift of $w$.
\end{definition}

By Rouquier canonicity and Corollary \ref{cor:minimalCxUnique}, $F_w$ depends (up to unique isomorphism) only on $w$, not on the choice of reduced expression giving the positive braid lift.


\begin{defn}\label{def:perverse}  A complex $D \in \KC(\SBim)$ is \emph{perverse}\footnote{Sometimes also called \emph{linear} or \emph{diagonal}.} if the chain bimodule $D^k$ is a direct sum of bimodules of the form $B_w(k)$.  \end{defn}
There can be no nonzero homotopies \revise{(that is, maps of homological degree $-1$ and graded degree $0$)} between perverse complexes, because $\SBim_n$ is mixed, see Proposition \ref{prop:soergelmixed}.  A perverse complex is necessarily minimal.  

Observe that, for each $s \in S$, the complex $F_s = F(\s)$ is perverse (as is $F_s\inv = F(\s\inv)$). The following crucial result was proven in \cite[Theorem 6.9]{EWHodge}.

\begin{thm}\label{thm:diagonalmiracle}
The minimal complexes $F_w$ are perverse, for all $w \in W$.
\end{thm}

As noted previously, an object $A$ of a Krull-Schimdt additive category $\AC$ is pinned down uniquely up to isomorphism by its class in the Grothendieck group.  However if $C\in \KC^b(\AC)$ is a complex, then the class of $[C]\in K_0(\AC)$ certainly does not determine the chain groups $C^k$ (for a counterexample take $C=B\oplus B[1]$ which always has zero Euler characteristic).  However, knowing the symbol of a bounded perverse complex $C$ in a \revadd{mixed} category does determine the chain objects uniquely (but not the differential)! 

Combining Theorem \ref{thm:diagonalmiracle} and the Kazhdan-Lusztig inversion formula \eqref{eq:KLinversion} gives a description of the chain bimodules of $F_w$ for every $w\in W$.  Namely, each occurence of $(-v)^i b_x$ in the right-hand side of \eqref{eq:KLinversion} contributes a summand $B_x(i)$ in homological degree $i$.   In particular, \eqref{eq:KLhalftwist} implies that the KL polynomials $h_{1,x}$ determine the summands appearing in $F_{w_0}$.

\begin{defn}\label{def:HTandFT} If $W$ is a finite Coxeter group with longest element $w_0$, let $\HT$ denote $F_{w_0}$, and let $\FT=\HT \ot \HT$. We call $\HT$ the \emph{half-twist},
and $\FT$ the \emph{full-twist}. In type $A_{n-1}$, we often write the half-twist as $\HT_n$ and the full-twist as $\FT_n$. We may also use the terms half-twist and full-twist to refer
to the corresponding elements of the braid group, or their images in the Hecke algebra; the meaning will be clear from context. \end{defn}

\begin{example} Let $s$ and $t$ be the simple reflections of $S_3$. Then \begin{equation} \label{eq:Hsts} H_{sts} = b_{sts} - v(b_{st}+b_{ts}) + v^2(b_s+b_t) - v^3.\end{equation} Thus the Rouquier complex $F_{sts} = \HT_3$ has the form
\[ F_{sts} = \Big( \un{B_{sts}(0)} \longrightarrow B_{st}(1) \oplus B_{ts}(1) \longrightarrow B_s(2) \oplus B_t(2) \longrightarrow R(3) \Big). \]
We derive \eqref{eq:Hsts} above. Every element of $S_3$ is smooth, meaning that $h_{x,y} = v^{\ell(y)-\ell(x)}$. Thus by Corollary \ref{cor:KLhalftwist} the coefficient of $b_x$ in $H_{sts}$ is $(-v)^{\ell(w_0) - \ell(x)}$. \end{example}

\begin{ex} \label{ex:S4HTetc} Let $\{s,t,u\}$ be the simple reflections of $S_4$. For all but two elements of $S_4$, one has $h_{1,x} = v^{\ell(x)}$. However, one has $h_{1,tsut} = {\color{red} v^2} + v^4$ and $h_{1,sutsu} = {\color{red} v^3} + v^5$. Note that $tsut = w_0(su)$ and $sutsu = w_0(t)$. Thus the Rouquier complex $F_{w_0} = \HT_4$ has the following form.
\begin{equation} \label{eq:Fstsuts} \HT_4 \simeq 
\begin{tikzpicture}[baseline=-.2em]
\tikzstyle{every node}=[font=\scriptsize]
\node at (0,0) {$\underline{B_{stsuts}(0)}$};
\node (y) at (.65,0) {};
\node (z) at (1.5,0) {};
\node at (2,.5) {$B_{tstut}(1)$};
\node at (2,0) {$B_{tutst}(1)$};
\node at (2,-.5) {$B_{sutsu}(1)$};
\node (a) at (2.5,0) {};
\node (b) at (3.5,0) {};
\node at (4,1.25) {$B_{tstu}(2)$};
\node at (4,.75) {$B_{utst}(2)$};
\node at (4,.25) {$B_{tuts}(2)$};
\node at (4,-.25) {$B_{stut}(2)$};
\node at (4,-.75) {$B_{tsut}(2)$};
\node at (4,-1.25) {${\color{red} B_{su}(2)}$};
\node (c) at (4.5,0) {};
\node (d) at (5.5,0) {};
\node at (6,1.5) {${\color{red} B_{t}(3)}$};
\node at (6,1) {$B_{stu}(3)$};
\node at (6,.5) {$B_{uts}(3)$};
\node at (6,0) {$B_{tst}(3)$};
\node at (6,-.5) {$B_{tut}(3)$};
\node at (6,-1) {$B_{tsu}(3)$};
\node at (6,-1.5) {$B_{sut}(3)$};
\node (e) at (6.5,0) {};
\node (f) at (7.5,0) {};
\node at (8,1) {$B_{st}(4)$};
\node at (8,.5) {$B_{ts}(4)$};
\node at (8,0) {$B_{ut}(4)$};
\node at (8,-.5) {$B_{tu}(4)$};
\node at (8,-1) {$B_{su}(4)$};
\node (g) at (8.5,0) {};
\node (h) at (9.5,0) {};
\node at (10,.5) {$B_{s}(5)$};
\node at (10,0) {$B_{t}(5)$};
\node  at (10,-.5) {$B_{u}(5)$};
\node (i) at (10.5,0) {};
\node (j) at (11.5,0) {};
\node at (12,0) {$B_1(6)$};
\path[->,>=stealth',shorten >=1pt,auto,node distance=1.8cm,
  thick]
(y) edge node[above] {} (z)
(a) edge node[above] {} (b)
(c) edge node[above] {} (d)
(e) edge node[above] {} (f)
(g) edge node[above] {} (h)
(i) edge node[above] {} (j);
\end{tikzpicture}
\end{equation}
The nontrivial Kazhdan-Lusztig polynomials mentioned above give rise to the ``additional" terms $B_{su}(2)$ and $B_t(3)$ appearing above. \end{ex}

\section{Abstract cell theory}
\label{sec:cells}

In this chapter we discuss cells in algebras and cells in monoidal categories. The eventual application will be to Hecke algebras and Soergel bimodules. In type $A$, the resulting cell theory is much simpler than the general case, having a nice combinatorial description, so the novice reader is encouraged to skip ahead and read \S\ref{sec:typeAcells} concurrently.

In \S\ref{subsec:twists} we introduce the theory of twist-like complexes, and study how they act on cells.

\subsection{Cells in algebras}
\label{subsec:algCells}

Let $\kbbm$ be a commutative ring and let $A$ be a $\kbbm$-algebra.  Let $W$ be some indexing set and $\{b_x\}_{x \in W}$ be a $\kbbm$-basis of $A$. For an arbitrary element $a \in A$, we say that $$b_i \babysumset a$$ if, writing $a$ as a linear combination of basis elements, the term $b_i$ appears with nonzero coefficient.   One can place several relations on the elements of this basis (or on the indexing set $I$), by saying that
\[ b_x \le_{L} b_y \quad \textrm{if} \quad b_x \babysumset m \cdot b_y \textrm{ for some } m \in A. \]
\[ b_x \le_{R} b_y \quad \textrm{if} \quad b_x \babysumset b_y \cdot n \textrm{ for some } n \in A. \]
\[ b_x \le_{LR} b_y \quad \textrm{if} \quad b_x \babysumset m \cdot b_y \cdot n \textrm{ for some  } m, n \in A. \]
If the structure coefficients of $A$ in the basis $\{b_x\}$ are non-negative integers, then the above relations are transitive.  If these relations are not transitive, we consider instead their transitive closures, and abusively use the same notation for the transitive closure.

The equivalence classes under these relations are known as \emph{left}, \emph{right}, and \emph{two-sided cells} respectively. We write $b_x \sim_L b_y$ if these two basis elements are in the same left cell, and similarly for $b_x \sim_R b_y$ and $b_x \sim_{LR} b_y$.  The set of left (or right, or two-sided) cells inherits a partial order from $\leq_L$ (resp.~$\leq_R$, resp.~$\leq_{LR}$).  Each two-sided cell is a union of left cells (or right cells).

One of the main reasons for the appearance of cells is that they determine a filtration of the algebra $A$ by ideals.

\begin{prop}\label{prop:cellFiltration}
Let $A$ be an algebra with basis $\{b_x\}_{x\in W}$.  For a given two-sided cell $\l$, let $I_{< \l}\subset A$ denote the subspace spanned by basis elements $b_y$ in (two-sided) cells strictly lower than $\lambda$.  Then $I_{< \l}$ is a two-sided ideal.
\end{prop}

\begin{proof}
Immediate from the definitions.
\end{proof}

Similarly, $I_{\le \l}$ is a two-sided ideal, as is $I_{\ngeq \l}$, and so forth. By considering left or right cells instead, one obtains families of left or right ideals in $A$. By
taking the subquotients of these filtrations, one obtains potentially interesting representations of $A$.

Note that the cell theory of an algebra depends heavily on the chosen basis. For instance, if $A=\HB$ is the Hecke algebra associated to a Coxeter system $(W,S)$, then the cell theory
associated to the standard basis $\{H_w\}$ is boring: since each $H_w$ is invertible, they all generate the unit ideal, so they are left, right, and two-sided equivalent. For the
Kazhdan-Lusztig basis, the resulting cell theory and filtration is far from boring. Let us mention some easy examples of cells.

\begin{ex} \label{ex:threecells} Inside $\HB$, $1=b_1$ is identity, so $b_x \le b_1$ for all $x$ (for all three relations). It turns out that $\{1\}$ is an equivalence class (under all
three relations), called the \emph{identity cell}. At the other extreme, $b_{w_0}$ spans a one-dimensional ideal in $\HB$, so that $\{w_0\}$ is an equivalence class (under all three
relations), called the \emph{longest cell}. For $S_3$, there is only one other 2-sided cell $\{b_s,b_t,b_{st},b_{ts}\}$, which splits into two left (resp. right) cells. In general, as
was shown by \revise{ Lusztig \cite[Proposition 3.8]{LuszSqInt},} all the non-identity elements of $W$ with a unique reduced expression (such as the simple reflections) form a two-sided cell, called the \emph{simple cell}\footnote{If $W$ is not irreducible, then the set of non-identity elements with a unique reduced expression split into a disjoint union of two-sided ``simple cells,'' one for each connected component of the Coxeter graph.}. \end{ex}

In summary, given a choice of basis of $A$, one has an associated cell theory \revcomment{removed: which may or may not be boring}.  Now we categorify this notion, replacing the algebra with a monoidal additive category. The (arbitrary) choice of basis becomes the (intrinsic) collection of indecomposable objects.

\subsection{Categorical cells}
\label{subsec:celldef}

 
\revcomment{There are some revisions to this section to try to make the point more clearly.}

\revise{Throughout this section let $\AC$ be a Krull-Schmidt additive category.} Recall that a full subcategory $\BC\subset \AC$ is determined by the objects of $\AC$ which are contained in $\BC$. Let $\ind \AC$ denote the indecomposable objects of $\AC$.

\begin{defn}
A full subcategory $\BC\subset \AC$ is \emph{essential} if it is closed under isomorphisms, and \emph{thick} if it is closed under taking direct sums and direct summands.
\end{defn}

\revcomment{Made some remarks into lemmas, more clearly stated.} The following lemmas are easy.

\begin{lemma}
There is a bijection between thick, essential full subcategories and subsets of $\ind \AC$, given by $\BC \mapsto \BC \cap \ind \AC$.
\end{lemma}

\begin{lemma} There is a bijection between thick, essential full subcategories of $\AC$ and ideals in $\AC$ which are generated by identity maps. This bijection sends $\BC$ to the ideal $I_{\BC}$ generated by the identity maps of objects in $\BC$. The inverse bijection sends an ideal $I$ to the full subcategory of objects whose identity maps are in $I$. \end{lemma}

\begin{proof} (Sketch) The identity map of $M \oplus N$ is contained in an ideal if and only if the identity maps of $M$ and $N$ are contained in that ideal. \end{proof}

\revcomment{More revamped remarks}

\begin{remark} Let $I_{\BC}$ be the ideal generated by the identity maps of a set $\BC$ of objects. Then $f \in I_{\BC}$ if and only if $f = \sum g_i$ is a linear combination of morphisms $g_i$, where each $g_i$ factors through some object in $\BC$.
	
The quotient category $\AC/I_{\BC}$ is typically denoted by $\AC/\BC$. The latter notation often refers to killing morphisms which factor through $\BC$, and these concepts agree by the previous paragraph.  \end{remark}


\revcomment{More reorganization from here} Now we add a monoidal structure, and focus on monoidal ideals (which are closed under both composition and tensor product). The corresponding full subcategories are called tensor ideals, a poor but standard choice of notation (as tensor ideals are not actually ideals, but full subcategories).

\begin{defn} Let $\AC$ be an Krull-Schmidt additive monoidal category. A \emph{left tensor ideal} $\IC \subset \AC$ is a full subcategory which is essential and thick, such that $M
\ot B \in \IC$ whenever $B \in \IC$. The definitions of a \emph{right tensor ideal} and a \emph{two-sided tensor ideal} are similar. \end{defn}

If the identity of $B$ is in an ideal, and the ideal is closed under tensor product on the left with arbitrary morphisms (i.e. it is a left-monoidal ideal), then the identity of $M
\ot B$ is also in the ideal. This leads to the following easy lemma.

\begin{lemma} An thick, essential full subcategory $\IC$ is a left tensor ideal if and only if the corresponding ideal $I_{\IC}$ is left-monoidal. Similar statements can be made
for right and two-sided tensor ideals. \end{lemma}

\begin{remark} If $\AC$ is graded then a tensor ideal (left, right, or two-sided) is automatically closed under grading shifts, because $\one(k) \ot B \cong B \ot \one(k) \cong
B(k)$. \end{remark}

\revcomment{Expansion here} Let $B$ be an object, and let $\AC \cdot B$ be the smallest left tensor ideal containing $B$. Let $I_{\AC \cdot B}$ be the smallest left-monoidal ideal
containing the identity map of $B$. Then $I_{\AC \cdot B}$ is generated by identity maps, not just as a left-monoidal ideal but as an ordinary ideal; after all, for any morphism $f \colon M \to N$, $f \ot \id_B$ is in the ideal generated by $\id_M \ot \id_B$. It is straightforward to verify that $I_{\AC \cdot B}$ is the ideal associated to the thick subcategory $\AC \cdot B$.
	
Given indecomposables $B,B'\in \AC$ we write $B \le_{L} B'$ if $B \in \AC \cdot B'$, or equivalently, if $\id_B \in I_{\AC \cdot B'}$. Unlike the previous section, this relation is
always transitive. Similarly, one defines $B \le_R B'$ if $B \in B' \AC$, and $B \le_{LR} B'$ if $B'\in \AC B\AC$. The corresponding equivalence classes on $\ind \AC$ are called
\emph{left cells}, \emph{right cells}, and \emph{two-sided cells} respectively. Left cells inherit the partial order $\le_{L}$, and so forth.

To any two-sided cell $\l$, we have a tensor ideal $\IC_{\le \l}$ and a monoidal ideal $I_{\le \l}$. We may set $I_{\le \l}$ to be the monoidal ideal generated by $\id_B$ for any
indecomposable $B$ in cell $\l$. For an indecomposable object $B'$ we have $\id_{B'} \in I_{\le \l}$ if and only if $B' \le_{LR} B$, by definition of cells. If $B'$ is in two-sided
cell $\mu$, then this is equivalent to $\mu \le \l$.

\revcomment{returning roughly to what was there before} The Grothendieck group $[\AC]$ is an algebra, \revise{which is free as a $\ZM$-module} with basis given by the set $\ind
\AC$ of indecomposables in $\AC$ up to isomorphism. Multiplication is given by $[M][N] = [M \ot N]$. Recall from \S\ref{subsec:SoergelCat} that an object in a Krull-Schmidt category is uniquely determined up to isomorphism by its symbol in the Grothendieck group. An implication is that the decomposition of tensor products into indecomposable objects is controlled by multiplication in the Grothendieck group. For sake of illustration, let $W$ be an indexing set for the isomorphism classes of indecomposable objects, named $\{B_w\}_{w \in W}$, and let $b_w = [B_w]$ denote the corresponding basis of $[\AC]$.

\begin{proposition}\label{prop:structureCoeffs}
If
\begin{equation} \label{eq:multdownstairs} b_w b_x = \sum c^y_{w,x} b_y \end{equation}
for integers $c^y_{w,x}$, then
\begin{equation} \label{eq:multupstairs} B_w \ot B_x \cong \bigoplus B_y^{\oplus c^y_{w,x}},\end{equation}
and in particular the coefficients $c^y_{w,x}$ are non-negative.\qed
\end{proposition}
If $\AC$ is graded with grading shift functor $(1)$ and $\{B_w\}_{w\in W}$ instead denotes a set of representatives of the indecomposables up to grading shift, then a similar statement is true, where now the structure coefficients $c^y_{w,x}$ are elements of $\Z[v,v\inv]$.

For this reason, the cell theory of $\AC$ and the cell theory of $[\AC]$ agree.



\subsection{The example of $\HT_3$ and $\FT_3$}
\label{subsec:HT3}

We wish to prove some abstract results about the interaction between cells and operators which act like ``twists.'' Eventually, this will be applied to the Hecke algebra and the Soergel
category, and illustrates why cell theory is so important for diagonalizing the full twist. So, to motivate these results, we briefly illustrate this phenomenon for $S_3$.

As mentioned in Example \ref{ex:threecells}, $S_3$ has only three two-sided cells. In order from largest to smallest with respect to $\le_{LR}$, they are: the {\color{brown} identity cell}, the {\color{red} simple cell}, and the {\color{blue} longest cell}. We color-code these cells in the complexes below.

Direct computation shows that
\begin{eqnarray}
\HT_3 \ot {\color{brown}\one} & \simeq & \begin{tikzpicture}[baseline=-.5em]
\node  (a) at (0,0) {$\underline{{\color{blue} B_{sts}}}(0)$};
\node at (2,.5) {${\color{red}B_{st}}(1)$};
\node (b) at (2,0) {$\oplus$};
\node at (2,-.5) {${\color{red}B_{ts}}(1)$};
\node at (4,.5) {${\color{red}B_{s}}(2)$};
\node (c) at (4,0) {$\oplus$};
\node at (4,-.5) {${\color{red}B_{t}}(2)$};
\node (d) at (6,0) {${\color{brown}\one}(3)$};
\path[->,>=stealth',shorten >=1pt,auto,node distance=1.8cm,
  thick]
(a) edge node[above] {} (b)
(b) edge node[above] {} (c)
(c) edge node[above] {} (d);
\end{tikzpicture}
\\ \label{eq:HTBs}
\HT_3\otimes {\color{red}B_s} & \simeq & \begin{tikzpicture}[baseline=-.5em]
\node (a) at (0,0) {$\underline{{\color{blue}B_{sts}}}(-1)$};
\node (b) at (2,0) {${\color{red}B_{ts}}(0)$};
\path[->,>=stealth',shorten >=1pt,auto,node distance=1.8cm,
  thick]
(a) edge node[above] {} (b);
\end{tikzpicture}
\\ \label{eq:HTBst}
\HT_3\otimes {\color{red}B_{st}} & \simeq & \begin{tikzpicture}[baseline=-.5em]
\node (a) at (0,0) {$\underline{{\color{blue}B_{sts}}}(-2)$};
\node (b) at (2,0) {${\color{red}B_{t}}(0)$};
\path[->,>=stealth',shorten >=1pt,auto,node distance=1.8cm,
  thick]
(a) edge node[above] {} (b);
\end{tikzpicture}
\\
\HT_3\otimes {\color{blue}B_{sts}}  & \simeq &  \underline{{\color{blue}B_{sts}}}(-3)
\end{eqnarray}
Equations \eqref{eq:HTBs} and \eqref{eq:HTBst} are also true after swapping $s$ and $t$.  Here is the analogous computation for the full twist:

\[
\FT_3 \ot {\color{brown}\one} \ \ \simeq  \ \ \begin{tikzpicture}[baseline=-.2em]
\tikzstyle{every node}=[font=\scriptsize]
\node (a) at (0,0) {$\underline{{\color{blue}B_{sts}}}(-3)$};
\node at (2,.25) {${\color{blue}B_{sts}}(-1)$};
\node at (2,-.25) {${\color{blue}B_{sts}}(-1)$};
\node (c) at (4,.75) {${\color{blue}B_{sts}}(1)$};
\node at (4,.25) {${\color{blue}B_{sts}}(1)$};
\node at (4,-.25) {${\color{red}B_{s}}(1)$};
\node at (4,-.75) {${\color{red}B_{t}}(1)$};
\node (d) at (6,.5) {${\color{blue}B_{sts}}(3)$};
\node at (6,0) {${\color{red}B_{st}}(2)$};
\node at (6,-.5) {${\color{red}B_{ts}}(2)$};
\node (e) at (8,.25) {${\color{red}B_{st}}(4)$};
\node at (8,-.25) {${\color{red}B_{ts}}(4)$};
\node (f) at (10,.25) {${\color{red}B_{s}}(5)$};
\node at (10,-.25) {${\color{red}B_{t}}(5)$};
\node (g) at (12,0) {${\color{brown}\one}(6)$};
\node (b1) at (1.6,0) {};
\node (b2) at (2.4,0) {};
\node (c1) at (3.6,0) {};
\node (c2) at (4.4,0) {};
\node (d1) at (5.6,0) {};
\node (d2) at (6.4,0) {};
\node (e1) at (7.6,0) {};
\node (e2) at (8.4,0) {};
\node (f1) at (9.6,0) {};
\node (f2) at (10.4,0) {};
\path[->,>=stealth',shorten >=1pt,auto,node distance=1.8cm,
  thick]
(a) edge node[above] {} (b1)
(b2) edge node[above] {} (c1)
(c2) edge node[above] {} (d1)
(d2) edge node[above] {} (e1)
(e2) edge node[above] {} (f1)
(f2) edge node[above] {} (g);
\end{tikzpicture}
\]
\begin{eqnarray}
\label{eq:FTBs}
\FT_3\otimes {\color{red}B_s} & \simeq & \begin{tikzpicture}[baseline=-.2em]
\node (a) at (0,0) {$\underline{{\color{blue}B_{sts}}}(-4)$};
\node (b) at (2.5,0) {${\color{blue}B_{sts}}(-2)$};
\node (c) at (5,0) {${\color{red}B_{s}}(0)$};
\path[->,>=stealth',shorten >=1pt,auto,node distance=1.8cm,
  thick]
(a) edge node[above] {} (b)
(b) edge node[above] {} (c);
\end{tikzpicture}
\\  \label{eq:FTBst}
\FT_3\otimes {\color{red}B_{st}} & \simeq & \begin{tikzpicture}[baseline=-.2em]
\node (a) at (0,0) {$\underline{{\color{blue}B_{sts}}}(-5)$};
\node (b) at (2.5,0) {${\color{blue}B_{sts}}(-1)$};
\node (c) at (5,0) {${\color{red}B_{ts}}(0)$};
\path[->,>=stealth',shorten >=1pt,auto,node distance=1.8cm,
  thick]
(a) edge node[above] {} (b)
(b) edge node[above] {} (c);
\end{tikzpicture}
\\
\FT_3\otimes {\color{blue}B_{sts}}  & \simeq &  {\color{blue}B_{sts}}(-6)
\end{eqnarray}
Again, equations \eqref{eq:FTBs} and \eqref{eq:FTBst} remain true after swappping $s$ and $t$.

We make the following observation: let $F$ denote either $\HT_3$ or $\FT_3$.  Then for each two-sided cell $\l$ there exists a number $\nb_F(\l)$ such that:
\begin{itemize}
\item if $B$ is in cell $\l$ then every occurence of $B$ in the minimal complex of $F$ is in homological degrees $\geq \nb_F(\l)$.
\item if $B$ is in cell $\l$ then the minimal complex of $F\otimes B$ is supported in homological degrees $\leq \nb_F(\l)$.
\end{itemize}
Furthermore, $\FT_3$ appears to act in a ``block upper triangular'' fashion with respect to the cell filtration.  Before moving on to more advanced cell theory techniques, let us abstract the above phenomena.

\subsection{Twists in general}
\label{subsec:twists}

Throughout this section let $\AC$ be an additive monoidal category with the Krull-Schmidt property (so that there is a well behaved notion of minimal complexes).

\begin{definition}\label{def:twistlike}  We say that an invertible \revcomment{removed: minimal} complex $F\in \KC^b(\AC)$ is \emph{twist-like} if for all $B \in \AC$, $F\otimes B\otimes F\inv$ is homotopy equivalent to an object of $\AC$ (that is, a complex supported in degree zero).  
\end{definition}

\begin{remark}\label{rmk:twistsAreTwists}
We expect that the half and full twists are twist-like (hence the name) and, in fact, 
\[
\HT\otimes B_x\otimes \HT\inv \simeq B_{\tau(x)}, \qquad\qquad \FT\otimes B_x\otimes \FT\inv \simeq B_x.
\]
This is well known in type $A$ (and we include a proof in  \S \ref{subsec:commutehalffull}), but we are not aware of any proof in the literature for other finite Coxeter groups.
\end{remark}

Now we define some numerical invariants of twist-like complexes $F$.

\begin{defn}\label{def:nb}
If $F$ is twist-like and $B\in \AC$ is any object, let $\nb_F(B)$ denote the smallest integer such that $F\otimes B$ is homotopy equivalent to a complex supported in degrees $\leq \nb_F(B)$.
\end{defn}

\begin{lemma}\label{lemma:tailProperty}
If $F$ is twist-like and $B'\leq_{LR} B$, then $\nb_F(B')\leq \nb_F(B)$.  In particular, $\nb_F(B)$ depends only on the two-sided cell containing $B$.
\end{lemma}

\begin{proof}
Suppose $B'$ is a direct summand of $M\otimes B\otimes N$ \revise{for some objects $M, N \in \AC$}.  Then in the homotopy category $F\otimes B'$ is a direct summand of
\[
F\otimes M\otimes B\otimes N \simeq (F\otimes M\otimes F\inv) \otimes  (F\otimes B)\otimes N. 
\]
The first factor $F\otimes M\otimes F\inv$ is homotopy equivalent to an object of $\AC$ by assumption.  The middle factor $F\otimes B$ is homotopy equivalent to complex in degrees $\leq \nb_F(B)$.  From this, the inequality $\nb_F(B')\leq \nb_F(M\otimes B\otimes N)\leq \nb_F(B)$ is clear.
\end{proof}

\begin{definition}\label{def:n}
If $F\in \KC^b(\AC)$ is twist-like and $\l$ is a two-sided cell of $\AC$, let $\nb_F(\l)$ denote $\nb_F(B)$ for any $B\in \AC$ in cell $\l$.  We refer to the mapping $\l\mapsto \nb_F(\l)$ as the \emph{homological spectrum} of $F$. 
\end{definition}

Now we relate the invariants $\nb_F(\l)$ to some internally defined numbers.

\revcomment{Made modifications to the definition below since we no longer assume $F$ is a minimal complex.  Also eliminated the confusing $\pb_F(\l)$ from this definition; it is not needed.}
\begin{definition}\label{def:m}
Let $\pi_{\geq \l}\colon \KC^b(\AC)\to \KC^b(\AC/\IC_{\not \geq \l})$ denote the functor which quotients by the ideal of cells not greater than or equal to $\l$.  

For any twist-like complex $F\in \KC^b(\AC)$, let $\mb_F(\l)$ denote the largest integer $k$ such that $\pi_{\geq \l}(F)$ is homotopy equivalent to a complex supported in homological degrees $\geq k$.

We say that a twist-like $F$ is \emph{increasing} if $\mb_F(\l)<\mb_F(\mu)$ whenever $\l<\mu$.

\end{definition}
\revcomment{ removed: the parenthetical (strictly) before increasing above, and the following paragraph: Note that $\mb_F(\l)\le\mb_F(\mu)$ whenever $\l < \mu$ by definition; a twist-like complex is strictly increasing if each of these inequalities is strict.  We will henceforth omit the adverb `strictly' when refering to increasing twist-like complexes.}

\revcomment{Added this remark:}
\begin{remark}
A twist-like $F$ is increasing iff for each cell $\l$, the quotient $\pi_{\geq \l}(F)$ is homotopy equivalent to a complex $X$ which is supported in homological degrees $\geq \mb_F(\l)$ and whose degree $\mb_F(\l)$ chain object is in cell exactly $\l$ (and not higher cells).
\end{remark}

\revcomment{Deleted a lemma on $\pb_F(\l)$ that is no longer used}

\begin{example}
From the computations in \S \ref{subsec:HT3}, we see that for the half-twist $F=\HT_3$ we have \revcomment{fixed mathematical typos}
\begin{itemize}
\item $\nb_F(\l)=\mb_F(\l)=3$ for the maximal cell $\l=(3)$.
\item $\nb_F(\l)=\mb_F(\l)=1$ for the simple cell $\l=(2,1)$.
\item $\nb_F(\l)=\mb_F(\l)=0$ for the minimal cell $\l=(1,1,1)$.
\end{itemize}
Thus $\HT_3$ is twist-like and increasing according to the definition.
\end{example}

\begin{example} Let $G = \HT_3\inv$, the \emph{negative half twist}. Then $\nb_G(\l)=\mb_G(\l)=0$ for all cells $\l$.  Thus, $\HT_3\inv$ is twist-like but not increasing.
\end{example}

\revcomment{Modifying this remark to be compatible with the new definition of $\mb_F(\l)$.}
\begin{remark} \label{rmk:frombelow} 
We could analogously define $\nb^\ast_F(\l)$ to be the largest $k$ such that $F\ot B$ is supported in degrees $\geq k$ (up to homotopy equivalence, for all $B$ in cell $\l$), and $\mb^\ast_F(\l)$ to be the smallest $k$ such that $\pi_{\geq \l}(F)$ is supported in degrees $\leq k$ (up to homotopy equivalence).  A twist-like complex is \emph{(strictly) decreasing} if $\mb_F^\ast(\l)>\mb_F^\ast(\mu)$ when $\l<\mu$. Then $\HT_3\inv$ is strictly decreasing, while $\HT_3$ is not.
%
\end{remark}
	
\begin{remark} We could also define statistics associated to tensor product with $F$ on the right, e.g. the maximal homological degree appearing in the minimal complex of $B \ot F$. We postpone discussion of these generalizations until later in this section. \end{remark}


\begin{lemma} \label{lem:orderedeasy}
Let $F$ be twist-like.  If $\l$ and $\mu$ are two-sided cells and $\l \le \mu$, then $\nb_F(\l) \le \nb_F(\mu)$ and $\mb_F(\l) \le \mb_F(\mu)$. \end{lemma}

\begin{proof} That $\nb_F(\l)\le \nb_F(\mu)$ follows immediately from Lemma \ref{lemma:tailProperty}. That $\mb_F(\l) \le \mb_F(\mu)$ follows from the definition of $\mb_F$, since any cell $\ge \mu$ is also $\ge \l$. \end{proof}

\begin{lemma}\label{lemma:nmineq}
Let $F,G\in \KC^b(\AC)$ be twist-like complexes, and $\l$ be a fixed two-sided cell.  We have
\begin{subequations}
\begin{equation}\label{eq:nAndm}
\nb_F(\l)\geq \mb_F(\l)
\end{equation}
\begin{equation}\label{eq:nFG}
\nb_{G\otimes F}(\l)\leq \nb_F(\l)+\nb_G(\l)
\end{equation}
\begin{equation}\label{eq:mFG}
\mb_{G\otimes F}(\l)\geq \mb_{F}(\l)+\mb_G(\l).
\end{equation}
\end{subequations}
\end{lemma}

\revcomment{The proof below has been revised to eliminate any references to minimal complexes.}
\begin{proof}
%
Fix $B$ in cell $\l$, and let \revise{$X\simeq F\ot B$ be a complex supported in homological degrees $\le \nb_F(\l)$ and in cells $\le \l$, which exists by definition of $\nb_F(\l)$.  Let $\pi_{\geq \l}$ be the quotient functor which kills cells $\not \geq \l$; it is a monoidal functor. From the definition of $\mb_F(\l)$ we have that $\pi_{\geq \l}(X)$ is homotopy equivalent to a complex $X'$ supported in homological degrees $\geq \mb_F(\l)$.  We know that $X'$ is not the zero complex: applying $\pi_{\geq \l}$ to the isomorphism $F\inv\ot X\simeq B$ gives $\pi_{\geq \l}(F\inv)\ot X'\simeq \pi_{\geq \l}(B)$, which is not zero since $B$ is in cell $\l$. The degree $k$ chain object of $X'$ is zero unless $\mb_F(\l)\leq k \leq \nb_F(\l)$. The inequality \eqref{eq:nAndm} follows.}


To show \eqref{eq:nFG}, consider $G \ot X$. This is the total complex of a bicomplex whose $i$-th column is $G \ot X^i[-i]$, where $X^i$ denotes the $i$-th chain object of $X$. Since
each $X^i$ is in cells $\le \l$, $G \ot X^i[-i]$ \revise{is homotopy equivalent to a complex} supported in degrees $\le i+\nb_G(\l)$ (here we have used Lemma \ref{lem:orderedeasy}).  Applying these homotopy equivalences to each column and noting that $X^i=0$ for $i>\nb_F(\l)$ shows that $G \ot X$ is homotopy equivalent to a complex supported in degrees $\le \nb_F(\l) + \nb_G(\l)$, as desired.

The inequality \eqref{eq:mFG} is relatively straightforward, using $\pi_{\geq \l}(F\ot G)\cong \pi_{\geq \l}(F)\otimes \pi_{\geq \l}(G)$. There exist complexes $Y,Z$ such that $\pi_{\geq \l}(F)\simeq Y$, supported in degrees $\geq \mb_F(\l)$ and $\pi_{\geq \l}(G)\simeq Z$, supported in degrees $\geq \mb_G(\l)$.  It follows that $\pi_{\geq \l}(F\ot G)\simeq Y\ot Z$, which is supported in degrees $\geq \mb_F(\l)+\mb_G(\l)$. 

\end{proof}

The property $\nb_F(\l)=\mb_F(\l)$ seems to be quite special, hence we give it a name.

\begin{definition}\label{def:lsharp}
We say that a twist-like complex $F\in \KC^b(\AC)$ is \emph{$\l$-sharp} if $\nb_F(\l)=\mb_F(\l)$.  We call $F$ \emph{sharp} if it is $\l$-sharp for all $\l$.
\end{definition}

\begin{lemma} \label{lem:sharptensorclosed} If $F$ and $G$ are $\l$-sharp then $\nb_{G \ot F}(\l) = \nb_F(\l) + \nb_G(\l)$. Sharp, $\l$-sharp (for any given $\l$), and sharp increasing complexes are closed under tensor product. 
\end{lemma}

\begin{proof} If $\nb_F(\l)=\mb_F(\l)$ and $\nb_G(\l)=\mb_G(\l)$, then the inequalities from Lemma \ref{lemma:nmineq} give:
\[
\mb_{G\otimes F}(\l) \leq \nb_{G\otimes F}(\l)\leq  \nb_F(\l) +\nb_G(\l) = \mb_F(\l) + \mb_G(\l) \leq \mb_{G\otimes F}(\l),
\]
which forces $\mb_{G\otimes F}(\l) = \nb_{G\otimes F}(\l) = \nb_F(\l) + \nb_G(\l)$. \revise{Thus $G \ot F$ is $\l$-sharp. Moreover, if $F$ and $G$ are both $\l$-sharp and $\mu$-sharp for $\l < \mu$, and $\mb_F(\l) < \mb_F(\mu)$ and $\mb_G(\l) < \mb_G(\mu)$, then
\[ \mb_{G\otimes F}(\l) = \mb_F(\l) + \mb_G(\l) < \mb_F(\mu) + \mb_G(\mu) = \mb_{G \ot F}(\mu). \]
Thus if $F$ and $G$ are sharp increasing, then so is $G \ot F$.}
\end{proof}

Let us summarize and improve upon some of the ideas of Lemma \ref{lemma:nmineq}, in the special case that $F$ is $\l$-sharp.

\revcomment{Again, modifying the statement below to eliminate references to minimal complexes.}
\begin{lemma}\label{lemma:tail1}
Suppose $F$ is $\l$-sharp.  If $B$ is an indecomposable object in cell $\l$, then 
$F\ot B$ is homotopy equivalent to a complex $X$ of the form
\[
X = \cdots \rightarrow X^{\nb_F(\l)-2}  \rightarrow X^{\nb_F(\l)-1}  \rightarrow X^{\nb_F(\l)},
\]
where \begin{itemize}
\item $X^k$ is in cells strictly less than $\l$ for $k< \nb_F(\l)$, 
\item $X^{\nb_F(\l)}$ has exactly one indecomposable summand in cell $\l$, and the remaining summands are in cells $< \l$.
\end{itemize}
\end{lemma}

\begin{proof} The first paragraph of the proof of Lemma \ref{lemma:nmineq} implies that $X^k$ is in cells $<\l$ for $k < \mb_F(\l) = \nb_F(\l)$, and that $X^{\nb_F(\l)}$ is in cells $\le
\l$. It remains to prove that $X^{\nb_F(\l)}$ has exactly one indecomposable summand in cell $\l$.

\revcomment{Decided to reorganize this paragraph, hopefully it is clearer.}
Consider the quotient category $\AC/\IC_{< \l}$, which is still a monoidal category since $\IC_{< \l}$ is a two-sided tensor ideal. The identity morphism of $B$ is not in $\IC_{< \l}$ by the definition of cell ideals, so $B$ descends to a nonzero object $\overline{B}$ in the quotient. The endomorphism algebra of $B$ in $\AC$ was local, and remains local after any further quotient, so $\overline{B}$ is indecomposable.  Any invertible additive functor must preserve nonzero indecomposable objects. So $\overline{F} \ot \overline{B} = \overline{X}$ is indecomposable in $\KC^b(\AC/\IC_{<\l})$. But $\overline{X}$ is supported in a single homological degree $\nb_F(\l)$, so the chain object $\overline{X^{\nb_F(\l)}}$ must be indecomposable. This implies that $X$ had at most one indecomposable summand in cell $\l$ and degree $\nb_F(\l)$ to begin with.
\end{proof}

\begin{example}
To illustrate the statement of Lemma \ref{lemma:tail1}, take  $F=\HT_3$ and $\l=(2,1)$, so that $\nb_F(\l)=1$.  Note that $\HT_3\otimes B_s \simeq \underline{B_{sts}}(-1)\rightarrow B_{ts}$. The term in homological degree $1$ is in the same cell as $B_s$, and every term in smaller homological degree is in strictly smaller cells.

Moreover, the fact that $\HT_3$ is sharp implies that $\FT_3$ has the same property, by the last statement of Lemma \ref{lem:sharptensorclosed}.
\end{example}

We can improve slightly on Lemma \ref{lemma:tail1} when $F$ is increasing, but for this we need some additional assumptions on $\AC$.

\begin{hypothesis}\label{hyp:rigidetc}
Let $\AC$ be a Krull-Schmidt additive monoidal category. For each indecomposable object $B$ in $\AC$ there is some $B^\vee$ in $\AC$ such that $(-) \ot B$ and $(-) \ot B^\vee$ are both left and right adjoint to each other. Moreover, $B^\vee$ is in the same cell as $B$. Finally, the cells are \emph{nondegenerate} in that, for each indecomposable $B$, there is some indecomposable $B'$ in the same cell as $B$ such that $B$ is isomorphic to a direct summand of $B \ot B'$. \end{hypothesis}

\begin{remark} \label{rmk:rigidmonoidal} \revcomment{Edits to this remark} These hypotheses hold for $\SBim$. The category $\SBim$ is pivotal, implying that an object $B^{\vee}$ exists for each $B$. Since $B_s^\vee \cong B_s$, one can deduce that $B_w^\vee \cong B_{w\inv}$. That $w$ and $w\inv$ are in the same cell holds by Proposition \ref{prop:someP}, (P14). The object $B'$ can be taken to be $B_d$, where $d$ is the distinguished involution in the same left cell as $B$, see Proposition \ref{prop:dacts}. \end{remark}
	
%


\begin{proposition}\label{prop:tailProperty}
Assume Hypothesis \ref{hyp:rigidetc}. Let $F$ be twist-like, sharp, and increasing. If $B$ is any indecomposable object in cell $\l$ and $X\simeq F\otimes B$ is the minimal complex, then
\begin{enumerate}
\item if $k<\nb_F(\l)$ then the chain object $X^k$ is in cells strictly less than $\l$.
\item if $k=\nb_F(\l)$ then the chain object $X^k$ is an indecomposable object in cell exactly $\l$.
\end{enumerate}
\end{proposition}


\begin{proof} The statement (1) is part of Lemma \ref{lemma:tail1}. Now we prove (2). By Lemma \ref{lemma:tail1}, if (2) fails then $X^{\nb_F(\l)}$ has a direct summand of the form $C$, where $C$ is an indecomposable object in cell $\mu < \l$. The inclusion of $C$ induces a chain map $C[-\nb_F(\l)]\rightarrow X$, since $X$ is zero in homological degrees $>\nb_F(\l)$. This will also be thought of as a map $j:C\rightarrow X[\nb_F(\l)]$. We claim that $j$ is null homotopic. If we can prove this, then this would imply that the component of the differential $X^{\nb_F(\l)-1}\rightarrow C$ is projection to a direct summand (with a splitting provided by the null-homotopy for $j$). Then $C$ can be cancelled by a Gaussian
elimination, contradicting the fact that $X$ is a minimal complex.

It remains to prove that $j:C\rightarrow X[\nb_F(\l)]$ is null-homotopic.  By hypothesis \ref{hyp:rigidetc} , there is an object $C'$ in cell $\mu$ such that $C$ is isomorphic to a direct summand of $C \otimes C'$.  Taking hom complexes, it follows that $\Homc_{\AC}(C,X)\simeq \Homc_{\AC}(C, F\otimes B)$ is isomorphic to a direct summand of
\[
\Homc_{\AC}(C\otimes C', F\otimes B)\cong \Homc_{\AC}(C, F\otimes B \ot (C')^\vee)
\]
Meanwhile, each indecomposable summand of $B \ot (C')^\vee$ is in a cell $\le \mu < \l$. Thus $F\otimes B \ot (C')^\vee$ is homotopy equivalent to a complex supported in degrees $\leq \nb_F(\mu)<\nb_F(\l)$. It follows that the homology of $\Homc_{\AC}(C, F\otimes B)$ is supported in homological degrees $<\nb_F(\l)$.  In particular $j$ is null-homotopic.  This completes the proof.
\end{proof}

\begin{definition}\label{def:head} Assume the hypotheses and notation of Proposition \ref{prop:tailProperty}. The maximal chain object $X^{\nb_F(\l)}[-\nb_F(\l)]$ (viewed as a complex
living in its usual homological degree) will be referred to as the \emph{head} of $F\otimes B$. The \emph{tail} of $F\otimes B$ is by definition the truncated complex $\cdots \rightarrow
X^{\nb_F(\l)-2}\rightarrow X^{\nb_F(\l)-1}$. \end{definition}

\begin{remark} When $F$ is sharp but not increasing, one may wish to define the head of $F \ot B$ as the unique indecomposable summand of $X^{\nb_F(\l)}$ living in cell $\l$, guaranteed by Lemma \ref{lemma:tail1}. However, this summand need not be canonically defined (though it is canonically defined in the quotient category $\AC / \IC_{< \l}$).  \revcomment{Removed the sentence: For simplicity we restrict to increasing  twist-like complexes.}
\end{remark}

Each of the results in this section has an analog where the homological degrees are reversed, as in Remark \ref{rmk:frombelow}. We call $F$ \emph{co-sharp} if $\nb^\ast_F(\l) =
\mb^\ast_F(\l)$ for all $\l$. As noted above $\HT_3$ is sharp and increasing, while $\HT_3\inv$ is co-sharp and decreasing.

\begin{lemma} If $F$ is sharp and $F\inv$ is co-sharp, then $\nb^\ast_{F\inv}(\l) = -\nb_F(\l)$. Thus if $F$ is sharp and increasing and $F\inv$ is co-sharp, then $F\inv$ is also
decreasing. \end{lemma}

\begin{proof} Fix $B$ in cell $\l$. By Lemma \ref{lemma:tail1} we see that, modulo lower cells, $F \ot B$ is a single indecomposable in cell $\l$ and homological degree $\nb_F(\l)$.
Thus, by the same argument, $F\inv \ot (F \ot B)$ is (modulo lower cells) a single indecomposable in cell $\l$ and homological degree $\nb_F(\l) + \nb^\ast_{F\inv}(\l)$. But since this
tensor product is just $B$ in homological degree $0$, we must have $\nb_F(\l) + \nb^\ast_{F\inv}(\l) = 0$. \end{proof}

\begin{remark} Recall that $\nb_F(\l)$ and $\nb_F^\ast(\l)$ were defined in terms of $F\otimes B$ when $B\in\AC$ is in cell $\l$.  We could also have made similar definitions using $B\otimes F$, which we temporarily denote $\nb_{F,R}(\l)$ and $\nb_{F,R}^\ast(\l)$.  If $\AC$ satisfies Hypothesis \ref{hyp:rigidetc}, then taking right adjoints gives a functor $B \mapsto B^\vee$, and this functor extends to a functor on $\KC(\AC)$ which reverses homological degree. Since the adjoint of an invertible operator is also its inverse, this functor sends $F$ to $F\inv$ and $F \ot B$ to $B^\vee \ot F\inv$. Consequently, $\nb_F(\l) = \nb^\ast_{F\inv,R}(\l)$ and $\nb_F^\ast(\l) = \nb_{F\inv,R}(\l)$.

Moreover, the category of Soergel bimodules has a covariant (!) autoequivalence which swaps the order of tensor products and preserves cells, namely $B_w\mapsto B_{w\inv}$. Applying this equivalence we see that $\nb_F(\l) = \nb_{F,R}(\l)$ for twist-like complexes in $\KC^b(\SBim_n)$.
\end{remark}

\subsection{Cell triangularity: prelude and conventions}
\label{subsec:celltriprelude}

We continue to assume that $\AC$ is an additive monoidal category with the Krull-Schmidt property.   

A sharp twist-like complex $F$ \revise{acts} in a block upper-triangular fashion with respect to cells. Namely, an indecomposable object $B$ in cell $\l$ is sent to $B'[-\nb_F(\l)]$ for
some other indecomposable object $B'$ in the same cell $\l$, modulo $\IC_{< \l}$. However, the indecomposables in a given cell can be permuted in an interesting way by $F$. 

\begin{example} Let $F = \HT_3$ and recall the examples in \S\ref{subsec:HT3}. For $B = \one$ we have $B' = \one(3)$; for $B = B_s$ we have $B' = B_{ts}(0)$; for $B = B_{ts}$ we have $B'
= B_s(0)$; for $B = B_{sts}$ we have $B' = B_{sts}(-3)$. \end{example}

\begin{example} Let $F = \FT_3$. Then $B'$ agrees with $B$ up to a grading shift which only depends on $\l$! \end{example}

\revcomment{From here, major revisions. Moved grading conventions here.} The situation for $\FT_3$ is much nicer than for $\HT_3$, because the way in which $\FT_3$ acts
on cells (in the associated graded) is encapsulated by a grading shift functor. In order to make general definitions (beyond the case of graded modules over a ring), let us briefly
discuss some formalities involving graded categories, taking some material from \cite[Proposition 11.2.1]{EMTW}. 

Let $\Gamma$ be an abelian group. Consider an additive category $\AC'$ where the morphism spaces are $\Gamma$-graded abelian groups, and composition is degree preserving (meaning $\deg(f\circ g) = \deg(f)+\deg(g)$). To this
category one can associate a new additive category $\AC$, where the objects are formal direct sums of formal shifts $A(i)$, where $A \in \AC'$ and $i \in \Gamma$, and such that
$\Hom_{\AC}(A(i),B(j))$ agrees with the degree $j-i$ part of the morphism space $\Hom_{\AC'}(A,B)$. There is a strict action of $\Gamma$ on $\AC$ by formal shift functors $(i)
\colon \AC \to \AC$. These shifts satisfy natural compatibilities by construction, e.g. one has strict equalities of functors $(i) \circ (j)= (i+j)$ and $(0)= \Id_{\AC}$. When $\AC'$ is monoidal, it is easy to equip $\AC$ with an monoidal structure for which $A(i) \ot B(j) := (A \ot B)(i+j)$. When $\AC$ arises by this construction, we say that $\AC$ is equipped with a strict $\Gamma$ action. Note that letting $\AC=\SBim$ with $\Gamma=\Z$ gives an
instance of the above situation.

It is the category $\AC$ we typically work with, so that $\Hom$ by default represents a space of degree zero morphisms. The original Hom space in $\AC'$ will be denoted as
\begin{equation}\label{eq:graded hom space}
\Hom^\Gamma_{\AC}(A,B):=\bigoplus_{i\in \Gamma}\Hom_{\AC}(A,B(i)).
\end{equation}

In this context, the functor $(i) \colon \AC \to \AC$ for $i \in \Gamma$ we will regard as an \emph{invertible scalar functor}. Any
direct sum of such functors will be called a scalar functor. When $\AC$ is monoidal, $A(i) \cong A \ot \one(i) \cong \one(i) \ot A$, so that $(i) \cong \one(i) \ot (-) \cong (-)
\ot \one(i)$ as functors. The object $\one(i)$ is called an \emph{invertible scalar object}. It satisfies the requirements to be a scalar object as specified in \cite[\S
6.1]{ElHog17a}.

For any additive category, the homotopy category $\KC(\AC)$ has a strict $\Z$ action coming from the homological shift. That is, there is a category with $\Z$-graded morphism spaces 
\[
\Hom^\Z_{\KC(\AC)}(C,D) = H^{\bullet}(\Homc_{\KC(\AC)}(C,D)).
\]
If $\AC$ has the structure of a strict $\Gamma$-action, then $\KC(\AC)$ has the structure of a strict $\Gamma\times \Z$ action, and a shift like $(i)[j]$ for $i \in \Gamma$ and $j \in \Z$ is a valid scalar functor. Thus we would write $\Hom^{\Gamma \times \Z}$ for the bigraded Hom space.

\subsection{Cell triangularity}
\label{subsec:celltri}

Now we resume our study of twist-like complexes, assuming in addition that $\AC$ has a strict $\Gamma$ action for an abelian group $\Gamma$.

\begin{defn} \label{def:twistedunitri}
\revadd{Let $F \in \KC^b(\AC)$ be twist-like, sharp, and increasing. We say that $F$ is \emph{cell triangular} if for each cell $\l$ there is a functor $\Sigma_\l= (a_\l)[b_\l]$, for some $a_{\l} \in \Gamma$ and $b_{\l} \in \Z$,} such that if $B$ is in cell $\l$ then the head of $F \ot B$ is isomorphic to $\Sigma_\l(B)$.
\end{defn}

\revadd{Note that in the above definition, $b_\lambda$ is forced to equal $-\nb_F(\lambda)$.}

\begin{defn}\label{def:lequiv}
Let $F$ be cell triangular. We say that a chain map $\a_\l \co \Sigma_\l(\one) \to F$ is a \emph{$\l$-equivalence} if $\Cone(\a_\l) \ot B$ is homotopy equivalent to a complex in cells $<\l$, \revise{for all $B$ in cell $\l$}.
\end{defn}

\revadd{Equivalently, $\a_\l$ is a $\l$-equivalence if $\Cone(\a_\l)\ot B$ becomes contractible in the cellular subquotient $\KC^b(I_{\leq \l}/I_{<\l})$.  In case $F=\FT_n\in \KC^b(\SBim)$ this is equivalent to the definition of $\l$-equivalence given in the introduction (Definition \ref{defn:introlequiv}).}

The $\l$-equivalence property can be reinterpreted by analyzing precisely what happens in homological degree $\nb_F(\lambda)$.  If $B$ is in cell $\l$ then $F\otimes B$ is homotopy equivalent to a minimal complex of the form
\begin{equation} \label{refereeX}
X := \cdots \rightarrow X^{\nb_F(\l)-2}\rightarrow X^{\nb_F(\l)-1}\rightarrow X^{\nb_F(\l)}.
\end{equation}
where $X^{\nb_F(\l)}[-\nb_F(\l)]=\Sigma_\l(B)$. \revise{Recall the definition of the head and tail of $X$ from Definition \ref{def:head}.} Let $\iota_B \co \Sigma_\l(B) \to X$ denote the inclusion of the head of $X$.

\begin{lemma} \label{lemma:onlyiota}
Let $F$ be cell triangular. Let $B\in \AC$ be an indecomposable object in cell $\l$. Then $\Hom(\Sigma_{\l}(B),F \ot B)$ (the morphism space in the homotopy category) is a nonzero cyclic right module for $\End(B)$. \end{lemma}

\begin{proof} 
Let $\KM := \End(B)$, a local ring, and let $X$ be the minimal complex of $F \ot B$ \revise{as in \eqref{refereeX}}. Note that $\KM \cong \End(\Sigma_{\l}(B))$ since $\Sigma_{\l}$ is an invertible functor, so $\Hom(\Sigma_{\l}(B),F \ot B)$ is a right $\KM$-module. It is enough to show that the isomorphic module $\Hom(\Sigma_{\l}(B),X)$ is cyclic, spanned by $\iota_B$ (as defined above). Now $\Sigma_\l(B)$ lives in homological degree $\nb_F(\l)$ and $X$ lives in homological degrees $\leq \nb_F(\l)$. Hence any chain map $\psi\in \Hom_{\Ch(\AC)}(\Sigma_\l(B), X)$ is supported in degree $\nb_F(\l)$ and has the same data as an endomorphism of $\Sigma_\l(B)$. Evidently, this space of chain maps is cyclic over $\KM$, spanned by $\iota_B$. If $\iota_B$ were null-homotopic then any homotopy would factor through the tail of $X$, which lives in cells $<\l$. However, the identity of $B$ is not in cells $< \l$, by the definition of cells \revise{(see
\S\ref{subsec:celldef})}. Thus the $\KM$-submodule of nulhomotopic maps is proper. \end{proof}

\begin{lemma} \label{lem:lequivcriterion} Let $F$ be cell triangular. A chain map $\a_\l : \Sigma_{\l}(\one) \to F$ is a $\l$-equivalence if and only if the degree $\nb_F(\l)$ component (which is the only nonzero component) of the chain map $\pi_B \circ (\a_\l \ot \Id_B)$ is an isomorphism, where $\pi_B$ is a homotopy equivalence from $F \ot B$ to its minimal complex $X$ \revise{(as in \eqref{refereeX})}, for all indecomposable objects $B$ in cell $\l$. When $\End(B)$ is a field, this is equivalent to $\a_\l \ot \Id_B \ne 0$.  \end{lemma}

\begin{proof} Tensoring the map $\a_\l$ with the identity map of $B$ gives a chain map
\[
\a_\l \ot \Id_B \co \Sigma_\l(\one)\otimes B\cong \Sigma_\l(B) \rightarrow F\otimes B,
\]
such that $\Cone(\a_\l \ot \Id_B) \cong \Cone(\a_\l)\otimes B$. By the previous lemma and its proof, the chain map $\pi_B \circ (\a_\l \ot \Id_B)$ is $\iota_B \cdot k$ for some $k \in \KM$. If $k$ is invertible, then the two copies of $\Sigma_\l(B)$ can be removed via Gaussian Elimination from $\Cone(\a_\l \ot \Id_B)$, leaving only the tail of $X$ which is supported in lower cells, as desired for a $\l$-equivalence. If not, $k$ is in the maximal ideal of $\KM$, and $\Cone(\pi_B \circ (\a_\l \ot \Id_B))$ is a minimal complex. In this case the cell $\l$ is there to stay, so $\a_{\l}$ is not a $\l$-equivalence.

When $\KM$ is a field, a nonzero cyclic right module is free. Thus $k$ is invertible if and only if it is nonzero, if and only if $\pi_B \circ (\a_\l \ot \Id_B)$ is nonzero. Since $\pi_B$ is a homotopy equivalence, this holds if and only if $\a_\l \ot \Id_B$ is nonzero.
\end{proof}

\revcomment{Revisions to this section mostly end here}

Note that if $\a_\l:\Sigma_\l(\one)\rightarrow F$ is a $\l$-equivalence and $\b$ is any map homotopic to $\a_\l$, then $\b$ is a $\l$-equivalence since $\Cone(\b)\cong \Cone(\a_\l)$.



	


If one has a family of $\l$-equivalences for each cell $\l$, then $F$ is very close to being categorically prediagonalizable \cite[Definition 6.13]{ElHog17a}. 

\begin{prop}\label{prop:lequivmeansdiag} Suppose that $\AC$ is a category with finitely many cells, and $F \in \KC^b(\AC)$ is positive twist-like, sharp, and cell triangular. If there exists a $\l$-equivalence $\a_\l$ for each cell $\l$, then
\[
\bigotimes_{\l}\Cone(\a_\l)\simeq 0
\]
for a particular ordering of the tensor factors (see the proof).
\end{prop}

\begin{remark}\label{rmk:obstructions} Recall that $\Cone(\a) \ot \Cone(\b)$ and $\Cone(\b) \ot \Cone(\a)$ need not be homotopy equivalent, so that the order on the tensor product matters. In \cite[\S 6.2]{ElHog17a}, the authors define certain obstructions, the vanishing of which guarantees that the cones tensor commute up to homotopy.  If the obstructions vanish, then the collection of maps $\{\a_\l\:|\: \text{$\l$ is a 2-sided cell}\}$ is said to be \emph{obstruction-free}.\end{remark}


\begin{cor} \label{cor:lequivmeansdiag} If the $\l$-equivalences of Proposition \ref{prop:lequivmeansdiag} are obstruction-free, then $F$ is categorically prediagonalizable. \end{cor}

\begin{proof}[Proof of Proposition \ref{prop:lequivmeansdiag}] Let $\PC$ denote the finite poset of two-sided cells of $\AC$.  We need only show that \[ T := \bigotimes_{\l \in \PC} \Cone(\a_\l) \simeq 0
\] for some ordering on the tensor product.
	
For any ideal $J$ in the poset $\PC$ (e.g. $J = \{\le \mu\}$ or $J = \{\ngeq \l\}$), let $\IC_J\subset \AC$ denote the two-sided tensor ideal whose indecomposables are those in the cells
of $J$. Let $\ess \KC^b(\IC_J)\subset \KC^b(\AC)$ denote the full subcategory of (complexes which are homotopy equivalent to) complexes whose terms are in $\IC_J$.

Tensoring with $\Cone(\a_\l)$ will preserve $\ess \KC^b(\IC_J)$ for all ideals $J$, and will send $\ess \KC^b(\IC_{\le \l})$ to $\ess \KC^b(\IC_{< \l})$. Thus if $J_1 \subsetneq J_2$ are two ideals with $J_2 = J_1 \cup \{\l\}$, then tensoring with $\Cone(\a_\l)$ will send $\ess\KC^b(\IC_{J_1})$ to $\ess\KC^b(\IC_{J_2})$.

Pick any filtration $\emptyset = J_0 \subset J_1 \subset \cdots \subset J_k = \PC$ by ideals, with $J_n \setminus J_{n-1} = \{\l_n\}$ a singleton, and consider the tensor product $T$
with the order given by this filtration, so that $\l_k$ appears on the far right and $\l_1$ on the far left. Then tensor product with $T$ will send $\KC^b(\AC) = \KC^b(\IC_{J_k})$ to
$\KC^b(\IC_{J_0}) = 0$. Thus $T \ot \one \simeq 0$, and therefore $T \simeq 0$, as desired. \end{proof}

\subsection{Recap}
\label{subsec:recap}

We conjecture in this paper that the full twists $\FT$ are categorically diagonalizable, and in particular are categorically prediagonalizable. Our method of proof in type $A$ will begin by proving that $\FT_n$ is positive twist-like, sharp, and cell triangular. This will require a deep dive into the cell theory of Hecke algebras, which we discuss in the following section for arbitrary Coxeter groups. Then we will construct $\l$-equivalences for each cell $\l$. This will be done by computing the space of all maps from $\Sigma(\one)$ to $\FT_n$ for
invertible scalars $\Sigma$, and picking one with satisfactory properties, which is the content of \S\ref{sec:constructing}. An analogous computation shows that a certain obstruction space is trivial, so the $\l$-equivalences are obstruction-free. Then we will apply Corollary \ref{cor:lequivmeansdiag} to prove prediagonalizability.

\section{The asymptotic Hecke algebra and the eigenvalues of the full twists}
\label{sec:cellsplus}

A general reference for Hecke algebras and their cell theory is Lusztig's book on Hecke algebras with unequal parameters \cite{LuszUnequal14}, which has been updated on the
arXiv to account for the results of \cite{EWHodge}.

Lusztig has developed a beautiful theory for studying Hecke algebras in terms of their ``asymptotic'' behavior. When the parameter $v$ of the Hecke algebra $\HB$ tends towards $0$, only the smallest power of $v$ in a polynomial should survive. By examining the smallest powers of $v$ which appear in certain key coefficients, Lusztig extracts interesting numerical information. The relationship between this asymptotic behavior and the cells of $\HB$ is encapsulated in Lusztig's famous conjectures P1-P15, see \cite[Section 14]{LuszUnequal14}. Lusztig proved (see \cite[Section 15]{LuszUnequal14}) that when the structure coefficients $c^z_{x,y}$ (to be defined below) and the Kazhdan-Lusztig polynomials $h_{y,w}$ are always positive (that is, in $\NM[v,v^{\inv}]$) then conjectures P1-P15 are true. For Weyl groups (such as the symmetric group), this positivity was known since the proof of the Kazhdan-Lusztig conjectures. For arbitrary Coxeter groups, this positivity is also known thanks to \cite{EWHodge}.\footnote{Conjectures P1-P15 were also made for Hecke algebras with unequal parameters, where they remain open in general.}

Again, the reader new to this theory should read \S\ref{sec:typeAcells} concurrently, to see these results in the more familiar situation of type $A$. In type $A$, the two-sided cells of
$\HB(S_n)$ are in bijection with partitions of $n$.

Let us note a major notational change with the literature. What is commonly called Lusztig's $\ab$-function we will be denoting by $\rbb$, and calling the $\rbb$-function. In type $A$,
$\rbb$ corresponds to the row function of a partition. We will eventually introduce other functions $\cbb$ and $\xbb$ derived from $\rbb$, which correspond in type $A$ to the column and
content functions of partitions. To keep notation consistent with type $A$, we have departed from tradition and used $\rbb$ for the general case.

\subsection{The $\rbb$-function}
\label{subsec:afunc}

Henceforth, we will be studying the cell theory (as defined in \S\ref{subsec:algCells}) of the KL basis of $\HB(W)$. Thus, one has three transitive relations $\le_{L}$, $\le_{R}$, and $\le_{LR}$ on
the set $W$, whose equivalence classes are called \emph{left}, \emph{right}, and \emph{two-sided cells} respectively. We write $x \sim_L y$ if $x \le_L y$ and $y \le_L x$, and so forth.
Unless otherwise specified, a \emph{cell} will refer to a two-sided cell.

\begin{lemma} For $x, y \in W$ one has $x \le_L y$ (resp. $x \le_{LR} y$) if and only if $x^{-1} \le_R y^{-1}$ (resp. $x^{-1} \le_{LR} y^{-1}$). \end{lemma}

\begin{proof} This follows immediately from the existence of an algebra antiautomorphism on $\HB$ which sends $H_x\mapsto H_{x\inv}$ and $b_x\mapsto b_{x^{-1}}$. \end{proof}
	
We now recall Lusztig's $\rbb$-function $W \to \NM$, see \cite[Section 13]{LuszUnequal14}.

\begin{defn}\label{def:a} For $x,y,z\in W$, let $c^z_{x,y} \in \Z[v,v^{\inv}]$ be the structure coefficients\footnote{Denoted $h_{x,y,z}$ in \cite[Section 13]{LuszUnequal14}.} of the
Kazhdan-Lusztig basis, that is, $$b_x b_y=\sum_{z}c^z_{x,y} b_z.$$ Note that the polynomials $c^z_{x,y}$ are preserved by swapping $v$ and $v^{\inv}$.

Define $\rbb(z) \in \NM$ such that $v^{-\rbb(z)}$ is the smallest power of $v$ (and $v^{+\rbb(z)}$ the largest power) appearing in $c^z_{x,y}$ for all $x,y \in W$. Let $t^z_{x,y}$\footnote{Denoted $\gamma_{x,y,z^{-1}}$ in \cite[Section 13]{LuszUnequal14}.} denote the coefficient of $v^{-\rbb(z)}$ inside $c^z_{x,y}$.  \end{defn}

The following facts were proven by Lusztig, see \cite[Section 14]{LuszUnequal14}. Note that $t^z_{x,y}$ is zero unless $x$ and $y$ succeed in minimizing the power of $v$ which appears in $c^z_{x,y}$.

\begin{prop} \label{prop:someP} The following properties hold. \begin{itemize}
\item (P14) For all $w$, $w \sim_{LR} w^{-1}$.
\item (P8) If $t^z_{x,y} \ne 0$, then $x \sim_R z$, $y \sim_L z$, and $x \sim_L y\inv$. In particular, they are all in the same two-sided cell.
\item (P7) For any $x, y, z$, $t^z_{x,y} = t^{x\inv}_{y,z\inv}$.
\item (P4) If $x \le_{LR} y$ then $\rbb(x) \ge \rbb(y)$. In particular, $\rbb(x) = \rbb(\l)$ only depends on the two-sided cell $\l$ of $x$.
\item (P11) If $\l \ne \mu$ are two-sided cells with $\l < \mu$, then $\rbb(\l) > \rbb(\mu)$.
\item (P9-10) If $x \le_L y$ and $\rbb(x) = \rbb(y)$ then $x \sim_L y$. Consequently, two comparable left cells in the same two-sided cell are equal. The same holds, replacing left with right.
\end{itemize} \end{prop}

Because of (P4), we will mostly think of $\rbb$ as a function on two-sided cells, rather than a function on $W$. Because of (P8), one can reinterpret the function $\rbb$ as follows: for a
two-sided cell $\l$, $v^{-\rbb(\l)}$ is the minimal power of $v$ appearing in $c^z_{x,y}$ for any $x, y, z$ in cell $\l$.

\begin{example} \label{ex:boundviolated}
Let $\{s,t,u\}$ be the simple reflections of $S_4$, \revise{as in Example \ref{ex:S4HTetc}.} For this example we focus on $x = tsut$ and $w = sutsu$, which happen to be the two involutions in $S_4$ which are not longest elements of parabolic subgroups, and also happen to be the only non-smooth elements.

The element $x$ lies in cell $\mu = (2,2)$, and it turns out that $\rbb(\mu)=2$. One has \begin{equation} \label{eq:x2} b_x b_x = [2]^2 b_x + [2]^2 b_{w_0} + b_{tst} + b_{tut} + b_{tstut}
+ b_{tutst}.\end{equation} Aside from $[2]^2 b_x$, all the other terms lie in lower cells. Thus $t^x_{x,x} = 1$, while $c^y_{x,x} = 0$ for other elements $y$ in this cell.

The element $w$ lies in the cell $\l = (2,1,1)$, and $\rbb(\l) = 3$. One has \begin{equation} \label{eq:w2} b_w b_w = [2]^3 b_w + [2]^4 b_{w_0}.\end{equation} The term $b_{w_0}$ lies in a
lower cell. Thus $t^w_{w,w} = 1$. Note that $c^{w_0}_{w,w}$ has minimal power $v^{-4}$, illustrating that $-\rbb(\l)$ is not the minimal power of $v$ appearing in $c^z_{x,y}$ for $x, y$
in cell $\l$, but also requires $z$ in $\l$. 
\end{example}

To emphasize a takeaway from this example, it is important to remember that the bound $-\rbb(\l)$ on the $v$-degrees of summands is only valid within the same cell, and that lower cells
may also reach or even exceed this bound.

The \emph{asympotic Hecke algebra} or \emph{J-ring} is morally obtained by multiplying KL basis elements, but only remembering the minimal $v$-degree which ``should" appear.

\begin{prop} \label{prop:Jassoc} Fix a two-sided cell $\l$. Let $J_{\l}$ be the algebra with basis $\{j_x\}_{x \in \l}$ and with multiplication $j_x j_y = \sum_{z \in \l} t^z_{x,y} j_z$. Then $J_\l$ is an associative algebra (possibly without unit when $W$ is infinite, see Remark \ref{rmk:Junit}). \end{prop}

\begin{remark} One could also define the ring $J$ with a basis $\{j_x\}_{x \in W}$ and multiplication as above. Because of (P8), this would just be the product of the rings $J_\l$ for each cell. \end{remark}

\begin{ex} Consider $S_3$ with simple reflections $\{s,t\}$, and its simple cell $\l$ consisting of $\{s,t, st, ts\}$. Then $\rbb(\l)=1$, as can be seen from the products $b_s b_s =
(v+v\inv)b_s$ and $b_s b_t = b_{st}$ (and similarly, swapping $s$ and $t$). The equation $b_sb_s=(v+v\inv)b_s$ implies that in the $J$ ring $j_s j_s = j_s$, while $b_s b_t=b_{st}$
implies $j_s j_t=0$. In particular $j_s$ and $j_t$ are orthogonal idempotents in the $J$-ring. The elements $b_{st}$ and $b_{ts}$ satisfy. \[ b_{st}b_{ts}=b_sb_tb_tb_s =
(v+v\inv)b_sb_tb_s=(v+v\inv)(b_{sts}+b_s). \] This implies that $j_{st} j_{ts} = j_s$. Hence the ring $J_\l$ is a $2\times2$ matrix algebra. In fact, for any cell $\l$ in type $A$, the
$J$-ring is a matrix algebra (see Lemma \ref{lem:howPQmult}), though this does not generalize to arbitrary Coxeter groups. \end{ex}

\subsection{The $\Delta$-function}
\label{subsec:Dfunc}

Now we recall Lusztig's Delta function, a function $W \to \NM$ which does not only depend on the cell.

\begin{defn}\label{def:D} For $x \in W$, define $\D(x) \in \NM$ such that $v^{\D(x)}$ is the smallest power of $v$ appearing in the Kazhdan-Lusztig polynomial $h_{1,x} \in \ZM[v]$.
\end{defn}	

Note that $\D(x) = \D(x\inv)$, because $h_{1,x} = h_{1\inv,x\inv} = h_{1,x\inv}$.

In the next section, we discuss the connection between $\D$ and the half twist. The first use of $\D(x)$ is to determine, by comparison with $\rbb(x)$, a special set of elements of $W$. Here are more facts proven by Lusztig, see \cite[Section 14]{LuszUnequal14}.

\begin{defn}\label{def:duflo}
Let $\DC\subset W$ denote the set of elements $x$ such that $\Delta(x)=\rb(x)$.  An element of $\DC$ is called a \emph{distinguished involution} or \emph{Duflo involution}.
\end{defn}

The proposition below implies that distinguished involutions are in fact involutions.  We often use the letter `$d$' to denote distinguished involutions.

\begin{prop}\label{prop:moreP} The following properties hold. \begin{itemize}
\item (P1) For all $x \in W$, $\D(x) \ge \rbb(x)$.
\item (P13,P6) If $d\in \DC$ then $d$ is an involution.  Each left cell (resp. right cell) contains a unique element of $\DC$. 
\item (P5 and positivity) If $d \in \DC$ then the coefficient of $v^{\D(d)}$ in $h_{1,d}$ is $1$.
\item (P3,P13,P5) If $d \in \DC$ then $t^d_{x,y} \ne 0$ if and only if $x = y\inv$ and $y \sim_L d$. Moreover, in this case $t^d_{y\inv,y} = 1$.
\end{itemize}
\end{prop}

In type $A$ each left cell contains a unique involution, so it must be distinguished, and thus $\DC$ is the set of all involutions.

\begin{ex} In any dihedral group with simple reflections $\{s,t\}$, there are three cells: the identity, longest, and simple cells, as in Example \ref{ex:threecells}. Within the simple
cell there are two left cells, one containing those elements whose unique reduced expression ends in $s$, and the other those ending in $t$. The simple cell $\l$ satisfies $\rbb(\l)=1$.
All elements are smooth, thus $\D(x) = \ell(x)$. Hence the simple cell contains two distinguished involutions of length $1$, namely $s$ and $t$. Outside of the special case of type
$A_2$, both of these left cells contain multiple involutions, but only $s$ and $t$ are distinguished. \end{ex}

We state a useful corollary which gives a practical way to compute $\rbb(\l)$.

\begin{cor}\label{cor:computinga}
If $\l$ is a two-sided cell in $W$, then $\rbb(\l)$ can be characterized as the minimal power of $v$ occuring in $c^d_{d,d}$, for any distinguished involution $d \in \DC \cap \l$.
\end{cor}


\begin{ex} \revcomment{minor revisions here to clarify the purpose of the example} Note that a ``converse'' to the above corollary does not hold: distinguished involutions are not the only elements $x$ (or even the only involutions) for which $c^x_{x,x}$ has minimal power $v^{-\rbb(\l)}$. For example, consider a dihedral group with simple reflections $\{s,t\}$ with $m_{st} \ge 6$. Then $j_{sts} j_{sts} = j_s + j_{sts} + j_{ststs}$. Even as we proceed to study and use the polynomials $c^x_{d,x}$ for distinguished involutions $d$, we emphasize that distinguished involutions $d$ are characterized instead by the behavior of $h_{1,d}$. \end{ex}

\subsection{The action of distinguished involutions}
\label{subsec:invact}

Combining the results above, Lusztig proves the following.

\begin{prop} \label{prop:dacts} If $d \in \DC$ and $x \in W$ satisfies $x \sim_{LR} d$, then $t^y_{d,x} = 0$ for all $y$ unless $x \sim_R d$, in which case $t^x_{d,x} = 1$ and
$t^y_{d,x} = 0$ for $y \ne x$. In other words, left multiplication by $j_d$ in the $J$-ring is the same as projection to those $j_x$ for which $x$ in the same right cell as $d$.
Similarly, right multiplication by $j_d$ is projection to the left cell of $d$. \end{prop}

\begin{proof} This is proven by combining (P7) and (P8) from Proposition \ref{prop:someP} with the last property from Proposition \ref{prop:moreP}. \end{proof}
	
\begin{remark} \label{rmk:Junit} It follows that if a two-sided cell $\l$ contains finitely many left cells, then the finite sum $\sum_{d \in \DC \cap \l} j_d$ is the identity element of $J_{\l}$. Even if
there are infinitely many left cells, $J_{\l}$ is locally unital for the orthogonal idempotents $\{j_d\}_{d \in \DC}$. \end{remark}	

A restatement of this proposition says that for $x$ and $d$ in cell $\l$, with $d$ a distinguished involution, we have
\begin{equation} \label{eq:dacts} v^{\rbb(d)} b_d b_x \equiv \d_{x \sim_R d} b_x + \HB^+_{\l} + \HB_{< \l}, \end{equation}
where $\HB_{< \l}$ denotes the ideal spanned by $b_y$ for $y$ in a strictly smaller cell than $\l$, and $\HB^+_\l$ denotes those linear combinations of $b_y$ for $y$ in cell $\l$ whose coefficients live in $v\Z[v]$.

A categorical analogue of this proposition is the following.

\begin{cor}\label{cor:weakJringCat}
Let $d\in W$ be a distinguished involution in cell $\l$. If $x \sim_R d$ then $B_x(0)$ is a direct summand of $B_d(\rbb(\l)) \ot B_x$ with multiplicity one, and all other summands in cell $\l$ have strictly positive grading shifts. If $x \sim_{LR} d$ but $x \nsim_R d$ then every summand of $B_d(\rbb(\l)) \ot B_x$ has a strictly positive grading shift. \end{cor}

\begin{proof} This follows directly from the decategorified statement. \end{proof}

\subsection{First consequences for the half twist}
\label{subsec:halftwistconsequences}

The following is a result of Lusztig, see \cite[Corollary 11.7]{LuszUnequal14}.

\begin{prop} \label{prop:w0oncells} Let $W$ be a finite Coxeter group with longest element $w_0$. Then $x \le_L y$ if and only if $w_0 x \ge_L w_0 y$ if and only if $x w_0 \ge_L y w_0$.
Similar statements hold replacing $L$ with $R$ or $LR$. Thus, multiplication by $w_0$ induces an order-reversing involution on the set of left cells, right cells, and two-sided cells.
\end{prop}

\begin{lemma} Given a two-sided cell $\l$, the two-sided cells $w_0 \l$ and $\l w_0$ are equal.\footnote{This statement is well-known, and was known to Lusztig in the early days of cell theory. We could not find an early reference in the literature, however. The proof below was told to the first author by Victor Ostrik.} Equivalently, two-sided cells are preserved by conjugation by the longest element. Moreover $w_0 \DC = \DC w_0$ as sets, or equivalently, $w_0 d w_0$ is distinguished whenever $d$ is. \end{lemma}

\begin{proof} A two-sided cell gives rise to a left-module of $\HB$, by taking the associated graded in the cell filtration of $\HB$. \revise{The modules associated to distinct two-sided cells will have different annihilators\footnote{The two-sided cell module for $\l$ is annihilated by $\HB_{\ngeq \l}$, since $\HB_{\ngeq \l} \cap \HB_{\le \l} \subset \HB_{< \l}$.}, so they are non-isomorphic representations of $\HB$.} Because of Lemma \ref{lem:ftIsCentral_decat}, twisting under conjugation by $H_{w_0}$ will send the module associated to $w_0 \l$ to the module associated to $\l w_0$. However, conjugation is an inner automorphism, so it preserves the isomorphism class of a left module. Thus $w_0 \l$ and $\l w_0$ are the same two-sided cell.\footnote{Note that distinct left cells may give rise to isomorphic left modules, so this proof does not show that conjugation by $w_0$ fixes left cells. In fact, it does not.}

The automorphism $\tau$ sends $H_x$ to $H_{w_0 x w_0}$, and similarly sends $b_x$ to $b_{w_0 x w_0}$. Thus $\tau$ preserves KL polynomials, structure coefficients, etcetera. It follows that $\Delta(\tau(x)) = \Delta(x)$ and $\rbb(\tau(x)) = \rbb(x)$. Thus $w_0 d w_0$ must be distinguished whenever $d$ is distinguished. \end{proof}


\begin{defn} \label{defn:c} For a two-sided cell $\l$, let $\l^t$ denote the two-sided cell $w_0 \l$. Let $\cbb(\l)$ denote $\rbb(\l^t)$. For any distinguished involution $d \in \DC$, we refer to $w_0 d$ as a \emph{$w_0$-twisted distinguished involution}\footnote{It need not be an involution.}. \end{defn}

\begin{lemma}\label{lemma:cineq}
If $\l<\mu$ then $\cbb(\l)<\cbb(\mu)$.
\end{lemma}
\begin{proof}
Follows from Proposition \ref{prop:w0oncells} and Proposition \ref{prop:someP} (P11).
\end{proof}

Here is a crucial consequence of this result for the half twist.

\begin{prop} \label{prop:htm} Let $\l$ be a two-sided cell. In the language of \S\ref{subsec:twists}, one has $\mb_{\HT}(\l) = \cbb(\l)$. Moreover, the summands of $\HT$ in homological degree $\cbb(\l)$ and cell $\l$ are precisely the $w_0$-twisted involutions \begin{equation} \bigoplus_{d \in \DC \cap \l^t} B_{w_0 d}(\cbb(\l)). \end{equation} \end{prop}
	
\begin{proof} By Theorem \ref{thm:diagonalmiracle}, $\HT$ is perverse. By \eqref{eq:KLhalftwist}, $B_x$ appears in $\HT$ with graded multiplicity (for both homological shift and bimodule grading shift) given by $h_{1,w_0 x}$, whose minimum power of $v$ is $\D(w_0 x)$ by definition. As $x$ ranges over the cell $\l$, $w_0 x$ will range over $\l^t$, and by (P1) of Proposition
\ref{prop:moreP} the minimum value of $\D$ in this cell is $\rbb(\l^t) = \cbb(\l)$. This proves that $\mb_{\HT}(\l) = \cbb(\l)$. Moreover (again see Proposition \ref{prop:moreP}), the
only elements of the cell $\l^t$ which minimize $\D$ are the distinguished involutions $d \in \DC \cap \l^t$, for which the coefficient of $v^{\cbb(\l)}$ in $h_{1,d}$ is $1$. \end{proof}

Thus, if the half twist were twist-like, then it would be increasing by Lemma \ref{lemma:cineq}.

\subsection{The Sch\"utzenberger involution and the eigenvalues of the full twist}
\label{subsec:schutz}
We now discuss Sch\"utzenberger duality, which is a surprising involutory operation which fixes the left cells. Here is a restatement of a theorem\footnote{This theorem of Mathas was
recently generalized to Hecke algebras with unequal parameters by Lusztig in \cite{LuszLongest15}. In type $A$, Mathas attributes this theorem to J.~J.~Graham's thesis.} of Mathas
\cite[Theorem 3.1]{Mathas96}.

\begin{defn} Let $\xbb(\l) = \cbb(\l) - \rbb(\l)$. We refer to this as the \emph{content} of a two-sided cell. \end{defn}

\begin{thm} \label{thm:schugeneral} If $W$ is any finite Coxeter group, then there is an involution $\Schu_L \co W \to W$ which preserves each left cell, so that if $y$ is in cell $\l$
then \begin{equation} \label{eq:schugeneral} H_{w_0} b_y \equiv (-1)^{\cbb(\l)} v^{\xbb(\l)} b_{\Schu_L(y)} + \HB_{<\l}. \end{equation} Moreover, $\Schu_L(y)$ is in the same right cell
as $w_0 y w_0$.
\end{thm}

\begin{remark} \revcomment{Moved all reference to $\Schu_R$ here.} Similarly, there is an involution $\Schu_R \co W \to W$ and an analogous formula which involves the right action of $H_{w_0}$ instead. The relationship between $\Schu_L$ and $\Schu_R$ is explored in \cite[Proposition 3.9]{Mathas96}; in particular, $\Schu_L(y\inv) = \Schu_R(y)\inv$. \end{remark}

In type $A$ an element is determined by its left and right cells, so $\Schu_L(y)$ is the unique element with $\Schu_L(y) \sim_L y$ and $\Schu_L(y) \sim_R w_0 y w_0$. In
other types, Mathas provides an additional condition on $\Schu_L$ which pins down the element $\Schu_L(y)$ uniquely; see Lemma \ref{lem:schuvstwist}. 
	
In type $A$, where the elements of a given left cell are parametrized by standard Young tableau with a given shape, $\Schu_L$ corresponds to a standard combinatorial involution on tableau called the \emph{Sch\"utzenberger involution}.  For this reason we call $\Schu_L$ the \emph{(generalized) Sch\"utzenberger involution}.

\begin{remark} If $\HT$ were sharp, then \eqref{eq:schugeneral} would be categorified in the best possible way: the minimal complex of $\HT \ot B_y$ would have most terms in cells $<
\l$, and there would be a unique term in cell $\l$, namely $B_{\Schu_L(y)}$ living in homological degree $\cbb(\l)$ with grading shift $\xbb(\l)$. This idea is pursued in the next
section. \end{remark}

%
%

Note the following crucial consequence of Theorem \ref{thm:schugeneral}.

\begin{cor}\label{cor:ftUpperTriangDecat} The operator of left multiplication by $H_{w_0}^2$ is upper-triangular with respect to the cell filtration, and on $\HB_{\le \l}/\HB_{< \l}$ it acts by the scalar\footnote{Of course, $2\cbb(\l)$ is even so the factor $(-1)^{2\cbb(\l)}$ is redundant. Nonetheless, we keep it around in anticipation of the fact that it will be categorified by a nontrivial homological shift.} $(-1)^{2\cbb(\l)} v^{2\xbb(\l)}$. \end{cor}

We have already seen that $H_{w_0}^2$ is central, hence acts diagonalizably on $\HB\otimes_{\Z[v,v\inv]}\Q(v)$.  The above pins down the eigenvalues exactly, and also illustrates the interaction with cell theory.

We conclude this section with another way to think about the Sch\"utzenberger involution. The lemma below was proven by Mathas as well (see \cite[Theorem 3.1]{Mathas96}) but we include a proof to emphasize again the same ideas used in Proposition \ref{prop:htm}.

\begin{lemma} \label{lem:schuvstwist} Fix $y \in W$ in two-sided cell $\l$. Then there is a unique distinguished involution $d$ in the opposite cell $\l^t$, and a unique $z$ in cell
$\l$, such that $t_{w_0 d, y}^z \ne 0$. In particular, $d$ is the distinguished involution in the same right cell as $y w_0$, $z$ is equal to $\Schu_L(y)$, and $t_{w_0 d, y}^z = 1$.
\end{lemma}

\begin{proof} Let $d$ be a distinguished involution in $\l^t$ and $z \in W$. If $t_{w_0 d, y}^z \ne 0$ then by (P8) of Proposition \ref{prop:someP} we have \begin{itemize} \item $w_0 d
\sim_R z$, \item $y \sim_L z$, and \item $w_0 d \sim_L y\inv$. \end{itemize} This last property is equivalent to $d w_0 \sim_R y$ and $d \sim_R y w_0$. Multiplying by $w_0$ on the left,
this is equivalent to $w_0 d \sim_R w_0 y w_0$. Combining this with the first property, we see that $z \sim_R w_0 y w_0$. In particular, $d$ is uniquely specified as above, and $z$ is at
least in the same right and left cell as $\Schu_L(y)$. (In type $A$ this already implies $z = \Schu_L(y)$.)

By the KL inversion formula \eqref{eq:KLhalftwist}, we can write $H_{w_0} b_y$ as
\begin{equation} H_{w_0} b_y = \sum_{x \in W} (-1)^{\ell(w_0) - \ell(w_0 x)} h_{1, x} b_{w_0 x} b_y = \sum_{x, z \in W} (-1)^{\ell(x)} h_{1, x} c_{w_0 x,y}^z b_z. \end{equation}
We simplified this formula by observing that $\ell(w_0) - \ell(w_0 x) = \ell(x)$. Let us consider this formula modulo $\HB_{< \l}$. By \eqref{eq:schugeneral}, we know that the answer must be $(-1)^{\cbb(\l)} v^{\xbb(\l)} b_{\Schu_L(y)}$.

All the terms with $w_0 x$ in cells $< \l$ disappear modulo $\HB_{< \l}$. In fact, if $w_0 x$ is in a cell incomparable to $\l$ then $b_{w_0 x} b_y$ is in cells $< \l$, so these terms
disappear as well. Hence we can assume that $w_0 x$ is in cells $\ge \l$. Therefore, $\D(x) \ge \rbb(x) \ge \rbb(\l^t) = \cbb(\l)$, with \revise{$\D(x) = \rbb(\l^t)$} precisely when $x$ is in cell $\l^t$ and
is a distinguished involution. Consequently, the minimal power of $v$ appearing in $h_{1,x}$ is $v^{\cbb(\l)}$, for all $x$ which are relevant modulo $\HB_{< \l}$.

All terms with $z$ in cells $< \l$ disappear modulo $\HB_{< \l}$. Moreover, every nonzero term has $z$ in cells $\le \l$, because $b_z$ is a summand of $b_{w_0 x} b_y$. So we can assume
$z \sim_{LR} y$. The minimal power of $v$ appearing in $c_{w_0 x,y}^z$ when $z \sim_{LR} y$ is $v^{-\rbb(\l)}$, with equality precisely when $t_{w_0 x,y}^z \ne 0$.

Hence, every coefficient of $b_z$ has power of $v$ greater than $v^{\cbb(\l) - \rbb(\l)} = v^{\xbb(\l)}$, and the terms which contribute to this particular degree are \begin{equation}
\sum_{z \in \l, \\ d \in \DC \cap \l^t} (-1)^{\ell(d)} t_{w_0 d, y}^z b_z. \end{equation} We have already observed in Lemma \ref{lem:schuvstwist} that there is a unique $d \in \DC \cap
\l^t$ for which $t_{w_0 d,y}^z$ can be nonzero, namely the distinguished involution in the same right cell as $y w_0$. Thus, the coefficient of any $b_z$ is simply $\pm t_{w_0 d,y}^z$
for this particular $d$.

Since the answer must be $(-1)^{\cbb(\l)} v^{\xbb(\l)} b_{\Schu_L(y)}$, we see that $t_{w_0 d,y}^{\Schu_L(y)} = \pm 1$, and all other $t_{w_0 d, y}^z$ are zero. We already know by positivity that all structure coefficients live in $\N[v,v\inv]$, so $t_{w_0 d,y}^{\Schu_L(y)} = 1$. \end{proof}

\begin{remark} It is not difficult to show that, for any distinguished involution $d$ in cell $\l^t$, $\ell(d)$ and $\rbb(\l^t)$ have the same parity, which avoids the need for
positivity in this last step. This uses the fact that the KL polynomials $h_{x,y}$ and structure coefficients $c_{x,y}^z$ are all parity (they have only even or odd powers of $v$).
\end{remark}

\begin{cor} \label{cor:schudistinguished} Let $e$ be a distinguished involution in cell $\l$. Then there is a unique distinguished involution $d$ in cell $\l^t$ such that $\Schu_L(e) = d
w_0$. Namely, $d$ is the distinguished involution in the same right cell as $w_0 e$. \end{cor}
	
\begin{proof} By the last property in Proposition \ref{prop:moreP}, if $t_{x,y}^e \ne 0$ then $y = x^{-1}$. Applying this result when $y = \Schu_L(e)$ and $x = w_0 d$ for the
distinguished involution $d$ in the same right cell as $\Schu_L(e) w_0$, Lemma \ref{lem:schuvstwist} says that $t_{x,y}^e \ne 0$, and thus $y = x^{-1}$. Thus $y = d w_0$. Note that
$\Schu_L(e)$ is in the same right cell as $w_0 e w_0$, so $\Schu_L(e) w_0$ is in the same right cell as $w_0 e$. \end{proof}

Thus the set of $w_0$-twisted involutions $w_0 \DC$ is the same as the set of Sch\"utzenberger-twisted involutions $\Schu_L(\DC)$, though this bijection reverses cells, i.e. $w_0 d$ for
$d \in \l^t$ is $\Schu_L(e)$ for $e \in \l$.

\subsection{Main conjectures}
\label{subsec:sharpconj}

\begin{conjecture} \label{conj:HTaction} For any finite Coxeter group, the complex $\HT \in \KC^b(\SBim)$ is twist-like and sharp. \end{conjecture}

We prove this conjecture in type $A$ in \S\ref{subsec:bounding}. Let us describe some of the first consequences of this conjecture.

\begin{prop}\label{prop:FTconsequence}
Suppose $\HT$ is twist-like and sharp. Then $\nb_{\HT}(\l)=\mb_{\HT}(\l)=\cbb(\l)$ and $\nb_{\FT}(\l)=\mb_{\FT}(\l)=2\cbb(\l)$ for all two-sided cells $\l$.  In particular both $\HT$ and $\FT$ are sharp and increasing.  Moreover:
\begin{enumerate}
\item If $B_x$ is in cell $\l$ then the head of $\HT\otimes B_x$ equals $B_{\Schu_L(x)}[-\cbb(\l)](\xbb(\l))$.
\item If $B_x$ is in cell $\l$ then the head of $\FT\otimes B_x$ equals $\Sigma_\l(B_{x})$, where
\[
\Sigma_\l=[-2\cbb(\l)](2\xbb(\l)).
\]
\item The terms in $\HT$ in cell $\l$ and minimal homological degree $\cbb(\l)$ are exactly
\[ \bigoplus_d B_{\Schu_L(d)}(\cbb(\l)), \] where $d$ ranges over the distinguished involutions in cell $\l$.
\item The terms in $\FT$ in cell $\l$ and minimal homological degree $2\cbb(\l)$ are exactly
\[ \bigoplus_d B_d(\cbb(\l) + \xbb(\l)), \] where $d$ ranges over the distinguished involutions in cell $\l$.
\end{enumerate}
\end{prop}

\begin{proof}
Assume that $\HT$ is twist-like and sharp.  In Proposition \ref{prop:htm} we showed that $\mb_{\HT}(\l)=\cbb(\l)$, hence the inequality $\mb_{\HT}(\l)<\mb_{\HT}(\mu)$ when $\l<\mu$ follows from Lemma \ref{lemma:cineq}.  Thus, $\HT$ is increasing.  Also $\nb_{\HT}(\l)=\mb_{\HT}(\l)=\cbb(\l)$ by sharpness.  Furthermore, by Lemma \ref{lem:sharptensorclosed} we can also deduce that $\FT$ is sharp and increasing, and $\mb_{\FT}(\l) =  2\cbb(\l)$.

By Proposition \ref{prop:tailProperty}, we know a lot about $\HT \ot B_x$ for any $x$ in cell $\l$. In particular, the only indecomposable summand of the minimal complex of $\HT \ot B_x$ in cell $\l$ is its head. \revise{The image of $\HT \ot B_x$ in the Grothendieck group is governed by \eqref{eq:schugeneral}, which also has exactly one summand in cell $\l$. So} we deduce that its head is as desired. Similar arguments describe the head of $\FT \ot B_x$.  This proves (1) and (2).

We have already proven in Proposition \ref{prop:htm} that the terms in cell $\l$ and minimal homological degree are precisely
\[
\bigoplus_d B_{w_0 d}(\cbb(\l)),
\]
where the sum is over distinguished involutions $d$ in the opposite cell $\l^t$.  This proves (3).

It remains to study the terms in $\FT$ in homological degree $2\cbb(\l)$ and cell $\l$. For the rest of this proof, let us fix $\l$ and set $\cbb = \cbb(\l)$, $\xbb = \xbb(\l)$, $\rbb = \rbb(\l)$.  If $C$ is a complex and $k$ is an integer, we let $C^k$ denote the $k$-th chain object of $C$.

Let us study $\HT\otimes \HT$, whose minimal complex is $\FT$.   Remark that each indecomposable summand of the chain object $\FT^m$ is a direct summand of $\HT^k\otimes \HT^{\ell}$ for some $k,\ell$ with $k+\ell=m$.   If $k<\cbb$ or $l<\cbb$ then $\HT^k\otimes \HT^l$ is a direct sum of indecomosables in cells less than or incomparable to $\l$.  Thus all terms of $\FT$ in cell $\l$ come from terms of the form $\HT^k\otimes \HT^{\ell}$ with $k,\ell\geq \cbb$.  In particular summands of $\FT^{2\cbb}$ in cell $\l$ are all summands of $\HT^{\cbb}\otimes \HT^{\cbb}$.

Let $Z:=\bigoplus_{d\in \DC\cap \l^t} B_{w_0d}(\cbb)$.  As noted above $\HT^{\cbb}$ is isomorphic to a direct sum of $Z$ and other bimodules in cells smaller than or incomparable to $\l$.  Consider the minimal complex of $\HT\otimes Z$.  First, recall that each $w_0d$ can be written uniquely as $w_0d=\Schu_L(e)$ for some $e\in \DC\cap \l$, by Corollary \ref{cor:schudistinguished}.  The head of $\HT\otimes B_{\Schu_L(e)}$ is isomorphic to $B_e[-\cbb](\xbb)$ since $\Schu_L(\Schu_L(e))=e$.  In particular the head is in homological degree $\cbb$.  It follows that each summand $B_{\Schu_L(e)}(\cbb)\subset Z$ ultimately contributes a summand $B_e(\cbb+\xbb)$ to $\FT$, in homological degree $2\cbb$.

It remains to prove that these terms survive to the minimal complex $\FT$ of $\HT \ot \HT$.  This can be done by proving that no isomorphic terms appear in adjacent homological degrees $2\cbb \pm 1$, so that no Gaussian elimination can cancel these terms.  The constraints $k, \ell \ge \cbb$ already imply that nothing in cell $\l$ appears in homological degree $2\cbb-1$. In homological degree $2\cbb+1$, either $k = \cbb$ and $\ell = \cbb + 1$ or vice versa; without loss of generality assume $k = \cbb$.  Every indecomposable summand of $\HT^{\cbb+1}$ is of the form $B_y(k+1)$ by perversity.   If $B_z(i)$ is a direct summand of $B_y(\cbb+1) \ot B_x(\cbb)$ in cell $\l$, then we must have $i\geq \cbb+1 + \cbb - \rbb = \cbb + \xbb +1$.  Hence nothing with grading shift $\cbb + \xbb$ can appear. This concludes the proof.
\end{proof}

The shift functors $\Sigma_\l=[-2\cbb(\l)](2\xbb(\l))$ therefore play the role of the ``eigenvalues'' of the full twist (modulo Conjecture \ref{conj:HTaction}), and $\FT$ is
cell triangular in the sense of Definition \ref{def:twistedunitri}. Thus it makes sense to look for $\l$-equivalences for each cell $\l$, see Definition \ref{def:lequiv}.

\begin{conjecture}[Eigenmap conjecture]\label{conj:eigenmap}
If $W$ is a finite Coxeter group and $\l$ is a two-sided cell in $W$ then there exists a map $\a_\l:\Sigma_\l(\one)\rightarrow \FT$ which is a $\l$-equivalence.  These maps are obstruction free in the sense of \cite[Definition 6.12]{ElHog16a}.
\end{conjecture}

We prove this for $W=S_n$ in \S \ref{sec:constructing}.

\subsection{Dot maps}
\label{subsec:dots}

Previously we used the numerics of the functions $\D$ and $\rbb$ to study complexes of Soergel bimodules, using perversity and \revadd{tensor product decomposition} rules. Now we discuss another categorical shadow of these numerical invariants, this time using the Soergel Hom formula. We need not assume $W$ is finite in this section.

\begin{proposition}\label{prop:dotspace} Let $x$ be in two-sided cell $\l$. The smallest $k$ such that there exists a nonzero map $R \rightarrow B_x(k)$ equals $\Delta(x)$. Within each
cell, this minimum equals $\rbb(\l)$, and is achieved only for the distinguished involutions. If $x$ is a distinguished involution, then $\Hom(R,B_x(\rbb(\l)))$ is 1-dimensional.
\end{proposition}

\begin{proof}
By the Soergel Hom formula, the graded rank of $\Homg(R,B_x)$ over $R$ is $h_{1,x}$, so its minimal degree is $\D(x)$. The remaining statements follow \revise{accordingly} from Proposition \ref{prop:moreP}.
\end{proof}


\begin{definition}\label{def:dots}
For each $d\in \DC$ in two-sided cell $\l$, choose $\xi_d: R\rightarrow B_{d}(\rbb(\l))$ which \revadd{spans} the corresponding Hom space.  Let $\xi_d^\ast$ denote the dual map $B_{d}\rightarrow R(\rbb(\l))$. We call these \emph{(generalized) dot maps} or \emph{$d$-dot maps}. Let $\tbarb{d}$, the \emph{(generalized) barbell} denote the composition $\xi_d^\ast \circ \xi_d$.  Then $\tbarb{d}$ is a degree $2\rbb(\l)$ element of $R=\Endg(R)$. 
\end{definition}

\begin{example}\label{ex:usualDot} Fix a simple reflection $s$, which is a distinguished involution in the simple cell $\l$ satisfying $\rbb(\l)=1$. We may choose $\xi_s:R\rightarrow
B_s(1)$ and $\xi_s^\ast:B_s\rightarrow R(1)$ to be the usual \emph{dot maps} in the diagrammatic calculus, see Remark \ref{rmk:dotsdefn}. Composing the two dots gives a morphism $R
\to R(2)$ which is multiplication by $\a_s$, and is pictorially represented by a \emph{barbell}. \end{example}

\begin{example}\label{ex:longestBarbell}
Let $W$ be finite. Then $w_0$ is in the longest cell $\l$ with $\rbb(\l) = \ell(w_0)$. Let $\ell = \ell(w_0)$.  Recall\footnote{\revise{This statement can be proven in a straightforward fashion using Soergel's theory of support filtrations \cite{Soer07}, but also follows quickly from \cite[Theorem 3]{WillSingular}.}} that $B_{w_0}=R\otimes_{R^W}R(\ell)$, and $R$ is a Frobenius extension over $R^W$ with Frobenius trace $\pa_{w_0}$.  We choose $\xi_w:R \rightarrow B_{w_0}(\ell)$ to be the map which sends
\[
1\mapsto \sum a_i \ot b_i,
\]
where the sum is over dual bases for $R$ over $R^W$ with respect to $\pa_{w_0}$. The dual map $B_{w_0}\rightarrow R(\ell)$ sends $f\otimes g\mapsto fg$.  Their composition $\tbarb{w_0}$ is the product of all the positive roots. See \cite[\S 3.6]{EThick} for more discussion of these maps. \end{example}

The importance of $d$-barbells will be highlighted in the next chapter. For now, we focus on $d$-dot maps.

\begin{lemma} \label{lem:notinlower} For $d$ in two-sided cell $\l$, the $d$-dot maps are not in the ideal $\IC_{< \l}$, or the ideal $\IC_{\ngeq \l}$. \end{lemma}

\begin{proof} We need to prove that $\xi_d$ is not a linear combination of maps which factor as $R \to B_z \to B_d$ for $z$ in lower cells. But the minimal degree of a map $R \to B_z$ is
$\D(z) \ge \rbb(z) \ge \rbb(\l)$, and the minimal degree of a map $B_z \to B_d$ is $+1$ (as for any two non-isomorphic indecomposable Soergel bimodules). Thus no map of degree $\rbb(\l)$
can factor through lower cells. Since $\xi_d$ does factor through $B_d$ in cell $\l$, it is in $\IC_{\le \l}$, so not being in $\IC_{< \l}$ is equivalent to not being in $\IC_{\ngeq
\l}$. \end{proof}

Now let us fix a finite Coxeter group $W$ and assume Conjectures \ref{conj:HTaction} and \ref{conj:eigenmap}, so that $\l$-equivalences $\a_\l \co \Sigma_\l (R)\to \FT$ exist. Since $\Sigma_\l(R)$ is just a shift of $R$ concentrated in a single homological degree, \revise{the chain map $\a_\l$ is concentrated in a single degree, where it is a bimodule map from $R(2\xbb(\l))$ to the} homological degree $2\cbb(\l)$ chain object of the full twist, which contains summands \[ \bigoplus_d B_d(\cbb(\l) + \xbb(\l)) \] for all $d \in \DC \cap \l$. Thus $\a_\l$ determines a collection of morphisms $R\rightarrow B_d(\cbb(\l) - \xbb(\l))$, which must be scalar multiples of the dot maps since $\cbb(\l)-\xbb(\l)=\rbb(\l)$.  We call these the \emph{$\xi_d$ components} of $\a_\l$.  \revcomment{Previously it said: Strictly speaking, these components are only well-defined up to unit multiple, since there is a choice in how $B_d(\cbb(\l) + \xbb(\l))$ sits as a direct summand of $\FT^{2\cbb(\l)}$.}

\revcomment{turned a paragraph of discussion into a lemma.}

\begin{lemma} The $\xi_d$ components of a chain map $\a_\l \co \Sigma_{\l}(R) \to \FT$ depend only on the homotopy class of $\a_\l$. \end{lemma}
	
\begin{proof}  If $\b-\a_\l = [d,h]$ then actually $\b-\a_\l = d\circ h$ (the differential on $\Sigma_\l(R)$ being zero), and $h$ is concentrated in one degree, where it is a bimodule morphism from a shift of $R$ to the homological degree $2\cbb(\l)-1$ part of $\FT$.  It follows that $h$ factors through terms in cells $\ngeq \l$, hence $d\circ h$ cannot affect the $\xi_d$ components by Lemma \ref{lem:notinlower}.\end{proof}

\begin{proposition}\label{prop:dotsFromEigenmaps}
Assume that $W$ is a finite Coxeter group and Conjectures \ref{conj:HTaction} and \ref{conj:eigenmap} hold. Let $x \in W$ be in the same right cell as $d \in \DC$. Then $\xi_d\otimes \Id_{B_x}$ is the inclusion of a direct summand while $\xi_d^\ast\otimes \Id_{B_x}$ is the projection onto a direct summand. The $\xi_d$ component of $\a_\l$ is an invertible\footnote{We work over a field, so we may as well say nonzero. When working over other base rings, it is not hard to prove this slightly stronger result.} scalar multiple of $\xi_d$.
\end{proposition}

\begin{proof}
Fix a $\l$-equivalence $\a_\l$, for the two-sided cell $\l$ containing $x$ and $d$. As before, let $\cbb = \cbb(\l)$, etcetera. By Proposition \ref{prop:FTconsequence} we know that the head of the minimal complex of $\FT \ot B_x$ is $\Sigma_\l(B_x)$ and that $\a_\l \ot \Id_{B_x}$ is an isomorphism from $\Sigma_\l(B_x)$ to this head. We claim that it is only the $\xi_d$ component of $\a_\l$ which could contribute to the isomorphism in $\a_\l \ot \Id_{B_x}$.

In homological degree $2\cbb$ the complex $\FT$ is a direct sum of $B_{d'}(\cbb+\xbb)$, where $d'$ ranges over the distinguished involutions in cell $\l$, together with objects in cells
$\ngeq \l$, by Proposition \ref{prop:FTconsequence}. Consequently, the only direct summands of $\FT \ot B_x$ in this homological degree and in cell $\l$ are summands $B_{d'} \ot B_x(\cbb
+ \xbb)$. The minimal grading shift which could occur in such a summand is $B_z(\cbb + \xbb - \rbb) = B_z(2\xbb)$, and it is realized only when $t_{d',x}^z \ne 0$. But by Proposition
\ref{prop:dacts}, this happens precisely when $d'$ is in the same right cell as $x$ (i.e. $d' = d$), and $z = x$. Hence the summand $B_x(2\xbb)$ which survives to the minimal complex of
$\FT \ot B_x$ must have arisen as a direct summand of $B_d(\cbb + \xbb) \ot B_x$.

If the $\xi_d$ component of $\a_\l$ is zero (or noninvertible), then so is the component of $\a_\l \ot \Id_{B_x} \co B_x[-2\cbb](2\xbb) \to \FT \ot B_x$ which goes to $B_d(\cbb+\xbb) \ot
B_x$ and its summand $B_x(2\xbb)$. Gaussian elimination from $\FT \ot B_x$ to its minimal complex can change the components of a chain map, but only by null-homotopic maps. Any homotopy
would factor through homological degrees $< 2\cbb$ in $\FT \ot B_x$, and therefore would factor through cells $< \l$. Thus it can not make a non-isomorphism into an isomorphism. This
contradicts Lemma \ref{lem:lequivcriterion}. Thus the $\xi_d$-component of $\a_\l$ is invertible, $\xi_d \ot \Id_{B_x}$ must be the inclusion of the $B_x$ summand.

Duality implies that $\xi_D^\ast$ is projection to a summand. \end{proof}

One interprets Proposition \ref{prop:dotsFromEigenmaps} as a strong categorification of the unit axiom in the $J$-ring. Decategorified, if $d$ is a distinguished involution then
$v^{\rbb(d)} b_d$ acts by the identity on its right cell, modulo positive powers of $v$ and lower two-sided cells. Corollary \ref{cor:weakJringCat} was a weak categorification of the
same statement. A strong categorification would give a natural transformation between the identity functor and the functor $B_d(\rbb(d)) \ot (-)$, which realizes the fact that
$B_d(\rbb(d))$ acts by the identity on its right cell (modulo positive shifts and lower cells). By Proposition \ref{prop:dotsFromEigenmaps}, the $d$-dot $\xi_d$ serves this role!

Moreover, in cell $\l$ and homological degree $2\cbb(\l)$, $\FT$ is (up to shift) just a categorification of the element $v^{\rbb(\l)} \sum_{d \in \DC \cap \l} b_d$ in $\HB$ which ``descends'' to the unit of the ring $J_\l$. The $\l$-equivalence $\a_\l$ combines the $d$-dots into a map which realizes the unit axiom of the $J$-ring.

\begin{remark} \label{rmk:syzygy}
We wish to point out an interesting interpretation of the above.  First, let us choose a two-sided cell $\l$, and let $F = \FT[2\cbb(\l)](-2\xbb(\l))$. Let's work modulo lower cells, i.e.~ in the category $\SBim_n/I_{<\l}$.  Modulo Conjecture \ref{conj:HTaction}, Proposition \ref{prop:FTconsequence} tells us that
\[
F = F^0\rightarrow F^1\rightarrow \cdots \rightarrow F^{m},
\]
where $F^0 = \bigoplus_{d\in \DC\cap \l} B_d(\rbb(\l))$ is the Soergel bimodule analogue of the unit in the $J$-ring.  The maximal homological degree happens to be $m=2\ell(w_0)-2\cbb(\l)$, but this won't be relevant.  We know that if $x\in W$ is in cell $\l$ then
\[
F^0\otimes B_x \cong B_x \oplus Y \qquad \text{modulo $I_{<\l}$}
\]
where $Y$ is a direct sum $B_y(k)$ with $y$ in cell $\l$ and $k>0$.  We also know (again modulo Conjecture \ref{conj:HTaction}) that
\[
F\otimes B_x \simeq B_x \qquad \text{modulo $I_{<\l}$}.
\]
That is to say, the terms in $F^0\otimes B_x$ with positive degree shifts are being cancelled by $F^1\otimes B_x$, and the surviving terms in $F^1\otimes B_x$ are cancelled by $F^2\otimes B_x$, and so on.  Thus, one can think of $F^0$ more accurately as a first approximation to the unit in the $J$-ring, $F^0\rightarrow F^1$ as a better approximation, and so on, until $F^\bullet$ itself actually behaves as the unit in the $J$-ring associated to $\l$.

Note that $(F^0\rightarrow \cdots \rightarrow F^i)\otimes B_x\simeq B_x \oplus Y^i[-i]$ where $Y^i\in \SBim$ ($0\leq i\leq m$); these yield potentially interesting new invariants $[Y^i]$ in the Hecke algebra.
\end{remark}


Let us conclude this section by motivating and discussing one additional application of Proposition \ref{prop:dotsFromEigenmaps}.

When $d = w_I$ is the longest element of a finite parabolic subgroup $W_I$, the Soergel bimodule $B_{w_I} \cong R \ot_{R^I} R(\ell(w_I))$ is a graded Frobenius algebra object of degree $\ell(w_I) = \rbb(w_I)$. This is because $R^I \subset R$ is a Frobenius extension of commutative graded rings. In this case, the dot maps $\xi$ and $\xi^\ast$ give rise to the unit and counit maps in the Frobenius algebra structure. The multiplication and comultiplication maps arise via the decomposition 
\[ B_{w_I} \ot B_{w_I} \cong B_{w_I}(-\ell(w_I)) \oplus \cdots \oplus B_{w_I}(\ell(w_I)),\]
namely, they are the projection to the minimal degree summand, and the inclusion from the maximal degree summand. All four of these maps are well-defined up to scalar, and some choice of scalars makes the Frobenius algebra axioms hold.

\begin{conj} For any Coxeter group and any distinguished involution $d$, $B_d$ has the structure of a graded Frobenius algebra object in $\SBim$ of degree $\rbb(d)$. The unit and counit are given (up to scalar) by $\xi_d$ and $\xi_d^\ast$ respectively. The multiplication (resp. comultiplication) map is given (up to scalar) by a projection to (resp. inclusion of) the minimal (resp. maximal) degree summand in
\[ B_d \ot B_d \cong B_d(-\rbb(\l)) \oplus \cdots \oplus B_d(\rbb(\l)) \oplus Z. \]
Here $Z$ represents those summands in cells $< \l$, which may in theory have more extreme degree shifts, but which we ignore when we use the words ``minimal'' and ``maximal.''
\end{conj}

\revcomment{Removed a promise about a follow-up paper.}
Proposition \ref{prop:dotsFromEigenmaps} is in essence a proof of the unit axiom for the Frobenius algebra.


\section{Cell theory in type $A$}
\label{sec:typeAcells}

In this chapter we explain how many of the previous constructions (cell theory, $\rbb$-function, distinguished involutions, Sch\"utzenberger involution) work in the special case of type
$A$. We go into more detail on examples to give more intuition to the novice reader. For a good introduction to this topic in type $A$, we recommend \cite{ArikiCells} or \cite{GeordiePQ}.

\subsection{Cells in type $A$}
\label{subsec:cellA}

Let $\PC(n)$ denote the set of partitions of $n$, and $\SYT(\l)$ denote the set of standard Young tableaux of shape $\l$, for $\l \in \PC(n)$. The following theorem is a major tool in
combinatorics; for a good reference we recommend \cite{FultonTab}.

\begin{thm}(Robinson-Schensted correspondence) There is an explicit bijection between elements of $S_n$ and triples $(P,Q,\l)$ where $\l \in \PC(n)$, and $P, Q \in \SYT(\l)$. The
algorithm which takes $w \in S_n$ and returns a triple $(P,Q,\l)$ is called Schensted's \emph{bumping algorithm}. \end{thm}

We will not discuss the bumping algorithm here, though we will state many consequences of it below. Since $\l$ is determined by $P$ and $Q$ we often just write $(P,Q)$ for the pair of
tableaux, with the understanding that $P$ and $Q$ have the same shape. We write $w(P,Q,\l)$ or $w(P,Q)$ for the corresponding element of $S_n$. For shorthand, we write the corresponding
indecomposable Soergel bimodule $B_{w(P,Q)}$ as $B_{P,Q}$, and the KL basis element as $b_{P,Q}$.

The following is a crucial theorem of Kazhdan-Lusztig\footnote{The reader will not find this theorem in this form inside \cite{KazLus79}, though the results of that paper imply it when one understands Knuth equivalence. One can find several expository resources aimed at making this theorem explicit from Kazhdan-Lusztig's work or proving it more simply, see \cite{ArikiCells} or \cite{GeordiePQ}.} \cite{KazLus79}, describing the cells in the Hecke algebra with respect to the KL basis of $\HB(S_n)$.  Since the KL basis corresponds to the indecomposable objects in $\SBim_n$ by Theorem \ref{thm:SC}, this also describes the cells in $\SBim_n$.

\begin{thm}\label{thm:soergelcells}
Let $w = w(P,Q,\l)$ and $w' = w(P',Q',\l')$. Then $w \sim_{LR} w'$ iff $\l = \l'$, $w \sim_L w'$ iff $Q=Q'$, and $w \sim_R w'$ iff $P=P'$. Thus, we associate the two-sided cells of $S_n$ with partitions $\l$, and the left (or right) cells with tableaux. Moreover, we have $\l \le_{LR} \mu$ if and only if $\l \le \mu$ in the dominance order.
\end{thm}

This theorem describes the left and right cells, and the two-sided cells, as well as the partial order on two-sided cells. The partial order on left cells is not as clear-cut, and even the following important result of Lusztig does not have an elementary proof (it follows from (P11) of Proposition \ref{prop:someP}).

\begin{prop} \label{prop:leftincomp} Two non-equal left cells in the same two-sided cell are incomparable. \end{prop}

\begin{remark} \revcomment{revised for brevity} The partial order on left cells is mysterious, and does not agree with the dominance order on standard tableaux. See \cite{Taskin} for a comparison of various partial orders on standard tableaux. \end{remark}
	

\begin{example}
Let $w\in S_n$ be in cell $\l\in \PC(n)$ and $v\in S_m$ be in cell $\mu\in \PC(m)$. With respect to the inclusion $S_n\times S_m\rightarrow S_{n+m}$, $(w,v)$ is in cell $\l+\mu$, which is the partition of $n+m$ given by $(\l+\mu)_i = \l_i+\mu_i$ (extending $\l$ or $\mu$ by zero if necessary).   This is easy to prove using the bumping algorithm.  The tableaux also add in a similar way, though we will not need this.
\end{example}

\begin{example}\label{ex:longesteltcell}
If $w_0\in S_n$ is the longest element, then $w_0$ is in cell $1^n$, the one-column paritition.  More generally, if $w$ is the longest element of a parabolic subgroup $S_{k_1}\times \cdots \times S_{k_r}\subset S_n$, then $w$ is in cell $\l$, where $\l\in\PC(n)$ is the partition whose column lengths coincide with the multiset $\{k_i\}$ (but in decreasing order).
\end{example}

An important special case occurs when the integers $k_i$ in the previous example are non-increasing.

\begin{ex} \label{ex:columnreading} Let $\l \in \PC(n)$, with column sizes $k_1 \ge k_2 \ge \ldots \ge k_r > 0$. Let $P_{\col}$ be the \emph{column-reading tableau}, obtained by placing
the numbers $1$ through $k_1$ in the first column, $k_1 + 1$ through $k_1 + k_2$ in the second column, and so forth. Then $w(P_{\col},P_{\col},\l)$ is the longest element $w_{\l}$ of the
parabolic subgroup $S_\l = S_{k_1} \times S_{k_2} \times \cdots \times S_{k_r} \subset S_n$. \end{ex}

Here is a key property of the Robinson-Schensted correspondence.
	
\begin{prop} \label{prop:involution} If $x = w(P,Q,\l)$, then $x\inv = w(Q,P,\l)$. In particular, $x$ is an involution if and only if $x = w(P,P,\l)$ for some $P \in \SYT(\l)$. \end{prop}

\subsection{Numerics of the $\rbb$-function in type $A$}
\label{subsec:longestnumerics}

The statements in this section are all well-known. We go into more detail than necessary, to aid the novice reader.

For any Coxeter group $W$, each left cell contains a unique distinguished involution by Proposition \ref{prop:moreP}. For the symmetric group, every left cell $(-,P,\l)$ contains exactly
one involution $w(P,P,\l)$. Thus in type $A$, $\DC$ is the set of all involutions. By Corollary \ref{cor:computinga}, to compute $\rbb(\l)$ we need only find an involution $d$ in cell $\l$
and compute the largest and smallest powers of $v$ appearing in the expansion of $b_d^2$. The easiest case to analyze will typically be when $d$ is the longest element of a parabolic
subgroup.

Now we recall a property of KL polynomials of longest elements.

\begin{defn} Let $W$ be a finite Coxeter group. Then $\pi(W)$ is its \emph{balanced Poincare polynomial}, which is $v^{-\ell(w_0)} \sum_{w \in W} v^{2 \ell(w)}$. \end{defn}
	
It is clear from the definition that the lowest degree of $v$ which appears in this polynomial is $v^{-\ell(w_0)}$, and the highest degree is $v^{+\ell(w_0)}$. When $W_I = S_{k_1} \times \cdots \times S_{k_r}$, $\pi(W_I)$ is a product of quantum factorials:
\begin{equation}
\pi(W_I) = [k_1]! [k_2]! \cdots [k_r]!.
\end{equation}

\begin{lemma} \label{lem:longestsquared} When $w_I$ is the longest element of a parabolic subgroup $W_I$, one has \begin{equation} \label{eq:longestsquared} b_{w_I} b_{w_I} = \pi(W_I)
b_{w_I}.\end{equation} More generally, if $x \in W$ and $sx < x$ for all $s \in I$, then \begin{equation} \label{eq:descentfactorial} b_{w_I} b_x = \pi(W_I) b_x. \end{equation}
\end{lemma}

\begin{proof} Let us sketch the proof of this well-known fact. One has 
\begin{equation} \label{eq:swq2}
b_s b_x = (v+v^{-1}) b_x
\end{equation}
whenever $sx < x$. This tells you how $H_s$ acts on $b_x$, and thereby how $H_w$ acts on $b_x$ for all $w \in W_I$. From the fact that $b_{w_I}$ is smooth, it is easy to deduce the lemma. \end{proof}

For the symmetric group $S_4$, every involution is the longest element of a parabolic subgroup except two: $x = tsut$ and $w = sutsu$. Computations for these elements were done in
Example \ref{ex:boundviolated}.

Now we introduce some statistics associated to partitions.

\begin{defn}\label{defn:xcr} We think of a partition $\l \in \PC(n)$ as a Young diagram, drawn in the ``English style'':
\[
\ig{.6}{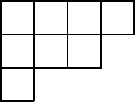}
\]
Suppose a box $\square$ is in the $i$-th column and $j$-th row.  Here columns and rows are counted left-to-right and top-to-bottom, starting at zero. We say that the box has \emph{column number} $\cbb(\square) = i$, \emph{row number} $\rbb(\square) = j$, and \emph{content} $\xbb(\square) = i-j$.

For a partition $\l$, we set $\cbb(\l) = \sum_{\square} \cbb(\square)$, $\rbb(\l) = \sum_\square \rbb(\square)$, and $\xbb(\l) = \sum_\square \xbb(\square)$. That is, the column number of a partition is the sum of the column numbers of each box, etcetera.  Equivalently, we have
\[
\rbb(\l) = \sum_i (i-1)\l_i, \qquad \qquad \cbb(\l) = \rbb(\l^t),\qquad\qquad \xbb(\l) = \cbb(\l)-\rbb(\l).
\]
\end{defn}

It is easy to observe that $\rbb(\l)$ is the length of $w_\l$, the longest element of the parabolic subgroup $S_\l$. It is also the length of $w_I$ for any conjugate parabolic subgroup $W_I$.

\begin{example} For the Young diagram $\l=(4,3,1)$ pictured above, one has $\cbb(\l) = 9$, $\rbb(\l) = 5$, and $\xbb(\l)=4$. The number $\rbb(\l)=5$ is the length of $w_{\l} = s_1 s_2 s_1 s_4 s_6 \in S_8$, \revise{see Example \ref{ex:columnreading}}.  \end{example}

\begin{example}
The two-row partition $\lambda=(n-1,1)$ has $\rbb(\lambda)=1$ and $\cbb(\lambda)=\binom{n-1}{2}$ for all $n$.
\end{example}

\begin{proposition}\label{prop:abEqualsrb} For a 2-sided cell $\l \in \PC(n)$, the row number $\rbb(\l)$ agrees with Lusztig's $\rbb$-function for the two-sided cell $\l$.\qed
\end{proposition}

\begin{proof} \revise{By Corollary \ref{cor:computinga} and Lemma \ref{lem:longestsquared} and the preceding discussion, we need only show that $\ell(w_{\l}) = \rbb(\l)$, which is straightforward. } \end{proof}

\begin{lemma} \label{lem:rowOrder}
If $\mu < \lambda$ in the dominance order then $\rbb(\mu) > \rbb(\l)$. In particular, if $\rbb(\l)=\rbb(\mu)$, then either $\mu=\l$ or $\mu$ and $\l$ are incomparable in the dominance order. \end{lemma}

\begin{proof} This is an easy exercise.  Alternatively, one can use (P4) from Proposition \ref{prop:someP}. \end{proof}

\begin{example} The partitions $\l=(3,1,1,1)$ and $\mu=(2,2,2)$ satisfy $\rbb(\l)=\rbb(\mu)$, and are incomparable. They also satisfy $\cbb(\l) = \cbb(\mu)$ and $\xbb(\l) = \xbb(\mu)$. \end{example}

\subsection{Multiplication in a cell}
\label{subsec:morestill}

\begin{notation} \label{not:HBP-} For $\l \in \PC(n)$, let $\HB_\l \subset \HB$ denote the $\Z[v,v\inv]$-span of $b_{P,Q}$, for $P,Q \in \SYT(\l)$. We may refine this span by fixing $Q$
and taking $\HB_{-,Q}$, the $\Z[v,v\inv]$-span of $b_{P,Q}$ for $P \in \SYT(\l)$. We define $\HB_{P,-}$ similarly. Following this analogy, $\HB_{P,Q}$ is just the $\Z[v,v\inv]$-span of
the element $b_{P,Q}$.
	
We define $\HB_\l^+ \subset \HB_\l$ as the $\Z[v]$-span of $v b_{P,Q}$, i.e. as those elements of $\HB_\l$ with strictly positive powers of $v$. For fixed $P \in \SYT(\l)$, we define
$\HB_{P,-}^+$ as $\HB_{P,-} \cap \HB_\l^+$. We define $\HB_{-,Q}^+$ similarly. \end{notation}

Let $P,Q,U,V \in \SYT(\l)$. The Definition \ref{def:a} of the $\rbb$-function gives the following equality:
\[ v^{\rbb(\l)} b_{P,Q} b_{U,V} \equiv \sum_{X,Y \in \SYT(\l)} t_{(P,Q),(U,V)}^{(X,Y)} b_{X,Y} + \HB^+_\l + \HB_{< \l}.\]
Most of the time, this coefficient $t$ is zero. For example, since $\HB_{-,V} + \HB_{< \l}$ forms a left ideal, we must have $Y = V$ to have a nonzero contribution to the sum. By similar arguments, $X = P$. Stronger still, property (P8) from Proposition \ref{prop:someP} implies that $t_{(P,Q),(U,V)}^{(X,Y)} = 0$ unless $X=P$, $Y=V$, and $Q=U$.

The last part of Proposition \ref{prop:moreP} implies that $t_{(P,Q),(Q,P)}^{(P,P)} = 1$, which treats the case when $P = V$. \revcomment{reorganized below, and slightly adapted the segues.} In fact, 
$t_{(P,Q),(Q,V)}^{(P,V)}=1$ for any $P,Q,V \in \SYT(\l)$, which we now justify. The Hecke algebra in type $A$ is a \emph{cellular algebra}\footnote{It is unfortunate that the word ``cell'' was appropriated by the literature for two similar but logically independent concepts: the cells described previously, and the cells in the cellularity theory of Graham and Lehrer. Thankfully, in this case, the Graham-Lehrer cells agree with the two-sided cells.} as defined by Graham and Lehrer \cite{GraLeh}. A major implication of this is that
\begin{equation} \label{eq:howPQmultcloser} b_{P,Q} b_{U,V} = \phi(Q,U) b_{P,V} + \HB_{< \l}. \end{equation}
Here, $\phi(Q,U)$, often called the \emph{cellular form}, is a function which takes $Q, U \in \SYT(\l)$ and returns a coefficient in $\Z[v,v\inv]$. In particular, this coefficient of $b_{P,V}$ does not depend on $P$ and $V$, only on $Q$ and $U$. This reduces the computation of $t_{(P,Q),(Q,V)}^{(P,V)}$ to the special case where $P = V$. For more on the cellular structure in type $A$, see \cite{GeordiePQ}.

\revadd{Thus we have justified the following well-known result.}

\begin{lemma} \label{lem:howPQmult} For any $P, Q, U, V \in \SYT(\l)$, we have
\begin{equation} \label{eq:howPQmult} v^{\rbb(\l)} b_{P,Q} b_{U,V} \equiv \delta_{Q,U} b_{P,V} + \HB^+_\l + \HB_{< \l}.\end{equation} \end{lemma}

Let $P,Q \in \SYT(\l)$, and let $a(P,Q)$ denote the most negative power of $v$ appearing with nonzero coefficient in $\phi(P,Q)$. The results above can be restated as follows: $a(P,Q)
\ge -\rbb(\l)$, with equality if and only if $P = Q$. Moreover, if $P=Q$, then the coefficient of $v^{-\rbb(\l)}$ is one.

\begin{remark} \label{rmk:notametric} A warning for the reader. The numbers $a(P,Q)$ for $P \ne Q$ are quite mysterious. A related number is $\D(P,Q)$, the minimal power of $v$ appearing
in $h_{1,w(P,Q)}$. It is tempting to view these as some measure of the distance between two standard tableaux, but this must be taken with a grain of salt. Set $d(P,Q) = \D(P,Q) -
\rbb(\l)$. It is true that $d(P,Q)$ takes positive values and is zero precisely when $P = Q$, but it fails\footnote{We were unable to find this result in the literature. Counterexamples
for $S_6$ were computed by Benjamin Young in response to our query.} to satisfy the triangle inequality, so it is not a metric. Note that $\D(P,Q)$ controls the ``asympotic'' distribution
of KL basis elements appearing in the half twist (e.g. in which cohomological degree a given bimodule appears for the first time in the half twist complex). The distribution of
indecomposables in the full twist is equally mysterious. \end{remark}

\subsection{Sch\"utzenberger in type $A$}
\label{subsec:more}

\begin{prop} \label{prop:Schutz}  There is an involution $P\mapsto P^\vee$ on the set of SYT \revise{of shape $\l$}, called the \emph{Sch\"utzenberger involution}, such that if $x = w(P,Q,\l)$ then
\begin{equation}
w_0 x = w((P^\vee)^t,Q^t,\l^t),
\end{equation}
\begin{equation}
x w_0 = w(P^t,(Q^\vee)^t,\l^t).
\end{equation}
The $w_0$-twisted involutions are therefore precisely the elements in $S_n$ of the form $w(P^\vee,P)$. \end{prop}

We let $\tau$ denote the Dynkin diagram automorphism of $S_n$, which can also be realized as conjugation by the longest element. The above implies that if $x = w(P,Q,\l)$, then $\tau(x) = w(P^\vee,Q^\vee,\l)$.

\begin{ex} In $S_4$, with simple reflections $\{s,t,u\}$, the $w_0$-twisted involutions are: $w_0$ in the longest cell; $tstut$, $tutst$ and $sutsu$ in the subminimal cell; $su$ and
$tsut$ in the middle cell; $stu$, $t$, and $uts$ in the simple cell; and the identity in the identity cell. Comparing this with \eqref{eq:Fstsuts}, we see that the $w_0$-twisted
involutions in cell $\l$ are precisely the terms in cell $\l$ which appear in homological degree $\cbb(\l)$ in the half twist. This was proven in general in Proposition \ref{prop:htm}.
\end{ex}

The Sch\"utzenberger involution has an explicit combinatorial description in terms of jeu-de-taquin\footnote{\revise{For several combinatorial approaches to this involution, see \cite{Gansner}. This involution was proven to agree with the involution on canonical bases in Hecke algebras and relatedly in quantum groups by Berenstein-Zelevinsky \cite{BerenZelev}.}}.  This explicit construction is not yet illuminating to the authors, so we will not recall it, and the reader is welcome to use Proposition \ref{prop:Schutz} as a definition of the Sch\"utzenberger involution.

To relate this construction with the operator $\Schu_L$ defined in \S\ref{subsec:schutz}, we have 
\begin{equation} \Schu_L(w(P,Q)) = w(P^\vee,Q). \end{equation} Consequently, \eqref{eq:schugeneral} becomes
\begin{equation} \label{eq:longesttimescell} H_{w_0} b_{(P,Q,\l)} \equiv (-1)^{\cbb(\l)} v^{\xbb(\l)} b_{(P^\vee, Q, \l)} + \HB_{< \l}. \end{equation}

\begin{ex} In $S_3$, we have already seen the categorification of \eqref{eq:longesttimescell} in \S\ref{subsec:HT3}. \end{ex}

\subsection{Dots in type $A$}
\label{subsec:dotsA}
Recall that every involution in $S_n$ is distinguished. \revcomment{Changed $T$ to $P$ throughout.}

\begin{defn} Let $\l \in \PC(n)$ and $P \in \SYT(\l)$.  If $w=w(P,P,\l)$ is an involution, then let $\xi_P = \xi_{w}:R\rightarrow B_w(\rbb(\l))$ and $\tbarb{P} = \tbarb{w}= \xi_P^\ast \circ \xi_P : R\rightarrow R(2\rbb(\l))$ be as in Definition \ref{def:dots}.
\end{defn}

\revise{Recall that $\xi_P$ is a basis vector in a one-dimensional Hom space, so it and $\tbarb{P}$ are well-defined up to invertible scalar.} For those $P$ such that $w(P,P,\l)$ is the longest element of a parabolic subgroup, these ``thick dots'' and ``thick barbells'' were already studied in type $A$ in \cite{EThick}.

For the longest element $w_0 \in S_n$, letting $R = \RM[x_1, \ldots, x_n]$ and $\ell = \ell(w_0)$, we set $\xi_{w_0}$ to be the map $R \to B_{w_0}(\ell) = R \ot_{R^{S_n}} R(2\ell)$ sending
\begin{equation}
1\mapsto \prod_{1 \le i<j \le n}(x_i\otimes 1 - 1\otimes x_j).
\end{equation}
The dual map $\xi_{w_0}^\ast \co B_{w_0}\rightarrow R(\ell)$ sends $1\otimes 1\mapsto 1$.  Their composition is the polynomial
\begin{equation} \label{eq:barbellA}
\tbarb{w_0}=\prod_{1 \le i<j \le n}(x_i -  x_j) \in R,
\end{equation}
which is the product of the positive roots.

\section{Twists in type $A$}
\label{sec:typeAtwists}

The eventual goal of this chapter is to prove Conjecture \ref{conj:HTaction} in type $A$.

\subsection{The braid group in type $A$}
\label{subsec:braid}

The braid group associated to $S_n$ is the usual braid group on $n$ strands. It has invertible generators $\{\s_i\}_{i=1}^{n-1}$ called \emph{(positive) crossings} corresponding to the
simple reflections $s_i \in S$. We may write $\s$ or $\s_s$ to denote one of these generators, the one corresponding to the simple reflection $s$. These generators satisfy the braid
relations, but do not satisfy $\s^2 = 1$.


Given an element $w \in S_n$, the \emph{positive lift} to the braid group is obtained by taking $\s_{i_1} \s_{i_2} \cdots \s_{i_d}$, whenever $\un{w} = (s_{i_1}, \ldots, s_{i_d})$ is a rex for $w$. This lift is independent of the choice of rex.

Let $\hT_n$ denote the positive lift of the longest element $w_0$ in $S_n$. Let $\fT_n$ denote $\hT_n^2$. These are called the \emph{half twist} and \emph{full twist} respectively. It is known that $\fT_n$ generates the center of $\Br_n$.

Specifically, the half twist is the braid $\hT_n=\s_1(\s_2\s_1)\cdots (\s_{n-1}\cdots\s_2\s_1)$ and the full twist is $\fT_n=\hT_n^2$.  Graphically, these may be pictured as follows:
\[
\hT_4 = \ig{1}{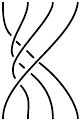},\qquad\qquad \fT_4 = \ig{1}{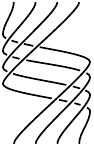}.
\]

The braid group admits a homomorphism to the symmetric group sending $\s_s$ to $s$, and admits a homomorphism to the (unit group of the) Hecke algebra sending $\s_s$ to $H_s$.

\subsection{The external product}
\label{subsec:external}

\revcomment{Revised this paragraph to remove reference to string diagrams}
We use $\sqcup$ to denote the \emph{external product} $S_i \times S_j \to S_{i+j}$, which sends simple reflections to simple
reflections in the obvious way. We also use $\sqcup$ for the corresponding map of braid groups; the braid $\beta_1 \sqcup \beta_2$ has diagram given by the horizontal concatenation of the braid diagrams for $\beta_1$ and $\beta_2$.

The external product on braid groups induces an external product of Hecke algebras. Letting $\HB_i$ denote $\HB(S_i)$, we have $\sqcup \co \HB_i \times \HB_j \to \HB_{i+j}$. The
following proposition is a straightforward consequence of the definitions.

\begin{proposition}
We have $H_{w\sqcup x} = H_w \sqcup H_x$ and $b_{w\sqcup x} = b_w\sqcup b_x$ for all $w\in S_i$ and all $x\in S_j$.
\end{proposition}

Similarly, we may let $R_i = \RM[x_1,\ldots,x_i]$. We have the familiar isomorphism $R_i\otimes_\RM R_j \cong R_{i+j}$. Thus, tensoring over $\RM$ gives rise to a bilinear functor
$\SBim_i \times \SBim_j \rightarrow \SBim_{i+j}$, which we continue to denote by $\sqcup$, and call the \emph{external product}. In the literature this functor is often denoted by
$\boxtimes$. We use the notation $\one_i$ for the monoidal identity inside $\SBim_i$. If $A$ is an object of $\SBim_i$, we may write $A$ or $(A \sqcup \one_j)$ for the image of $A$ inside $\SBim_{i+j}$.

The external product is monoidal: if $A,A'\in \SBim_i$ and $B,B'\in \SBim_j$ then $(A\sqcup B)\otimes (A'\sqcup B')\cong (A\otimes A')\sqcup (B\otimes B')$. Proving this is an easy
exercise, and holds more generally for the external tensor product of any algebras over a field. \revise{An implication is that if $A \in \SBim_i$ and $B \in \SBim_j$, then 
\begin{equation}\label{sqcupinterchange} (A \sqcup \one_j) \ot (\one_i \sqcup B) \cong (\one_i \sqcup B) \ot
(A \sqcup \one_j).\end{equation} If one were to abuse notation and identify $A$ with its image $A \sqcup R_j$ in $\SBim_{i+j}$, and similarly for $B$ and $R_i \sqcup B$, then we can abbreviate this to $A \ot B \cong B \ot A$. In other words, objects in the images of $\SBim_i$ and $\SBim_j$ tensor commute with each other (functorially, even).}

\begin{remark} When working with Soergel bimodules in finite characteristic, it is better to use the diagrammatically defined category from \cite{EWsoergelCalc} than the actual category of bimodules. Nonetheless, all the properties of the external product that we use can be proven using easy diagrammatic arguments. \end{remark}

\begin{remark} \revcomment{Revised this remark} The reader familiar with diagrammatics may think of $\otimes$ as being the horizontal concatenation of (diagrams representing) morphisms. This should not be confused with $\sqcup$, despite both being associated with horizontal concatenation. These concatenation operations are ``monoidal structures at different categorical levels.'' \end{remark}

\subsection{Bimodule sliding}
\label{subsec:sliding}

Throughout out the remainder of this chapter we will be illustrating certain statements diagrammatically.  Each of the diagrams is meant to indicate a complex of Soergel bimodules.  For instance, the Rouquier complex associated to a braid $\b$ will simply be drawn using usual pictures for braids.  If $w_0\in S_n$ is the longest element, then $B_{w_0}$ will be drawn by $n$-strands merging into one, and then splitting back into $n$.  For example $B_{w_0}\in \SBim_4$ would be pictured as
\[
B_{w_0} = \eqig{.9}{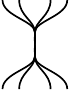}.
\]
The external tensor product $\sqcup$ is drawn by placing diagrams side-by-side, and $\otimes$ is indicated by vertical concatenation.

\begin{lemma} \label{lem:BFcommute} $F(\s_i) \ot B_j \cong B_j \ot F(\s_i)$ whenever $|i-j|>1$. \end{lemma}
This identity is pictured schematically below (in case $i<j$):\vskip4pt
\begin{equation} \label{eq:BFcommute}
\eqig{.9}{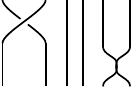} \ \ \ \simeq \ \ \ \eqig{.9}{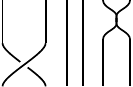}.
\end{equation}\vskip4pt

\begin{proof} This follows from the fact that $\sqcup$ is a monoidal functor $\SBim_m\times \SBim_n\rightarrow \SBim_{m+n}$, and hence induces a monoidal functor
$\KC^b(\SBim_m) \times \KC^b(\SBim_n) \to \KC^b(\SBim_{m+n})$. \revadd{Said another way, this is just a consequence of \eqref{sqcupinterchange}.} \end{proof}

\begin{lemma}\label{lem:slidingBim}
$F(\sigma_i\sigma_{i+1})\otimes B_i \simeq   B_{i+1}\otimes F(\sigma_i\sigma_{i+1})$.
\end{lemma}
This is pictured schematically below:
\[
\eqig{.8}{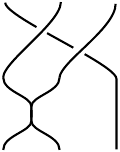} \ \ \ \simeq \ \ \ \eqig{.8}{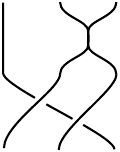}.
\]

\begin{proof}
Let $s = s_i$ and $t = s_{i+1}$. A straightforward computation shows that
\[
F(\sigma_s\sigma_t)\otimes B_s \simeq \Big(\underline{B_{sts}}\rightarrow B_{ts}(1)\Big) \simeq B_t\otimes F(\sigma_s\sigma_t).
\]
The differential is the composition $B_{sts}\subset B_sB_tB_s\rightarrow B_t B_s = B_{ts}$, where the dot map is applied to the \revise{first} tensor factor $B_s$.
\end{proof}

\begin{lemma}\label{lem:comparingfunctors}
Let $\FC_1,\FC_2$ be additive graded functors from $\SBim_n$ to $\AC$ for some additive graded category $\AC$ which has unique direct sum decompositions.  If $\FC_1(B)\cong \FC_2(B)$ for every Bott-Samelson bimodule $B\in \SBim_n$, then $\FC_1(B)\cong \FC_2(B)$ for every $B\in \SBim_n$.
\end{lemma}

\begin{proof} It is enough to prove the statement for $B = B_w$ indecomposable, which we prove by induction on the Bruhat order. For $\ell(w) \le 1$, $B_w$ is a Bott-Samelson, so there
is nothing to prove. So fix $w$ with length $\ge 2$, and assume that $\FC_1(B_x) \cong \FC_2(B_x)$ for all $x < w$. Choose a reduced expression $\un{w}$ for $w$. By Theorem
\ref{thm:SCT}, $\BS(\un{w})$ has a unique summand of the form $B_w$, and the remaining summands have the form $B_x(k)$ for $x < w$ and $k \in \Z$. Since $\FC_1(\BS(\un{w})) \cong
\FC_2(\BS(\un{w}))$ and $\FC_1(B_x(k)) \cong \FC_2(B_x(k))$ for all $x < w$, one can cancel summands and deduce that $\FC_1(B_w) \cong \FC_2(B_w)$. \end{proof}

We will apply this lemma shortly to functors from $\SBim_n$ to $\KC^b(\SBim_n)$. We note for this purpose that $\KC^b(\SBim_n)$ is Krull-Schmidt (as is the bounded homotopy category of any Hom-finite Krull-Schmidt additive category, see \cite[Theorem 3.4]{Schnurer-pp11} and \cite[Theorem 6.1]{Shah23}).

\begin{remark}
We warn the reader that there is no statement of naturality in Lemma \ref{lem:comparingfunctors}, so that $\FC_1,\FC_2$ may be non-isomorphic as functors, even if $\FC_1(B)\cong \FC_2(B)$ for all objects $B\in \SBim_n$.
\end{remark}

\begin{proposition}\label{prop:slidingBimodules}
For any $B\in \SBim_k$, we have the following isomorphism in $\KC^b(\SBim_{k+1})$:
\begin{equation}\label{eq:slidingBimodules}
F(\sigma_1\cdots \sigma_k) \otimes (B\sqcup \one_1)\simeq (\one_1\sqcup B)\otimes F(\sigma_1\cdots \sigma_k).
\end{equation}
\end{proposition}

This is pictured schematically below:
\[
\begin{minipage}{1.3in}
\labellist
\small
\pinlabel $B$ at 34 19
\endlabellist
\begin{center}\includegraphics[scale=.8]{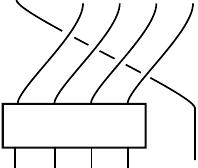}\end{center}
\end{minipage}
=
\begin{minipage}{1.3in}
\labellist
\small
\pinlabel $B$ at 58 59
\endlabellist
\begin{center}\includegraphics[scale=.8]{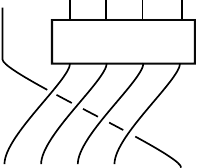}\end{center}
\end{minipage}
\]

\begin{proof} For each $B\in\SBim_k$, let $\FC_1(B)$ and $\FC_2(B)$ denote the left and right-hand sides of \eqref{eq:slidingBimodules}, respectively. Using Lemma
\ref{lem:comparingfunctors} it suffices to show that $\FC_1(B) \cong \FC_2(B)$ when $B$ is a Bott-Samelson bimodule. The result for Bott-Samelson bimodules follows quickly from the
result for $B_s$. This in turn follows from repeated use of Lemmas \ref{lem:BFcommute} and \ref{lem:slidingBim}. \end{proof}

\begin{remark} \label{rmk:provingconjfunctorial} \revcomment{reshuffled details from a comment in a later section to here, altered to discuss recent work} The homotopy equivalence \eqref{eq:slidingBimodules} is natural with respect to $B \in \SBim_k$. That is, \eqref{eq:slidingBimodules} can be upgraded to an isomorphism of functors $\SBim_k \to \KC^b(\SBim_{k+1})$. This was proven in recent work of Mackaay-Miemietz-Vaz, see \cite[Lemmas 4.17 through 4.22]{MacMieVaz}. \end{remark}
	
%
	

\subsection{Conjugation by twists}
\label{subsec:commutehalffull}

\begin{defn} As in \S\ref{subsec:Hw0}, we let $\tau$ denote the involutory automorphism of the symmetric group $S_n$ coming from the Dynkin diagram automorphism, defined by
$\tau(s_i)=s_{n-i}$. Equivalently $\tau(x) = w_0xw_0$ for all $x\in S_n$. We let $\tau$ act on $R$ by $\tau(x_i) = x_{n+1-i}$. We also let $\tau$ denote the
corresponding automorphism of $\SBim_n$ and of $\KC^b(\SBim_n)$. This automorphism sends $B_{s_i}$ to $B_{s_{n-i}}$, and acts on morphisms by the \emph{color swap}: given a diagram (a
decorated colored graph), one replaces the color $s_i$ with the color $s_{n-i}$. \end{defn}

Note that $\tau(F(\s_i)) = F(\s_{n-i}) = F(\tau(\s_i))$. Also note that $\hT_n \s_i = \tau(\s_i) \hT_n$ in the braid group. Thus $\hT_n \b = \tau(\b) \hT_n$ for any braid $\b$.
Consequently, Rouquier canonicity gives an isomorphism \begin{equation} \HT_n \ot F(\b) \simeq \tau(F(\b)) \ot \HT_n\end{equation} for any braid $\b$. Similarly, one has an isomorphism
\begin{equation} \FT_n \ot F(\b) \simeq F(\b) \ot \FT_n, \end{equation} as noted in \S\ref{subsec:RouqCanon}.
	
Now we prove the same result, replacing $F(\b)$ by a complex supported in a single degree.

\begin{corollary}\label{cor:HTcommutes}
We have
\begin{equation}\label{eq:HTconjugation}
\HT_n\otimes B\simeq \tau(B)\otimes \HT_n
\end{equation}
for all $B\in \SBim_n$.
\end{corollary}

\begin{proof} By Lemma \ref{lem:comparingfunctors} it suffices to prove \eqref{eq:HTconjugation} for the Bott-Samelson bimodules, and again, just for $B_s$. This follows easily from
Lemma \ref{lem:BabsorbsF} and Proposition \ref{prop:slidingBimodules}. \end{proof}

\begin{corollary}\label{cor:FTcommutes}
We have $\FT_n\otimes B\simeq B\otimes \FT_n$ for all $B\in \SBim_n$.\qed
\end{corollary}

A consequence of these results is that conjugation by $\HT_n$ and $\FT_n$ are both endofunctors of $\SBim_n$; they preserve complexes which are supported in a single degree.  In the language of \S \ref{subsec:twists}, $\HT_n$ and $\FT_n$ are twist-like.

\subsection{Conjugation by twists, functorially}
\label{subsec:commutehalffulltrue}

\revcomment{Many edits throughout this section.}

We do not use the results of this section, but include them to fend off some common misunderstandings.

The previous section gives the misleading impression that the functors $B\mapsto \HT_n \ot B \ot \HT_n\inv$ and $B\mapsto \tau(B)$ from $\SBim_n$ to $\SBim_n$ should be isomorphic.  Instead, $\HT_n \ot (-) \ot \HT_n\inv$ is isomorphic to $\tau'$, which is obtained by composing $\tau$ with an automorphism of $\SBim_n$ which fixes objects and multiplies certain morphisms by signs\footnote{In diagrammatic language, the end-dot and splitting trivalent vertices should be multiplied by a sign, while the other morphisms are fixed. This leads to the barbell being multiplied by a sign, consistent with the fact that $\tau(x_i - x_{i-1})$ is a negative root, not a positive root.}. The homotopy equivalence in Lemma \ref{lem:BabsorbsF} is not functorial, but involves twisting morphisms with signs, which leads to this result. Since $\tau'$ is an involution, conjugation by the full twist is isomorphic to the identity functor on $\SBim_n$.

\begin{remark} \label{rmk:provingcentral} We prove the aforementioned properties of conjugation by $\HT_n$ in a manuscript in preparation. \end{remark}
	
\begin{remark} \label{rmk:HTnotCentral} Let $W$ be a finite Coxeter group with no diagram automorphisms.  Then the half-twist $H_{w_0}$ is central in the Hecke algebra $\HB(W)$. However, the corresponding Rouquier complex $\HT$ does not functorially commute with all objects, due to signs which appear on morphisms. \end{remark}

The fact that $\FT_n$ commutes functorially (up to homotopy) with all objects in $\SBim_n$ does not imply that it commutes functorially with complexes in $\KC^b(\SBim_n)$! This too we plan to address in a manuscript in preparation, where we construct a functorial isomorphism
\begin{equation} \label{eq:FTreallycommutes} \FT_n\otimes C\simeq C\otimes \FT_n \end{equation}
for all $C \in \KC^b(\SBim_n)$. In other words, $\FT_n$ is an object in the Drinfeld center of $\KC^b(\SBim_n)$.

In this paper we only use the isomorphism \eqref{eq:FTreallycommutes} when $C$ is the Rouquier complex of a braid, or when $C$ is concentrated in a single degree; both of these special
cases were proven in the previous section.


\subsection{Bounding the action of the half twist}
\label{subsec:bounding}

\begin{thm} \label{thm:HTactionA} The half twist $\HT_n$ in type $A$ is twist-like, increasing, and sharp. In other words, conjecture \ref{conj:HTaction} holds in type $A$. \end{thm}

\begin{proof} We have shown that $\HT_n$ is twist-like in \S\ref{subsec:commutehalffull}. We have proven that $\mb_{\HT_n}(\l) = \cbb(\l)$ in Proposition \ref{prop:htm}. What remains to be proven is that $\nb_{\HT_n}(\l) = \cbb(\l)$ for all $\l \in \PC(n)$. Note that $\nb_{\HT_n} \ge \mb_{\HT_n}$ by Lemma \ref{lemma:nmineq}, so we need only prove $\nb_{\HT_n} \le \cbb(\l)$. By Lemma \ref{lemma:tailProperty}, $\nb_{\HT_n}(B)$ only depends on the two-sided cell of an indecomposable object $B$, so it suffices to choose a single indecomposable object $B$ in each cell $\l$, and prove that $\HT_n \ot B$ is supported in homological degrees $\le \cbb(\l)$. Thus, the theorem is deduced from Proposition \ref{prop:HTsupport}. \end{proof}

\begin{prop}\label{prop:HTsupport}
For $\l\in \PC(n)$, let $w_\l$ denote the longest element of the parabolic subgroup $S_\l = S_{k_1}\times \cdots \times S_{k_r}$ where $k_1\geq \cdots \geq k_r\geq 1$ are the column lengths of $\l$. Then $\HT_n\otimes B_{w_\l}$ is homotopy equivalent to a complex in homological degrees $\le \cbb(\l)$.
\end{prop}

Before discussing the proof, let us state the implications of Theorem \ref{thm:HTactionA}.

\begin{theorem}\label{thm:halftwistprops}
In case $W=S_n$, we have $\nb_{\HT}(\l)=\mb_{\HT}(\l)=\cbb(\l)$ for every partition $\l\in \PC(n)$.  If $P,Q\in \SYT(\l)$, then the head of $\HT_n\ot B_{P,Q}$ is isomorphic to $B_{P^\vee,Q}[-\cbb(\l)](\xbb(\lambda))$.  
\end{theorem}

\begin{proof} We can apply Proposition \ref{prop:tailProperty} to determine the shape of $\HT_n \ot B_{P,Q}$. \revise{In particular, modulo cells less than $\l$, $\HT_n \ot
B_{P,Q}$ is isomorphic to its head, which is an indecomposable object in homological degree $\cbb(\l)$. When decategorified, according to \eqref{eq:longesttimescell}, $\hT_n
b_{P,Q} \equiv (-1)^{\cbb(\l)} v^{\xbb(\l)} b_{P^\vee,Q}$ modulo lower cells. These two statements are only compatible if the head is as claimed.} \end{proof}

As an immediate corollary we have the following result on the full-twists.

\begin{corollary}\label{cor:fulltwistprops}
We have $\nb_{\FT}(\l)=\mb_{\FT}(\l)=2\cbb(\l)$.  If $P,Q\in \SYT(\l)$, then the head of $\FT_n\otimes B_{P,Q}$ is isomorphic to $B_{P,Q}[2\cbb(\l)](2\xbb(\lambda))$.
\end{corollary}

Now we get to the proof of Proposition \ref{prop:HTsupport}. In this section we reduce the proposition to a key lemma. In the next section, we prove this lemma. 

Let $w_k$ denote the longest element of $S_k$. Let $T$ and $T^\vee$ be the tableaux of shape $(2,1^{k-1})$ such that
\[
w_k \sqcup 1_1 = w(T,T) \ \ \ \ \ \ \ \ 1_1\sqcup w_k = w(T^\vee,T^\vee).
\]
In $T$, the box in the second column is labeled $k+1$, while in $T^\vee$ this box has label $2$. Note that $T^\vee$ is the Sch\"utzenberger dual of $T$.

\begin{lemma}\label{lem:thickCrossing}
Using the notation of the previous paragraph, we have
\begin{equation}\label{eq:thickcrossing}
F(\sigma_1\cdots\sigma_k)\otimes(B_{w_k}\sqcup \one_1) \ \ \simeq \ \ (\underline{B_{w_{k+1}}}\rightarrow B_{w(T^\vee,T)}(1))
\end{equation}
where $w_{k+1}\in S_{k+1}$ is the longest element. In particular, it is supported in homological degrees $0$ and $1$.
\end{lemma}

The left hand side of \eqref{eq:thickcrossing} is the complex pictured below (when $k=4$):
\begin{equation}\label{eq:forkslideish}
\eqig{.9}{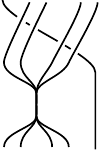}\ \ \simeq \ \ \eqig{.9}{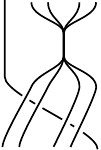},
\end{equation}
where the equivalence of these two complexes is given by bimodule sliding \revise{(Proposition \ref{prop:slidingBimodules})}.

\begin{remark} \revcomment{remark adjusted} This lemma is closely related to results about Rouquier complexes for singular Soergel bimodules. The expert can compare the two-term crossing above with the complex in \cite[Corollary 5.7]{WW09} for $n=k$ and $m=1$, which describes the crossing of a strand of thickness $k$ with a strand of thickness $1$. In theory one can relate this complex with the one in \eqref{eq:thickcrossing} using ``fork-sliding'', which is pictured\footnote{Note that Queffelec-Rose discuss fork-sliding in a different context, nor do they provide a complete proof. For a version of fork-sliding yet another context, see work of Cautis \cite[Lemma 5.2]{CautisClasp}.} in \cite[Equation (4.3)]{QR16}. Singular Rouquier complexes and fork-sliding are not sufficiently established in the literature at this time (and this approach would require a great deal of additional background as well). Instead of relying on fork-lore, we provide a self-contained proof below. \end{remark}
	

The rest of this section is a proof of Proposition \ref{prop:HTsupport}, given Lemma \ref{lem:thickCrossing}.

For each pair of integers $k,\ell$, the \emph{cabled crossing} $x(k,\ell)$ is the element of $S_{k+\ell}$ which crosses the first $k$ strands over the last $\ell$, without permuting either block of strands. In other words, $x(k,\ell)$ is the shortest length element in the coset $w_{k+\ell} (S_k \times S_\ell)$. Let $X(k,\ell)$ denote the Rouquier complex associated to the positive braid lift of $x(k,\ell)$. Explicitly,
\begin{equation}
X(k,\ell) =  F\Big((\sigma_k\cdots\sigma_{k+\ell-1})(\sigma_{k-1}\cdots \sigma_{k+\ell-2})\cdots (\sigma_1\cdots \sigma_\ell)\Big).
\end{equation}
We may shorten this to
\begin{equation}
X(k,\ell) =  F_{k, k+1, \ldots, k+\ell-1} \ot \cdots \ot F_{2,3,\ldots,k+1} \otimes  F_{1,2,\ldots,k}.
\end{equation}
For example
\[
X(3,2) = \eqig{.8}{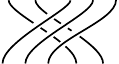}.
\]

Let $\GB(k,\ell)$ denote the complex
\begin{equation}
\GB(k,\ell):= X(k,\ell)\otimes (B_{w_k}\sqcup \one_{\ell})
\end{equation}
where $w_k\in S_k$ is the longest element.  For example
\[
\GB(3,2)= \eqig{.8}{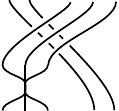}.
\]

\begin{lemma}\label{lem:Gkl}
The complex $\GB(k,\ell)$ is homotopy equivalent to a complex supported in homological degrees between $0$ and $\ell$.
\end{lemma}

\begin{proof}
The main trick in this proof is the fact that, for any $m \ge 1$, $B_{w_k}^{\ot m}$ is just a direct sum of many copies of $B_{w_k}$ (with grading shifts). Thus, if $X(k,\ell) \ot (B_{w_k}^{\ot m} \sqcup \one_\ell)$ can be bounded in homological degree, then so can $\GB(k,\ell)$.

From Lemma \ref{lem:thickCrossing} we know that $\GB(k,1)$ is homotopy equivalent to a complex supported in homological degrees 0,1. This handles the $\ell = 1$ case.

In case $\ell=2$, we claim that
\begin{equation}\label{eq:foobar1}
X(k,2) \ot (B_{w_k}^{\ot 2} \sqcup \one_2) \simeq F_{2,3,\ldots,k+1} \otimes (\one_1\sqcup B_{w_k}\sqcup \one_1)\otimes F_{1,2,\ldots,k}\otimes(B_{w_k}\sqcup \one_2),
\end{equation} where we have applied Lemma \ref{lem:BFcommute} to slide one $B_{w_k}$ past $F_{1, 2, \ldots, k}$. But \eqref{eq:foobar1} can clearly be rewritten as 
\begin{equation}
(\GB(k,1)\sqcup \one_1)\otimes (\one_1\sqcup \GB(k,1)),
\end{equation}
pictured for $k=4$ below:
\[
\eqig{.8}{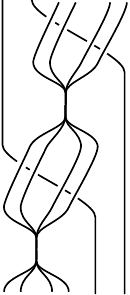}.
\]
By the $\ell = 1$ case, this is homotopy equivalent to a complex supported in homological degrees between 0 and 2. Since $G(k,2)$ is a direct summand of the complex in \eqref{eq:foobar1}, we have proven the $\ell = 2$ case.

In general $\GB(k,\ell)$ is a direct summand of \revcomment{replaced erroneous $B_z$ with $B_{w_k}$.}
\[
\bigotimes_{i=0}^{\ell-1} F_{i+1,\ldots,i+k}\otimes(\one_i\sqcup B_{w_k}\sqcup\one_{\ell-i}) = \bigotimes_{i=0}^{\ell-1} \one_{i}\sqcup \GB(k,1)\sqcup \one_{\ell-i-1},
\]
and application of Lemma \ref{lem:thickCrossing} to each occurence of $\GB(k,1)$ completes the proof.
\end{proof}

Now let $x(k_1, k_2, \ldots, k_r)$ denote the \emph{cabled half twist} in $S_n$, where $n = \sum k_i$. Explicitly, the cabled half twist is the shortest length element in the coset $w_n (S_{k_1} \times \cdots \times S_{k_r})$. When $k_1 \ge k_2 \ge \ldots \ge k_r \ge 1$, and $\l$ is the partition with these column sizes, we let $x_\l$ denote this cabled half twist, and $X_\l$ denote the Rouquier complex associated to its positive lift.

For instance if $\l$ is the partition with column lengths $3,2,2$, then
\[
X_{\l} = \eqig{.8}{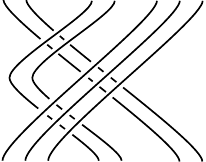}.
\]

\begin{lemma}\label{lem:thickHT}
We have
\begin{equation}
X_\l \ot B_{w_\l} \simeq (\GB(k_{r},0)\sqcup \one_{k_1+\cdots+k_{r-1}})\otimes (\GB(k_{r-1},k_r)\sqcup \one_{k_1+\cdots+k_{r-2}})\cdots \otimes\GB(k_1,k_2+\cdots+k_r)
\end{equation}
\end{lemma}
\begin{proof}
This is straightforward, and best illustrated by example.  For instance, when $\l=(3,3,2)$ we have $(k_1,k_2,k_3)=(3,2,2)$, and
\begin{equation}
\eqig{.8}{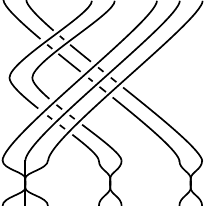} \  \ \ \simeq \ \ \   \eqig{.8}{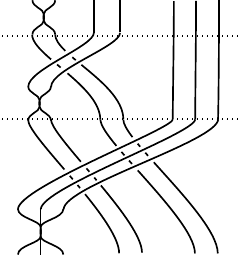}.
\end{equation}
The picture on the left denotes $X_\l \otimes B_{w_\l}$, and the picture on the right is $(\GB(2,0)\sqcup \one_{5})\otimes (\GB(2,2)\sqcup \one_{3})\otimes \GB(3,4)$.  These complexes are homotopy equivalent by bimodule sliding (Proposition \ref{prop:slidingBimodules}).
\end{proof}

Recall that $k_1,\ldots,k_r$ are the column lengths of $\l$, so that
\[
(k_2+\cdots+k_n)+(k_3+\cdots+k_n)+\cdots + k_n = \cbb(\l).
\]
Consequently, the following corollary is immediate from Lemma \ref{lem:thickHT} and Lemma \ref{lem:Gkl}.

\begin{cor} \label{cor:thickHT}
The complex $X_\l \ot B_{w_\l}$ is supported in homological degrees between $0$ and $\cbb(\l)$. \qed\end{cor}

\begin{proof}[Proof of Proposition \ref{prop:HTsupport}]
We observe that $w_0 = x_\l w_\l$ with $\ell(w_0) = \ell(x_\l) + \ell(w_\l)$, so that $\HT_n \simeq X_\l \ot F_{w_\l}$. Using Lemma \ref{lem:BabsorbsF}, $F_{w_\l} \ot B_{w_\l} \simeq B_{w_\l}(-\ell(w_\l))$, and no homological shift is created. Therefore,
\begin{equation} \HT_n \ot B_{w_\l} \simeq X_\l \ot B_{w_\l}(-\ell(w_\l)), \end{equation}
which is supported in homological degrees between $0$ and $\cbb(\l)$ by Corollary \ref{cor:thickHT}. This proves Proposition \ref{prop:HTsupport}, given Lemma \ref{lem:thickCrossing}.
\end{proof}

\subsection{Proof of the lemma}
\label{subsec:boundinglemma}

To prove Lemma \ref{lem:thickCrossing} we must prove
\begin{equation} \label{eq:tC}
F(\sigma_1\cdots\sigma_k)\otimes(B_{w_k}\sqcup \one_1) \ \ \simeq \ \ (\underline{B_{w_{k+1}}}\rightarrow B_{w(T^\vee,T)}(1))
\end{equation}
Here, $w_k \in S_k$ is the longest element, and $T, T^\vee$ are tableau describing the involutions $w_k \sqcup 1_1$ and $1_1 \sqcup w_k$ respectively.

\revcomment{This section has been mostly rewritten, and new notation is introduced which modifies the later part of the section}

We begin with some computations in the Hecke algebra. For $x \in W$ define 
\begin{equation} \sigma_x := \sum_{y \le x} v^{\ell(x) - \ell(y)} H_y. \end{equation}
This element is not always self-dual. When $\sigma_x$ is self-dual, then $\sigma_x = b_x$, in which case we say that $x$ is \emph{smooth}\footnote{Sometimes this is called \emph{rationally smooth} in the literature, because it happens if and only if the corresponding Schubert variety is rationally smooth.} (i.e. all the Kazhdan-Lusztig polynomials $h_{y,x}$ are trivial).

Let $z_{k+1} = w_k \sqcup 1_1$, and for $1 \le i \le k$ let
\begin{equation} z_i = s_i s_{i+1} \cdots s_k z_{k+1}. \end{equation}
Then $z_i$ is a maximal length coset representative for an element of $S_{k+1} / (S_k \times S_1)$, such that $z_i(k+1) = i$. In particular, $z_1 = w_{k+1}$, and $z_2 = w(T^\vee,T)$.

\begin{lemma} The elements $z_i$ are smooth, for $1 \le i \le k+1$. Moreover, $b_{s_k} b_{z_{k+1}} = b_{z_{k}}$, and $b_{s_i} b_{z_{i+1}} = b_{z_i} + b_{z_{i+2}}$ for all $1 \le i \le k-1$. \end{lemma}
	
\begin{proof} It is well-known that the longest element of any parabolic subgroup is smooth, implying that $z_{k+1} = w_K \sqcup 1_1$ is smooth. We now prove the result by descending induction on $i$. (Indeed, this gives another proof that $w_{k+1}$ is smooth, so one need not even rely on the well-known result.)

Recall that for a simple reflection $s$ and an element $y \in W$ we have
\begin{equation} b_s \cdot H_y = \begin{cases} H_{sy} + v H_y & \text{ if } sy>y, \\ H_{sy} + v^{-1} H_y & \text{ if } sy < y. \end{cases} \end{equation}
Since no element $y \le z_{k+1}$ satisfies $s_k y < y$, we easily compute that 
\[ b_{s_k} \sigma_{z_{k+1}} = \sigma_{z_k}.\]
Since the left-hand side is self-dual, so is the right-hand side, and thus $b_{z_k} = \sigma_{z_k}$ is smooth.
	
For $1 \le i \le k-1$, it is an exercise to verify that those elements $y < z_{i+1}$ for which $sy<y$ are precisely those elements satisfying $y < z_{i+2}$. It is an exercise to use these facts to verify that
\begin{equation} b_{s_i} \cdot \sigma_{z_{i+1}} = \sigma_{z_i} + \sigma_{z_{i+2}}. \end{equation}
Since every term other than $\sigma_{z_i}$ in this equation is self-dual (by induction), we deduce that $\sigma_{z_i}$ is self-dual. Consequently $b_{z_i} = \sigma_{z_i}$ is smooth. \end{proof}

As a consequence of this lemma, we have
\begin{equation} \label{eq:tensordecompforz}
B_{s_k} \ot B_{z_{k+1}} \cong B_{z_{k}}, \qquad B_{s_i} \ot B_{z_{i+1}} \cong B_{z_i} \oplus B_{z_{i+2}} \quad \text{ for } 1 \le i \le k-1. \end{equation}

Let us explain how one constructs the projection map $B_{s_i} \ot B_{z_{i+1}} \to B_{z_{i+2}}$, which we will call a \emph{pitchfork projection}.

Fix simple reflections $s$ and $t$ with $m_{st} = 3$. Fix an element $z \in W$ with $tz> z$ and $sz<z$. There is a degree $+1$ dot map $B_t \to R$, which when tensored with
$\id_{B_z}$ produces a degree $+1$ map $B_t B_z \to B_z$ which we also call a dot map. There is a one-dimensional space of degree $+1$ morphisms $B_{tz} \to B_z$, which must be
spanned by the inclusion $B_{tz} \to B_t B_z$ followed by the dot map $B_t B_z \to B_z$.

Because $B_s B_z \cong B_z(-1) \oplus B_z(1)$ (see \eqref{eq:sdown}), there is a morphism of degree $-1$ from $B_s B_z$ to $B_z$ (projection to the first factor), and a degree $-1$ map from $B_z$ to $B_s B_z$ (inclusion from the second factor). Thus there is a degree zero morphism $B_s B_t B_z \to B_z$, obtained as the composition
\begin{equation} B_s B_t B_z \to B_s B_z \to B_z \end{equation}
of the degree $+1$ dot map $B_t \to R$ (tensored with identity maps), and the aforementioned degree $-1$ map. Precomposed with the inclusion map $B_{tz} \to B_t B_z$, this is the pitchfork projection $B_s B_{tz} \to B_z$. When precomposed with the degree $+1$ dot map $R \to B_s$ (tensored with $\id_{B_{tz}}$), we obtain the dot map $B_{tz} \to B_z$. Consequently, the pitchfork is nonzero.

There is also a degree zero map in the other direction, constructed similarly: the pitchfork inclusion.

\begin{lemma} \label{lem:pitchforks} Setting $s = s_i$ and $t = s_{i+1}$ and $z = z_{i+2}$, the pitchfork projections and inclusions constructed above are degree zero maps between $B_{s_i} \ot B_{z_{i+1}}$ and $B_{z_{i+2}}$, which compose to be a non-zero scalar multiple of the identity map of $B_{z_{i+2}}$. \end{lemma}
	
\begin{proof} Above we have already demonstrated that these maps are nonzero. When $B_z$ appears with multiplicity $1$ as a direct summand of $B_s B_{tz}$, there is a one-dimensional space of degree $0$ morphisms $B_s B_{tz} \to B_z$. Thus the pitchfork map spans that space, and must agree with any projection map up to nonzero scalar. Similarly, the pitchfork inlcusion is a valid inclusion map. Their composition is therefore a non-zero scalar multiple of the identity map $\id_{B_z}$. \end{proof}

\begin{proof}[Proof of Lemma \ref{lem:thickCrossing}]
We are interested in the complex
\begin{equation} F_{s_1} F_{s_2} \cdots F_{s_k} B_{z_{k+1}}.\end{equation}

By \eqref{eq:tensordecompforz}, $F_{s_k} B_z$ is the two-term
complex \[\Big( \un{B_{z_k}(0)} \to B_{z_{k+1}}(1) \Big). \] Under the isomorphism $B_{z_k} \cong B_{s_k} B_{z_{k+1}}$, the differential is just the dot map $B_{s_k} \to R(1)$ applied to the first tensor factor.

By Lemma \ref{lem:BabsorbsF} we have $F_{s_{k-1}} B_{z_{k+1}}(1) \cong B_{z_{k+1}}(0)$. Thus $F_{s_{k-1}} F_{s_k} B_{z_{k+1}}$ has the form
\begin{equation} \label{eq:foobar22}
\Big( \un{B_{s_{k-1}} B_{z_k} (0)} \to B_{z_{k+1}}(0) \oplus B_{z_k}(1) \Big).
\end{equation}
By \eqref{eq:tensordecompforz}, $B_{s_{k-1}} B_{z_k} \cong B_{z_{k-1}} \oplus B_{z_{k+1}}$. Keeping careful track of the differential in \eqref{eq:foobar22}, the map from $B_{s_{k-1}} B_{z_k}(0) \to B_{z_{k+1}}(0)$ is precisely the pitchfork projection! Thus, one can apply Gaussian elimination to remove the $B_z(0)$ summand from both terms. The result is a two-term complex
\begin{equation}
F_{s_{k-1}} F_{s_k} B_{z_{k+1}} \simeq \Big( \un{B_{z_{k-1}}(0)} \to B_{z_k}(1) \Big).
\end{equation}
Keeping track of the differential is another exercise: it is the inclusion of the direct summand $B_{z_{k-1}}(0)\rightarrow B_{s_{k-1}}B_{s_k}B_{z_{k+1}}(0)$ followed by the dot map $B_{s_{k-1}}B_{s_k}B_{z_{k+1}}(0)\rightarrow B_{s_k}B_{z_{k+1}}(1)$ followed by the projection map $B_{s_k}B_{z_{k+1}}(1) \to B_{z_k}(1)$. This final projection map is just the identity map, but in later stages of the induction it will be a non-trivial projection.

Assume by (descending) induction that
\begin{equation}\label{eq:twoTermCx}
F_{s_i\cdots s_{k-1}s_k}B_{z_{k+1}}\simeq \Big(\un{B_{z_i}} \rightarrow B_{z_{i+1}}(1)\Big).
\end{equation}
where the differential is induced from the dot map $B_{s_i}\rightarrow R(1)$ on the Bott-Samelson bimodules of which these indecomposables are direct summands. We wish to prove the same, replacing $i$ with $i-1$. Once again, we tensor the above two-term complex with $F_{s_{i-1}}$ and use \eqref{eq:tensordecompforz} to decompose the result.  The argument above, involving the pitchfork projection, works almost verbatim prove that the result is homotopy equivalent to 
\[ \Big(\un{B_{z_{i-1}}} \rightarrow B_{z_i}(1)\Big).\]
Thus \eqref{eq:twoTermCx} is true for all $i\in\{1,\ldots,k-1\}$ by induction.

When $i=1$ this gives the desired result.
\end{proof}

\begin{remark} Addendum: Inspired by this Lemma, in more recent work of the first author \cite[Chapters 11.1-11.3]{EGaitsgory}, the idempotents projecting to each indecomposable
bimodule $B_{z_i}$ inside the tensor product $B_{s_i} \cdots B_{s_k} B_{z_{k+1}}$ were computed. It was shown that for all the pitchforks in this inductive decomposition, the
scalar multiple from Lemma \ref{lem:pitchforks} is equal to $(-1) \cdot \id_{B_{z_i}}$. As a consequence, Lemma \ref{lem:thickCrossing} holds integrally, i.e. it applies to Soergel
bimodules in any characteristic. \end{remark}

\section{The eigenmaps}
\label{sec:constructing}

In this chapter we will prove Conjecture \ref{conj:eigenmap} on the existence of special maps $\a_\l:\Sigma_\l(\one)\rightarrow \FT$  in type $A$. 

\subsection{Reminder on $\l$-equivalences, left versus right}
\label{subsec:lambdaEquiv}

Recall the result of Corollary \ref{cor:fulltwistprops}, which states that
\begin{equation}\label{eq:FTwithtail}
\FT_n\ot B \simeq (\text{tail}\to \shift{\l}{B})
\end{equation}
for each $B\in \SBim_n$ in cell $\l$, where the tail consists of terms in cells $<\l$ and homological degrees $<2\cbb(\l)$.  In particular, the inclusion of the term in maximal homological degree is a chain map \[ \iota_B:\shift{\l}{B} \to \FT_n\ot B\] whose mapping cone is homotopy equivalent to the  tail. Our goal is to show that this inclusion map $\shift{\l}{B}\to \FT_n\ot B$ is induced from a map $\shift{\l}{\one}\to \FT_n$ after tensoring with $B$. Before proving this, we develop the theory of such maps.

We have already discussed such maps in \S\ref{subsec:celltri}, where we called them $\l$-equivalences. We restate the definition here, in order to be more precise about left versus right actions.  Recall that $\FT_n \ot B \simeq B \ot \FT_n$, so that there is also a corresponding inclusion map \[\iota'_B \co \shift{\l}{B} \to B \ot \FT_n.\]

\begin{definition}\label{def:lambdaMap} Choose a partition $\l \in \PC(n)$. A chain map $\a \co \shift{\l}{\one} \to \FT_n$ is said to be a \emph{left $\l$-equivalence} if $\Cone(\a_\l)
\ot B$ is homotopy equivalent to a complex in cells strictly less than $\l$, for each $B\in \SBim_n$ in cell $\l$. It is a \emph{right $\l$-equivalence} if $B \ot \Cone(\a_\l)$  is homotopy equivalent to a complex in cells strictly less than $\l$. \end{definition}

It is easy to see one can check whether $\a$ is a (left or right) $\l$-equivalence by checking its defining condition only when $B$ is indecomposable. Assuming that $\End(B)$ is a field for all indecomposable $B$ in cell $\l$ (which is true by Theorem \ref{thm:SC}), Lemma \ref{lem:lequivcriterion} proved that $\a$ is a left $\l$-equivalence if and only if $\a \ot \Id_B$ is a nonzero scalar multiple of $\iota_B$, if and only if $\a \ot \Id_B$ is not null-homotopic. Similarly, it is a right $\l$-equivalence if and only if $\Id_B \ot \a$ is a nonzero scalar multiple of $\iota'_B$.

\begin{ex}
Consider the one-row partition $\lambda=(n)$. The only indecomposable $B$ in cell $\l$ is $\one$. Then $\colsum(\l)=\binom{n}{2}$, $\rowsum(\l)=0$, and the inclusion of the maximal degree chain bimodule $\a \co \one[-n(n-1)](n(n-1)) \rightarrow \FT$ is a (left or right) $\l$-equivalence. This is because $\Cone(\a)$ is homotopy equivalent to a complex which has no summands isomorphic to shifted copies of $\one$, and hence it lies in lower cells.
\end{ex}

\revcomment{Edits to this example to make it connect with [AH17] more directly.}
\begin{ex}\label{ex:one-col-eigen}
Let $\lambda$ be the one-column partition.  Then $B_{w_0}$ is the unique indecomposable bimodule in cell $\l$, and we have $\colsum(\l) = 0$, $\rowsum(\l)=-\binom{n}{2}$.  Moreover, $\l$ is the minimal cell, so a map $\a\colon \one(-\binom{n}{2}) \to \FT$ is a left (resp.~right) $\l$-equivalence iff $\Cone(\a)\otimes B_{w_0}\simeq 0$ (resp.~$B_{w_0}\otimes \Cone(\a)\simeq 0$).  

In \cite[Corollary 2.15]{AbHog17} it is shown that if $Z\in \KC^b(\SBim_n)$ is acyclic then $Z\ot B_{w_0}\simeq 0\simeq B_{w_0}\ot Z$.  Thus if $\a\colon \one(-\binom{n}{2}) \to \FT$ is a quasi-isomorphism, then $\Cone(\a)$ is acyclic, hence $\a$ is both a left and right $\l$-equivalence.  A specific construction of this quasi-isomorphism $\a$ appears in Proposition 3.12 of \cite{AbHog17}. Note that $\a$ will in general not be a homotopy equivalence.
%

\end{ex}

We now prove that the notions of left and right $\l$-equivalence agree, so after this section we will simply write $\l$-equivalence.

\begin{lemma}\label{lemma:lequivsymmetry}
Let $B_{T,U}\in \SBim_n$ be in cell $\l$, and let $\a:\shift{\l}{\one}\to \FT_n$ be a chain map.  Then the following are equivalent:
\begin{enumerate}\setlength{\itemsep}{3pt}
\item $\Cone(\a)\ot B_{T,U}$ is in cells $<\l$ up to homotopy.
\item $B_{U,T}\ot \Cone(\a)$ is in cells $<\l$ up to homotopy.
\end{enumerate}
Consequently, $\a$ is a left $\l$-equivalence if and only if it is a right $\l$-equivalence.
\end{lemma}

\begin{proof}
Consider the tensor product
\begin{equation}\label{eq:lequivsym}
X = B_{U,T} \ot \Cone(\a) \ot B_{T,U}.
\end{equation}
If either (1) or (2) holds, then the complex $X$ must live in cells $< \l$ up to homotopy. Suppose that (2) holds but (1) fails. Then, by \revise{Lemma \ref{lem:lequivcriterion}}, we know that $\a \ot \Id_{B_{T,U}} \simeq 0$, and $\Cone(\a) \ot B_{T,U}\cong \shift{\l}{B_{T,U}}[1]\oplus (\FT_n\otimes B_{T,U})$.  Tensoring on the left with $B_{U,T}$, we see that the complex \eqref{eq:lequivsym} has a direct summand of the form $B_{U,T} \ot B_{T,U}$, which in turn has a direct summand of the form $B_{U,U}$ (up to shifts).  This contradicts the fact that $X$ is homotopy equivalent to a complex in cells $< \l$. A similar contradiction arises if we assume (1) holds but (2) fails. \end{proof}

\begin{remark} A similar proof will show that left and right $\l$-equivalence are equivalent, for any two-sided cell in a finite Coxeter group, assuming that $\FT$ is twist-like, sharp and increasing. \end{remark}

\subsection{The Specht module and barbells}
\label{subsec:specht}

\revcomment{Changed $T$ to $P$ throughout.}

Now we need a slight detour, which will eventually lead to Corollary \ref{cor:sufficientcondition}, a simple criterion for the existence of a $\l$-equivalence.

Fix $n \ge 2$, and let $R = \RM[x_1, \ldots, x_n]$ with its action of $S_n$. All partitions $\l$ will have $n$ boxes. For a standard tableau $P$, recall the definition of the $P$-barbell
$\tbarb{P}$ from \S \ref{subsec:dotsA}.

\begin{defn} \label{def:barbellSpecht} Let $\SM_\l'$ be the $\RM$-span of polynomials $\tbarb{P}\in R$, as $P$ ranges over all $P\in \SYT(\l)$.\end{defn}

Meanwhile, here is a definition dating back to Specht \cite{Specht35}, see also \cite{Peel75}.

\begin{defn}\label{def:specht}
For $P$ any tableau of shape $\l$, with entries $\{1, 2, \ldots, n\}$ but not necessarily standard,  define $$g_P := \prod (x_i - x_j) \in R$$ where the product is over pairs of indices
$1\leq i, j\leq n$ with $i$ above $j$ in the same column of $P$. Let $\SM_\l$ denote the $\RM$-span of the $g_P$, as $P$ ranges over all tableaux of shape $\l$. This is the \emph{Specht module}\footnote{The Specht module has various other constructions, but this is the original construction due to Specht.} of $\l$.
\end{defn}

Observe that $P$ is uniquely determined by $g_P$. Also note that $w(g_P) = g_{w(P)}$ for $w \in S_n$, where $w$ acts on $P$ by permuting the box labels. Hence $\SM_\l$ is naturally an $S_n$-representation. The following is a theorem of Specht.

\begin{thm} The $S_n$-representation $\SM_\l$ is irreducible, and it has a basis given by $g_P$ for standard Young tableaux $P \in \SYT(\l)$. \end{thm}

The main result of this section is the following.

\begin{thm}\label{thm:specht}
We have $\SM_\l'=\SM_\l$ as subspaces of $R$.
\end{thm}

The proof follows from a series of Lemmas.  First, we warn the reader that $\tbarb{P}\neq g_P$ in general, as the following example shows.

\begin{ex} \label{ex:fPnotgP} Let $\l = (2,2)$ and let $P$ be the unique standard tableau with $g_P=(x_1-x_3)(x_2-x_4)$. One has $w(P,P) = tsut \in S_4$, and one can compute that $\tbarb{P}$ is $(x_1-x_4)(x_3 - x_2)\neq g_P$. \end{ex}
	
However, there is one case where $\tbarb{P}$ and $g_P$ automatically agree.

\begin{lemma} When $P = P_{\col}$ is the column-reading tableau, one has $\tbarb{P} = g_P$. \end{lemma}
	
\begin{proof} As noted in \eqref{eq:barbellA}, $\tbarb{P}$ is the product of the positive roots in the parabolic subgroup $S_\l$. Clearly, so is $g_P$. \end{proof}

For each $\l$, let $I_{\leq\l}\subset R$ denote the homogeneous ideal generated by polynomials which factor through Soergel bimodules in cells $\le \l$.

\begin{lemma}\label{lemma:specht1}
The ideal $I_{\leq\l}\subset R$ is supported in degrees $\geq 2\rbb(\l)$.  The degree $2\rbb(\l)$ component equals $\SM_\l'$.
\end{lemma}

\begin{proof}
In Proposition \ref{prop:dotspace}, we proved that if $w\in S_n$ is in  cell $\l$, then $\Homg(\one,B_w)$ is supported in degrees $\geq \rbb(\l)$.  Further, the degree $\rbb(\l)$ component is one dimensional if $w$ is an involution and zero otherwise.   Consequently, any endomorphism of $\one$ which factors through an indecomposable $B_w$ in cell $\l$ will have degree at least $2\rbb(\l)$, with equality only if $w$ is an involution.  It follows that the degree $2\rbb(\l)$ component of $I_{\leq\l}$ is spanned by the $w$-barbells $\tbarb{w}$.

Moreover, if $\mu<\l$ then any endomorphism of $\one$ factoring through an indecomposable in cell $\mu$ will have degree at least $2\rbb(\mu)$, which is strictly larger than $2\rbb(\l)$. Thus the ideal $I_{\leq \l}$ is supported in degrees $\geq 2\rb(\l)$.  This completes the proof.
\end{proof}

\begin{lemma}\label{lemma:specht2} The ideal $I_{\le \l}$ is closed under the action of $W$. Consequently, $\SM_{\l}'$ is an $S_n$-representation. \end{lemma}
	
\begin{proof}  Let $\IC_{\leq \l}$ denote the ideal of morphisms in $\SBim_n$ factoring through bimodules in cells $\leq \l$.   The ideal $\IC_{\le \l}$ is closed under tensor products with arbitrary morphisms in $\SBim_n$.  Pick a simple reflection $s$. If $f\in R$ is in $I_{\le \l}$, then so is the morphism $c_s(f)$ given by
\[
R \to B_s \otimes B_s \buildrel m_f \over \longrightarrow B_s \otimes B_s \to R,
\]
where $m_f$ denotes middle multiplication by $f$, and the remaining maps are the units and counits of biadjunction for $B_s$.   Here, middle multiplication by $f$ is the endomorphism of $B\otimes_R B'$ sending
\[
b\otimes b'\mapsto bf\otimes b' = b\otimes fb',
\]
where $B$ and $B'$ are $(R,R)$ bimodules and $\otimes = \otimes_R$, as usual.

Diagrammatically, $c_s(f)$ corresponds to placing $f$ inside a circle colored $s$. A simple computation involving \cite[Equation 5.2]{EWsoergelCalc} will imply that $c_s(f)=f-s(f)$. So $f - s(f) \in I_{\le \l}$, implying that $s(f) \in I_{\le \l}$. This proves that $I_{\le \l}$ is closed under the action of $S_n$. \end{proof}

\begin{proof}[Proof of Theorem \ref{thm:specht}]
We know that $\SM_\l$ is defined as the span of the $S_n$-orbit of $g_{P_{\col}}$, that $g_{P_{\col}} = \tbarb{P_{\col}}$, and that $\SM_\l'$ is closed under the action of $S_n$. Hence $\SM_\l \subset \SM_\l'$. On the other hand, we know that the dimension of $\SM_\l'$ is at most the size of $\SYT(\l)$, which is the dimension of $\SM_\l$. We conclude that $\SM_\l'=\SM_\l$ by a dimension count.
\end{proof}

\begin{remark}
For any two-sided cell $\l$ in any Coxeter group $W$, one can define a vector space $\SM_\l'$ as the span of $\tbarb{d}$ over all distinguished involutions in $\l$. The proofs above work verbatim to show that $\SM_\l'$ is closed under the action of $W$.  As far as we are aware, these analogues $\SM_\l'$ of Specht modules have not appeared elsewhere in the literature! We conjecture that the $d$-barbells form a basis for $\SM_\l'$.
\end{remark}

To avoid primes we write $\SM_\l$ instead of $\SM_\l'$ below, but the crucial point is that it is spanned by barbells.

\subsection{A sufficient condition for existence of $\l$-equivalences}
\label{subsec:lambdaEquivsufficient}

For each partition $\l$ recall the grading shift functor $\Sigma_\l=[-2\cbb(\l)](2\xbb(\l))$.

\begin{lemma} If $\a\co\shift{\l}{\one}\to \FT_n$ is a chain map and $\tbarb{P} \cdot \a$ is not null-homotopic, then $\Id_B \ot \a$ and $\a \ot \Id_B$ are not null-homotopic, where $B = B_{P,P}$. \end{lemma}

\begin{proof}  The reader may wish to recall the trivialities of \S\ref{subsec:trivial}, in particular that the left and right actions of $R$ on $\Hom^{\Z\times \Z}(\one,C)$ coincide for all complexes $C$.   Then since $\tbarb{P}$ is the composition of maps
\[
\one\to B_{P,P}(\rbb(\l)) \to \one(2\rbb(\l)),
\]
it follows that $\tbarb{P} \cdot \a=\tbarb{P}\otimes \a = \a\otimes \tbarb{P}$ factors through $\Id_B \ot \a$ and through $\a \ot \Id_B$.   If $\Id_B\otimes \a$ or $\a\otimes \Id_B$ were null-homotopic, then so would be $\tbarb{P}\cdot \a$.\end{proof}
	
\begin{cor} Suppose that a chain map $\a\co\shift{\l}{\one}\to \FT_n$ satisfies the condition that $\tbarb{P} \cdot \a$ is not null-homotopic, for all $P \in \SYT(\l)$. Then $\a$ is a
$\l$-equivalence. \end{cor}

\begin{proof} This follows from the above lemma and Lemma \ref{lem:lequivcriterion}.\end{proof}

Recall the action of the braid group on $\Homgg(\one, \FT_n)$, given by conjugation, as discussed in \S\ref{subsec:conjugate}.  When we refer to the conjugate of a map $\a \co \one \to \FT_n$, we refer to the map $\psi_\b(\a)$ for some braid $\b$.

\begin{lemma}\label{lem:criterion} Fix $\l \in \PC(n)$ and a chain map $\a\co\shift{\l}{\one}\to \FT_n$. Suppose there is a nonzero polynomial $g \in \SM_\l$ such that $g \cdot \a$ is not null-homotopic. Then one can find a linear combination $\a'$ of conjugates of $\a$, such that $\tbarb{P} \cdot \a'$ is not null-homotopic for all $P \in \SYT(\l)$. \end{lemma}
	
\begin{proof} We first claim that for each nonzero polynomial $f \in \SM_\l$ there exists a linear combination $\gamma$ of conjugates of $\a$ such that $f\cdot \g\not\simeq 0$. Indeed, let
\[
X := \{f \in \SM_\l \:|\: f\cdot \psi_\b(\a) \simeq 0 \textrm{ for all } \b\in \Br_n\}.
\]
Clearly $X$ is a subspace of $\SM_\l$, and it is proper because $g \notin X$.  We wish to show that $X$ is preserved by $S_n$, in which case it must be $0$ by the irreducibility of $\SM_\l$. However, we have seen in \eqref{eq:conjugation1} that
\[
\psi_\b(f\cdot \g)\simeq w(f) \cdot \psi_\b(\g)
\]
for any polynomial $f$, map $\g$, and braid $\b$, where $w$ is the element of the symmetric group corresponding to $\b$. Thus, if $f \cdot \psi_\b^{-1} (\gamma) \simeq 0$ then $w(f) \cdot \g \simeq 0$, for a braid lifting $w$. In particular, if $f \in X$ then $w(f) \in X$. Thus we conclude that $X = 0$.

Now let $V\subset \Homgg(\one,\FT_n)$ be the vector space spanned by the classes of the conjugates $\psi_\b(\a)$, as $\b$ ranges over all elements in $\Br_n$. For each $P\in \SYT(\l)$, let
$V_P\subset V$ denote the subspace consisting of those homotopy classes $[\g]$ such that $\tbarb{P}\cdot \g \simeq 0$. If $V_P = V$ then $\tbarb{P} \in X$, a contradiction. Thus each $V_P$ is a proper subspace of $V$. Since there are only finitely many such $P$, the union of the $V_P$ over $P \in \SYT(\l)$ can not be all of $V$. Choosing any element $\a' \in V\setminus
\bigcup_{P\in \SYT(\l)}V_P$, we see that $\tbarb{P} \cdot \a' \not\simeq 0$ for all $P \in \SYT(\l)$. \end{proof}

\begin{cor} \label{cor:sufficientcondition} If there is a map $\a \co \shift{\l}{\one}\to \FT_n$ and a polynomial $g \in \SM_\l$ such that $g \cdot \a$ is not null-homotopic, then some linear combination $\a'$ of braid conjugates of $\a$ is a $\l$-equivalence.\qed \end{cor}

\subsection{Addendum: the braid group action on the full twist}
\label{subsec:symgrp}

During the unusually long preparation of this work, the second author and Eugene Gorsky posted a preprint \cite{GorHog17} which can be used to strengthen and simplify Corollary \ref{cor:sufficientcondition}. We felt it was worth including this result, even though it comes from a different geological time period than the rest of this paper.

Briefly stated, the result proves that the braid group action on $\Homgg(\one,\FT_n)$ factors through the symmetric group, and that the important part of any $\l$-equivalence lives in a copy of the sign representation. As a consequence, one need not take a linear combination $\a'$ of conjugates in the statement of Corollary \ref{cor:sufficientcondition}; the original map $\a$ will be a $\l$-equivalence. We believe this is worth noting for future work, although it does not make any appreciable difference to the arguments in this paper, so we continue to use Corollary \ref{cor:sufficientcondition} in the rest of the paper.


\begin{definition}
Let $\Q[x_1,\ldots,x_n,y_1,\ldots,y_n]$ be a bigraded ring with $\deg(x_i) = (2,0)$ and $\deg(y_i)=(-2,2)$, and let $S_n$ act on $\Q[x_1,\ldots,x_n,y_1,\ldots,y_n]$ by permuting both sets of variables. Let $A\subset \Q[x_1,\ldots,x_n,y_1,\ldots,y_n]$ denote the subspace of polynomials which are anti-symmetric with respect to the $S_n$ action.  Let $I\subset \Q[x_1,\ldots,x_n,y_1,\ldots,y_n]$ denote\footnote{The ideal $I$ is denoted $J_n$ in \cite{GorHog17}.} the ideal generated by $A$.
\end{definition}

Recall again from \S\ref{subsec:conjugate} the conjugation action of the braid group on $\Homgg(\one, \FT_n)$. As proven in Lemma \ref{lem:howtoconj}, the action of
polynomials and braids together yields an action of $\Q[\Br_n]\ltimes \Q[x_1,\ldots,x_n]$ on $\Homgg(\one,\FT_n)$, where the action of the braid group on polynomials factors through the symmetric group. Meanwhile, $\Homgg(\FT_n,\one)$ is also a module over $\Q[\Br_n]\ltimes \Q[x_1,\ldots,x_n]$, and one can prove (say, by \cite{LibWil}) that $\Homgg(\FT_n,\one) \cong R$ as $R$-bimodules. So the action of $\Q[\Br_n] \ltimes \Q[x_1, \ldots, x_n]$ on $\Homgg(\FT_n,\one)$ is determined by the action of $\Br_n$ on the one-dimensional space corresponding to $1 \in R$. Let this one-dimensional representation of $\Br_n$ be temporarily denoted $V$. Being a one-dimensional representation, $V$ is invertible in the sense that there exists $V\inv$ such that $V\otimes_\Q V\inv \cong V\inv\otimes_\Q V\cong \Q$, the trivial representation of $\Q[\Br_n]$.

\begin{theorem}\label{thm:ftsymgrp}
The braid group action on $\Homgg(\one,\FT_n) \ot_\Q V\inv$ factors through the symmetric group, and
 \[
 \Homgg(\one,\FT_n) \ot_\Q V\inv \cong I / (y_1,\ldots,y_n)I
 \]
as bigraded $\Q[S_n]\ltimes \Q[x_1,\ldots,x_n]$-modules.
\end{theorem}

\begin{remark} It is expected that $V$ is the trivial representation, though this was not proven in \cite{GorHog17}. In our proofs below, it is irrelevant what $V$ is.   \end{remark}


\revcomment{Removed a remark here.}

\begin{example}
In case $n=2$, $A$ is the $\Q[x_1,x_2,y_1,y_2]^{S_2}$-submodule of $\Q[x_1,x_2,y_1,y_2]$ generated by $\a_1 = x_1-x_2$ and $\a_2 = y_1-y_2$.  Then $I/(y_1,y_2)I$ is the $\Q[x_1,x_2]$-module generated by $\a_1,\a_2$ modulo $(x_1-x_2)\a_2 = 0$, since \[ (x_1 - x_2) \a_2 = (y_1 - y_2) \a_1 \in (y_1,y_2)I.\] These two generators have bidegrees $(2,0)$ and $(-2,2)$, corresponding to maps $\one(-2)[0] \to \FT_2$ and $\one(2)[-2] \to \FT_2$.
\end{example}

\begin{cor}\label{cor:FTmodMaxl}
Let $\mg\subset \Q[x_1,\ldots,x_n]$ denote the maximal graded ideal (generated by $x_1,\ldots,x_n$).   Then modulo $\mg$ every $\a\in \Homgg(\one,\FT_n)$ spans (the zero subspace or) a copy of the sign representation. \qed
\end{cor}

\begin{lemma}\label{lemma:alphaNotInM}
\revcomment{Changed $\b$ in this proof to $\xi$.} Let $\xi \in \mg \Homgg(\one,\FT_n)$ Then $\a\in \Homgg(\one,\FT_n)$ is a $\l$-equivalence if and only if $\a + \xi$ is. In particular, any $\l$-equivalence is not in $\mg \Homgg(\one,\FT_n)$.  
\end{lemma}

\begin{proof}
This is true for degree reasons.   To be precise, recall that $\Sigma_\l(\one)$ is supported in homological degree $2\cbb(\l)$, hence the only nonzero component of $\a$ is $\a^{2\cbb(\l)}$.  In cell $\l$ in this homological degree, $\FT_n$ is a direct sum of $B_d(2\xbb(\l)+\rbb(\l))$, as $d$ ranges over the involutions in cell $\l$.  In fact each component of the composition
\[
\a^{2\cbb(\l)}: \one(2\xbb(\l)) \rightarrow \FT_n^{2\cbb(\l)}\twoheadrightarrow \bigoplus_d B_d(2\xbb(\l)+\rbb(\l))
\]
is a nonzero multiple of the $d$-dot map (\S \ref{subsec:dots}).  The $d$-dots are living in the smallest possible degree, hence any map with smaller bimodule degree must have zero components mapping into $B_d$.  Thus, any morphism $\xi \in \mg \Homgg(\one,\FT_n)$ will factor through cells smaller than or incomparable to $\l$. Consequently, if $B$ is any indecomposable in cell $\l$, the coefficient of the identity of $\Sigma_\l(B)$ in $\xi \ot \Id_B$ is zero, because this morphism lives in cells $< \l$. In particular, $\a$ and $\a + \xi$ have the same coefficient for the identity of $\Sigma_\l(B)$.
\end{proof}

\begin{lemma}\label{lemma:primesnotnecessary}
\revise{Let $\a$ be a chain map $\a \co \shift{\l}{\one}\to \FT_n$.} If some element $\a'$ of $\Q[S_n]\cdot \a\subset \Homgg(\one,\FT_n)$ is a $\l$-equivalence, then $\a$ is already a $\l$-equivalence.
\end{lemma}

\begin{proof} By Lemma \ref{lemma:alphaNotInM}, $\a' \notin \mg\Homgg(\one,\FT_n)$, and so $\a\notin \mg\Homgg(\one,\FT_n)$ since $\mg\Homgg(\one,\FT_n)$ is preserved under the action of
$S_n$. Moreover, by Corollary \ref{cor:FTmodMaxl} $\a$ and $\a'$ are colinear (and nonzero) modulo $\mg\Homgg(\one,\FT_n)$. Any invertible scalar multiple of a $\l$-equivalence is a
$\l$-equivalence, so replacing $\a'$ with a scalar multiple, we can assume $\a$ and $\a'$ are equal modulo $\mg\Homgg(\one,\FT_n)$. But then by Lemma \ref{lemma:alphaNotInM}, $\a$ is a
$\l$-equivalence since $\a'$ is. \end{proof}


Thus we have our improvement upon Corollary \ref{cor:sufficientcondition}.

\begin{cor} \label{cor:sufficientcondition2} If there is a map $\a \co \shift{\l}{\one}\to \FT_n$ and a polynomial $g \in \SM_\l$ such that $g \cdot \a$ is not null-homotopic, then $\a$
is a $\l$-equivalence.\qed \end{cor}

\subsection{Homology of the full twist}
\label{subsec:FTHHH}

\revcomment{reordering of this section}

To each complex $C\in \KC^b(\SBim_n)$ we have a triply-graded vector space $\HHH(C)=\bigoplus_{i,j,k\in \Z}\HHH^{ijk}(C)$ where
\[
\HHH^{ijk}(C)=H^k(\Ext_{(R,R)}^j(\one,C(i)),
\]
and the $\Ext$ groups are taken in the category of graded $(R,R)$ bimodules. That is, one applies the Hochschild cohomology functor $\Ext_{R,R}(\one,-)$ to each homological degree (obtaining a complex of bigraded vector spaces), and then takes the cohomology of this complex (obtaining a triply graded vector space). Here, $k$ is the homological grading, $j$ is the Hochschild grading, and $i$ is the bimodule grading. We are most interested in the bigraded space
\[
\HHH^{\bullet,0,\bullet} =: \Homgg(\one,C),
\]
whose degree $(i,k)$ component consists of chain maps $\one\to C(i)[k]$ modulo homotopy.

For each braid $\b\in \Br_n$, let $\hat{\b}$ denote braid closure, i.e.~ the oriented link in $\R^3$ obtained by connecting the top and bottom boundary points of $\b$ in a planar fashion. 
\begin{remark} When $C=F(\b)$ is the Rouquier complex attached to a braid, Khovanov showed \cite{Kh07} that the triply-graded vector $\HHH(C)$ is isomorphic to the Khovanov-Rozansky homology of $\hat{\b}$ up to overall shift.  In particular, classes in $\Homgg(\one,\FT_n)$ correspond to classes in the Khovanov-Rozansky homology of the $(n,n)$ torus link with Hochschild degree zero.
\end{remark}

The entire triply graded homology $\HHH(\FT_n)$ was computed in the authors' earlier work \cite{ElHog16a}.  Let us recall how this computation is accomplished. \revcomment{additional elaboration begins here, continues for a few paragraphs} We begin by recalling notation relating to twisted complexes, also called convolutions.


If $X$ is a complex with differential $d_X$ and $\e$ is a degree 1 endomorphism of $X$ such that $(d_X+\e)^2=0$, then we write $\tw_\e(X)$ for the complex $X$ with ``twisted differential'' $d_X+\e$.  If $\xi\colon C\to A[1]$ is a chain map then we write $(C\buildrel\xi\over\rightarrow A)$ for the complex $\tw_\e(C\oplus A)$ in which $\e=\smMatrix{0 & 0\\ \xi& 0}$.   Note that $B:=(C\to A)$ fits into a distinguished triangle of the form
\[
A\rightarrow B\rightarrow C\buildrel \d[1]\over \rightarrow A[1].
\]

More generally let $I$ be a finite poset and suppose we have an $I$-indexed family of complexes $\{A_i\}_{i\in I}$ and $B=\tw_\e\Big(\bigoplus_{i\in I} A_i\Big)$.  Assume that the component of $\e$ from $A_i$ to $A_j$ is zero unless $i<j$, so that $\e$ is represented by a strictly lower triangular matrix with respect to the partial order on $I$.  In this case we say that $B$ is a \emph{convolution} of the $A_i$, and we write $B = \tilde{\bigoplus}_{i\in I} A_i$, with $\e$ hidden from the notation. The convolution $B$ has an $I$-indexed filtration, whose subquotients are the original complexes $A_i$.


\begin{notation} \label{not:JMbraid} Let $y_n\in \Br_n$ denote the \emph{Jucys-Murphy braids}, defined by $y_1=\sigma_1^2$ and $y_{n+1} = \sigma_n y_n\sigma_n$.  Let $Y_n$ denote the Rouquier complex associated to $y_n$. \end{notation}

\begin{notation} To streamline the results below, we abbreviate shift functors by writing $q=(-2)$ and $t=(2)[-2]$. \end{notation}

\begin{proposition}[\cite{Hog15}]\label{prop:KnProps}
There exist complexes $\KB_n\in \KC^b(\SBim_n)$ such that
\begin{subequations}
\begin{equation} \label{eq:K1is}
\KB_1 = (q \RM[x_1]\buildrel x_1\over \longrightarrow \underline{\RM[x_1]})
\end{equation}
\begin{equation} \label{eq:KkillsBw}
\KB_n\ot B_w\simeq 0 \simeq B_w\ot \KB_n \ \ \ \ \ \ \text{ for each $w\neq 1\in S_n$}
\end{equation}
\begin{equation}\label{eq:homOneToK}
\Homgg(\one_n,\KB_n) \cong t^{\binom{n}{2}} \RM[x_1,\ldots,x_n] / (x_i-x_j)_{1\leq i<j\leq n}
\end{equation}
\begin{equation} \label{eq:Ktriangle}
(\KB_{n-1}\sqcup \one_1) Y_n \simeq \Big(\KB_{n} \longrightarrow  q\KB_{n-1}\sqcup \one_1\Big).
\end{equation}
\begin{equation}\label{eq:RouquierAbsorbing}
\KB_n\ot F(\b)\simeq t^{e/2}\KB_n\simeq \KB_n\ot  F(\b) \ \ \ \ \ \ \ \  \text{ for all $\b\in \Br_n$}
\end{equation}
\end{subequations}
where $e=e(\b)$ is the \emph{braid exponent} or \emph{writhe} (that is, the signed number of crossings in a diagram for $\b$).
\end{proposition}

\revcomment{This part is abbreviated.}


\revise{The complex $\KB_n$ categorifies the (renormalized) Young symmetrizer associated to the sign representation, as evidenced by \eqref{eq:KkillsBw} and \eqref{eq:RouquierAbsorbing}.} Using the complexes $\KB_k$ for $k \le n$, we can give a very useful expression for the Rouquier complex $\FT_n$.

\begin{definition}\label{def:shuffleAndDv}
Each sequence $v \in \{0, 1\}^n$ with $k$ zeroes determines a permutation $\pi_v\in S_n$---which we will call a \emph{shuffle}---that sends $1,\ldots,k$ to the zeroes of $v$ and $k+1,\ldots,n$ to the ones of $v$ in an order-preserving fashion.  Let $\b_v$ denote the positive braid lift of $\pi_v$ and $\omega(\b_v)$ the positive braid lift of $\pi_v\inv$.  Then set $D_v := F(\b_v)\ot(\KB_k\sqcup \FT_{n-k})\ot F(\omega(\b_v))$.
\end{definition}
For example, here is $D_{01001010}$, which occurs (up to shift) in our expression for $\FT_8$.
 \[
 D_{01001010}=
\begin{minipage}{1.1in}
\labellist
\small
\pinlabel $\FT_3$ at 80 40
\pinlabel $\KB_5$ at 30 40
\endlabellist
\begin{center}\includegraphics[scale=1]{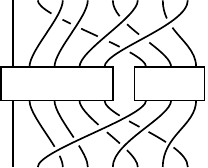}\end{center}
\end{minipage}
 \]
Inside $v$, the zeroes indicate which strands are connected to $\KB_k$, and the ones indicate which are connected to the full-twist $\FT_{n-k}$.  We warn the reader that the conventions for $D_v$ differ by a symmetry from those in \cite{ElHog16a}.  The main result of \cite{ElHog16a} uses the above properties of the complexes $\KB_n$ to prove the following results.

\begin{theorem}\label{thm:ftcomputation}
\revise{Given binary sequences $v,w\in \{0,1\}^n$, write $v<w$ for the usual lexicographic ordering, i.e.~there is an index $j\in \{1,\ldots,n\}$ such that $v_j<w_j$ and $v_i=w_i$ for $1\leq i<j$.  Let $I\subset \{0,1\}^n$ be the subset consisting of sequences with $v_1=0$, endowed with the opposite of the lexicogrphic ordering.  Finally, let $|v|=\sum_i v_i$ be the number of ones in $v$.}  Then we have
\begin{equation}\label{eq:ftexpression}
\FT_n \simeq  \widetilde{\bigoplus}_{v\in I} q^{|v|}D_v
\end{equation}
\end{theorem}

\begin{theorem}\label{thm:parity}
For each $v\in \{0,1\}^n$ the bigraded space of homs $\Homgg(\one,D_v)$ is supported in even homological degrees.  \qed
\end{theorem}

The Poincar\'e polynomial of these spaces are given by an explicit recursion. We do not need the explicit formulae here.  One very important consequence of the parity of these hom spaces is the following isomorphism.  Before stating, observe that the homotopy equivalence \eqref{eq:ftexpression} endows $\Homgg(\one,\FT_n)$ with an $I$-indexed filtration (in which the subspace $\FC^v(\Hom(\one,\FT_n))$ is spanned by homology classes with respresentatives living in summands $\Hom(\one,q^{|w|}D_w)$ with $w\geq v$ in $I$.

\revcomment{Revisions from here to end of section, to remove any errors in associated-graded-fu, and corresponding change in notation in the next section.}

\begin{corollary}\label{cor:FThoms}
We have isomorphisms
\begin{eqnarray}
\Homgg(\one,\FT_n) &\cong& \operatorname{AssGr}(\Homgg(\one,\FT_n))
\label{eq:homOneFT}\\
&\cong& \bigoplus_{v\in I} \Homgg(\one,q^{|v|} D_v)\label{eq:homOneFT 2}
\end{eqnarray}
of bigraded \revise{vector spaces over $\R$, respecting the evident $I$-filtrations.}
\end{corollary}

\begin{proof}
Let $\Homc^{\Z}(\one,\FT_n)$ denote the chain complex of bihomogeneous bimodule maps, whose homology is $\Homgg(\one,\FT_n)$.  Theorem \ref{thm:ftcomputation} gives a finite filtration on $\FT_n$, which gives a finite filtration on $\Homgg(\one,\FT_n)$.  The spectral sequence which computes the homology has $E_2$ page equal to $\bigoplus_v \Homgg(\one,q^{|v|}D_v)$ and  $E_\infty$ page equal to $\operatorname{AssGr}(\Homgg(\one,\FT_n))$.  The $E_2$ page is supported in even homological degrees by Theorem \ref{thm:parity}.  Thus, the differentials $d_r$ must be zero for $r\geq 3$, since $d_r$ increases homological degree by 1.  This shows that the spectral sequence satisfies $E_2 = E_\infty$ which gives the isomorphism \eqref{eq:homOneFT 2}. Since we work over a field, there is no issue with the extension problem in passing from $E_\infty$ to $\Homgg(\one,\FT_n)$.
\end{proof}


To restate the result slightly, the filtration on $\FT_n$ given in \eqref{eq:ftexpression} leads to a filtration on the $R$-module $\Homgg(\one,\FT_n)$, with submodules
\[ F_{\ge v} := \Homgg(\one,\widetilde{\bigoplus}_{w \ge v} q^{|w|} D_w), \]
whose subquotients $F_{\ge v}/F_{> v}$ are isomorphic as $R$-modules to $\Homgg(\one,q^{|v|} D_v)$. Given an element $\b\in \Homgg(\one,q^{|v|}D_v)$, we say that $\a \in F_{\ge v}$ is a \emph{lift} of $\b$ if it descends to $\b$ in the quotient $F_{\ge v}/F_{> v}$. The non-canonical isomorphism between a filtered vector space its associated graded vector space is constructed by choosing lifts of a basis of the associated graded. However, \eqref{eq:homOneFT} is not an isomorphism of $R$-modules! For example, if $\a$ is a lift of $\b$, and $f \in R$, it may be the case that $f \ot \b = 0$ whereas $f \ot \a$ is nonzero, but lives within $F_{> v}$. Nevertheless, we have the following result. 

\begin{lemma} \label{lem:fvsgrf} Let $f \in R$ and $\b\in \Homgg(\one,q^{|v|}D_v)$ be given, and choose any lift $\a \in F_{\ge v}$ (with notation as in the previous paragraph).  If $f\ot \b$ is nonzero, then $f \ot \a$ is nonzero. \end{lemma}

\begin{proof} The quotient map $F_{\ge v} \to F_{\ge v}/F_{> v}$ is an $R$-module map, so evidently $f \ot \a$ descends to $f \ot \b$. Since the latter is nonzero, so too must be the former. \end{proof}

\subsection{The eigenmap theorem}
\label{subsec:eigenmaptheorem}

Let $\l=(\l_1,\ldots,\l_r)$ be given.  Let $n_1 = n - \l_1 = \l_2 + \cdots + \l_r$.  Note that the right-hand side of \eqref{eq:homOneFT} has a unique summand corresponding to the 01-sequence $(0^{\l_1}1^{n_1}) = (00\cdots 011\cdots 1)$, having the form
\[
\Homgg(\one, q^{n_1}\KB_{\l_1}\sqcup \FT_{n_1}) \cong q^{n_1}\Homgg(\one, \KB_{\l_1})\ot_\RM \Homgg(\one,\FT_{n_1}).
\]
In the second isomorphism we used the fact that for any Soergel bimodules $B_i,B_i'\in \SBim_{m_i}$, $i=1,2$, we have
\[
\Hom_{\SBim_{m_1+m_2}}(B_1\sqcup B_2, B_1'\sqcup B_2')\cong \Hom_{\SBim_{m_1}}(B_1,B_1')\otimes_\RM \Hom_{\SBim_{m_2}}(B_2,B_2').
\]
\revcomment{removed a sentence citing diagrammatic calculus}

Iterating this for $\Homgg(\one,\FT_{n_1})$, we see that $\Homgg(\one,\FT_n)$ has a unique summand of the form \revcomment{slight revisions here for overfill too}
\[
q^{\rbb(\l)}\Homgg(\one,\KB_{\l_1})\ot_\RM\cdots \ot_\RM\Homgg(\one,\KB_{\l_r}),
 \]
where $n_i:=\l_{i+1}+\l_{i+2}+\cdots+\l_r$, and $n_r := 0$. This summand is isomorphic to
\begin{equation}\label{eq:specialSummand}
\Homgg(\one, q^{n_1} \KB_{\l_1}\sqcup \cdots \sqcup  q^{n_r}\KB_{\l_r}).
\end{equation}

Now, for our given partition $\l\in\PC(n)$, let $P_{\row}$ be the row-reading tableau of shape $\l$, and let $Z_\l$ denote the set of pairs of indices $i<j$ which appear in the same row of $P_{\row}$.    Several applications of \eqref{eq:homOneToK} tell us that \eqref{eq:specialSummand} is isomorphic to
\begin{equation}\label{eq:specialSummand2}
 q^{\rbb(\l)}t^{\cbb(\l)} \RM[x_1,\ldots,x_n]/(x_i-x_j\:|\: (i,j)\in Z_\l),
\end{equation}
where we have used the fact that  $n_1+\cdots+n_r=\rbb(\l)$ and $\binom{\l_1}{2}+\cdots+\binom{\l_r}{2}=\cbb(\l)$.  This module has a distinguished generator, namely 1, which gives rise to a distinguished element of \eqref{eq:specialSummand} \revadd{which we temporarily denote by $\b_{\l}$}.

\revcomment{minor revisions below until the theorem}

\begin{definition}\label{def:alphaprime}
Let $\a_\l \in \Homgg(\one,\FT_n)$ denote \revise{any lift of the distinguished element $\b_{\l}$ (see the discussion before Lemma \ref{lem:fvsgrf}).}
\end{definition}
\revadd{We reiterate that $\a_\l$ is not unique, as the choice of lift (or the choice of isomorphism \eqref{eq:homOneFT}) is not canonical. By construction $\a_\l$ has degree $(\rbb(\l),\cbb(\l))$ when regarded as an element of $\Homgg(\one,\FT_n)$.  In the sequel will regard $\a_\l$ instead as a degree zero map $\one(-2\rbb(\l))[-\cbb(\l)]\to \FT_n$.}


Now we make the following crucial observation.

\begin{lemma}\label{lemma:eigenmapNotKilled}   Recall the polynomial $g_{P_{\row}}$ from as in Definition \ref{def:specht}.  That is,
\[
g_{P_{\row}} :=\prod (x_i-x_j)
\]
where the product is over pairs $i<j$ such that $i,j$ are in the same column of $P_{\row}$.  Then $g_{P_{\row}} \cdot \a_\l$ is not null-homotopic. \end{lemma}

\begin{proof}
There is no pair $i < j$ where the numbers $i, j$ appear in both the same row and the same column of $P_{\row}$. Consequently, if $i$ and $j$ are in the same column, then $(i,j) \notin Z_\l$, and $(x_i - x_j)$ does not act by zero on the $R$-module in \eqref{eq:specialSummand2}.  Thus
$g_{P_{\row}}\ot \b_\l\neq 0$.   Consequently, by Lemma \ref{lem:fvsgrf}, $g_{P_{\row}} \ot \a_\l$ is also nonzero.
%
\end{proof}

\begin{remark}  Note that $f \ot \b_\l \simeq 0$ for various polynomials $f \in \SM_\l$, including $\tbarb{P}$ for $P \in \SYT(\l)$. However, this does not imply that $f \ot \a_\l \simeq 0$, so it is not clear from this argument whether or not $\a_\l$ is a $\l$-equivalence.  In fact, it is a $\l$-equivalence, see \S\ref{subsec:symgrp}.\end{remark}


\begin{theorem}\label{thm:lambdaMaps}
For each partition $\l \in \PC(n)$ there exists a $\l$-equivalence $\a_\l: \shift{\l}{\one}\to \FT_n$. Furthermore, the collection of maps $\{\a_\l\}$ is obstruction-free in the sense of \cite{ElHog17a}. 
\end{theorem}

Forgiving the abuse of notation, the map $\a_\l$ of Theorem \ref{thm:lambdaMaps} may be assumed to be a linear combination of braid conjugates of the map $\a_\l$ constructed in
Definition \ref{def:alphaprime}.

\begin{proof}
Everything but the last statement follows from Lemma \ref{lemma:eigenmapNotKilled} and Corollary \ref{cor:sufficientcondition}.  The obstructions referred to in the statement are all degree $-1$ classes in $\Homgg(\Sigma_\l\circ \Sigma_\mu(\one), \FT_n^{\otimes 2})$.  Since $\Sigma_\l\circ \Sigma_\mu$ involves an even homological degree shift for each $\l,\mu\vdash n$, the obstructions correspond to classes in $\Homgg(\one,\FT_n^{\otimes 2})$ with odd homological degree.  However, $\Homgg(\one,\FT_n^{\otimes 2})$ is supported in even homological degrees by \cite[Theorem 3.18]{Hog17b}.
\end{proof}

\begin{corollary}
Let $\KB_\l:=\bigotimes_{\mu\neq \l}\Cone(\a_\mu)$ \revise{be the tensor product of all eigencones save for $\Cone(\a_{\l})$. The order of the tensor product is irrelevant by \eqref{eq:conescommute}. For all cells $\l$} we have
\begin{subequations}
\begin{equation} \label{eq:conescommute}
\Cone(\a_\l)\otimes \Cone(\a_\mu)\simeq \Cone(\a_{\mu})\otimes \Cone(\a_\l) \revise{\text{ for all } \mu \ne \l}.
\end{equation}
\begin{equation}\label{eq:conesquared}
\Cone(\a_\l)^{\otimes 2} \simeq (\shift{\l}{\one}[1] \oplus F) \otimes \Cone(\a_\l).
\end{equation}
\begin{equation} \label{eq:itsprediag} \bigotimes_{\l\in\PC(n)} \Cone(\a_\l) \simeq 0 \qquad \textrm{for all orderings of the factors}. \end{equation}
\begin{equation}\label{eq:quasiidemp}
\KB_\l^{\otimes 2} \simeq \left(\bigotimes_{\mu\neq \l} \revise{\shift{\l}{\one}\oplus \shift{\mu}{\one}[1]}\right)\otimes  \KB_\l.
\end{equation}
\begin{equation}\label{eq:KisEigenobject}
\KB_\l\otimes \Cone(\a_\l)\simeq 0 \simeq \Cone(\a_\l)\otimes \KB_\l.
\end{equation}
\end{subequations}
In particular, $\FT_n$ is categorically prediagonalizable.
\end{corollary}

\begin{proof}
The commuting of cones \eqref{eq:conescommute} follows because the maps $\a_\l$ are obstruction-free, see \cite[\S 6.2]{ElHog17a}. Given this, \eqref{eq:itsprediag} follows from Proposition \ref{prop:lequivmeansdiag}.  Equation \eqref{eq:KisEigenobject} follows from \eqref{eq:itsprediag}. Equation \eqref{eq:conesquared} is proven in \cite[\S A.2]{ElHog17a}, and \eqref{eq:quasiidemp} is an immediate consequence.
\end{proof}

\section{Relative cell theory and relative diagonalization in type $A$}
\label{sec:diag}

In the last section we constructed maps $\a_\l$ for each two-sided cell of $S_n$, and proved that $\FT_n$ was categorically prediagonalizable. However, we can not directly conclude that
$\FT_n$ is categorically diagonalizable. The issue is that there may exist distinct partitions $\l$ with the same shift $\Sigma_\l$, and which precludes an application of the categorical
diagonalization theorem \cite[Theorem 8.1]{ElHog17a}. To circumvent this problem, it is necessary to use the more technical relative diagonalization theorem \cite[Theorem 8.2]{ElHog17a}
(relative to an inductively constructed diagonalization of $\FT_{n-1}$).

\subsection{Diagonalization of $\FT_n$}

The purpose of this section is to state and understand some basic consequences of our main theorem.

\begin{theorem}\label{thm:typeAdiag}
There exist complexes $\PB_\l\in \KC^-(\SBim_n)$ indexed by partitions $\l$ of $n$, such that $\{(\PB_\l, \a_{\l})\}_{\PC(n)^{\op}}$ is a diagonalization of $\FT_n$, where $\a_\l$ are the maps from Theorem \ref{thm:lambdaMaps}.  In particular, we have a decomposition of identity 
\begin{equation}\label{eq:resofId}
\one \simeq \Big(\bigoplus_{\l} \PB_\l, \sum_{\mu,\l} d_{\mu\l}\Big)
\end{equation}
with $\PB_\l\otimes \Cone(\a_\l) \simeq \Cone(\a_\l)\otimes \PB_\l$.  Furthermore, $\PB_\l$ interacts with the cell theory of $\SBim_n$ in the following way:
\begin{enumerate}\setlength{\itemsep}{2pt}
\item the chain bimodules of $\PB_\l$ are in cells $\leq \l$.
\item $\PB_\l$ annihilates bimodules in cells $\not\geq \l$.
\end{enumerate}
\end{theorem}

The proof is an induction on $n$, carried out in \S \ref{subsec:proof}. But first we develop some consequences which will be helpful in the induction step. \revadd{We assume Theorem
\ref{thm:typeAdiag} for the rest of this section.}

If $(\Omega,\leq)$ is a poset, recall that a subset $I\subset \Omega$ is an \emph{ideal} if $i\leq j$ and $j\in I$ implies $i\in I$.  A subset $K\subset \Omega$ is \emph{convex} if $i\leq j\leq k$ and $i,k\in K$ implies $j\in K$.  Any convex set $K$ can be written as $K=J\setminus I$ for some poset ideals $I\subset J\subset \Omega$.

\begin{definition}\label{def:PK}
For each convex set $K\subset \PC(n)$, we have subquotient complex
\[
\PB_K:= \Big(\bigoplus_{\l\in K}\PB_\l, \sum_{\mu,\l\in K} d_{\mu\l}\Big).
\]
In particular $\PB_{\PC(n)}$ denotes the right-hand side of \eqref{eq:resofId}. 
\end{definition}

The differential $d_{\mu\l}\in \Homc^{1}(\PB_\l,\PB_\mu)$ is nonzero only if $\l\geq \mu$.  Thus, if $I\subset \PC(n)$ is a poset ideal, then $\PB_I$ is a subcomplex of $\PB_{\PC(n)}$, and the inclusion $\e:\PB_I\rightarrow \PB_{\PC(n)}\simeq \one$ satisfies $\Cone(\e)\simeq \PB_{I^c}$ where $I^c = \PC(n)\setminus I$.  This gives us a distinguished triangle
\[
\PB_I\rightarrow \one\rightarrow \PB_{I^c}\rightarrow \PB_I[1].
\]
The orthogonality of the $\PB_\l$ with respect to $\otimes$ implies that $\PB_I \otimes \PB_{I^c}\simeq 0 \simeq \PB_{I^c} \ot \PB_I$, hence $\PB_I$ has the structure of a counital idempotent in $\KC^-(\SBim_n)$ with complementary unital idempotent $\PB_{I^c}$.  The same argument shows more generally that if $I\subset J\subset \PC(n)$ are poset ideals and $K=J\setminus I$, then we have a distinguished triangle
\[
\PB_I\rightarrow \PB_J\rightarrow \PB_K\rightarrow \PB_I[1],\qquad\qquad \PB_I\otimes \PB_K\simeq 0 \simeq \PB_K\otimes\PB_I.
\]
That is, $\PB_I$ is a counital idempotent relative to $\PB_J$, with relative complement $\PB_{J\setminus I}$.   We collect some important properties of the $\PB_K$ below.

\begin{lemma}\label{lemma:PK}
The idempotents $\PB_K$, indexed by convex subsets $K\subset \PC(n)$, satisfy
\begin{enumerate}
\item $\PB_{\PC(n)}\simeq \one$ and $\PB_\emptyset = 0$.
\item $\PB_K\otimes \PB_L\simeq \PB_{K\cap L}$.
\end{enumerate}
In particular each $\PB_K$ is idempotent with respect to $\otimes$, up to homotopy.
\end{lemma}
\revcomment{We have added a proof of this.}
\begin{proof}
Statement (1) is clear. Now, let $K,L\subset \PC(n)$ be convex sets.  Observe that $\PB_K\ot \PB_L = \widetilde{\bigoplus}_{\l\in L}\PB_K\ot \PB_\l$.  If $\l\not\in K$ then $\PB_K\ot \PB_\l\simeq 0$. Contracting these complexes does not affect the differential among the surviving terms, by the following argument.

Let $\e$ be the twist which when added to the differential of $\bigoplus_{\l\in L}\PB_K\ot \PB_\l$ yields $\PB_K\ot \PB_L$.  This twisting element is strictly decreasing with respect to the dominance order on partitions, so contracting the terms with $\l\not\in K$ can only affect the the component of the differential from two surviving terms $\PB_K\ot \PB_{\mu_1}$ to $\PB_K\ot \PB_{\mu_2}$ (with $\mu_i\in K$) if there exists $\l\in I$ such that $\mu_1>\l>\mu_2$.  But this is not possible since $K$ is convex.  This gives us that
\begin{equation}\label{eq:PK PL}
\PB_K\ot \PB_L\simeq \PB_K\ot \PB_{K\cap L}
\end{equation}
In the special case $L=\PC(n)$ this shows that $\PB_K\simeq \PB_K\ot \PB_K$, using (1).

By symmetry we also have  $\PB_K\ot \PB_L\simeq \PB_{K\cap L}\ot \PB_L$. Applying this with $K \cap L$ replacing $L$, we deduce
\[
\PB_K\ot \PB_{K\cap L}\simeq \PB_{K\cap L}\ot \PB_{K\cap L}.
\]
Combining this with \eqref{eq:PK PL} and idempotence of $\PB_{K\cap L}$ (established earlier in the proof) proves (2).
\end{proof}

The counital idempotents $\PB_I$, indexed by poset ideals $I\subset \PC(n)$, are often easier to deal with in practice. 
\begin{proposition}\label{prop:counitalChar}
If $I\subset \PC(n)$ is a poset ideal then $\PB_I$ satisfies
\begin{enumerate}
\item $\PB_I$ is homotopy equivalent to a complex whose chain objects are in cells $\l\in I$.
\item there exists a chain map $\e:\PB_I\rightarrow \one$ such that $\Cone(\e)\otimes B\simeq 0 \simeq B\otimes \Cone(\e)$ for all $B$ in cell $\l$ with $\l\in I$.
\end{enumerate}
Furthermore,
\begin{enumerate}
\item[(3)] a complex $C\in \KC^-(\SBim_n)$ satisfies $\PB_I\otimes C\simeq C$ if and only if $C$ is homotopy equivalent to a complex whose chain bimodules lie in cells $\l \in I$.
\item[(4)] properties (1) and (2) characterize the pair $(\PB_I,\e)$ up to canonical isomorphism in $\KC^-(\SBim_n)$.
\end{enumerate}
\end{proposition}
\begin{proof}
Property (1) follows from property (1) of Theorem \ref{thm:typeAdiag}, given that $\PB_I$ is a convolution of complexes $\PB_\l$ with $\l\in I$.

We choose $\e:\PB_I\rightarrow \one$ to be the first map in a distinguished triangle $\PB_I\rightarrow \one\rightarrow \PB_{I^c}\rightarrow \PB_I[1]$, so that $\PB_{I^c} \simeq \Cone(\e)$.  If $B$ is a bimodule in cell $\l\in I$, then $\PB_{I^c}\otimes B\simeq 0\simeq B\otimes \PB_{I^c}$ by property (2) of Theorem \ref{thm:typeAdiag}.  Here we are using the fact that $\PB_{I^c}$ is a convolution of $\PB_\l$ with $\l\in I^c$, and no element of $I$ can be larger or equal to an element of $I^c$.  This proves (2).

Now, let $C\in \KC^-(\SBim_n)$ be given.  Then each chain bimodule of $\PB_I\otimes C$ in a cell $\l\in I$, since the same is true of $\PB_I$.  This proves the ``only if'' direction of (3).  Conversely, if each summand of each chain bimodule of $C$ is in a cell $\l\in I$, then $\PB_{I^c}\otimes C\simeq 0\simeq C\otimes \PB_{I^c}$ by (2). The existence of the distinguished triangle $\PB_I\rightarrow \one\rightarrow \PB_{I^c}\rightarrow \PB_I[1]$ then implies that $\PB_I\otimes C\simeq C\simeq C\otimes \PB_I$.  Finally, the uniqueness statement (4) follows by elementary properties of categorical idempotents, see Theorem 4.28 in \cite{Hog17a}. 
\end{proof}

\begin{lemma}\label{lemma:projcentral}
Each $\PB_K$ is central in $\KC^-(\SBim_n)$.
\end{lemma}
\begin{proof}
Any convex set $K$ can be expressed as $K=J\setminus I$ for some poset ideals $I\subset J\subset \PC(n)$.  In this case $\PB_K$ can be expressed as
\[
\PB_K \simeq \PB_J\otimes \PB_{I^c}.
\]
A counital idempotent is central if and only if its unital complement is central (Theorem 4.23 in \cite{Hog17a}) hence we are reduced immediately to the case when $K=I$ is a poset ideal. 
This, in turn, follows from general arguments.  Let $\EB$ be a unital or counital idempotent in any triangulated monoidal category.  Let $\AC \EB\subset \AC$ denote the full subcategory on the set of objects $A$ with $A\otimes \EB\simeq A$.  Then $\EB$ is central in $\AC$ if and only if $\AC \EB$ is a two-sided tensor ideal.   The lemma now follows directly from \revadd{statement (3)} of Proposition \ref{prop:counitalChar}.
\end{proof}

Note that the singleton $K=\{\l\}$ is automatically a convex subset of $\PC(n)$, hence $\PB_\l =\PB_{\{\l\}}$ is central as well.


\subsection{Adding a strand}

Suppose that we have proven Theorem \ref{thm:typeAdiag} for $n-1$. Then for any $\mu \in \PC(n-1)$, one has an idempotent complex $\PB_{\mu}$ with certain properties. The key to our inductive argument will be understanding how the complex $(\PB_{\mu} \sqcup \one_1)$ acts on $\KC^b(\SBim_n)$. We begin with some elementary observations.

\begin{lemma}[Projector sliding]\label{lemma:projsliding1}
We have $F(\sigma_1\sigma_2\cdots\sigma_{n-1})\otimes (\PB_K\sqcup \one_1)\simeq (\one_1\sqcup \PB_K)\otimes F(\sigma_1\sigma_2\cdots \sigma_{n-1})$ inside $\KC^b(\SBim_{n})$, for any convex subset $K$ of $\PC(n-1)$.
\end{lemma}
We refer to the the statement of this lemma colloquially as ``projectors slide past strands.''

\begin{proof}
\revcomment{Removed redundant first paragraph.}

Denote $\XB:=F(\sigma_1\sigma_2\cdots \sigma_{n-1})$.  Abuse notation by writing $\PB_K = \PB_K\sqcup \one_1$ and $\PB_K'=\one_1\sqcup \PB_K$.  We also will omit the tensor product symbol for aesthetic reasons.  We may write $K=J\setminus I$ for some poset ideals $I\subset J\subset \PC(n-1)$, and $\PB_K \simeq \PB_J\otimes \PB_{I^c}$. Thus, the desired equivalence $\XB \PB_K\simeq \PB_K'\XB$ reduces immediately to the case when $K=I$ or $K=I^c$.

If $I\subset \PC(n-1)$ is an ideal, then there is a distinguished triangle
\[
\PB_I\rightarrow \one \rightarrow \PB_{I^c}\rightarrow \PB_I[1],
\]
where the chain bimodules of $\PB_{I}$ are in cells $\l\in I$ and $\PB_{I^c}$ kills all bimodules of $\SBim_{n-1}$ in cells  $\l\in I$.  This yields a distinguished triangle
\[
\PB_I' \XB \PB_I\rightarrow \PB_I' \XB \rightarrow \PB_I'\XB \PB_{I^c}\rightarrow \PB_I'\XB \PB_I[1].
\]
The third term is contractible since each chain bimodule of $\PB_{I}'$ is in cells $\l\in I$, slides past $\XB$ by Proposition \ref{prop:slidingBimodules}, and is annihilated by $\PB_{I^c}$.  This implies that
\[
\PB_I'\XB\PB_I \simeq \PB_I'\XB.
\]
A similar argument shows that
\[
\PB_I'\XB\PB_I \simeq \XB\PB_I,
\]
from which it follows that $\XB\PB_I\simeq \PB_I'\XB$, as claimed.  A similar argument shows that $\XB\PB_{I^c}\simeq \PB_{I^c}\XB$.   This completes the proof.
\end{proof}

\begin{lemma}\label{lemma:projectorSliding}
We have $(\PB_{\mu}\sqcup \one_1)\otimes \FT_n \simeq \FT_n\otimes (\PB_{\mu} \sqcup \one_1)$.
\end{lemma}
\begin{proof}
This follows immediately from Lemmas \ref{lemma:projcentral} and \ref{lemma:projsliding1}.
\end{proof}

\begin{corollary} \label{cor:projectorconecommute}
$\PB_{\mu}\sqcup \one_1$ commutes with $\Cone(\a_\l)$ for all $\l \in \PC(n)$ and all $\mu \in \PC(n-1)$.
\end{corollary}
\begin{proof}
See Lemma 6.19 in \cite{ElHog17a} and its proof.
\end{proof}

Now that these basic properties are in place, we must ask the harder questions. The idempotent $\PB_\mu$ interacts with the cell theory of $\SBim_{n-1}$ in a particular way; how does the
idempotent $\PB_\mu \sqcup \one_1$ interact with the cell theory of $\SBim_n$? We will eventually show that tensoring with $\PB_\mu \sqcup \one_1$ will allow us to restrict our attention
to those $\l \in \PC(n)$ which are obtained from $\mu$ by adding a box, which is a totally ordered set. For this purpose, we now digress and discuss relative cell theory.

\subsection{Parabolic subgroups and subtableaux}
\label{subsec:subtableaux}

The Hecke algebra $\HB(S_k)\cong \HB(S_k\times S_1^{\ell})$ embeds into $\HB(S_{k+\ell})$ using the external product. The goal of this section and the next will be to generalize \eqref{eq:dacts}, taking the operator $v^{\rbb(\l)} b_{T,T} \in \HB(S_k)$ and seeing how it acts on $\HB(S_{k+\ell})$.

Let $W_I$ be a parabolic subgroup of $W$. Each right coset of $W_I$ has a unique element of minimal length; let $Y_I$ denote this set of minimal length right coset representatives.
Similarly, let $X_I$ denote the set of minimal length left coset representatives. Each $w \in W$ has a unique representation as $uy$ for $u \in W_I$ and $y \in Y_I$, or as $xt$ for $t
\in W_I$ and $x \in X_I$, and moreover $\ell(w) = \ell(u) + \ell(y) = \ell(t) + \ell(x)$.

Our first goal is to relate the language of cosets to the language of tableaux. Given a tableau $P$ with $n$ boxes and $k<n$, let $P^k$ denote the tableau with $k$ boxes obtained by
remembering the first $k$ boxes of $P$.

\begin{lemma} \label{lem:cosetvtab} Let $W = S_{k+\ell}$ and $W_I = S_k=S_k\times S_1\times \cdots \times S_1$. Let $w = uy = xt$ for $w \in W$, $x \in X_I$, $y \in Y_I$, and $t, u \in W_I$. Suppose that $w$ corresponds to
$(P,Q)$ under the Robinson-Schensted correspondence. Then $u$ corresponds to $(P^k,A)$ for some tableau $A$, and $t$ corresponds to $(B,Q^k)$ for some tableau $B$. \end{lemma}

\begin{proof} \revcomment{Merged remark with proof.} The result for right cosets and the result for left cosets imply each other, under taking inverses. Thus, we need prove only the second statement, that $t$ corresponds to
$(B,Q^k)$. Because of the asymmetry between $P$-symbols and $Q$-symbols in the Robinson-Schensted algorithm, we found left cosets easier to use than right cosets.

The permutation $t \in S_k$ is easily described as the permutation which puts $\{1, \ldots, k\}$ in
the same order as $\{w(1), \ldots, w(k)\}$. More precisely, let $Z$ be the image of the set $\{1, \ldots, k\}$ under $w$, and let $\xi \co Z \to \{1, \ldots, k\}$ denote the unique
order-preserving bijection. Then $x(i) = \xi \circ w(i)$ for all $1\leq i\leq k$.

In the Robinson-Schensted algorithm, one constructs the tableau-pair $(P,Q)$ associated to $w$ one box at a time, by inserting $w(1)$, then $w(2)$, etcetera. At the $k$-th step in this
process, one has a pair $(P_k, Q_k)$ where $P_k$ is not a standard tableau (its entries are $Z$, not $\{1, \ldots, k\}$), while $Q_k$ is standard and agrees with $Q^k$. The only
difference between performing the first $k$ steps of this algorithm for $w$ and for $t$ is that the boxes in $P_k$ are relabeled via $\xi$; the recording tableau $Q_k$ is unchanged, and
only depends on the relative order of $w(1), \ldots, w(k)$. Hence $Q^k$ agrees with the recording tableau for $t$, as desired. \end{proof}

\begin{remark} In particular, this proof also implies that the tableau $B$ can be described in a straightforward way from $(P,Q)$ as well: by ``unbumping'' the extraneous boxes of $P$ in
the order determined by $Q$ to recover the non-standard tableau $P_k$, and then relabeling it via $\xi$. \end{remark}

\begin{notation} For a tableau $T$ with shape $\l$, we write $\shape(T) = \l$. For $k \le n$ and $w \in S_n$, write $\shape_{L,k}(w) = \shape(P^k)$, where $w$ corresponds to $(P,Q)$. Write $\shape_{R,k}(w) = \shape(Q^k)$. These are the shapes of $u$ and $t$ respectively. \end{notation}
	
In the opposite direction, there is an inclusion map $S_k \to S_{k+\ell}$ which sends $w \mapsto w \sqcup 1_{\ell}$. If $w$ corresponds to $(P,Q)$, then $w \sqcup 1_{\ell}$ corresponds to $(U,V)$, where $U$ is obtained from $P$ by adding the boxes $k+1, k+2, \ldots, k+\ell$ to the first row, and likewise for $V$ and $Q$.

\subsection{Relative KL cells and the relative action of involutions}
\label{subsec:relativecells}

Meinolf Geck in \cite{GeckRelative} laid down the framework for the study of relative cells in $\HB(W)$, relative to a parabolic subgroup $W_I$. We continue to use the notation $X_I$ and $Y_I$ from the previous section.

Geck defines a new partial order $\le_{L,I}$ on $W$, or on the set of KL basis elements $\{b_w\}_{w \in W}$, which is analogous to the left cell order $\le_L$ (see
\S\ref{subsec:algCells}) but for left multiplication by elements of the subalgebra $\HB(W_I) \subset \HB(W)$. Namely\footnote{Again, one takes the transitive closure of this relation.}, \[ b_x \le_{L,I} b_y \quad \textrm{if} \quad b_x \babysumset h \cdot
b_y \textrm{ for some } h \in \HB(W_I). \] One defines $\le_{R,I}$ and $\le_{LR,I}$ analogously.

We can equip $W_I$ with its usual partial order $\le_L$; this will not be confused with the corresponding order on $W$, since we are only interested in $\le_{L,I}$ in this section. We may equip $X_I$ and $Y_I$ with the partial order $\le$ induced from the Bruhat order. We state results for left multiplication by $W_I$ (hence right cosets); the results on the other side are analogous.

\begin{lemma} \label{lem:Geck1}(\cite[Proposition 4.4]{GeckRelative}) For $t,u \in W_I$ and $x,y \in Y_I$ one has
\begin{equation} \label{eq:Geck1} tx \le_{L,I} uy \quad \implies \quad t \le_{LR} u \quad \textrm{and} \quad x \le y. \end{equation}
Thus, if $tx \sim_{L,I} uy$ then $t \sim_{LR} u$ and $x = y$. \end{lemma}

One also has an analog of Proposition \ref{prop:leftincomp}: that left cells within a given two-sided cell are incomparable.

\begin{lemma} \label{lem:Geck2} (\cite[Theorem 4.8]{GeckRelative}) For $t,u \in W_I$ and $x,y \in Y_I$, if $tx \le_{L,I} uy$ and $t \sim_{LR} u$ then $x = y$ and $t \sim_L u$.
\end{lemma}

We will not use Lemma \ref{lem:Geck2} directly, but it gives an idea of how Geck's relative cells function. The key result we will use from Geck is the following. As before, let $c_{x,y}^z$ denote\footnote{This is denoted $h_{x,y,z}$ in \cite{GeckRelative}.} the coefficient of $b_z$ in the product $b_x b_y$.

\begin{lemma} \label{lem:Geck3} (\cite[Lemma 4.7]{GeckRelative}) For $t,u,w \in W_I$ and $x, y \in Y_I$ one has \begin{enumerate}
\item If $c_{w, ty}^{ux} \ne 0$ then $u \le_{LR} t$ and $x \le y$. (This is a restatement of Lemma \ref{lem:Geck1}.)
\item If $x=y$ then $c_{w,ty}^{uy} = c_{w,t}^{u}$.
\item Assume that $t \sim_{LR} u$ inside $W_I$, both living in cell $\mu$. The coefficient of $v^{-a}$ in $c_{w,ty}^{ux}$ is zero for any $a > \rbb(\mu)$, and if the coefficient of $v^{-\rbb(\mu)}$ is nonzero, then $x=y$. \footnote{The statement of \cite[Lemma 4.7]{GeckRelative} is slightly weaker than this, but the proof suffices to show this stronger result. Geck also considers positive exponents of $v$ rather than negative exponents, but these polynomials are self-dual.}
\end{enumerate} \end{lemma}

Let us rephrase this result using the language of tableau, in the special case when $W = S_n$ for $n = k + \ell$, and $W_I = S_k=S_k\times (S_1)^\ell$. First, let us introduce notation.

\begin{notation} \label{not:relideal} Let $k \le n$ and $\mu \in \PC(k)$. Let $\HB_{L,<\mu} \subset \HB = \HB(S_n)$ denote the span of those $b_w$ for which $\shape_{L,k}(w) < \mu$.
Define $\HB_{L,\le \mu}$ similarly. Let $\HB_{L,\mu}$ denote the span of those $b_w$ for which $\shape_{L,k}(w) = \mu$. Let $\HB^+_{L,\mu}$ denote the $\Z[v]$-span of $vb_w$ where
$\shape_{L,k}(w) = \mu$. \end{notation}

Note that, while $\HB_{L,<\mu}$ is not a left ideal, it is preserved under the left action of the subalgebra $\HB(S_k)$, thanks to Lemma \ref{lem:Geck1}.

\begin{prop} \label{prop:relaction}  Fix $k < n$. Let  $w \in S_n$, and $w = ty$ for $t \in W_I$ and $y \in Y_I$. Suppose $w = w(P,Q,\l)$, and $t = w(P^k,A,\mu)$, so that $\mu = \shape_{L,k}(w)$. Fix any other element $z \in W_I$ with shape $\mu$, so that $z = w(U,V,\mu)$. Let $g = w(U,A,\mu) \in S_k$. Then
\begin{equation} v^{\rbb(\mu)} b_z b_w = \delta_{V,P^k} b_{gy} + \HB^+_{L,\mu} + \HB_{L,< \mu}.\end{equation}
\end{prop}

\begin{proof} By definition,
\begin{equation} b_z b_{ty} = \sum_{u \in W_I, x \in Y_I} c_{z,ty}^{ux} b_{ux}. \end{equation}
By part (1) of Lemma \ref{lem:Geck3}, the only terms appearing in this sum will have $u \le_{LR} t$ and $x \le y$. In particular, $\shape(u) = \rho \in \PC(k)$ with $\rho \le \mu$. Let us ignore those terms with $\rho < \mu$, and consider only those with $\rho = \mu$, i.e. $u \sim_{LR} t$.  Then no $b_{ux}$ appears with $v^{-a}$ for $a > \rbb(\mu)$, and $v^{-\rbb(\mu)}$ will appear only if $x = y$. Moreover, when $x = y$, $c_{z,ty}^{ux} = c_{z,t}^{u}$. But $b_z b_t = b_{U,V} b_{P^k,A} = \phi(V,P^k) b_{U,A}$, so that $c_{z,t}^{u} = 0$ unless $u = w(U,A,\mu) = g$. Now we can use \eqref{eq:howPQmult} to say that a coefficient with $v^{-\rbb(\mu)}$ will appear if and only if $V = P^k$.
\end{proof}

This leads to our relative version of \eqref{eq:dacts}.

\begin{cor} \label{cor:relactioninvolution} Fix $k < n$, and let $w \in S_n$ correspond to $(P,Q)$. Let $P^k$ have shape $\mu$, and fix another tableau $V \in \SYT(\mu)$. Then one has
\begin{equation} \label{eq:dactsrel} v^{\rbb(\mu)} (b_{V,V} \sqcup 1_{n-k}) b_w = \delta_{V,P^k} b_w + \HB^+_{L,\mu} + \HB_{L,< \mu}. \end{equation} \end{cor}
	
Multiplication formulas for KL basis elements in the Hecke algebra give rise to decompositions of tensor products in the category of Soergel bimodules.

\begin{cor}\label{cor:relativeUnit}
Suppose $w=w(P,Q,\l) \in S_n$, and fix $k < n$.  Let $\mu = \shape_{L,k}(w)$. Then $B_w$ is isomorphic to a direct summand of $(B_{P^k,P^k}(\rbb(\mu)) \sqcup \one_{n-k})\otimes B_w$. \qed
\end{cor}

\begin{proof} This follows immediately from \eqref{eq:dactsrel}. \end{proof}

Finally, we need one other result of Geck, whose proof uses separate ideas (mixing the KL basis and the standard basis) which we choose not to recall.

\begin{lemma} \label{lem:itsarightideal}
The subspace $\HB_{L,\le \mu}$ from Notation \ref{not:relideal} is a right ideal in $\HB(W)$. \end{lemma}

\begin{proof} This subspace appears in \cite[Corollary 3.4]{GeckInduction}, where it is shown to agree with a space Geck calls $\MC$. In \cite[Lemma 2.2]{GeckInduction}, it is proven that $\MC$ is a right\footnote{Geck is using left cosets rather than right cosets, so he obtains a left ideal.} ideal. \end{proof}

The consequence in the category of Soergel bimodules is that those indecomposables $B_w$ where $\shape_{L,k}(w) \le \mu$ (i.e. the categorification of $\HB_{L,\le
\mu}$) form a right tensor ideal. 

\subsection{Implications for induced projectors}
\label{subsec:implications}

Assuming that Theorem \ref{thm:typeAdiag} is proven for $S_k$, for any $\mu \in \PC(k)$ we refer to the idempotent $\PB_\mu \sqcup \one_{n-k}$ as an induced projector. Now we use the
results of the previous section to study induced projectors and their interaction with cell theory. 

\begin{lemma} \label{lem:relprojkills} Fix $k < n$. Assume Theorem \ref{thm:typeAdiag} is proven for $k$, and fix $\mu \in \PC(k)$. Then \begin{enumerate}
	\item the chain bimodules of $\PB_\mu \sqcup \one_{n-k}$ are direct sums of shifts of $B_w$ for $w \in S_n$ with $\shape_{L,k}(w) \le \mu$.
	\item $\PB_\mu \sqcup \one_{n-k}$ annihilates bimodules $B_w$ for $w \in S_n$ with $\shape_{L,k}(w) \ngeq \mu$.
\end{enumerate}
\end{lemma}

\begin{proof} We know by Theorem \ref{thm:typeAdiag} that $\PB_\mu$ annihilates indecomposable objects of $\SBim_{k}$ in cells $\ngeq \mu$. Let $w = w(P,Q,\l)$ with $\shape_{L,k}(\l) \ngeq \mu$. Thus
\begin{equation} (\PB_\mu \sqcup \one_{n-k}) \ot (B_{P^k,P^k} \sqcup \one_{n-k}) \ot B_w \simeq 0, \end{equation}
because the shape of $P^k$ is $\ngeq \mu$. However, $B_w$ is a summand of $(B_{P^k,P^k} \sqcup \one_{n-k}) \ot B_w$, by Corollary \ref{cor:relativeUnit}, so $(\PB_\mu \sqcup \one_{n-k}) \ot B_w \simeq 0$.

Now, every chain object of $\PB_\mu$ is in cells $\le \mu$ in $\SBim_k$. In particular, every indecomposable summand of a chain object of $\PB_\mu \sqcup \one_{n-k}$ is $B_w$ for $w \in S_k \subset S_n$ in cell $\nu \le \mu$. This is a much stronger result than merely asserting that $\shape_{L,k}(w) \le \mu$. Of course, $\nu = \shape_{L,k}(w)$ for such $w$.  \end{proof}

Now we prove the key lemma for our inductive proof of Theorem \ref{thm:typeAdiag}. We set up notation. Assume Theorem \ref{thm:typeAdiag} is proven for $n-1$, and fix $\mu \in \PC(n-1)$. Let $\l^1 \le \l^2 \le \ldots \le \l^r$ be the partitions obtained from $\mu$ by adding a box, ordered via the dominance order, so that $\l^r$ is obtained by adding a box to the first row, and $\l^1$ is obtained by adding a box to the first column. (Here we have used superscripts to avoid confusion with the parts of a partition.) We write $\l \supset \mu$ to indicate that $\l = \l^i$ for some $i$. Recall that the maps $\a_{\l^i}$ have been defined already in Theorem \ref{thm:lambdaMaps}. 

\begin{lemma} \label{lem:relprojcones} Assume Theorem \ref{thm:typeAdiag} is proven for $n-1$, and fix $\mu \in \PC(n-1)$. With notation as in the previous paragraph, the tensor product
\begin{equation}\label{eq:reltensorprod}
(\PB_{\mu}\sqcup \one_1) \otimes \bigotimes_{\l\supset \mu} \Cone(\a_\l) \simeq 0
\end{equation}
is contractible. \end{lemma}

\begin{proof}
\revcomment{moved paragraph before into proof.} Let us write $\PB$ for $(\PB_{\mu}\sqcup \one_1)$ during this proof. Since the cones commute, this tensor product does not depend on the ordering of the factors. Note that $\PB$ is an idempotent, and commutes with each of these cones by Corollary \ref{cor:projectorconecommute}.

\begin{subequations}
Let \begin{equation} \XB^r = \PB \ot \Cone(\a_{\l^r}), \end{equation} \begin{equation} \XB^{r-1} = \PB \ot \Cone(\a_{\l^r}) \ot \Cone(\a_{\l^{r-1}}), \end{equation} etcetera, so that $\XB^1$ is just the tensor product in \eqref{eq:reltensorprod}.
\end{subequations}

Given a collection of complexes $C_1,\ldots,C_r\in \KC^-(\SBim_n)$, let $\ip{C_1,\ldots,C_r} \subset \KC^-(\SBim_n)$ denote smallest full subcategory containing the $C_i$ and closed under grading shifts $(\pm 1)$, $[\pm 1]$, and locally finite convolutions.   \revadd{
Here, a convolution $\tilde{\bigoplus}_{i\in I} A_i$ is called \emph{locally finite} if the direct sum on chain objects is finite in each homological degree.}
 Note that $\ip{B_w\:|\: \shape_{L,n-1}(w)\leq \mu}$ is a right tensor ideal in $\KC^-(\SBim_n)$ by Lemma \ref{lem:itsarightideal}, and $\PB\in \ip{B_w\:|\: \shape_{L,n-1}(w)\leq \mu}$ by Lemma \ref{lem:relprojkills}.  Thus every complex under consideration for the remainder of the proof will be in $\ip{B_w\:|\: w\in \shape_{L,n-1}(w)\leq \mu}$.

We prove by descending induction that
\begin{equation}\label{eq:thegoal1}
\XB^i\in \ip{\PB \ot B_w\:|\: \text{ $\shape_{L,n-1}(w)\leq \mu$ and $w$ is in some cell $< \l^i$}}.
\end{equation}
First observe that $\PB$ annihilates $B_w$ when $\shape_{L,k}(w) \ngeq \mu$, hence we need only consider those $w$ with $\sh_{L,n-1}=\mu$.  But this implies that $w$ is in cell $\l^j$ for some $j$, since removing a box yields $\mu$. Then $\l^j<\l^i$ implies $j<i$.  Thus \eqref{eq:thegoal1} is equivalent to
\begin{equation}\label{eq:thegoal2}
\XB^i\in \ip{\PB \ot B_w\:|\: \text{ $\shape_{L,n-1}(w) = \mu$ and $w$ is in some cell $\l^j$ for $j<i$}}
\end{equation}
We remark that if $i=1$, then this reduces to $\XB^1\simeq 0$, which is what we wish to prove.

For any indecomposable bimodule $B \in \SBim_{n-1}$ in cell $\mu$, $B \sqcup \one_1$ is in cell $\l^r$. Thus the chain objects of $\PB$ lie in cells $\leq
\l^r$. Tensoring with $\Cone(\a_{\l^r})$ takes objects in cells $\le \l^r$ to complexes in cells $< \l^r$. This proves that \eqref{eq:thegoal1} (equivalently \eqref{eq:thegoal2}) is satisfied for $\XB^r$.

Assume by induction that \eqref{eq:thegoal2} is satisfied by $\XB^{i+1}$.  Then
\[
\XB^i\in\ip{\PB\otimes B_w\otimes \Cone(\a_{\l^i})\:|\: \text{ $w$ is in some cell $\l^j$ for $j \le i$}}.
\]
Tensoring with $\Cone(\a_{\l^i})$ takes complexes  in cells $\le \l^i$ to complexes in cells $< \l^i$. We deduce that
\[
\XB^i\in\ip{\PB\otimes B_w\:|\: \text{$w$ is in some cell $<\l^j$ for $j \le i$}}.
\]
As remarked above, $\ip{\PB\otimes B_w\:|\: \shape_{L,n-1}(w)\leq \mu}$  is a right tensor ideal, hence \eqref{eq:thegoal1} (equivalently \eqref{eq:thegoal2}) is satisfied for $\XB^i$. 

Thus, $\XB^1$ satisfies \eqref{eq:thegoal1} by induction. There are no elements $w$ in cells $< \l^1$ for which $\shape_{L,k}(w) = \mu$, hence $\XB^1 \simeq 0$.  \end{proof}

\subsection{Proof of Theorem \ref{thm:typeAdiag}}
\label{subsec:proof}

The proof is by induction on $n$.  The base case $n=1$ is trivial.  Assume by induction that we have constructed $\PB_{\mu}\in \KC^-(\SBim_{n-1})$ as in the statement.  In particular,
\[
\one_{n-1}\simeq \bigoplus_{\mu \in \PC(n-1)} \PB_{\mu} \qquad \text{ with twisted differential}.
\]
The twisted differential respects the opposite of the usual dominance order on partitions, meaning that the component of the differential from $\PB_\mu$ to $\PB_{\mu'}$ is zero unless $\mu\geq \mu'$.

Let us match the notation of the Relative Diagonalization Theorem \cite[Theorem 8.2]{ElHog17a}. Let $F = \FT_n$. Let $\XC$ (resp $\YC$) denote the set $\PC(n-1)$ (resp. $\PC(n)$) with its opposite poset structure. For $\mu \in \XC$, let $\YC_\mu \subset \YC$ denote the set $\{ \l \in \PC(n) \mid \l \supset \mu\}$. Then, Lemma \ref{lem:relprojcones} together with Lemma \ref{lemma:projectorSliding} and Theorem \ref{thm:lambdaMaps} give all the assumptions of the Relative Diagonalization Theorem.

Then \cite[Theorem 8.2 and following remarks]{ElHog17a} implies the following theorem.

\begin{theorem} \label{thm:relDiag1less} There is a diagonalization $\{(\PB_\l, \a_\l)\}$ of $\FT_n$ indexed by $\YC$, which is built by reassociating an idempotent decomposition $\{\PB_{\mu,\l}\}$ indexed by $\XC \times \YC$. We have
\begin{enumerate}\setlength{\itemsep}{2pt}
\item $\PB_\l \simeq \bigoplus_{\mu \in \XC} \PB_{\mu,\l}$ with twisted differential.
\item $\PB_{\mu,\l}\simeq 0$ unless $\mu\subset \l$.
\item If $\mu\subset \l$, then $\PB_{\mu,\l}$ is constructed as a locally finite convolution of shifts of complexes of the form
\begin{equation} \label{eq:doubleidemp}
(\PB_{\mu}\sqcup \one_1)\otimes \bigotimes_{\nu}\Cone(\a_\nu)
\end{equation}
where $\nu$ ranges over partitions containing $\mu$ but not equal to $\l$.
\end{enumerate}
Furthermore, the relative diagonalization is always \emph{relatively tight} in the sense that
\begin{enumerate}
\item[(4)] if $C\in \KC^-(\SBim_n)$ satisfies \revise{$(\PB_{\mu}\sqcup \one_1)\otimes C \simeq C$} and $\Cone(\a_\l)\otimes C\simeq 0$ for $\mu\subset \l$, then $\PB_{\mu,\l}\otimes C\simeq C$ (and conversely).
\end{enumerate}
\end{theorem}

\begin{lemma} \label{lem:doubleidempcell}
$\PB_{\mu,\l}$ lives in cells $\leq \l$ and kills complexes in cells $\not\geq \l$.
\end{lemma}

\begin{proof}
We need only prove the same result for the tensor product in \eqref{eq:doubleidemp}, which we temporarily call $\XB_{\mu,\l}$

Let us resume the notation of the proof of Lemma \ref{lem:relprojcones}.  In particular, $\PB:=\PB_\mu\sqcup \one_1$.  Let $\l = \l^i$. By \eqref{eq:thegoal2},
\[ \XB^{i+1} \in \ip{\PB \ot B_w\:|\: \text{ $w$ is in some cell $\le \l$}}. \] In particular, this implies that
\[ \XB^{i+1} \in \ip{B_w\:|\: \text{ $w$ is in some cell $\le \l$}}. \] Since $\XB_{\mu,\l}$ is obtained by taking $\XB^{i+1}$ and tensoring with additional cones, it is also in cells $\le \l$, as desired.

Because $\XB_{\mu,\l} \simeq \XB_{\mu,\l} \ot \PB$, we need to prove that $\XB_{\mu,\l} \ot \PB \ot B_w \simeq 0$ whenever $w$ is in cell $\nu$ with $\nu \ngeq \l$. In fact, we can assume that
$\shape_{L,k}(w) \le \mu$. This is because any indecomposable summand in $\PB \ot B_w$ is $B_z$ with $\shape_{L,k}(z) \le \mu$, by Lemma \ref{lem:relprojkills}, and $\PB \ot B_w \simeq
\PB \ot \PB \ot B_w$. But by Lemma \ref{lem:relprojkills}, we may assume $\shape_{L,k}(w) \ge \mu$ or else $\PB \ot B_w \simeq 0$. Hence $\shape_{L,k}(w) = \mu$, meaning that $\nu =
\l^j$ for some $j$. Then $\nu \ngeq \l$ implies $j < i$. The same inductive argument as Lemma \ref{lem:relprojcones} will prove that the cones $\Cone(\a_{\l^k})$ for $k \le j$ will
suffice to kill $\PB \ot B_w$, meaning that $\XB_{\mu,\l} \ot \PB \ot B_w \simeq 0$ as well. \end{proof}

Finally, since $\PB_\l$ is obtained by associating together the complexes $\PB_{\mu,\l}$ with $\mu\subset \l$, it shares the properties stated in Lemma \ref{lem:doubleidempcell}. This
completes the proof of Theorem \ref{thm:typeAdiag}.

\subsection{The primitive idempotents}
\label{subsec:primitives}

\begin{definition} \label{defn:PS}
Let $S=(\l_S^{1},\ldots,\l_S^{n})\in \PC(1)\times \cdots \times \PC(n)$ be a sequence of partitions with $\l_S^{i}\in\PC(i)$.  Set
\[
\PB_S := \bigotimes_{i=1}^n (\PB_{\l_S^i}\sqcup \one_{n-i}).
\]
Note that $\PB_{\l_S^1}=\one$, hence can be removed from this tensor product with no effect.
\end{definition}

Each $\PB_{\l_S^i}$ is central in $\KC^b(\SBim_i)$, so while $(\PB_{\l_S^i}\sqcup \one_{n-i})$ is not central, it commutes with $(\PB_{\l_S^j}\sqcup \one_{n-j})$ for $j<i$. Hence the entire family commutes, and the ordering of the above tensor product is irrelevant up to homotopy equivalence.

Let $\SYT(n)$ denote the set of standard tableaux with $n$ boxes. Given $T \in \SYT(n)$, we have already defined the tableau $T^k \in \SYT(k)$ in
\S\ref{subsec:subtableaux}. Let $\l_T^k = \shape(T^k)$. Hence $T$ gives rise to a sequence of partitions $(\l_T^1, \l_T^2, \ldots, \l_T^n)$, which determines $T$ uniquely. Using this
injective map $\SYT(n) \to \PC(1) \times \cdots \times \PC(n)$, we define $\PB_T$ as above. Conversely, a sequence of partitions corresponds to a tableau if and only if $\l^i \subset
\l^{i+1}$ for all $i$. Note that we do not yet assume in Definition \ref{defn:PS} that $S$ corresponds to a standard tableau.

\begin{definition}
Given $S,S'\in \PC(1)\times \cdots \times \PC(n)$ we write $S\leq T$ if $\l_S^k\leq \l_{S'}^k$ in the dominance order, for all $1\leq k\leq n$.  We refer to the induced partial order on $\SYT(n)$ as the dominance order on standard tableaux.
\end{definition}

Note that the dominance order on standard tableaux does not agree with the Kazhdan-Lusztig left cell order (or right cell order), when one identifies a tableau $T$ with the corresponding
left cell $(-,T)$.

\begin{example} The first example where the dominance order on tableaux disagrees with the KL left cell order is in $S_3$, for the tableau $ ((1,3),(2)) < ((1,2),(3)) $. These tableaux have the same shape, but comparable left cells can not exist in the same two-sided cell. \end{example}

\begin{example} The first example where the dominance order on tableaux disagrees with the KL left cell order for left cells in distinct two-sided cells is in $S_4$, for the tableaux $P
= ((1,4),(2),(3)) < ((1,2,3),(4)) = Q$. The left cell of $Q$ is generated by $b_u$, and any element $b_w$ in this cell will have $u$ in its right descent set. The left cell of $P$ is
generated by $b_{sts}$, and any element $b_w$ in this cell will have $s$ and $t$ in its right descent set. Therefore, these left cells intersect only in $b_{w_0}$, which is in a strictly
smaller left cell than both of them. Hence the left cells are incomparable. \end{example}

\begin{remark} In \cite{Taskin} various partial orders on $\SYT(n)$ are compared, including the cell order. None of the orders in that paper is the dominance order we consider here; the dominance order in \cite{Taskin} also involves comparing skew tableau within a tableau. No two tableaux of the same shape can be comparable in any of the orders in \cite{Taskin}, while it is important that they be comparable here. \end{remark}

\begin{theorem}\label{thm:PTprops}
We have:
\begin{enumerate}\setlength{\itemsep}{2pt}
\item $\PB_S\not \simeq 0$ if and only if $S$ is a standard tableau.
\item If $T$ is a standard tableau, then
\begin{enumerate}\setlength{\itemsep}{2pt}
\item  $\PB_T$ is constructed from bimodules $B_{P,Q}$ with $P,Q\leq T$.
\item  $B_{P,Q}\otimes \PB_T \simeq 0$ unless $Q\geq T$.
 \item $\PB_T\otimes B_{P,Q}\simeq 0$ unless $P\geq T$.
 \end{enumerate}
 \item There is a decomposition of identity $\one\simeq (\bigoplus_T \PB_T, d)$ in which the differential respects the reverse of the dominance order.
 \item A complex $C\in \KC^-(\SBim_n)$ satisfies $\Cone(\a_{T^k})\otimes C\simeq 0$ for all $1\leq k\leq n$ if and only if $\PB_T\otimes C\simeq C$ (and similarly for tensoring with $C$ on the left).  Thus $\PB_T$ projects onto the joint $(\a_{T^1},\ldots,\a_{T^n})$-eigencategory of $\FT_1,\ldots,\FT_n$.
\end{enumerate}
\end{theorem}

Statement (3) deserves a bit more explanation.  The claim is that there is a differential $d$ on $\bigoplus_T \PB_T$ such that the component $d_{T,U}$ from $\PB_U$ to $\PB_T$ vanishes unless $U\geq T$, $d_{T,T}$ is the given differential on $\PB_T$, and the resulting complex is homotopy equivalent to the monoidal identity.

\begin{proof} 
Recall the relative idempotents $\PB_{\mu,\l}\in \KC^-(\SBim_n)$, indexed by pairs $(\mu,\l)\in \PC(n-1)\times \PC(n)$ with $\mu\subset \l$ (\S \ref{subsec:proof}).  Then $\PB_{\mu}\sqcup \one_1$ is a convolution of idempotents of the form $\PB_{\mu,\l}$ where $\l$ ranges over the partitions of $n$ containing $\mu$, while $\PB_\l$ is a convolution of idempotents $\PB_{\mu,\l}$ where $\mu$ ranges over the paritions of $n-1$ contained in $\l$.  If $\l'\in \PC(n)$ does not contain $\mu\in \PC(n-1)$, then $(\PB_{\mu}\sqcup \one_1) \otimes \PB_{\l'}$ is contractible, since it is a convolution of complexes $\PB_{\mu,\l}\otimes \PB_{\mu',\l'}$ with $\mu\neq \mu'$, and these idempotents are orthogonal.  This proves the ``only if'' direction of (1).

Statement (2) is an immediate consequence of Lemma \ref{lem:relprojkills}.  Finally, statement (3) follows from tensoring together the decompositions of identity $\one\simeq (\bigoplus_{\mu\in \PC(k)} \PB_{\mu}\sqcup \one_{n-k},d)$ for $1\leq k\leq n$, and contracting the terms which are contractible via the ``only if'' direction of (1).

Fix $T \in \SYT(n)$. If there were a complex $C \in \KC^-(\SBim_n)$ such that $\PB_U \otimes C \simeq 0$ for all $U \ne T$, then one must have $\PB_T \otimes C \simeq C$, using the decomposition of identity in (3).  If $C\not \simeq 0$, then this would imply $\PB_T\not\simeq 0$, giving the ``if" direction of (1). The bounded quasi-idempotent complex $\KB_T$, constructed in the next section, will serve as such a complex $C$.

Statement (4) is an implication of statement (4) in Theorem \ref{thm:relDiag1less}.
\end{proof}

The observations made following the statement of Theorem \ref{thm:typeAdiag} have analogues here as well: for each convex subset of $K\subset \SYT(n)$ we may define a subquotient idempotent $\PB_K$, just as in Definition \ref{def:PK}.  If $I\subset \SYT(n)$ is a poset ideal, then $\PB_I$ is has the structure of a counital idempotent with complementary unital idempotent $\PB_{I^c}$, where $I^c=\SYT(n)\setminus I$.  We have the following characterization of the counital idempotents constructed in this way.  We state it only for $\PB_{\leq T}$, to avoid a potentially confusing conflict of notation with Proposition \ref{prop:counitalChar}.

\begin{proposition}\label{prop:PTleqchar}
Let $T\subset \SYT(n)$ be given. The idempotent $\PB_{\leq T}$ satisfies:
\begin{enumerate}
\item $\PB_{\leq T}$ is constructed from bimodules $B_{P,Q}$ with $P,Q\leq T$.
\item there is a chain map $\e:\PB_{\leq T}\rightarrow \one$ such that $\Cone(\e)\otimes B_{P,Q}\simeq 0$ when $P\leq T$ and $B_{P,Q}\otimes \Cone(\e)\simeq 0$ when $Q\leq T$.
\end{enumerate}
Furthermore,
\begin{enumerate}
\item[(3)] a complex $C\in \KC^-(\SBim_n)$ satisfies $\PB_{\leq T}\otimes C\simeq C$ (resp.~$C\otimes \PB_{\leq T}\simeq C$) if and only if $C\in \ip{B_{P,Q}\:|\: P\leq T}$ (resp.~$C\in \ip{B_{P,Q}\:|\: Q\leq T}$).
\item[(4)] properties (1) and (2) characterize the pair $(\PB_{\leq T},\e)$ up to canonical isomorphism in $\KC^-(\SBim_n)$.
\end{enumerate}
\end{proposition}

\begin{proof}
Essentially identical to the proof of Proposition \ref{prop:counitalChar}. 
\end{proof}
 
The idempotent $\PB_T$ itself can be characterized as the ``difference'' between $\PB_{\leq T}$ and $\PB_{<T}$ (see \S 5.2 of \cite{Hog17a}).

We now say a few words on the construction of $\PB_T$.  First, if $\mu,\l\in \PC(n)$ are comparable, recall the complexes $\CB_{\mu,\l}=\CB_{\a_\mu,\a_\l}(\FT_n)\in \KC^-(\SBim_n)$ from \S 7.2 of \cite{ElHog17a}.  We remind the reader that $\l$ and $\mu$ denote scalar objects /functors in \cite{ElHog17a}, and a shifted copy of $\one$ in $\KC^-(\SBim_n)$ is \emph{small} if it is supported in strictly negative homological degrees.  We continue to let $\l,\mu$ denote partitions, and $\Sigma_\l,\Sigma_\mu$ denote the corresponding scalar objects (shifts of the identity).  If we have partitions $\mu<\l$ in the dominance order then the shift $\Sigma_\l\Sigma_\mu\inv$ is small. \revadd{The infinite direct sum $\bigoplus_{k\geq 0} (\Sigma_\l\Sigma_\mu\inv)^k$ is well-defined in $\KC^-(\SBim_n)$, and tensoring with it preserves $\KC^-(\SBim_n)$ (without introducing infinite direct sums in any given homological degree). We have
\begin{equation} \label{eq:Clmreminder}
\CB_{\l,\mu} = \Sigma_\mu \Cone(\a_\l)\otimes \bigoplus_{k\geq 0} (\Sigma_\l\Sigma_\mu\inv)^k \qquad\text{ with twisted differential}.
\end{equation}
Then $\CB_{\l,\mu}$ is a categorification of $\frac{\ft_n - v^{2 \xbb(\l)}}{v^{2 \xbb(\mu)} - v^{2 \xbb(\l)}}$, which kills the eigenspace associated to $\l$ while preserving that associated to $\mu$. Meanwhile, $\CB_{\mu,\l}$ does not have the symmetric description (since $\bigoplus (\Sigma_{\mu} \Sigma_{\l}^{-1})^k \notin \KC^-(\SBim_n)$), but is another similarly-constructed infinite convolution in $\KC^-(\SBim_n)$.}

The important properties of the complexes $\CB_{\mu,\l}$ are the existence of distinguished triangles
\begin{subequations}
\begin{equation}
\CB_{\mu,\l}\rightarrow \one\rightarrow \CB_{\l,\mu}\rightarrow \CB_{\mu,\l}[1],\qquad\qquad (\mu<\l)
\end{equation}
\begin{equation}\label{eq:ctriang}
\Sigma_\l \CB_{\l,\mu}\rightarrow \Sigma_\mu \CB_{\l,\mu} \rightarrow \Cone(\a_{\l}) \rightarrow \Sigma_\l \CB_{\l,\mu}[1].
\end{equation}
\end{subequations}

The following is a direct application of Lemma 8.1 in \cite{ElHog17a}.
\begin{lemma}\label{lemma:PTinduction}
Fix $T\in \SYT(n)$ be a tableaux with shape $\l\in \PC(n)$, and let $U = T^{n-1} \in \SYT(n-1)$ with shape $\mu$.  Then
\[
\PB_T \simeq (\PB_U\sqcup \one_1)\otimes \bigotimes_{\nu\neq \l} \CB_{\nu,\l}.
\]
Here the tensor product is over all partitions $\nu\in \PC(n)$ with $\mu\subset \nu$ and $\nu\neq \l$.\qed
\end{lemma}
We remark that the complexes $\CB_{\nu,\l}$ commute with one another, given the vanishing of the obstructions in the family of eigenmaps $\{\a_\l\}_{\l\in \PC(n)}$ and Proposition A.17 in \cite{ElHog17a}.

This relationship between $\PB_T$ and $\PB_U$ may be viewed as a categorification of a well known relationship between $p_T$ and $p_U$ in the Hecke algebra $\HB_n$.  See \cite[Equation (11)]{IMO}.

\subsection{Finite quasi-idempotents}
\label{subsec:quasiidemp}

Now we consider a collection of finite complexes $\KB_T$ which are closely related to the idempotents $\PB_T$. They are defined in analogy with Lemma \ref{lemma:PTinduction} which (given this definition) will imply that $\PB_T$ is a locally finite convolution constructed from shifted copies of $\KB_T$.

\begin{definition}\label{def:KT}
Define complexes $\KB_T\in \KC^b(\SBim_n)$, $T\in \SYT(n)$, inductively by
\begin{subequations}
\begin{equation}\label{eq:K1}
\KB_{\square}=\one_1,
\end{equation}
\begin{equation}
\KB_T := (\KB_U\sqcup \one_1)\otimes \bigotimes_{\nu\neq \l} \Cone(\a_{\nu}).
\end{equation}
\end{subequations}
Here, $\shape(T) = \l\in \PC(n)$, $U = T^{n-1} \in \SYT(n-1)$, and $\nu$ runs over the partitions of $n$ such that $\mu\subset \nu$ and $\nu\neq \l$. 
\end{definition}

One should think that $\KB_T$ is the tensor product of all the $\Cone(\a_\nu)$ for partitions one could have taken but did not, along the path in the Young lattice which constructs $T$. That is,
\begin{equation} \label{eq:KTdef} \KB_T \simeq \bigotimes_{k=2}^n \bigotimes_{\substack{\nu \in \PC(k) \\ \nu \supset \l_T^{k-1} \\ \nu \ne \l_T^k}} \Cone(\a_{\nu}). \end{equation}

\begin{example} Let $T = {\: \Yvcentermath1 \young(12,34)}$. Then \[ \KB_T \cong \Cone(\a_{\yoo}) \ot \Cone(\a_{\yh}) \ot \Cone(\a_{\yho}) \ot \Cone(\a_{\ytoo}).\] \end{example}

\begin{lemma}\label{lemma:KabsorbsP}
We have $\PB_T\otimes \KB_T\simeq \KB_T\simeq \KB_T\otimes \PB_T$.
\end{lemma}
\begin{proof}
We focus on the equivalence $\KB_T\simeq \KB_T\otimes \PB_T$; the other is similar.  By (4) of Theorem \ref{thm:PTprops} we must show that $\KB_T\otimes \Cone(\a_{T^k})\simeq 0$ for all $1\leq k\leq n$. This we prove by induction \revise{ on $k$. The case $k=1$ is trivial, as $\KB_T = \PB_T = \one$. For a given $k$ the value of $n \ge k$ is irrelevant, so we may as well assume $n=k$.}  Let $U = T^{n-1}$. Then $\PB_U\sqcup\one_1$ commutes with $\Cone(\a_\l)$ for all $\l\in \PC(n)$ by Corollary \ref{cor:projectorconecommute}. From the definition of $\KB_T$ and the fact that $\PB_U\sqcup \one_1$ is idempotent it follows that
\begin{equation}
\KB_T \otimes (\PB_U\sqcup \one_1)\simeq \KB_T.
\end{equation}
Thus $\KB_T\otimes(-)$ annihilates $\Cone(\a_{T^k})$ for $1\leq k\leq n-1$, since the same is true of $\PB_U\sqcup \one_1$.  The fact that $\KB_T\otimes \Cone(\a_{T^n})\simeq 0$ follows from Lemma \ref{lem:relprojcones}.  This completes the proof.
\end{proof}

\begin{cor} \label{cor:KBtoo} Part (2) of Theorem \ref{thm:PTprops} also applies with $\KB_T$ replacing $\PB_T$. \qed \end{cor}

\begin{cor} If $U \ne T$ in $\SYT(n)$ then $\PB_U \ot \KB_T \simeq 0 \simeq \KB_T \ot \PB_U$. \qed \end{cor}

\begin{example}
If $T\in \SYT(n)$ is the unique one-row tableau, then $\PB_T=\PB_n$ is the categorified Jones-Wenzl idempotent and $\KB_T=\KB_n$ is the associated finite quasi-idempotent, both constructed in \cite{Hog15} (see also \cite{ElHog16a}).  Both of these complexes annihilate $B_w$ for $w\neq 1$.
\end{example}

\begin{example} \label{ex:Konecolumn}
If $T\in \SYT(n)$ is the unique one-column tableau, $\PB_T=\PB_{1^n}$ is the idempotent considered in \cite{AbHog17}.  The finite quasi-idempotent $\KB_T$ is homotopy equivalent to the $n-1$ dimensional Koszul complex associated to the action of $\a_1\otimes 1-1\otimes \a_1,\ldots,\a_{n-1}\otimes 1 -1\otimes \a_{n-1}$ on $B_{w_0}\in \SBim_n$, where $\a_i=x_i-x_{i+1}$.  We don't need or prove this fact, though it can be deduced from results in \cite{AbHog17}.  In any case, both $\PB_{1^n}$ and $\KB_{1^n}$ are built from $B_{w_0}$.
\end{example}

We will soon prove statements about $\PB_T$ in the Grothendieck group. This is subtle, because one must choose an appropriate subcategory of $\KC^-(\SBim_n)$ to get a non-zero Grothendieck group. However, $\KB_T$ is a bounded complex, so it has a well-defined image in the triangulated Grothendieck group $[\KC^b(\SBim_n)] \cong \HB_n$, which we now discuss.

Let us quickly recall how quasi-idempotents and idempotents are constructed when diagonalizing an operator. If $f$ is a linear operator with eigenvalues $\{\k_i\}_{i=0}^r$, and $\prod_i (f - \k_i) = 0$, then the projection to the $\k_i$ eigenspace is given by the formula
\begin{equation} p_i = \prod_{j \ne i} \frac{f - \k_j}{\k_i - \k_j}. \end{equation}
In particular, the element
\begin{equation} k_i = \prod_{j \ne i} (f - \k_j) \end{equation}
is a quasi-idempotent, and the scalar $\g = \prod_{j \ne i} (\k_i - \k_j)$ gives the proportionality between $k_i$ and $p_i$.

For example, \cite[Equation (11)]{IMO} uses this idea to inductively construct projections $p_T$ (which they denote $E_T$). One can inductively construct $k_T$ analogously, using only
the numerators in \cite[Equation (11)]{IMO}, and it is the obvious decategorification of Definition \ref{def:KT}. In particular, $[\KB_T] \mapsto k_T$ under the identification of the
Grothendieck group of $\KC^b(\SBim_n)$ with $\HB_n$.

\begin{cor} For all $T \in \SYT(n)$, $\KB_T$ is not contractible, and its image under the isomorphism of Grothendieck groups $[\KC^b(\SBim_n)]\cong \HB$ is a scalar multiple of $p_T$. \qed \end{cor}

\subsection{Grothendieck group considerations}
\label{subsec:groth}

If $C\in \KC(\SBim_n)$ is a minimal complex of Soergel bimodules, let $\supp(C)\subset \Z\times \Z$ denote the set of $(i,j)$ such that the chain bimodule $C^i$ has a summand of the form $B_w(-j)$.  Note that $\supp(C)[k](l) = (-k,-l)+\supp(C)$. If $C$ is not a minimal complex, we let $\supp(C):=\supp(C')$, where $C'\simeq C$ is minimal.  Since minimal complexes are unique up to isomorphism, this is well-defined.  Observe that
\[
\supp(\Sigma_\l)=\{(-2\cbb(\l),2\xbb(\l))\}.
\]

Let $\l,\nu\in \PC(n)$ be a pair of distinct partitions such that there exists $\mu\in \PC(n-1)$ with $\mu\subset \l$ and $\mu\subset \nu$.  In other words, $\l$ and $\nu$ become equal after deleting a single box.  Then $\l$ and $\nu$ are automatically comparable. Assume $\l>\nu$.  The set of such pairs $(\l,\nu)$ will be denoted $\QC(n)$.  If $(\l,\nu)\in \QC(n)$ then the box $\l\setminus \mu$ has larger content and is in a larger column than the box $\nu\setminus \mu$.   Thus, the vector
\[
\vb_{\l,\nu}:= (-\cbb(\l),\xbb(\l))-(-\cbb(\nu),\xbb(\nu))\in\Q\times \Q
\]
lies in the interior of the second quadrant.  Let $D(n)\subset \Q\times \Q$ denote the smallest convex subset containing $\Q_{\geq 0}\vb_{\l,\nu}$ for each $(\l,\nu)\in \QC(k)$ with $2\leq k\leq n$.  Note that $D(n)$ is an ``angle shaped'' region, bounded by two rays from the origin.

For each pair of subsets $E,E'\subset \Q\times \Q$, we have $E+E':=\{\eb+\eb'\:|\: \eb\in E, \eb'\in E'\}$ as usual.  For a subset $D\subset \Q\times \Q$ we write $E\subset \OC(D)$ if there is a finite set $E'$ such that $E\subset E'+D$.  Let $\KC^{\angle}(\SBim_n)\subset\KC^-(\SBim_n)$ denote the full subcategory consisting of complexes $C$ such that $\supp(C)\subset \OC(D(n))$.

Note that $\KC^{\angle}(\SBim_n)$ contains $\KC^b(\SBim_n)$ as a full subcategory.

\begin{lemma}
The category $\KC^{\angle}(\SBim_n)$ is closed under tensor products, direct sums, mapping cones, and shifts.
\end{lemma}
\begin{proof}
Closure under tensor products follows from the (easy) fact that $\supp(C\otimes C')\subset \OC(\supp(C)+\supp(C'))$ for all $C,C'\in \KC^-(\SBim_n)$ and $D(n)$ is closed under addition.  The other properties are clear.
\end{proof}

\begin{lemma}
The idempotents $\PB_T$ and $\PB_\l$ are in $\KC^{\angle}(\SBim_n)$ for all $T\in \SYT(n)$ and all $\l\in \PC(n)$.
\end{lemma}
\begin{proof}
Since $\PB_\l$ is a finite convolution involving the complexes $\PB_T$, it suffices to show that $\PB_T\in \KC^{\angle}(\SBim_n)$.  Recall the construction from Lemma \ref{lemma:PTinduction}, which states that $\PB_T$ is a tensor product of $(\PB_U \sqcup \one_1)$ and of $\CB_{\nu,\l},\CB_{\l,\nu}\in \KC^-(\SBim_n)$ for various $(\l,\nu)\in\QC(k)$.  Thus, it suffices to show each of these is in $\KC^{\angle}(\SBim_n)$. By induction, $\PB_U \in \KC^{\angle}(\SBim_n)$. There is a distinguished triangle \eqref{eq:ctriang} relating $\CB_{\nu,\l}$,$\CB_{\l,\nu}$, and $\one$, so it suffices to show that $\CB_{\l,\nu}\in\KC^{\angle}(\SBim_n)$ for all $(\l,\nu)\in \QC(k)$ with $2\leq k\leq n$.

\revcomment{Since the construction was moved earlier, recalled it here} The construction of $\CB_{\l,\nu}$ for $\l > \nu$ is given in \eqref{eq:Clmreminder}. It follows easily that $\CB_{\l,\nu}\in \KC^{\angle}(\SBim_n)$.  This completes the proof. \end{proof}


\begin{lemma} \label{lem:grothPvK}
For each $T\in \SYT(n)$ there exists a scalar $\g_T\in \Z[v,v\inv]$ such that $p_T = (\g_T)\inv k_T\in \HB_n^{\Q(v)}$ and $[\PB_T] = (\g_T)\inv [\KB_T]$ in $[\KC^{\angle}(\SBim_n)]$.
\end{lemma}

\begin{proof}
Fix a tableau $T\in \SYT(n)$.  Define a scalar $\g_T$ inductively by $\g_{\square}=1$ and 
\[
\g_T =  \g_U \cdot \prod_{\nu\neq \l}(v^{2\xbb(\l)} - v^{2\xbb(\nu)})
\]
where $\l$ is the shape of $T$, $U=T^{n-1}$, and $\nu$ ranges over the partitions of $n$ with $\mu\subset \nu$ and $\nu\neq \l$. 

Let $\nu^k\in \PC(k)$ be a partition.   Then we have the eigenmap $\a_{\nu^k}:\one_n\rightarrow \FT_k\sqcup \one_{n-k}$ and its mapping cone:
\[
(\Sigma_{\nu^k}[-1]\rightarrow \FT_k\sqcup \one_{n-k}).
\]
Tensoring (say, on the right) with $\PB_T$ gives a complex
\[
(\Sigma_{\nu^k}\PB_T[-1]\rightarrow \Sigma_{\l_T^k}\PB_T),
\]
which can be viewed as the (the shift $\Sigma_{\l_T^k}$ applied to the) mapping cone on an endomorphism of $\PB_T$ of degree $\Sigma_{\nu^k}\Sigma_{\l_T^k}\inv$.  This endomorphism will be denoted by $\frac{\a_{\nu^k}}{\a_{\l_T^k}}\in \Endgg(\PB_T)$, as in \cite[\S 7.4]{ElHog17a}.

To avoid copious notation, we will let the partitions $\nu^k \in \PC(k)$ below come from the same set as in the tensor product of \eqref{eq:KTdef}. We have
\begin{eqnarray*}
\KB_T &\simeq & \KB_T\otimes \PB_T\\
& =& \left(\bigotimes_{k=2}^n \bigotimes_{\nu^k} \Cone(\a_{\nu^k})\right)\otimes \PB_T\\
& \simeq & \left(\bigotimes_{k=2}^n \bigotimes_{\nu^k} \Cone(\a_{\nu^k})\otimes \PB_T\right)\\
& \simeq & \bigotimes_{k=2}^n \bigotimes_{\nu^k} \Sigma_{\l_T^k}\Cone\left(\frac{\a_{\nu^k}}{\a_{\l_T^k}}\right).\\
\end{eqnarray*}
In the first equivalence we used Lemma \ref{lemma:KabsorbsP}.  The second is by \eqref{eq:KTdef}.  The third uses the fact that $\PB_T$ is idempotent and commutes with each $\Cone(\a_{\nu^k})$, and the last equivalence follows from the comments made above.  It follows from this that
\[
[\KB_T]=\g_T [\PB_T].
\]
\end{proof}

\begin{proposition}\label{prop:groth}
The completed Grothendieck group \cite{AchStr} of $\KC^{\angle}(\SBim_n)$ is isomorphic to $\HB_n\otimes_{\Z[v,v\inv]} \Z[v]\llbracket v\inv \rrbracket$.  Under this isomorphism the class of $\PB_T$ gets sent to $p_T$.
\end{proposition}

\begin{proof}
In the completed Grothendieck group, the class of a complex is determined by its Euler characteristic, viewed as a Laurent series with coefficients in the usual Grothendieck group $[\KC^b(\SBim_n)]\cong \HB_n$.  This gives the first statement, modulo standard details.  The second statement follows from Lemma \ref{lem:grothPvK}.
\end{proof}

\section{Examples}
\label{sec:examples}

\revcomment{shortened this bit.}

In this section we give examples of eigenmaps, finite quasi-idempotents, and projectors related to $\FT_n$ for $n \le 3$.  For shorthand, we write $\CB_\l$ to denote $\Cone(\a_\l)$.  

\subsection{Two strands}

\revcomment{new intro paragraphs. added polynomial notation to many complexes.}

In type $A_1$ (with simple reflection $s$) morphisms are simple enough that we have the luxury of describing them in multiple ways. We can describe them directly using operations on polynomials. We can also describe them compactly using the diagrammatic calculus of \cite{EKho}, for readers familiar with this notation. When we display the differentials in a complex, we typically put the diagrammatic notation for the map above the arrow, and the polynomial notation below. We set $R = \R[x_1, x_2]$.

We also introduce the following useful shorthand:
\begin{equation}\label{eq:triangleCalc}
\ig{1}{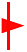} = \frac{1}{2}\left(\ig{1}{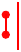} - \ig{1}{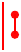}\right) = \frac{1}{2}(\a_s \ot 1 - 1 \ot \a_s).
\end{equation}
This morphism lives in the endomorphism ring of $B_s = R \ot_{R^s} R(1)$, which can be identified with $R \ot_{R^s} R$. Here we write $\a_s$ for the simple root associated to the simple reflection $s$. Because of potential conflicts with notation for eigenmaps, and for brevity, we write $\beta = \beta_s := \frac{1}{2} \a_s$. We write $\beta^{\ell} = \beta \ot 1$ for left-multiplication by $\beta$, and $\beta^r = 1 \ot \beta$ for right-multiplication. Thus
\[ \ig{1}{triang_red} = \beta^{\ell} - \beta^r. \]

\begin{remark} While $\beta$ has a denominator and is not defined integrally, the operators $\beta^{\ell} \pm \beta^r$ are defined over $\Z[x_1, x_2]$; for example, $\beta{^\ell} - \beta^r = x_1 \ot 1 - 1 \ot x_2$.\end{remark}

Recall that
\begin{equation} \FT_2 = \left(\begin{tikzpicture}[baseline=-.5em]
\node  (a) at (0,0) {$\un{B_s}(-1)$};
\node (b) at (2.5,0) {$B_s(1)$};
\node (c) at (5,0) {$R(2)$};
\path[->,>=stealth',shorten >=1pt,auto,node distance=1.8cm,thick]
	(a) edge node[above] {$\ig{.8}{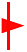}$} node[below] {$\beta^{\ell} - \beta^r$} (b)
	(b) edge node[above] {$\ig{.8}{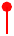}$} node[below] {$\mu$} (c);
\end{tikzpicture} \right)
\label{eq:FT2} \end{equation}
where the underline indicates homological degree zero. Here $\mu$ is the multiplication map sending $f \ot g \mapsto fg$ for $f,g \in R$. Thus $\mu(\beta^\ell) = \mu(\beta^r) = \beta$. Our two eigenmaps are given as follows.

\begin{equation} \FT_2 = \left(\begin{tikzpicture}[baseline=-.5em]
\node  (a) at (0,0) {$\un{B_s}(-1)$};
\node (b) at (2.5,0) {$B_s(1)$};
\node (c) at (5,0) {$R(2)$};
\node (d) at (0,-2) {$R(-2)[0]$};
\node  (e) at (5,-2) {$R(2)[-2]$};
\path[->,>=stealth',shorten >=1pt,auto,node distance=1.8cm,thick]
	(a) edge node[above] {$\ig{.8}{shorttriang_red}$} node[below] {$\beta^{\ell} - \beta^r$} (b)
	(b) edge node[above] {$\ig{.8}{counit_red}$} node[below] {$\mu$} (c)
	(d) edge node[right] {$\a_{\yoo}$} (a)
	(d) edge node[left] {$\ig{.8}{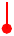}$} (a)
	(e) edge node[right] {$\a_{\yt}$} (c)
	(e) edge node[left] {$\Id$} (c);
\end{tikzpicture} \right)
\label{eq:FT2_withmaps} \end{equation}
\revadd{In polynomial notation the map $\a_{\yoo}$ is the coproduct map $\Delta$ sending $1 \mapsto \beta^\ell + \beta^r$. Note that $(\beta^\ell)^2 = (\beta^r)^2$ since $\beta^2$ is $s$-invariant.}

Thus
\begin{equation} \CB_{\yoo} \simeq \left( R(-2) \longrightarrow \un{B_s}(-1) \longrightarrow B_s(1) \longrightarrow R(2) \right), \end{equation}
\begin{equation} \CB_{\yt} \simeq \left( \un{B_s}(-1) \longrightarrow B_s(1) \right). \end{equation}
It is not hard to confirm directly that $\CB_{\yoo} \ot B_s \simeq 0 \simeq B_s \ot \CB_{\yoo}$ \revise{(indeed one may also see this by showing $\CB_{\yoo}$ is acyclic; see Example \ref{ex:one-col-eigen})}. Since $\CB_{\yt}$ is built from $B_s$, we see that $\CB_{\yoo} \ot \CB_{\yt} \simeq 0\simeq \CB_{\yt}\otimes \CB_{\yoo}$.   It is also easy to verify that $\CB_{\yoo}$ and $\CB_{\yt}$ are quasi-idempotent.

With only two eigenvalues, we have
\begin{equation} \KB_{\sytoct}\cong \KB_{\yoo} = \CB_{\yt}, \end{equation}
\begin{equation}\KB_{\sytot} \cong  \KB_{\yt} = \CB_{\yoo}. \end{equation}

One obtains $\PB_{\yt}$ by gluing together infinitely many copies of $\KB_{\yoo}$ as in the picture below.
\[
\PB_{\yt} = \left(\begin{tikzcd}[column sep = 1.4em]
&&& R(-4) \arrow[r,"\ig{.8}{unit_red}"] \arrow[rd,"-\Id"] & B_{s}(-3) \arrow[r,"\ig{.8}{shorttriang_red}"] & B_{s}(-1) \arrow[r,"\ig{.8}{counit_red}"] & \underline{R}\\
& R(-8) \arrow[r,"\ig{.8}{unit_red}"] \arrow[rd,"-\Id"] & B_{s}(-7) \arrow[r,"\ig{.8}{shorttriang_red}"] & B_{s}(-5)\arrow[r,"\ig{.8}{counit_red}"]  & R(-4) &&\\
\cdots \arrow[r,"\ig{.8}{shorttriang_red}"] &  B_{s} \arrow[r,"\ig{.8}{counit_red}"]  & R(-8) &&&&
\end{tikzcd}\right),
\]
which is homotopy equivalent to
\[
\PB_{\yt} \simeq \left(\begin{tikzcd}
\cdots \arrow{r}{\ig{.8}{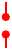}}[swap]{\beta^\ell + \beta^r}  & B_s(-7)   \arrow{r}{\ig{.8}{shorttriang_red}}[swap]{\beta^\ell - \beta^r}  & B_{s}(-5) \arrow{r}{\ig{.8}{twodots_red}}[swap]{\beta^\ell + \beta^r}  & B_{s}(-3) \arrow{r}{\ig{.8}{shorttriang_red}}[swap]{\beta^\ell - \beta^r}  & B_{s}(-1) \arrow{r}{\ig{.8}{counit_red}}[swap]{\mu} & \underline{R}
\end{tikzcd}\right).
\]
The complementary idempotent $\PB_{\yt}$ is obtained as below. 
\[
\PB_{\yoo} = \left(\begin{tikzcd}
\cdots   \arrow{r}{\ig{.8}{twodots_red}}[swap]{\beta^\ell + \beta^r}  & B_s(-7)   \arrow{r}{\ig{.8}{shorttriang_red}}[swap]{\beta^\ell - \beta^r} & B_{s}(-5)  \arrow{r}{\ig{.8}{twodots_red}}[swap]{\beta^\ell + \beta^r}  & B_{s}(-3) \arrow{r}{\ig{.8}{shorttriang_red}}[swap]{\beta^\ell - \beta^r} & \underline{B_{s}(-1)}
\end{tikzcd}\right).
\]
This is clearly a locally finite convolution constructed from shifted copies of $\KB_{\yt}$.

\subsection{Three strands}

\revcomment{modified this paragraph} Beyond $A_1$ it is no longer fruitful to use polynomial notation. We now study type $A_2$ with simple reflections $s$ and $t$. We adopt the diagrammatic calculus of \cite{EKho} and the ``thick'' calculus from \cite{EThick}. We denote $s$ by the color red and $t$ by the color blue in the diagrams below. The color purple represents $B_{sts}$, and morphisms involving it are explained in \cite[\S 4]{EThick}.

\subsubsection{The eigenmaps}
The minimal complex of $\FT_3$, with its eigenmaps, is pictured below.
\begin{equation}
\begin{tikzpicture}[baseline=-.2em]
\tikzstyle{every node}=[font=\scriptsize]
\node (a) at (0,0) {$\underline{{B_{sts}}}(-3)$};
\node at (2,.25) {$B_{sts}(-1)$};
\node at (2,-.25) {$B_{sts}(-1)$};
\node (c) at (4,.75) {$B_{sts}(1)$};
\node at (4,.25) {$B_{sts}(1)$};
\node at (4,-.25) {$B_{s}(1)$};
\node (yt) at (4,-.75) {$B_{t}(1)$};
\node (d) at (6,.5) {$B_{sts}(3)$};
\node at (6,0) {$B_{st}(2)$};
\node at (6,-.5) {$B_{ts}(2)$};
\node (e) at (8,.25) {$B_{st}(4)$};
\node at (8,-.25) {$B_{ts}(4)$};
\node (f) at (10,.25) {$B_{s}(5)$};
\node at (10,-.25) {$B_{t}(5)$};
\node (g) at (12,0) {$\one(6)$};
\node (x) at (0,-4) {$\one(-6)[0]$};
\node (y) at (4,-4) {$\one(0)[-2]$};
\node (z) at (12,-4) {$\one(6)[-6]$};
\node (b1) at (1.6,0) {};
\node (b2) at (2.4,0) {};
\node (c1) at (3.6,0) {};
\node (c2) at (4.4,0) {};
\node (d1) at (5.6,0) {};
\node (d2) at (6.4,0) {};
\node (e1) at (7.6,0) {};
\node (e2) at (8.4,0) {};
\node (f1) at (9.6,0) {};
\node (f2) at (10.4,0) {};
\path[->,>=stealth',shorten >=1pt,auto,node distance=1.8cm,
  thick]
(a) edge node[above] {} (b1)
(b2) edge node[above] {} (c1)
(c2) edge node[above] {} (d1)
(d2) edge node[above] {} (e1)
(e2) edge node[above] {} (f1)
(f2) edge node[above] {} (g)
(x) edge node[right] {$\a_{\yooo}$} (a)
(y) edge node[right] {$\a_{\yto}$} (yt)
(z) edge node[right] {$\a_{\yh}$} (g)
(x) edge node[left] {$\ig{.8}{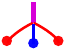}$} (a)
(y) edge node[left] {$\sqmatrix{0\\0\\ \ \ \ig{.8}{unit_red} \ \  \\ \ig{.8}{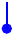}}$} (yt)
(z) edge node[left] {$\Id$} (g);
\end{tikzpicture}
\end{equation}
The components $d^k:C^k\rightarrow C^{k+1}$ of the differential in this compex is given by the following matrices:
\begin{equation} \label{thosearesomedifferentials}
d^0 =\sqmatrix{-\ig{.8}{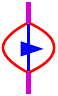}\\ \ig{.8}{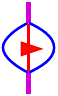}},\qquad
d^1 = \sqmatrix{\ig{.8}{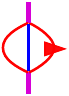} & \ig{.8}{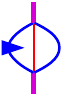} \\ \ig{.8}{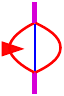} & \ig{.8}{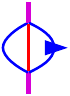} \\ \ig{.8}{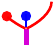} & 0 \\ 0 & \ig{.8}{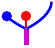}}, \qquad
d^2 = \sqmatrix{\ig{.8}{lefttriang_red} & -\ig{.8}{righttriang_red} & 0 & 0 \\ \ig{.8}{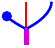} & 0 & -\ig{.8}{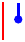} & \ig{.8}{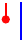}\\ 0 & \ig{.8}{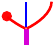} & \ig{.8}{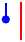} & -\ig{.8}{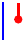}},
\end{equation}
\[
d^3 = \sqmatrix{\ig{.8}{sts_to_st} &  -\ig{.8}{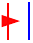} +  \ig{.8}{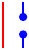} & \ig{.8}{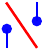} \\ \ig{.8}{sts_to_ts} &  \ig{.8}{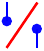} &  \ig{.8}{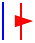} +  \ig{.8}{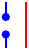} },\qquad
d^4 = \sqmatrix{-\ig{.8}{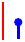} & \ig{.8}{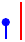} \\ \ig{.8}{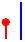} & -\ig{.8}{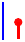}},\qquad
d^5 = \sqmatrix{\ig{.8}{counit_red} & \ig{.8}{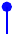}}.
\]

Let us indicate why each of the above maps $\a_\l$ is actually a $\l$-equivalence. 

For $\l=\yooo$, recall that the homology of a Rouquier complex $F(\b)$ (in the usual sense for complexes of $R$-bimodules) depends only on the permutation represented by $\b$ and the number of crossings in $\b$ (this is easy to see using the technology of standard bimodules).  In particular the homology $\FT_3$ is isomorphic to $R(-6)$, and in fact the map $\a_{\yooo}$ is a quasi-isomorphism of complexes of $R$-bimodules.  From this there are various approaches (in addition to direct computation) which one can use to quickly deduce that $\CB_{\yooo} \ot B_{w_0} \simeq 0$, and therefore $\a_{\yooo}$ is a $\yooo$-equivalence. One was done in \cite{AbHog17}.  Another
uses standard filtrations and their interaction with the functor $(-) \ot B_{w_0}$, together with work of Libedinsky and Williamson \cite{LibWil}.

It is obvious that $\a_{\yh}$ is a $\yh$-equivalence, since its cone is homotopy equivalent to a complex with no summands of the form $R(a)[b]$.

Now, consider $\Cone(\a_{\yto})$.  Note that the bimodules in cells $\neq \yooo$ form a subcomplex of the form
\begin{equation} \Phi = 
\begin{tikzpicture}[baseline=-.2em]
\tikzstyle{every node}=[font=\scriptsize]
\node (b) at (2,0) {$\one(0)$};
\node (c) at (4,.25) {$B_{s}(1)$};
\node at (4,-.25) {$B_{t}(1)$};
\node (d) at (6,.25) {$B_{st}(2)$};
\node at (6,-.25) {$B_{ts}(2)$};
\node (e) at (8,.25) {$B_{st}(4)$};
\node at (8,-.25) {$B_{ts}(4)$};
\node (f) at (10,.25) {$B_{s}(5)$};
\node at (10,-.25) {$B_{t}(5)$};
\node (g) at (12,0) {$\one(6)$};
\node (b1) at (1.6,0) {};
\node (b2) at (2.4,0) {};
\node (c1) at (3.6,0) {};
\node (c2) at (4.4,0) {};
\node (d1) at (5.6,0) {};
\node (d2) at (6.4,0) {};
\node (e1) at (7.6,0) {};
\node (e2) at (8.4,0) {};
\node (f1) at (9.6,0) {};
\node (f2) at (10.4,0) {};
\path[->,>=stealth',shorten >=1pt,auto,node distance=1.8cm,
  thick]
(b2) edge node[above] {} (c1)
(c2) edge node[above] {} (d1)
(d2) edge node[above] {} (e1)
(e2) edge node[above] {} (f1)
(f2) edge node[above] {} (g);
\end{tikzpicture}
\end{equation}
(The first term is in homological degree 1.)
	
This lovely complex $\Phi$ has the property that
\begin{equation} \label{eq:PhiBs}\Phi \ot B_s \cong \left( B_{w_0}(2) \longrightarrow B_{w_0}(4) \right), \end{equation}
\begin{equation} \label{eq:PhiBt}\Phi \ot B_t \cong \left( B_{w_0}(2) \longrightarrow B_{w_0}(4) \right). \end{equation} 
(The first term is in homological degree 3.)
From this it follows that $\CB_{\yto}$ sends the simple cell to the longest cell, so $\a_{\yto}$ is a $\yto$-equivalence.  We remark that while \eqref{eq:PhiBs} and \eqref{eq:PhiBt} look similar, the differentials are different and the complexes are non-isomorphic.

\begin{remark} The reader can observe from the differentials above that $\FT_3$ has a filtration $0 \subset F_{\yh} \subset F_{\yto} \subset F_{\yooo} = \FT_3$, where $F_{\l}$ is built
from shifts of objects $B_w$ in cells $\ge \l$. Moreover, the eigenmap $\a_\l$ maps into the subcomplex $F_\l$. The complex $\Phi$ above was the cone of the map $\Sigma_{\yto} \to
F_{\yto}$.

That such a filtration should exist is not obvious. For example, it is not clear that the components of the differential in $\FT_3$ mapping from $B_s(1) \oplus B_t(1)$ in homological
degree $2$ to $B_{sts}(3)$ in homological degree $3$ should be zero.

One can ask whether, for any finite Coxeter group, the full twist $\FT$ has a filtration as above, by subcomplexes associated with two-sided cells. We prove this for dihedral groups in
the sequel. \end{remark}

\begin{remark} Minimal complexes are well-defined up to isomorphism of chain complexes. The minimal complex $\FT_3$ does have non-trivial automorphisms, the most interesting of which come from the non-trivial maps $B_{st}(2) \to B_{sts}(3)$ and $B_{ts}(2) \to B_{sts}(3)$ in homological degree $3$. A subtle point is that an isomorphic minimal complex may give rise to a different filtration $0 \subset F'_{\yh} \subset F'_{\yto} \subset F'_{\yooo} \cong \FT_3$, where $F'_{\yto} \ncong F_{\yto}$ as complexes in $\KC^b(\SBim_n)$ (though they are isomorphic modulo $\SBim_{< \yto}$). \end{remark}

\subsubsection{The quasi-idempotents}

We have four tableaux for $n=3$, hence four complexes $\KB_T$:
\[ \KB_{\sytoth} = \CB_{\yoo} \ot \CB_{\yto},\qquad \KB_{\sytotch} = \CB_{\yoo} \ot \CB_{\yh},\]
\[ \KB_{\sytohct} = \CB_{\yt} \ot \CB_{\yooo},\qquad \KB_{\sytoctch} = \CB_{\yt} \ot \CB_{\yto}. \]
Let's record for posterity the involutions corresponding to these tableaux: $1, t,s,sts$, in this order.

We have found explicit forms of their minimal complexes which we include below.

The minimal complex of $\KB_{\sytoth}$ is 
\begin{equation}\label{eq:K123}
\begin{tikzpicture}[baseline=-.2em]
\tikzstyle{every node}=[font=\scriptsize]
\node (a) at (0,0) {$\one(-2)$};
\node at (1.7,.25) {$B_{s}(-1)$};
\node at (1.7,-.25) {$B_{t}(-1)$};
\node (b) at (1.7,0) {\ \ \ \ \ \ \ \ \ \  };
\node at (3.4,.25) {$B_{ts}(0)$};
\node at (3.4,-.25) {$B_{st}(0)$};
\node (c) at (3.4,0) {\ \ \ \ \ \ \ \ \ \  };
\node at (5.1,.75) {$B_{sts}(1)$};
\node at (5.1,.25) {$B_{st}(2)$};
\node at (5.1,-.25) {$B_{ts}(2)$};
\node at (5.1,-.75) {$\one(2)$};
\node (d) at (5.1,0) {\ \ \ \ \ \ \ \ \ \  };
\node at (6.8,1.25) {$B_{sts}(3)$};
\node at (6.8,.75) {$B_{sts}(3)$};
\node at (6.8,.25) {$B_{s}(3)$};
\node at (6.8,-.25) {$B_{s}(3)$};
\node at (6.8,-.75) {$B_{t}(3)$};
\node at (6.8,-1.25) {$B_{t}(3)$};
\node (e) at (6.8,0) {\ \ \ \ \ \ \ \ \ \ };
\node at (8.5,.75) {$B_{sts}(5)$};
\node at (8.5,.25) {$B_{st}(4)$};
\node at (8.5,-.25) {$B_{ts}(4)$};
\node at (8.5,-.75) {$\one(4)$};
\node (f) at (8.5,0) {\ \ \ \ \ \ \ \ \ \  };
\node at (10.2,.25) {$B_{ts}(6)$};
\node at (10.2,-.25) {$B_{st}(6)$};
\node (g) at (10.2,0) {\ \ \ \ \ \ \ \ \ \ };
\node at (11.9,.25) {$B_{s}(7)$};
\node at (11.9,-.25) {$B_{t}(7)$};
\node (h) at (11.9,0) {\ \ \ \ \ \ \ \ \ \  };
\node (i) at (13.6,0) {$\one(8)$};
\path[->,>=stealth',shorten >=1pt,auto,node distance=1.8cm,
  thick]
(a) edge node{} (b)
(b) edge node{} (c)
(c) edge node{} (d)
(d) edge node{} (e)
(e) edge node{} (f)
(f) edge node{} (g)
(g) edge node{} (h)
(h) edge node{} (i);
\end{tikzpicture}.
\end{equation}
Here, the left-most term is in homological degree 0.

The minimal complex of $\KB_{\sytotch}$ is
\begin{equation}\label{eq:K12c3}
\begin{tikzpicture}[baseline=-.2em]
\tikzstyle{every node}=[font=\scriptsize]
\node (b) at (1.7,0) {$B_{t}(-1)$};
\node at (3.4,.25) {$B_{ts}(0)$};
\node at (3.4,-.25) {$B_{st}(0)$};
\node (c) at (3.4,0) {\ \ \ \ \ \ \ \ \ \  };
\node at (5.1,.75) {$B_{sts}(1)$};
\node at (5.1,.25) {$B_{st}(2)$};
\node at (5.1,-.25) {$B_{ts}(2)$};
\node at (5.1,-.75) {$B_s(1)$};
\node (d) at (5.1,0) {\ \ \ \ \ \ \ \ \ \  };
\node at (6.8,1.25) {$B_{sts}(3)$};
\node at (6.8,.75) {$B_{sts}(3)$};
\node at (6.8,.25) {$B_{s}(3)$};
\node at (6.8,-.25) {$B_{s}(3)$};
\node at (6.8,-.75) {$B_{t}(3)$};
\node at (6.8,-1.25) {$B_{t}(3)$};
\node (e) at (6.8,0) {\ \ \ \ \ \ \ \ \ \ };
\node at (8.5,.75) {$B_{sts}(5)$};
\node at (8.5,.25) {$B_{st}(4)$};
\node at (8.5,-.25) {$B_{ts}(4)$};
\node at (8.5,-.75) {$B_s(5)$};
\node (f) at (8.5,0) {\ \ \ \ \ \ \ \ \ \  };
\node at (10.2,.25) {$B_{ts}(6)$};
\node at (10.2,-.25) {$B_{st}(6)$};
\node (g) at (10.2,0) {\ \ \ \ \ \ \ \ \ \ };
\node (h) at (11.9,0) {$B_{t}(7)$};
\path[->,>=stealth',shorten >=1pt,auto,node distance=1.8cm,
  thick]
(b) edge node{} (c)
(c) edge node{} (d)
(d) edge node{} (e)
(e) edge node{} (f)
(f) edge node{} (g)
(g) edge node{} (h);
\end{tikzpicture}
\end{equation}
The left-most term is in homological degree 1.

The minimal complex of $\KB_{\sytohct}$ is
\begin{equation}\label{eq:K13c2}
\begin{tikzpicture}[baseline=-.2em]
\node (b) at (0,0) {$B_{s}(-7)$};
\node at (2.5,.25) {$B_{sts}(-5)$};
\node at (2.5,-.25) {$B_{s}(-5)$};
\node (c) at (2.5,0) {\ \ \ \ \ \ \ \ \ \  \ \ \ \  };
\node at (5,.25) {$B_{sts}(-3)$};
\node at (5,-.25) {$B_{sts}(-3)$};
\node (d) at (5,0) {\ \ \ \ \ \ \ \ \ \   \ \ \ \ };
\node at (7.5,.25) {$B_{sts}(-1)$};
\node at (7.5,-.25) {$B_{s}(-1)$};
\node (e) at (7.5,0) {\ \ \ \ \ \ \ \ \ \  \ \ \ \ };
\node (f) at (10,0) {$B_s(1)$};
\path[->,>=stealth',shorten >=1pt,auto,node distance=1.8cm,
  thick]
(b) edge node{} (c)
(c) edge node{} (d)
(d) edge node{} (e)
(e) edge node{} (f);
\end{tikzpicture}.
\end{equation}
The left-most term is in homological degree $-1$.

Finally, the minimal complex of $\KB_{\sytoctch}$ is
\begin{equation}\label{eq:K1c2c3}
\begin{tikzpicture}[baseline=-.2em]
\node (c) at (0,0) {$B_{sts}(-5)$};
\node at (2.5,.25) {$B_{sts}(-3)$};
\node at (2.5,-.25) {$B_{st}(-3)$};
\node (d) at (2.5,0) {\ \ \ \ \ \ \ \ \ \  \ \ \ \  };
\node (e) at (5,0) {$B_{sts}(-1)$};
\path[->,>=stealth',shorten >=1pt,auto,node distance=1.8cm,
  thick]
(c) edge node{} (d)
(d) edge node{} (e);
\end{tikzpicture}.
\end{equation}
The left-most nonzero term is in homological degree 0.

Let us illustrate Corollary \ref{cor:KBtoo}, and give some remarks on the computations above.  For $\KB_{\sytoth}$, note that $\CB_{\yoo}\otimes (-)$ annihilates $B_w$ for $w=s, st, sts$.  Thus, to compute $\CB_{\yoo}\otimes \CB_{\yto}$ we may as well tensor $\CB_{\yoo}$ with $\Phi$, contract the contractible terms, and we obtain a convolution of the form (omitting shifts)
\begin{equation}
\KB_{\sytoth}  \ \ \simeq \ \ \left(\CB_{\yoo} \rightarrow \CB_{\yoo}  B_t\rightarrow \CB_{\yoo}B_{ts}\rightarrow \CB_{\yoo} B_{ts}\rightarrow \CB_{\yoo} B_t\rightarrow \CB_{\yoo}\right).
\end{equation}
Expanding $\CB_{\yoo}$, keeping track of differentials, and peforming a pair of Gaussian eliminations yields \eqref{eq:K123}. It is easy to prove that $\KB_{\sytoth} \ot (-)$ and $(-) \ot \KB_{\sytoth}$ annihilates any $B_w$ for $w \ne 1$, as tensoring with $\CB_{\yto}$ takes any such $B_w$ to a complex built from $B_{sts}$, and $\CB_{\yoo}$ kills $B_{sts}$.

We have a similar convolution description
\begin{equation}
\KB_{\sytotch}  \ \ \simeq \ \ \left(\CB_{\yoo}  B_t \rightarrow \CB_{\yoo}B_{ts}\rightarrow \CB_{\yoo} B_{ts}\rightarrow \CB_{\yoo} B_t\right),
\end{equation}
which yields \eqref{eq:K12c3} after expanding and simplifying. In fact, this gives rise to the formula
\begin{equation}
\KB_{\sytotch}  \ \ \simeq \ \ \CB_{\yoo} \ot B_t \ot \CB_{\yoo}.
\end{equation}
Now it is clear that $\KB_{\sytotch} \ot (-)$ will kill $B_w$ whenever $sw<w$, and $(-) \ot \KB_{\sytotch}$ will kill $B_w$ whenever $ws < w$.

Now, recall that $\CB_{\yt}=(B_s(-1)\rightarrow B_s(1))$.  Given that $B_s\otimes \FT_3\simeq (B_{sts}(-4)\rightarrow B_{sts}(-2)\rightarrow B_s(0))$ (see \eqref{eq:FTBs}) it follows that
\[
\CB_{\yt}\otimes \FT_3\simeq \Tot\left(
\begin{tikzpicture}[baseline=-2.5em]
\node (a) at (0,0) {$B_{sts}(-5)$};
\node (b) at (0,-2) {$B_{sts}(-3)$};
\node (c) at (2.5,0) {$B_{sts}(-3)$};
\node (d) at (2.5,-2) {$B_{sts}(-1)$};
\node (e) at (5,0) {$B_{s}(-2)$};
\node (f) at (5,-2) {$B_{s}(0)$};
\path[->,>=stealth',shorten >=1pt,auto,node distance=1.8cm,
  thick]
(a) edge node{} (c)
(c) edge node{} (e)
(b) edge node{} (d)
(d) edge node{} (f)
(a) edge node{} (b)
(c) edge node{} (d)
(e) edge node{} (f);
\end{tikzpicture}
\right).
\]
Then $\KB_{\sytohct}$ and $\KB_{\sytoctch}$ can both be expressed as mapping cones of maps
\[
\left(
\begin{tikzpicture}[baseline=-2.5em]
\node (a) at (0,0) {$B_{s}(-1)$};
\node (b) at (0,-2) {$B_{s}(1)$};
\path[->,>=stealth',shorten >=1pt,auto,node distance=1.8cm,
  thick]
(a) edge node{} (b);
\end{tikzpicture}
\right)(\text{shift})
\ \ \ \rightarrow \ \ \ \left(
\begin{tikzpicture}[baseline=-2.5em]
\node (a) at (0,0) {$B_{sts}(-5)$};
\node (b) at (0,-2) {$B_{sts}(-3)$};
\node (c) at (2.5,0) {$B_{sts}(-3)$};
\node (d) at (2.5,-2) {$B_{sts}(-1)$};
\node (e) at (5,0) {$B_{s}(-2)$};
\node (f) at (5,-2) {$B_{s}(0)$};
\path[->,>=stealth',shorten >=1pt,auto,node distance=1.8cm,
  thick]
(a) edge node{} (c)
(c) edge node{} (e)
(b) edge node{} (d)
(d) edge node{} (f)
(a) edge node{} (b)
(c) edge node{} (d)
(e) edge node{} (f);
\end{tikzpicture}
\right).
\]
From these we obtain \eqref{eq:K13c2} and \eqref{eq:K1c2c3}.

It is clear that $\KB_{\sytohct}$ kills $B_{sts}$ under tensor product on the right and left, since $\CB_{\yooo}$ does. Thus we have seen for all four tableaux that $\KB_T \ot B_{P,Q} \simeq 0$ for $P \ngeq T$, and $B_{P,Q} \ot \KB_T \simeq 0$ for $Q \ngeq T$.

\bibliographystyle{alpha}
\bibliography{bib}
\end{document}